\documentclass[10pt]{amsart}
\usepackage[english]{babel}
\usepackage[cp1251]{inputenc}
\usepackage{amsmath}
\usepackage{amssymb}
\usepackage{amsfonts}

\usepackage[linktocpage=true, colorlinks=true, linkcolor=blue, citecolor=blue, urlcolor=blue]{hyperref}

\begin{document}
\renewcommand{\refname}{References}
\renewcommand\contentsname{Contents}

\thispagestyle{empty}

\title[Mean-Square Approximation of Iterated Ito and 
Stratonovich Stochastic Integrals]
{Mean-Square Approximation of Iterated Ito and 
Stratonovich Stochastic Integrals of Multiplicities 1 to 6 
from the Taylor--Ito and 
Taylor--Stratonovich 
Expansions Using Legendre Polynomials}
\author[D.F. Kuznetsov]{Dmitriy F. Kuznetsov}
\address{Dmitriy Feliksovich Kuznetsov
\newline\hphantom{iii} Peter the Great Saint-Petersburg Polytechnic University,
\newline\hphantom{iii} Polytechnicheskaya ul., 29,
\newline\hphantom{iii} 195251, Saint-Petersburg, Russia}%
\email{sde\_kuznetsov@inbox.ru}

\maketitle {\small
\begin{quote}
\vspace{13mm}
\noindent{\sc Abstract.} 
The article is devoted to the practical material on
expansions and mean-square approximations of specific iterated
Ito and Stratonovich stochastic integrals of multiplicities 1 to 6
with respect to components of the multidimensional Wiener process
on the base of the method of generalized multiple Fourier series.
More precisely, we used the multiple Fourier--Legendre series
converging in the sense of norm in the space $L_2([t, T]^k)$ 
$(k=1,\ldots,6)$ for approximation of iterated Ito and 
Stratonovich stochastic integrals. The considered iterated Ito and 
Stratonovich stochastic integrals are part of the stochastic 
Taylor expansions (Taylor--Ito and Taylor--Stratonovich expansions).
Therefore, the results of the article can be useful for 
the construction of high-order strong numerical 
methods for Ito stochastic differential equations.
Expansions of iterated Ito and Stratonovich stochastic 
integrals of multiplicities 1 to 6
using Legendre polynomials 
are derived. The convergence with probability 1 of the mentioned  
method of generalized multiple Fourier series is proved for iterated
Ito stochastic integrals of arbitrary multiplicity $k$ $(k\in\mathbb{N})$
for the cases of multiple Fourier--Legendre series and
multiple trigonometric Fourier series.
\medskip
\end{quote}
}

\vspace{12mm}


\setlength{\baselineskip}{2.2em}

\tableofcontents

\setlength{\baselineskip}{1.2em}


\vspace{5mm}

\section{Introduction}

\vspace{5mm}

Let $(\Omega,$ ${\rm F},$ ${\sf P})$ be a complete probability space, let 
$\{{\rm F}_t, t\in[0,T]\}$ be a nondecreasing right-continous 
family of $\sigma$-algebras of ${\rm F},$
and let ${\bf f}_t$ be a standard $m$-dimensional Wiener stochastic 
process, which is
${\rm F}_t$-measurable for any $t\in[0, T].$ We assume that the components
${\bf f}_{t}^{(i)}$ $(i=1,\ldots,m)$ of this process are independent. 
Consider
an Ito stochastic differential equation (SDE) 
in the integral form

\vspace{-1mm}
\begin{equation}
\label{1.5.2}
{\bf x}_t={\bf x}_0+\int\limits_0^t {\bf a}({\bf x}_{\tau},\tau)d\tau+
\int\limits_0^t B({\bf x}_{\tau},\tau)d{\bf f}_{\tau},\ \ \
{\bf x}_0={\bf x}(0,\omega).
\end{equation}

\vspace{2mm}
\noindent
Here ${\bf x}_t$ is some $n$-dimensional stochastic process 
satisfying to the equation (\ref{1.5.2}). 
The nonrandom functions ${\bf a}: \mathbb{R}^n\times[0, T]\to\mathbb{R}^n$,
$B: \mathbb{R}^n\times[0, T]\to\mathbb{R}^{n\times m}$
guarantee the existence and uniqueness up to stochastic equivalence 
of a solution
of (\ref{1.5.2}) \cite{1982}. The second integral on 
the right-hand side of (\ref{1.5.2}) is 
interpreted as an Ito stochastic integral.
Let ${\bf x}_0$ be an $n$-dimensional random variable, which is 
${\rm F}_0$-measurable and 
${\sf M}\{\left|{\bf x}_0\right|^2\}<\infty$ 
(${\sf M}$ denotes a mathematical expectation).
We assume that
${\bf x}_0$ and ${\bf f}_t-{\bf f}_0$ are independent when $t>0.$

It is well known that one of the effective approaches 
to the numerical integration of 
Ito SDEs is an approach based on the Taylor--Ito and 
Taylor--Stratonovich expansions
\cite{1995}-\cite{2004}. The most important feature of such 
expansions is a presence in them of the so-called iterated
Ito and Stratonovich stochastic integrals, which play the key 
role for solving the 
problem of numerical integration of Ito SDEs and have the 
following form

\vspace{-1mm}
\begin{equation}
\label{ito}
J[\psi^{(k)}]_{T,t}=\int\limits_t^T\psi_k(t_k) \ldots \int\limits_t^{t_{2}}
\psi_1(t_1) d{\bf w}_{t_1}^{(i_1)}\ldots
d{\bf w}_{t_k}^{(i_k)},
\end{equation}

\vspace{2mm}
\begin{equation}
\label{str}
J^{*}[\psi^{(k)}]_{T,t}=
\int\limits_t^{*T}\psi_k(t_k) \ldots \int\limits_t^{*t_{2}}
\psi_1(t_1) d{\bf w}_{t_1}^{(i_1)}\ldots
d{\bf w}_{t_k}^{(i_k)},
\end{equation}

\vspace{1mm}
\noindent
where $\psi_1(\tau),\ldots,\psi_k(\tau)$ are  
nonrandom functions 
on $[t,T],$ ${\bf w}_{\tau}^{(i)}={\bf f}_{\tau}^{(i)}$
for $i=1,\ldots,m$ and
${\bf w}_{\tau}^{(0)}=\tau,$

\vspace{-3mm}
$$
\int\limits\ \hbox{and}\ \int\limits^{*}
$$ 

\vspace{3mm}
\noindent
denote Ito and 
Stratonovich stochastic integrals,
respectively; $i_1,\ldots,i_k = 0, 1,\ldots,m.$

Note that $\psi_l(\tau)\equiv 1$ $(l=1,\ldots,k)$ and
$i_1,\ldots,i_k = 0, 1,\ldots,m$ in  
\cite{1995}-\cite{2004}. At the same time 
$\psi_l(\tau)\equiv (t-\tau)^{q_l}$ ($l=1,\ldots,k$; 
$q_1,\ldots,q_k=0, 1, 2,\ldots $) and $i_1,\ldots,i_k = 1,\ldots,m$ in
\cite{1}-\cite{18}.

Effective solution 
of the problem of
mean-square approximation for collections 
of iterated Ito and Stratonovich stochastic integrals
(\ref{ito}) and (\ref{str})
composes the subject of this article.

\vspace{5mm}

\section{Theorems on Expansions of Iterated
Ito and Stratonovich Stochastic Itegrals}

\vspace{5mm}

Let us consider the effective
approach to expansion of the iterated Ito
stochastic integrals (\ref{ito}) 
\cite{3} (2006), \cite{3a}-\cite{17} (the so-called
method of generalized
multiple Fourier series). 
Sometimes these stochastic integrals are referred 
to in the literature as multiple stochastic integrals 
(see, for example, \cite{1995}).

The idea of this method is as follows: the iterated Ito stochastic 
integral (\ref{ito}) 
of multiplicity $k$ is represented as the multiple stochastic 
integral from the certain discontinuous nonrandom function of $k$ variables 
defined on the hypercube $[t, T]^k$, where $[t, T]$ is the interval of 
integration of the iterated Ito stochastic integral. Then, 
the indicated 
nonrandom function is expanded in the hypercube $[t, T]^k$
into the generalized 
multiple Fourier series converging 
in the mean-square sense
in the space 
$L_2([t,T]^k)$. After a number of nontrivial transformations we come 
(see Theorems 1, 2 below) to the 
mean-square convergening expansion of the iterated Ito stochastic 
integral into the multiple 
series of products
of standard  Gaussian random 
variables. The coefficients of this 
series are the coefficients of 
generalized multiple Fourier series for the mentioned nonrandom function 
of $k$ variables, which can be calculated using the explicit formula 
regardless of the multiplicity $k$ of the iterated Ito stochastic integral
(\ref{ito}).

Suppose that every $\psi_l(\tau)$ $(l=1,\ldots,k)$ is a continuous 
nonrandom function on $[t, T]$
(the case $\psi_1(\tau),\ldots,\psi_k(\tau)\in L_2([t, T])$
will be considered in Theorem~2).

Define the following function on the hypercube $[t, T]^k$

\vspace{-1mm}
\begin{equation}
\label{ppp}
K(t_1,\ldots,t_k)=
\begin{cases}
\psi_1(t_1)\ldots \psi_k(t_k)\ \ \hbox{for}\ \ t_1<\ldots<t_k\cr\cr\cr
0\ \ \hbox{otherwise}
\end{cases},\ \ t_1,\ldots,t_k\in[t, T],\ \ k\ge 2,
\end{equation}

\vspace{3mm}
\noindent
and 
$K(t_1)\equiv\psi_1(t_1)$ for $t_1\in[t, T].$

Suppose that $\{\phi_j(x)\}_{j=0}^{\infty}$
is a complete orthonormal system of functions in the space
$L_2([t, T])$.

The function $K(t_1,\ldots,t_k)$ is piecewise continuous in the 
hypercube $[t, T]^k.$
At this situation it is well known that the generalized 
multiple Fourier series 
of $K(t_1,\ldots,t_k)\in L_2([t, T]^k)$ is converging 
to $K(t_1,\ldots,t_k)$ in the hypercube $[t, T]^k$ in 
the mean-square sense, i.e.

\vspace{1mm}
$$
\hbox{\vtop{\offinterlineskip\halign{
\hfil#\hfil\cr
{\rm lim}\cr
$\stackrel{}{{}_{p_1,\ldots,p_k\to \infty}}$\cr
}} }\Biggl\Vert
K(t_1,\ldots,t_k)-
\sum_{j_1=0}^{p_1}\ldots \sum_{j_k=0}^{p_k}
C_{j_k\ldots j_1}\prod_{l=1}^{k} \phi_{j_l}(t_l)\Biggr\Vert_{L_2([t,T]^k)}=0,
$$

\vspace{4mm}
\noindent
where
\begin{equation}
\label{ppppa}
C_{j_k\ldots j_1}=\int\limits_{[t,T]^k}
K(t_1,\ldots,t_k)\prod_{l=1}^{k}\phi_{j_l}(t_l)dt_1\ldots dt_k
\end{equation}

\vspace{4mm}
\noindent
is the Fourier coefficient,
$$
\left\Vert f\right\Vert_{L_2([t,T]^k)}=\left(\int\limits_{[t,T]^k}
f^2(t_1,\ldots,t_k)dt_1\ldots dt_k\right)^{1/2}.
$$

\vspace{4mm}

Consider the partition $\{\tau_j\}_{j=0}^N$ of the interval
$[t,T]$ such that

\begin{equation}
\label{1111}
t=\tau_0<\ldots <\tau_N=T,\ \ \
\Delta_N=
\hbox{\vtop{\offinterlineskip\halign{
\hfil#\hfil\cr
{\rm max}\cr
$\stackrel{}{{}_{0\le j\le N-1}}$\cr
}} }\Delta\tau_j\to 0\ \ \hbox{if}\ \ N\to \infty,\ \ \
\Delta\tau_j=\tau_{j+1}-\tau_j.
\end{equation}

\vspace{4mm}                                                                   

{\bf Theorem 1} \cite{3} (2006), \cite{3a}-\cite{17}, \cite{xxx}-\cite{new-art-1-xxy}. 
{\it Suppose that
every $\psi_l(\tau)$ $(l=1,\ldots, k)$ is a continuous nonrandom function on 
$[t, T]$ and
$\{\phi_j(x)\}_{j=0}^{\infty}$ is a complete orthonormal system  
of continuous functions in the space $L_2([t,T]).$ Then

\vspace{1mm}
$$
J[\psi^{(k)}]_{T,t}\  =\ 
\hbox{\vtop{\offinterlineskip\halign{
\hfil#\hfil\cr
{\rm l.i.m.}\cr
$\stackrel{}{{}_{p_1,\ldots,p_k\to \infty}}$\cr
}} }\sum_{j_1=0}^{p_1}\ldots\sum_{j_k=0}^{p_k}
C_{j_k\ldots j_1}\Biggl(
\prod_{l=1}^k\zeta_{j_l}^{(i_l)}\ -
\Biggr.
$$

\vspace{4mm}
\begin{equation}
\label{tyyy}
-\ \Biggl.
\hbox{\vtop{\offinterlineskip\halign{
\hfil#\hfil\cr
{\rm l.i.m.}\cr
$\stackrel{}{{}_{N\to \infty}}$\cr
}} }\sum_{(l_1,\ldots,l_k)\in {\rm G}_k}
\phi_{j_{1}}(\tau_{l_1})
\Delta{\bf w}_{\tau_{l_1}}^{(i_1)}\ldots
\phi_{j_{k}}(\tau_{l_k})
\Delta{\bf w}_{\tau_{l_k}}^{(i_k)}\Biggr),
\end{equation}

\vspace{6mm}
\noindent
where $J[\psi^{(k)}]_{T,t}$ is defined by {\rm (\ref{ito}),}

\vspace{-1mm}
$$
{\rm G}_k={\rm H}_k\backslash{\rm L}_k,\ \ \
{\rm H}_k=\{(l_1,\ldots,l_k):\ l_1,\ldots,l_k=0,\ 1,\ldots,N-1\},
$$

$$
{\rm L}_k=\{(l_1,\ldots,l_k):\ l_1,\ldots,l_k=0,\ 1,\ldots,N-1;\
l_g\ne l_r\ (g\ne r);\ g, r=1,\ldots,k\},
$$

\vspace{3mm}
\noindent
${\rm l.i.m.}$ is a limit in the mean-square sense$,$
$i_1,\ldots,i_k=0,1,\ldots,m,$ 

\vspace{-1mm}
\begin{equation}
\label{rr23}
\zeta_{j}^{(i)}=
\int\limits_t^T \phi_{j}(s) d{\bf w}_s^{(i)}
\end{equation} 

\vspace{2mm}
\noindent
are independent standard Gaussian random variables
for various
$i$ or $j$ {\rm(}if $i\ne 0${\rm),}
$C_{j_k\ldots j_1}$ is the Fourier coefficient {\rm(\ref{ppppa}),}
$\Delta{\bf w}_{\tau_{j}}^{(i)}=
{\bf w}_{\tau_{j+1}}^{(i)}-{\bf w}_{\tau_{j}}^{(i)}$
$(i=0, 1,\ldots,m),$
$\left\{\tau_{j}\right\}_{j=0}^{N}$ is a partition of
the interval $[t, T],$ which satisfies the condition {\rm (\ref{1111})}.
}

\vspace{2mm}

It was shown 
that Theorem 1 is valid for convergence 
in the mean of degree $2n$ ($n\in \mathbb{N}$) 
\cite{10a} (Sect.~1.1.9, 1.11, 1.12), \cite{11} (Sect.~6, 15, 16)
and for convergence with probablity 1 (w.~p.~1) \cite{10a} (Sect.~1.7.2), \cite{9999}, \cite{new-new-2}.
Moreover, the complete orthonormal systems of Haar and 
Rademacher--Walsh functions in $L_2([t,T])$ 
can also be applied in Theorem 1
\cite{3}-\cite{11}.
The generalization of Theorem 1 for 
complete orthonormal with weigth $r(t_1)\ldots r(t_k)\ge 0$ systems
of functions in the space $L_2([t,T]^k)$ can be found in 
\cite{10}-\cite{11}, \cite{15b}, \cite{xxx}.
Another modification of Theorem 1 and Theorem~2 (see below) is connected with
the approximation of iterated stochastic integrals with respect 
to the infinite-dimensional $Q$-Wiener process 
\cite{10a}-\cite{10axx1} (Chapter 7), \cite{15f}, \cite{new-1}-\cite{new-4}.

In order to evaluate the significance of Theorem 1 for practice we will
demonstrate its transformed particular cases for 
$k=1,\ldots,6$ \cite{3}-\cite{17}

\vspace{1mm}
\begin{equation}
\label{a1}
J[\psi^{(1)}]_{T,t}
=\hbox{\vtop{\offinterlineskip\halign{
\hfil#\hfil\cr
{\rm l.i.m.}\cr
$\stackrel{}{{}_{p_1\to \infty}}$\cr
}} }\sum_{j_1=0}^{p_1}
C_{j_1}\zeta_{j_1}^{(i_1)},
\end{equation}

\vspace{4mm}
\begin{equation}
\label{leto5001}
J[\psi^{(2)}]_{T,t}
=\hbox{\vtop{\offinterlineskip\halign{
\hfil#\hfil\cr
{\rm l.i.m.}\cr
$\stackrel{}{{}_{p_1,p_2\to \infty}}$\cr
}} }\sum_{j_1=0}^{p_1}\sum_{j_2=0}^{p_2}
C_{j_2j_1}\Biggl(\zeta_{j_1}^{(i_1)}\zeta_{j_2}^{(i_2)}
-{\bf 1}_{\{i_1=i_2\ne 0\}}
{\bf 1}_{\{j_1=j_2\}}\Biggr),
\end{equation}

\vspace{8mm}
$$
J[\psi^{(3)}]_{T,t}=
\hbox{\vtop{\offinterlineskip\halign{
\hfil#\hfil\cr
{\rm l.i.m.}\cr
$\stackrel{}{{}_{p_1,\ldots,p_3\to \infty}}$\cr
}} }\sum_{j_1=0}^{p_1}\sum_{j_2=0}^{p_2}\sum_{j_3=0}^{p_3}
C_{j_3j_2j_1}\Biggl(
\zeta_{j_1}^{(i_1)}\zeta_{j_2}^{(i_2)}\zeta_{j_3}^{(i_3)}
-\Biggr.
$$

\vspace{1mm}
\begin{equation}
\label{leto5002}
\Biggl.-{\bf 1}_{\{i_1=i_2\ne 0\}}
{\bf 1}_{\{j_1=j_2\}}
\zeta_{j_3}^{(i_3)}
-{\bf 1}_{\{i_2=i_3\ne 0\}}
{\bf 1}_{\{j_2=j_3\}}
\zeta_{j_1}^{(i_1)}-
{\bf 1}_{\{i_1=i_3\ne 0\}}
{\bf 1}_{\{j_1=j_3\}}
\zeta_{j_2}^{(i_2)}\Biggr),
\end{equation}

\vspace{8mm}
$$
J[\psi^{(4)}]_{T,t}
=
\hbox{\vtop{\offinterlineskip\halign{
\hfil#\hfil\cr
{\rm l.i.m.}\cr
$\stackrel{}{{}_{p_1,\ldots,p_4\to \infty}}$\cr
}} }\sum_{j_1=0}^{p_1}\ldots\sum_{j_4=0}^{p_4}
C_{j_4\ldots j_1}\Biggl(
\prod_{l=1}^4\zeta_{j_l}^{(i_l)}
\Biggr.
-
$$
$$
-
{\bf 1}_{\{i_1=i_2\ne 0\}}
{\bf 1}_{\{j_1=j_2\}}
\zeta_{j_3}^{(i_3)}
\zeta_{j_4}^{(i_4)}
-
{\bf 1}_{\{i_1=i_3\ne 0\}}
{\bf 1}_{\{j_1=j_3\}}
\zeta_{j_2}^{(i_2)}
\zeta_{j_4}^{(i_4)}-
$$
$$
-
{\bf 1}_{\{i_1=i_4\ne 0\}}
{\bf 1}_{\{j_1=j_4\}}
\zeta_{j_2}^{(i_2)}
\zeta_{j_3}^{(i_3)}
-
{\bf 1}_{\{i_2=i_3\ne 0\}}
{\bf 1}_{\{j_2=j_3\}}
\zeta_{j_1}^{(i_1)}
\zeta_{j_4}^{(i_4)}-
$$
$$
-
{\bf 1}_{\{i_2=i_4\ne 0\}}
{\bf 1}_{\{j_2=j_4\}}
\zeta_{j_1}^{(i_1)}
\zeta_{j_3}^{(i_3)}
-
{\bf 1}_{\{i_3=i_4\ne 0\}}
{\bf 1}_{\{j_3=j_4\}}
\zeta_{j_1}^{(i_1)}
\zeta_{j_2}^{(i_2)}+
$$
$$
+
{\bf 1}_{\{i_1=i_2\ne 0\}}
{\bf 1}_{\{j_1=j_2\}}
{\bf 1}_{\{i_3=i_4\ne 0\}}
{\bf 1}_{\{j_3=j_4\}}
+
$$
$$
+
{\bf 1}_{\{i_1=i_3\ne 0\}}
{\bf 1}_{\{j_1=j_3\}}
{\bf 1}_{\{i_2=i_4\ne 0\}}
{\bf 1}_{\{j_2=j_4\}}+
$$
\begin{equation}
\label{leto5003}
+\Biggl.
{\bf 1}_{\{i_1=i_4\ne 0\}}
{\bf 1}_{\{j_1=j_4\}}
{\bf 1}_{\{i_2=i_3\ne 0\}}
{\bf 1}_{\{j_2=j_3\}}\Biggr),
\end{equation}

\vspace{8mm}

$$
J[\psi^{(5)}]_{T,t}
=\hbox{\vtop{\offinterlineskip\halign{
\hfil#\hfil\cr
{\rm l.i.m.}\cr
$\stackrel{}{{}_{p_1,\ldots,p_5\to \infty}}$\cr
}} }\sum_{j_1=0}^{p_1}\ldots\sum_{j_5=0}^{p_5}
C_{j_5\ldots j_1}\Biggl(
\prod_{l=1}^5\zeta_{j_l}^{(i_l)}
-\Biggr.
$$
$$
-
{\bf 1}_{\{i_1=i_2\ne 0\}}
{\bf 1}_{\{j_1=j_2\}}
\zeta_{j_3}^{(i_3)}
\zeta_{j_4}^{(i_4)}
\zeta_{j_5}^{(i_5)}-
{\bf 1}_{\{i_1=i_3\ne 0\}}
{\bf 1}_{\{j_1=j_3\}}
\zeta_{j_2}^{(i_2)}
\zeta_{j_4}^{(i_4)}
\zeta_{j_5}^{(i_5)}-
$$
$$
-
{\bf 1}_{\{i_1=i_4\ne 0\}}
{\bf 1}_{\{j_1=j_4\}}
\zeta_{j_2}^{(i_2)}
\zeta_{j_3}^{(i_3)}
\zeta_{j_5}^{(i_5)}-
{\bf 1}_{\{i_1=i_5\ne 0\}}
{\bf 1}_{\{j_1=j_5\}}
\zeta_{j_2}^{(i_2)}
\zeta_{j_3}^{(i_3)}
\zeta_{j_4}^{(i_4)}-
$$
$$
-
{\bf 1}_{\{i_2=i_3\ne 0\}}
{\bf 1}_{\{j_2=j_3\}}
\zeta_{j_1}^{(i_1)}
\zeta_{j_4}^{(i_4)}
\zeta_{j_5}^{(i_5)}-
{\bf 1}_{\{i_2=i_4\ne 0\}}
{\bf 1}_{\{j_2=j_4\}}
\zeta_{j_1}^{(i_1)}
\zeta_{j_3}^{(i_3)}
\zeta_{j_5}^{(i_5)}-
$$
$$
-
{\bf 1}_{\{i_2=i_5\ne 0\}}
{\bf 1}_{\{j_2=j_5\}}
\zeta_{j_1}^{(i_1)}
\zeta_{j_3}^{(i_3)}
\zeta_{j_4}^{(i_4)}
-{\bf 1}_{\{i_3=i_4\ne 0\}}
{\bf 1}_{\{j_3=j_4\}}
\zeta_{j_1}^{(i_1)}
\zeta_{j_2}^{(i_2)}
\zeta_{j_5}^{(i_5)}-
$$
$$
-
{\bf 1}_{\{i_3=i_5\ne 0\}}
{\bf 1}_{\{j_3=j_5\}}
\zeta_{j_1}^{(i_1)}
\zeta_{j_2}^{(i_2)}
\zeta_{j_4}^{(i_4)}
-{\bf 1}_{\{i_4=i_5\ne 0\}}
{\bf 1}_{\{j_4=j_5\}}
\zeta_{j_1}^{(i_1)}
\zeta_{j_2}^{(i_2)}
\zeta_{j_3}^{(i_3)}+
$$
$$
+
{\bf 1}_{\{i_1=i_2\ne 0\}}
{\bf 1}_{\{j_1=j_2\}}
{\bf 1}_{\{i_3=i_4\ne 0\}}
{\bf 1}_{\{j_3=j_4\}}\zeta_{j_5}^{(i_5)}+
{\bf 1}_{\{i_1=i_2\ne 0\}}
{\bf 1}_{\{j_1=j_2\}}
{\bf 1}_{\{i_3=i_5\ne 0\}}
{\bf 1}_{\{j_3=j_5\}}\zeta_{j_4}^{(i_4)}+
$$
$$
+
{\bf 1}_{\{i_1=i_2\ne 0\}}
{\bf 1}_{\{j_1=j_2\}}
{\bf 1}_{\{i_4=i_5\ne 0\}}
{\bf 1}_{\{j_4=j_5\}}\zeta_{j_3}^{(i_3)}+
{\bf 1}_{\{i_1=i_3\ne 0\}}
{\bf 1}_{\{j_1=j_3\}}
{\bf 1}_{\{i_2=i_4\ne 0\}}
{\bf 1}_{\{j_2=j_4\}}\zeta_{j_5}^{(i_5)}+
$$
$$
+
{\bf 1}_{\{i_1=i_3\ne 0\}}
{\bf 1}_{\{j_1=j_3\}}
{\bf 1}_{\{i_2=i_5\ne 0\}}
{\bf 1}_{\{j_2=j_5\}}\zeta_{j_4}^{(i_4)}+
{\bf 1}_{\{i_1=i_3\ne 0\}}
{\bf 1}_{\{j_1=j_3\}}
{\bf 1}_{\{i_4=i_5\ne 0\}}
{\bf 1}_{\{j_4=j_5\}}\zeta_{j_2}^{(i_2)}+
$$
$$
+
{\bf 1}_{\{i_1=i_4\ne 0\}}
{\bf 1}_{\{j_1=j_4\}}
{\bf 1}_{\{i_2=i_3\ne 0\}}
{\bf 1}_{\{j_2=j_3\}}\zeta_{j_5}^{(i_5)}+
{\bf 1}_{\{i_1=i_4\ne 0\}}
{\bf 1}_{\{j_1=j_4\}}
{\bf 1}_{\{i_2=i_5\ne 0\}}
{\bf 1}_{\{j_2=j_5\}}\zeta_{j_3}^{(i_3)}+
$$
$$
+
{\bf 1}_{\{i_1=i_4\ne 0\}}
{\bf 1}_{\{j_1=j_4\}}
{\bf 1}_{\{i_3=i_5\ne 0\}}
{\bf 1}_{\{j_3=j_5\}}\zeta_{j_2}^{(i_2)}+
{\bf 1}_{\{i_1=i_5\ne 0\}}
{\bf 1}_{\{j_1=j_5\}}
{\bf 1}_{\{i_2=i_3\ne 0\}}
{\bf 1}_{\{j_2=j_3\}}\zeta_{j_4}^{(i_4)}+
$$
$$
+
{\bf 1}_{\{i_1=i_5\ne 0\}}
{\bf 1}_{\{j_1=j_5\}}
{\bf 1}_{\{i_2=i_4\ne 0\}}
{\bf 1}_{\{j_2=j_4\}}\zeta_{j_3}^{(i_3)}+
{\bf 1}_{\{i_1=i_5\ne 0\}}
{\bf 1}_{\{j_1=j_5\}}
{\bf 1}_{\{i_3=i_4\ne 0\}}
{\bf 1}_{\{j_3=j_4\}}\zeta_{j_2}^{(i_2)}+
$$
$$
+
{\bf 1}_{\{i_2=i_3\ne 0\}}
{\bf 1}_{\{j_2=j_3\}}
{\bf 1}_{\{i_4=i_5\ne 0\}}
{\bf 1}_{\{j_4=j_5\}}\zeta_{j_1}^{(i_1)}+
{\bf 1}_{\{i_2=i_4\ne 0\}}
{\bf 1}_{\{j_2=j_4\}}
{\bf 1}_{\{i_3=i_5\ne 0\}}
{\bf 1}_{\{j_3=j_5\}}\zeta_{j_1}^{(i_1)}+
$$
\begin{equation}
\label{a5}
+\Biggl.
{\bf 1}_{\{i_2=i_5\ne 0\}}
{\bf 1}_{\{j_2=j_5\}}
{\bf 1}_{\{i_3=i_4\ne 0\}}
{\bf 1}_{\{j_3=j_4\}}\zeta_{j_1}^{(i_1)}\Biggr),
\end{equation}

\vspace{9mm}

$$
J[\psi^{(6)}]_{T,t}
=\hbox{\vtop{\offinterlineskip\halign{
\hfil#\hfil\cr
{\rm l.i.m.}\cr
$\stackrel{}{{}_{p_1,\ldots,p_6\to \infty}}$\cr
}} }\sum_{j_1=0}^{p_1}\ldots\sum_{j_6=0}^{p_6}
C_{j_6\ldots j_1}\Biggl(
\prod_{l=1}^6
\zeta_{j_l}^{(i_l)}
-\Biggr.
$$
$$
-
{\bf 1}_{\{i_1=i_6\ne 0\}}
{\bf 1}_{\{j_1=j_6\}}
\zeta_{j_2}^{(i_2)}
\zeta_{j_3}^{(i_3)}
\zeta_{j_4}^{(i_4)}
\zeta_{j_5}^{(i_5)}-
{\bf 1}_{\{i_2=i_6\ne 0\}}
{\bf 1}_{\{j_2=j_6\}}
\zeta_{j_1}^{(i_1)}
\zeta_{j_3}^{(i_3)}
\zeta_{j_4}^{(i_4)}
\zeta_{j_5}^{(i_5)}-
$$
$$
-
{\bf 1}_{\{i_3=i_6\ne 0\}}
{\bf 1}_{\{j_3=j_6\}}
\zeta_{j_1}^{(i_1)}
\zeta_{j_2}^{(i_2)}
\zeta_{j_4}^{(i_4)}
\zeta_{j_5}^{(i_5)}-
{\bf 1}_{\{i_4=i_6\ne 0\}}
{\bf 1}_{\{j_4=j_6\}}
\zeta_{j_1}^{(i_1)}
\zeta_{j_2}^{(i_2)}
\zeta_{j_3}^{(i_3)}
\zeta_{j_5}^{(i_5)}-
$$
$$
-
{\bf 1}_{\{i_5=i_6\ne 0\}}
{\bf 1}_{\{j_5=j_6\}}
\zeta_{j_1}^{(i_1)}
\zeta_{j_2}^{(i_2)}
\zeta_{j_3}^{(i_3)}
\zeta_{j_4}^{(i_4)}-
{\bf 1}_{\{i_1=i_2\ne 0\}}
{\bf 1}_{\{j_1=j_2\}}
\zeta_{j_3}^{(i_3)}
\zeta_{j_4}^{(i_4)}
\zeta_{j_5}^{(i_5)}
\zeta_{j_6}^{(i_6)}-
$$
$$
-
{\bf 1}_{\{i_1=i_3\ne 0\}}
{\bf 1}_{\{j_1=j_3\}}
\zeta_{j_2}^{(i_2)}
\zeta_{j_4}^{(i_4)}
\zeta_{j_5}^{(i_5)}
\zeta_{j_6}^{(i_6)}-
{\bf 1}_{\{i_1=i_4\ne 0\}}
{\bf 1}_{\{j_1=j_4\}}
\zeta_{j_2}^{(i_2)}
\zeta_{j_3}^{(i_3)}
\zeta_{j_5}^{(i_5)}
\zeta_{j_6}^{(i_6)}-
$$
$$
-
{\bf 1}_{\{i_1=i_5\ne 0\}}
{\bf 1}_{\{j_1=j_5\}}
\zeta_{j_2}^{(i_2)}
\zeta_{j_3}^{(i_3)}
\zeta_{j_4}^{(i_4)}
\zeta_{j_6}^{(i_6)}-
{\bf 1}_{\{i_2=i_3\ne 0\}}
{\bf 1}_{\{j_2=j_3\}}
\zeta_{j_1}^{(i_1)}
\zeta_{j_4}^{(i_4)}
\zeta_{j_5}^{(i_5)}
\zeta_{j_6}^{(i_6)}-
$$
$$
-
{\bf 1}_{\{i_2=i_4\ne 0\}}
{\bf 1}_{\{j_2=j_4\}}
\zeta_{j_1}^{(i_1)}
\zeta_{j_3}^{(i_3)}
\zeta_{j_5}^{(i_5)}
\zeta_{j_6}^{(i_6)}-
{\bf 1}_{\{i_2=i_5\ne 0\}}
{\bf 1}_{\{j_2=j_5\}}
\zeta_{j_1}^{(i_1)}
\zeta_{j_3}^{(i_3)}
\zeta_{j_4}^{(i_4)}
\zeta_{j_6}^{(i_6)}-
$$
$$
-
{\bf 1}_{\{i_3=i_4\ne 0\}}
{\bf 1}_{\{j_3=j_4\}}
\zeta_{j_1}^{(i_1)}
\zeta_{j_2}^{(i_2)}
\zeta_{j_5}^{(i_5)}
\zeta_{j_6}^{(i_6)}-
{\bf 1}_{\{i_3=i_5\ne 0\}}
{\bf 1}_{\{j_3=j_5\}}
\zeta_{j_1}^{(i_1)}
\zeta_{j_2}^{(i_2)}
\zeta_{j_4}^{(i_4)}
\zeta_{j_6}^{(i_6)}-
$$
$$
-
{\bf 1}_{\{i_4=i_5\ne 0\}}
{\bf 1}_{\{j_4=j_5\}}
\zeta_{j_1}^{(i_1)}
\zeta_{j_2}^{(i_2)}
\zeta_{j_3}^{(i_3)}
\zeta_{j_6}^{(i_6)}+
$$
$$
+
{\bf 1}_{\{i_1=i_2\ne 0\}}
{\bf 1}_{\{j_1=j_2\}}
{\bf 1}_{\{i_3=i_4\ne 0\}}
{\bf 1}_{\{j_3=j_4\}}
\zeta_{j_5}^{(i_5)}
\zeta_{j_6}^{(i_6)}+
{\bf 1}_{\{i_1=i_2\ne 0\}}
{\bf 1}_{\{j_1=j_2\}}
{\bf 1}_{\{i_3=i_5\ne 0\}}
{\bf 1}_{\{j_3=j_5\}}
\zeta_{j_4}^{(i_4)}
\zeta_{j_6}^{(i_6)}+
$$
$$
+
{\bf 1}_{\{i_1=i_2\ne 0\}}
{\bf 1}_{\{j_1=j_2\}}
{\bf 1}_{\{i_4=i_5\ne 0\}}
{\bf 1}_{\{j_4=j_5\}}
\zeta_{j_3}^{(i_3)}
\zeta_{j_6}^{(i_6)}
+
{\bf 1}_{\{i_1=i_3\ne 0\}}
{\bf 1}_{\{j_1=j_3\}}
{\bf 1}_{\{i_2=i_4\ne 0\}}
{\bf 1}_{\{j_2=j_4\}}
\zeta_{j_5}^{(i_5)}
\zeta_{j_6}^{(i_6)}+
$$
$$
+
{\bf 1}_{\{i_1=i_3\ne 0\}}
{\bf 1}_{\{j_1=j_3\}}
{\bf 1}_{\{i_2=i_5\ne 0\}}
{\bf 1}_{\{j_2=j_5\}}
\zeta_{j_4}^{(i_4)}
\zeta_{j_6}^{(i_6)}
+{\bf 1}_{\{i_1=i_3\ne 0\}}
{\bf 1}_{\{j_1=j_3\}}
{\bf 1}_{\{i_4=i_5\ne 0\}}
{\bf 1}_{\{j_4=j_5\}}
\zeta_{j_2}^{(i_2)}
\zeta_{j_6}^{(i_6)}+
$$
$$
+
{\bf 1}_{\{i_1=i_4\ne 0\}}
{\bf 1}_{\{j_1=j_4\}}
{\bf 1}_{\{i_2=i_3\ne 0\}}
{\bf 1}_{\{j_2=j_3\}}
\zeta_{j_5}^{(i_5)}
\zeta_{j_6}^{(i_6)}
+
{\bf 1}_{\{i_1=i_4\ne 0\}}
{\bf 1}_{\{j_1=j_4\}}
{\bf 1}_{\{i_2=i_5\ne 0\}}
{\bf 1}_{\{j_2=j_5\}}
\zeta_{j_3}^{(i_3)}
\zeta_{j_6}^{(i_6)}+
$$
$$
+
{\bf 1}_{\{i_1=i_4\ne 0\}}
{\bf 1}_{\{j_1=j_4\}}
{\bf 1}_{\{i_3=i_5\ne 0\}}
{\bf 1}_{\{j_3=j_5\}}
\zeta_{j_2}^{(i_2)}
\zeta_{j_6}^{(i_6)}
+
{\bf 1}_{\{i_1=i_5\ne 0\}}
{\bf 1}_{\{j_1=j_5\}}
{\bf 1}_{\{i_2=i_3\ne 0\}}
{\bf 1}_{\{j_2=j_3\}}
\zeta_{j_4}^{(i_4)}
\zeta_{j_6}^{(i_6)}+
$$
$$
+
{\bf 1}_{\{i_1=i_5\ne 0\}}
{\bf 1}_{\{j_1=j_5\}}
{\bf 1}_{\{i_2=i_4\ne 0\}}
{\bf 1}_{\{j_2=j_4\}}
\zeta_{j_3}^{(i_3)}
\zeta_{j_6}^{(i_6)}
+
{\bf 1}_{\{i_1=i_5\ne 0\}}
{\bf 1}_{\{j_1=j_5\}}
{\bf 1}_{\{i_3=i_4\ne 0\}}
{\bf 1}_{\{j_3=j_4\}}
\zeta_{j_2}^{(i_2)}
\zeta_{j_6}^{(i_6)}+
$$
$$
+
{\bf 1}_{\{i_2=i_3\ne 0\}}
{\bf 1}_{\{j_2=j_3\}}
{\bf 1}_{\{i_4=i_5\ne 0\}}
{\bf 1}_{\{j_4=j_5\}}
\zeta_{j_1}^{(i_1)}
\zeta_{j_6}^{(i_6)}
+
{\bf 1}_{\{i_2=i_4\ne 0\}}
{\bf 1}_{\{j_2=j_4\}}
{\bf 1}_{\{i_3=i_5\ne 0\}}
{\bf 1}_{\{j_3=j_5\}}
\zeta_{j_1}^{(i_1)}
\zeta_{j_6}^{(i_6)}+
$$
$$
+
{\bf 1}_{\{i_2=i_5\ne 0\}}
{\bf 1}_{\{j_2=j_5\}}
{\bf 1}_{\{i_3=i_4\ne 0\}}
{\bf 1}_{\{j_3=j_4\}}
\zeta_{j_1}^{(i_1)}
\zeta_{j_6}^{(i_6)}
+
{\bf 1}_{\{i_6=i_1\ne 0\}}
{\bf 1}_{\{j_6=j_1\}}
{\bf 1}_{\{i_3=i_4\ne 0\}}
{\bf 1}_{\{j_3=j_4\}}
\zeta_{j_2}^{(i_2)}
\zeta_{j_5}^{(i_5)}+
$$
$$
+
{\bf 1}_{\{i_6=i_1\ne 0\}}
{\bf 1}_{\{j_6=j_1\}}
{\bf 1}_{\{i_3=i_5\ne 0\}}
{\bf 1}_{\{j_3=j_5\}}
\zeta_{j_2}^{(i_2)}
\zeta_{j_4}^{(i_4)}
+
{\bf 1}_{\{i_6=i_1\ne 0\}}
{\bf 1}_{\{j_6=j_1\}}
{\bf 1}_{\{i_2=i_5\ne 0\}}
{\bf 1}_{\{j_2=j_5\}}
\zeta_{j_3}^{(i_3)}
\zeta_{j_4}^{(i_4)}+
$$
$$
+
{\bf 1}_{\{i_6=i_1\ne 0\}}
{\bf 1}_{\{j_6=j_1\}}
{\bf 1}_{\{i_2=i_4\ne 0\}}
{\bf 1}_{\{j_2=j_4\}}
\zeta_{j_3}^{(i_3)}
\zeta_{j_5}^{(i_5)}
+
{\bf 1}_{\{i_6=i_1\ne 0\}}
{\bf 1}_{\{j_6=j_1\}}
{\bf 1}_{\{i_4=i_5\ne 0\}}
{\bf 1}_{\{j_4=j_5\}}
\zeta_{j_2}^{(i_2)}
\zeta_{j_3}^{(i_3)}+
$$
$$
+
{\bf 1}_{\{i_6=i_1\ne 0\}}
{\bf 1}_{\{j_6=j_1\}}
{\bf 1}_{\{i_2=i_3\ne 0\}}
{\bf 1}_{\{j_2=j_3\}}
\zeta_{j_4}^{(i_4)}
\zeta_{j_5}^{(i_5)}
+
{\bf 1}_{\{i_6=i_2\ne 0\}}
{\bf 1}_{\{j_6=j_2\}}
{\bf 1}_{\{i_3=i_5\ne 0\}}
{\bf 1}_{\{j_3=j_5\}}
\zeta_{j_1}^{(i_1)}
\zeta_{j_4}^{(i_4)}+
$$
$$
+
{\bf 1}_{\{i_6=i_2\ne 0\}}
{\bf 1}_{\{j_6=j_2\}}
{\bf 1}_{\{i_4=i_5\ne 0\}}
{\bf 1}_{\{j_4=j_5\}}
\zeta_{j_1}^{(i_1)}
\zeta_{j_3}^{(i_3)}
+
{\bf 1}_{\{i_6=i_2\ne 0\}}
{\bf 1}_{\{j_6=j_2\}}
{\bf 1}_{\{i_3=i_4\ne 0\}}
{\bf 1}_{\{j_3=j_4\}}
\zeta_{j_1}^{(i_1)}
\zeta_{j_5}^{(i_5)}+
$$
$$
+
{\bf 1}_{\{i_6=i_2\ne 0\}}
{\bf 1}_{\{j_6=j_2\}}
{\bf 1}_{\{i_1=i_5\ne 0\}}
{\bf 1}_{\{j_1=j_5\}}
\zeta_{j_3}^{(i_3)}
\zeta_{j_4}^{(i_4)}
+
{\bf 1}_{\{i_6=i_2\ne 0\}}
{\bf 1}_{\{j_6=j_2\}}
{\bf 1}_{\{i_1=i_4\ne 0\}}
{\bf 1}_{\{j_1=j_4\}}
\zeta_{j_3}^{(i_3)}
\zeta_{j_5}^{(i_5)}+
$$
$$
+
{\bf 1}_{\{i_6=i_2\ne 0\}}
{\bf 1}_{\{j_6=j_2\}}
{\bf 1}_{\{i_1=i_3\ne 0\}}
{\bf 1}_{\{j_1=j_3\}}
\zeta_{j_4}^{(i_4)}
\zeta_{j_5}^{(i_5)}
+
{\bf 1}_{\{i_6=i_3\ne 0\}}
{\bf 1}_{\{j_6=j_3\}}
{\bf 1}_{\{i_2=i_5\ne 0\}}
{\bf 1}_{\{j_2=j_5\}}
\zeta_{j_1}^{(i_1)}
\zeta_{j_4}^{(i_4)}+
$$
$$
+
{\bf 1}_{\{i_6=i_3\ne 0\}}
{\bf 1}_{\{j_6=j_3\}}
{\bf 1}_{\{i_4=i_5\ne 0\}}
{\bf 1}_{\{j_4=j_5\}}
\zeta_{j_1}^{(i_1)}
\zeta_{j_2}^{(i_2)}
+
{\bf 1}_{\{i_6=i_3\ne 0\}}
{\bf 1}_{\{j_6=j_3\}}
{\bf 1}_{\{i_2=i_4\ne 0\}}
{\bf 1}_{\{j_2=j_4\}}
\zeta_{j_1}^{(i_1)}
\zeta_{j_5}^{(i_5)}+
$$
$$
+
{\bf 1}_{\{i_6=i_3\ne 0\}}
{\bf 1}_{\{j_6=j_3\}}
{\bf 1}_{\{i_1=i_5\ne 0\}}
{\bf 1}_{\{j_1=j_5\}}
\zeta_{j_2}^{(i_2)}
\zeta_{j_4}^{(i_4)}
+
{\bf 1}_{\{i_6=i_3\ne 0\}}
{\bf 1}_{\{j_6=j_3\}}
{\bf 1}_{\{i_1=i_4\ne 0\}}
{\bf 1}_{\{j_1=j_4\}}
\zeta_{j_2}^{(i_2)}
\zeta_{j_5}^{(i_5)}+
$$
$$
+
{\bf 1}_{\{i_6=i_3\ne 0\}}
{\bf 1}_{\{j_6=j_3\}}
{\bf 1}_{\{i_1=i_2\ne 0\}}
{\bf 1}_{\{j_1=j_2\}}
\zeta_{j_4}^{(i_4)}
\zeta_{j_5}^{(i_5)}
+
{\bf 1}_{\{i_6=i_4\ne 0\}}
{\bf 1}_{\{j_6=j_4\}}
{\bf 1}_{\{i_3=i_5\ne 0\}}
{\bf 1}_{\{j_3=j_5\}}
\zeta_{j_1}^{(i_1)}
\zeta_{j_2}^{(i_2)}+
$$
$$
+
{\bf 1}_{\{i_6=i_4\ne 0\}}
{\bf 1}_{\{j_6=j_4\}}
{\bf 1}_{\{i_2=i_5\ne 0\}}
{\bf 1}_{\{j_2=j_5\}}
\zeta_{j_1}^{(i_1)}
\zeta_{j_3}^{(i_3)}
+
{\bf 1}_{\{i_6=i_4\ne 0\}}
{\bf 1}_{\{j_6=j_4\}}
{\bf 1}_{\{i_2=i_3\ne 0\}}
{\bf 1}_{\{j_2=j_3\}}
\zeta_{j_1}^{(i_1)}
\zeta_{j_5}^{(i_5)}+
$$
$$
+
{\bf 1}_{\{i_6=i_4\ne 0\}}
{\bf 1}_{\{j_6=j_4\}}
{\bf 1}_{\{i_1=i_5\ne 0\}}
{\bf 1}_{\{j_1=j_5\}}
\zeta_{j_2}^{(i_2)}
\zeta_{j_3}^{(i_3)}
+
{\bf 1}_{\{i_6=i_4\ne 0\}}
{\bf 1}_{\{j_6=j_4\}}
{\bf 1}_{\{i_1=i_3\ne 0\}}
{\bf 1}_{\{j_1=j_3\}}
\zeta_{j_2}^{(i_2)}
\zeta_{j_5}^{(i_5)}+
$$
$$
+
{\bf 1}_{\{i_6=i_4\ne 0\}}
{\bf 1}_{\{j_6=j_4\}}
{\bf 1}_{\{i_1=i_2\ne 0\}}
{\bf 1}_{\{j_1=j_2\}}
\zeta_{j_3}^{(i_3)}
\zeta_{j_5}^{(i_5)}
+
{\bf 1}_{\{i_6=i_5\ne 0\}}
{\bf 1}_{\{j_6=j_5\}}
{\bf 1}_{\{i_3=i_4\ne 0\}}
{\bf 1}_{\{j_3=j_4\}}
\zeta_{j_1}^{(i_1)}
\zeta_{j_2}^{(i_2)}+
$$
$$
+
{\bf 1}_{\{i_6=i_5\ne 0\}}
{\bf 1}_{\{j_6=j_5\}}
{\bf 1}_{\{i_2=i_4\ne 0\}}
{\bf 1}_{\{j_2=j_4\}}
\zeta_{j_1}^{(i_1)}
\zeta_{j_3}^{(i_3)}
+
{\bf 1}_{\{i_6=i_5\ne 0\}}
{\bf 1}_{\{j_6=j_5\}}
{\bf 1}_{\{i_2=i_3\ne 0\}}
{\bf 1}_{\{j_2=j_3\}}
\zeta_{j_1}^{(i_1)}
\zeta_{j_4}^{(i_4)}+
$$
$$
+
{\bf 1}_{\{i_6=i_5\ne 0\}}
{\bf 1}_{\{j_6=j_5\}}
{\bf 1}_{\{i_1=i_4\ne 0\}}
{\bf 1}_{\{j_1=j_4\}}
\zeta_{j_2}^{(i_2)}
\zeta_{j_3}^{(i_3)}
+
{\bf 1}_{\{i_6=i_5\ne 0\}}
{\bf 1}_{\{j_6=j_5\}}
{\bf 1}_{\{i_1=i_3\ne 0\}}
{\bf 1}_{\{j_1=j_3\}}
\zeta_{j_2}^{(i_2)}
\zeta_{j_4}^{(i_4)}+
$$
$$
+
{\bf 1}_{\{i_6=i_5\ne 0\}}
{\bf 1}_{\{j_6=j_5\}}
{\bf 1}_{\{i_1=i_2\ne 0\}}
{\bf 1}_{\{j_1=j_2\}}
\zeta_{j_3}^{(i_3)}
\zeta_{j_4}^{(i_4)}-
$$
$$
-
{\bf 1}_{\{i_6=i_1\ne 0\}}
{\bf 1}_{\{j_6=j_1\}}
{\bf 1}_{\{i_2=i_5\ne 0\}}
{\bf 1}_{\{j_2=j_5\}}
{\bf 1}_{\{i_3=i_4\ne 0\}}
{\bf 1}_{\{j_3=j_4\}}-
$$
$$
-
{\bf 1}_{\{i_6=i_1\ne 0\}}
{\bf 1}_{\{j_6=j_1\}}
{\bf 1}_{\{i_2=i_4\ne 0\}}
{\bf 1}_{\{j_2=j_4\}}
{\bf 1}_{\{i_3=i_5\ne 0\}}
{\bf 1}_{\{j_3=j_5\}}-
$$
$$
-
{\bf 1}_{\{i_6=i_1\ne 0\}}
{\bf 1}_{\{j_6=j_1\}}
{\bf 1}_{\{i_2=i_3\ne 0\}}
{\bf 1}_{\{j_2=j_3\}}
{\bf 1}_{\{i_4=i_5\ne 0\}}
{\bf 1}_{\{j_4=j_5\}}-
$$
$$
-
{\bf 1}_{\{i_6=i_2\ne 0\}}
{\bf 1}_{\{j_6=j_2\}}
{\bf 1}_{\{i_1=i_5\ne 0\}}
{\bf 1}_{\{j_1=j_5\}}
{\bf 1}_{\{i_3=i_4\ne 0\}}
{\bf 1}_{\{j_3=j_4\}}-
$$
$$
-
{\bf 1}_{\{i_6=i_2\ne 0\}}
{\bf 1}_{\{j_6=j_2\}}
{\bf 1}_{\{i_1=i_4\ne 0\}}
{\bf 1}_{\{j_1=j_4\}}
{\bf 1}_{\{i_3=i_5\ne 0\}}
{\bf 1}_{\{j_3=j_5\}}-
$$
$$
-
{\bf 1}_{\{i_6=i_2\ne 0\}}
{\bf 1}_{\{j_6=j_2\}}
{\bf 1}_{\{i_1=i_3\ne 0\}}
{\bf 1}_{\{j_1=j_3\}}
{\bf 1}_{\{i_4=i_5\ne 0\}}
{\bf 1}_{\{j_4=j_5\}}-
$$
$$
-
{\bf 1}_{\{i_6=i_3\ne 0\}}
{\bf 1}_{\{j_6=j_3\}}
{\bf 1}_{\{i_1=i_5\ne 0\}}
{\bf 1}_{\{j_1=j_5\}}
{\bf 1}_{\{i_2=i_4\ne 0\}}
{\bf 1}_{\{j_2=j_4\}}-
$$
$$
-
{\bf 1}_{\{i_6=i_3\ne 0\}}
{\bf 1}_{\{j_6=j_3\}}
{\bf 1}_{\{i_1=i_4\ne 0\}}
{\bf 1}_{\{j_1=j_4\}}
{\bf 1}_{\{i_2=i_5\ne 0\}}
{\bf 1}_{\{j_2=j_5\}}-
$$
$$
-
{\bf 1}_{\{i_3=i_6\ne 0\}}
{\bf 1}_{\{j_3=j_6\}}
{\bf 1}_{\{i_1=i_2\ne 0\}}
{\bf 1}_{\{j_1=j_2\}}
{\bf 1}_{\{i_4=i_5\ne 0\}}
{\bf 1}_{\{j_4=j_5\}}-
$$
$$
-
{\bf 1}_{\{i_6=i_4\ne 0\}}
{\bf 1}_{\{j_6=j_4\}}
{\bf 1}_{\{i_1=i_5\ne 0\}}
{\bf 1}_{\{j_1=j_5\}}
{\bf 1}_{\{i_2=i_3\ne 0\}}
{\bf 1}_{\{j_2=j_3\}}-
$$
$$
-
{\bf 1}_{\{i_6=i_4\ne 0\}}
{\bf 1}_{\{j_6=j_4\}}
{\bf 1}_{\{i_1=i_3\ne 0\}}
{\bf 1}_{\{j_1=j_3\}}
{\bf 1}_{\{i_2=i_5\ne 0\}}
{\bf 1}_{\{j_2=j_5\}}-
$$
$$
-
{\bf 1}_{\{i_6=i_4\ne 0\}}
{\bf 1}_{\{j_6=j_4\}}
{\bf 1}_{\{i_1=i_2\ne 0\}}
{\bf 1}_{\{j_1=j_2\}}
{\bf 1}_{\{i_3=i_5\ne 0\}}
{\bf 1}_{\{j_3=j_5\}}-
$$
$$
-
{\bf 1}_{\{i_6=i_5\ne 0\}}
{\bf 1}_{\{j_6=j_5\}}
{\bf 1}_{\{i_1=i_4\ne 0\}}
{\bf 1}_{\{j_1=j_4\}}
{\bf 1}_{\{i_2=i_3\ne 0\}}
{\bf 1}_{\{j_2=j_3\}}-
$$
$$
-
{\bf 1}_{\{i_6=i_5\ne 0\}}
{\bf 1}_{\{j_6=j_5\}}
{\bf 1}_{\{i_1=i_2\ne 0\}}
{\bf 1}_{\{j_1=j_2\}}
{\bf 1}_{\{i_3=i_4\ne 0\}}
{\bf 1}_{\{j_3=j_4\}}-
$$
\begin{equation}
\label{a6}
\Biggl.-
{\bf 1}_{\{i_6=i_5\ne 0\}}
{\bf 1}_{\{j_6=j_5\}}
{\bf 1}_{\{i_1=i_3\ne 0\}}
{\bf 1}_{\{j_1=j_3\}}
{\bf 1}_{\{i_2=i_4\ne 0\}}
{\bf 1}_{\{j_2=j_4\}}\Biggr),
\end{equation}

\vspace{6mm}
\noindent
where ${\bf 1}_A$ is the indicator of the set $A$.

Thus, we obtain the following useful possibilities
of the method of generalized multiple Fourier series.

1. There is an explicit formula (see (\ref{ppppa})) for calculation 
of expansion coefficients 
of the iterated Ito stochastic integral with any
fixed multiplicity $k$. 

2. We have new possibilities for exact calculation of the mean-square 
approximation error for iterated Ito stochastic integrals
(see Theorem 3 below).

3. Since the used multiple Fourier series is a generalized in the sense
that it is built using various complete orthonormal
systems of functions in the space $L_2([t, T])$, we have new possibilities 
for approximation --- we can 
use not only the trigonometric functions as in \cite{1995}-\cite{2004}
but the Legendre polynomials.

4. As it turned out (see below), it is more convenient to work 
with Legendre polynomials for constructing approximations of iterated 
stochastic integrals. We can choose different numbers $q$ 
(see Sect.~4) for
approximations of different iterated Ito stochastic integrals.
This is impossible for approximations based on the approach
from \cite{1995}-\cite{2004}.
Approximations based on Legendre polynomials are much simpler 
than approximations based on trigonometric functions
(see (\ref{4002}), (\ref{4003}), 
(\ref{444}), (\ref{1970}) below).

5. The approach from \cite{1995}-\cite{2004}, \cite{KPS}-\cite{Zapad-9} 
leads to 
iterated series (iterated application of the operation
of limit transition) in contrast with multiple 
series from Theorem 1 (operation of limit transition is
implemented only once) starting at least from the 
second or third multiplicity of iterated stochastic integrals.
Multiple series are more convenient for approximation than the iterated ones, 
since partial sums of multiple series converge for any possible case of  
convergence to infinity of their upper limits of summation 
(let us denote them as $p_1,\ldots, p_k$). 
For example,
when $p_1=\ldots=p_k=p\to\infty$. 
For iterated series, the condition $p_1=\ldots=p_k=p\to\infty$ obviously 
does not guarantee the convergence of this series.
However, 
in \cite{1995}
(Sect.~5.8, pp.~202--204), \cite{KPS} (pp.~82-84),
\cite{Zapad-2} (pp.~438-439),  
\cite{Zapad-9} (pp.~263-264) 
the authors use (without rigorous proof)
the condition $p_1=p_2=p_3=p\to\infty$
within the frames of the mentioned approach \cite{1995}-\cite{2004},
\cite{KPS}-\cite{Zapad-9} 
based on the Karhunen--Loeve expansion of the Brownian bridge
process \cite{1988} together with the Wong--Zakai approximation
\cite{W-Z-1}-\cite{Watanabe} (see discussion in Sect.~8 for detail).

6. In a number of works of the author 
\cite{10a} (Chapter~2), \cite{12}, \cite{15a}, \cite{arxiv-4}. 
Theorem 1 has been adapted for the iterated  Stratonovich stochastic integrals
(\ref{str}) of multiplicities 1 to 8.

For further consideration, let us 
consider the generalization of formulas (\ref{a1})--(\ref{a6})                 
for the case of an arbitrary multiplicity $k$ $(k\in\mathbb{N})$ of 
the iterated Ito stochastic integral $J[\psi^{(k)}]_{T,t}$ defined by (\ref{ito}).
In order to do this, let us
introduce some notations. 
Consider the unordered
set $\{1, 2, \ldots, k\}$ 
and separate it into two parts:
the first part consists of $r$ unordered 
pairs (sequence order of these pairs is also unimportant) and the 
second one consists of the 
remaining $k-2r$ numbers.
So, we have

\begin{equation}
\label{leto5007}
(\{
\underbrace{\{g_1, g_2\}, \ldots, 
\{g_{2r-1}, g_{2r}\}}_{\small{\hbox{part 1}}}
\},
\{\underbrace{q_1, \ldots, q_{k-2r}}_{\small{\hbox{part 2}}}
\}),
\end{equation}

\vspace{4mm}
\noindent
where 

\vspace{-2mm}
$$
\{g_1, g_2, \ldots, 
g_{2r-1}, g_{2r}, q_1, \ldots, q_{k-2r}\}=\{1, 2, \ldots, k\},
$$

\vspace{4mm}
\noindent
braces   
mean an unordered 
set, and pa\-ren\-the\-ses mean an ordered set.

We will say that (\ref{leto5007}) is a partition 
and consider the sum with respect to all possible
partitions

\begin{equation}
\label{leto5008}
\sum_{\stackrel{(\{\{g_1, g_2\}, \ldots, 
\{g_{2r-1}, g_{2r}\}\}, \{q_1, \ldots, q_{k-2r}\})}
{{}_{\{g_1, g_2, \ldots, 
g_{2r-1}, g_{2r}, q_1, \ldots, q_{k-2r}\}=\{1, 2, \ldots, k\}}}}
a_{g_1 g_2, \ldots, 
g_{2r-1} g_{2r}, q_1 \ldots q_{k-2r}}.
\end{equation}

\vspace{4mm}

Below there are several examples of sums in the form (\ref{leto5008})

\vspace{2mm}
$$
\sum_{\stackrel{(\{g_1, g_2\})}{{}_{\{g_1, g_2\}=\{1, 2\}}}}
a_{g_1 g_2}=a_{12},
$$

\vspace{3mm}
$$
\sum_{\stackrel{(\{\{g_1, g_2\}, \{g_3, g_4\}\})}
{{}_{\{g_1, g_2, g_3, g_4\}=\{1, 2, 3, 4\}}}}
a_{g_1 g_2 g_3 g_4}=a_{1234} + a_{1324} + a_{2314},
$$

\vspace{3mm}
$$
\sum_{\stackrel{(\{g_1, g_2\}, \{q_1, q_{2}\})}
{{}_{\{g_1, g_2, q_1, q_{2}\}=\{1, 2, 3, 4\}}}}
a_{g_1 g_2, q_1 q_{2}}=
$$

$$
=a_{12,34}+a_{13,24}+a_{14,23}
+a_{23,14}+a_{24,13}+a_{34,12},
$$

\vspace{3mm}
$$
\sum_{\stackrel{(\{g_1, g_2\}, \{q_1, q_{2}, q_3\})}
{{}_{\{g_1, g_2, q_1, q_{2}, q_3\}=\{1, 2, 3, 4, 5\}}}}
a_{g_1 g_2, q_1 q_{2}q_3}
=
$$

$$
=a_{12,345}+a_{13,245}+a_{14,235}
+a_{15,234}+a_{23,145}+a_{24,135}+
$$
$$
+a_{25,134}+a_{34,125}+a_{35,124}+a_{45,123},
$$

\vspace{4mm}
$$
\sum_{\stackrel{(\{\{g_1, g_2\}, \{g_3, g_{4}\}\}, \{q_1\})}
{{}_{\{g_1, g_2, g_3, g_{4}, q_1\}=\{1, 2, 3, 4, 5\}}}}
a_{g_1 g_2, g_3 g_{4},q_1}
=
$$

$$
=
a_{12,34,5}+a_{13,24,5}+a_{14,23,5}+
a_{12,35,4}+a_{13,25,4}+a_{15,23,4}+
$$
$$
+a_{12,54,3}+a_{15,24,3}+a_{14,25,3}+a_{15,34,2}+a_{13,54,2}+a_{14,53,2}+
$$
$$
+
a_{52,34,1}+a_{53,24,1}+a_{54,23,1}.
$$

\vspace{5mm}

Now we can write (\ref{tyyy}) as

\vspace{1mm}

$$
J[\psi^{(k)}]_{T,t}=
\hbox{\vtop{\offinterlineskip\halign{
\hfil#\hfil\cr
{\rm l.i.m.}\cr
$\stackrel{}{{}_{p_1,\ldots,p_k\to \infty}}$\cr
}} }
\sum\limits_{j_1=0}^{p_1}\ldots
\sum\limits_{j_k=0}^{p_k}
C_{j_k\ldots j_1}\Biggl(
\prod_{l=1}^k\zeta_{j_l}^{(i_l)}+\sum\limits_{r=1}^{[k/2]}
(-1)^r \times
\Biggr.
$$

\vspace{3mm}
\begin{equation}
\label{leto6000hh}
\times
\sum_{\stackrel{(\{\{g_1, g_2\}, \ldots, 
\{g_{2r-1}, g_{2r}\}\}, \{q_1, \ldots, q_{k-2r}\})}
{{}_{\{g_1, g_2, \ldots, 
g_{2r-1}, g_{2r}, q_1, \ldots, q_{k-2r}\}=\{1, 2, \ldots, k\}}}}
\prod\limits_{s=1}^r
{\bf 1}_{\{i_{g_{{}_{2s-1}}}=~i_{g_{{}_{2s}}}\ne 0\}}
\Biggl.{\bf 1}_{\{j_{g_{{}_{2s-1}}}=~j_{g_{{}_{2s}}}\}}
\prod_{l=1}^{k-2r}\zeta_{j_{q_l}}^{(i_{q_l})}\Biggr),
\end{equation}

\vspace{5mm}
\noindent
where $[x]$ is an integer part of a real number $x;$
another notations are the same as in Theorem {\bf 1}.

\vspace{2mm}

In particular, from (\ref{leto6000hh}) for $k=5$ we obtain

\vspace{3mm}

$$
J[\psi^{(5)}]_{T,t}=
\hbox{\vtop{\offinterlineskip\halign{
\hfil#\hfil\cr
{\rm l.i.m.}\cr
$\stackrel{}{{}_{p_1,\ldots,p_5\to \infty}}$\cr
}} }\sum_{j_1=0}^{p_1}\ldots\sum_{j_5=0}^{p_5}
C_{j_5\ldots j_1}\Biggl(
\prod_{l=1}^5\zeta_{j_l}^{(i_l)}-\Biggr.
$$

\vspace{2mm}
$$
-
\sum\limits_{\stackrel{(\{g_1, g_2\}, \{q_1, q_{2}, q_3\})}
{{}_{\{g_1, g_2, q_{1}, q_{2}, q_3\}=\{1, 2, 3, 4, 5\}}}}
{\bf 1}_{\{i_{g_{{}_{1}}}=~i_{g_{{}_{2}}}\ne 0\}}
{\bf 1}_{\{j_{g_{{}_{1}}}=~j_{g_{{}_{2}}}\}}
\prod_{l=1}^{3}\zeta_{j_{q_l}}^{(i_{q_l})}+
$$

\vspace{2mm}
$$
+
\sum_{\stackrel{(\{\{g_1, g_2\}, 
\{g_{3}, g_{4}\}\}, \{q_1\})}
{{}_{\{g_1, g_2, g_{3}, g_{4}, q_1\}=\{1, 2, 3, 4, 5\}}}}
{\bf 1}_{\{i_{g_{{}_{1}}}=~i_{g_{{}_{2}}}\ne 0\}}
{\bf 1}_{\{j_{g_{{}_{1}}}=~j_{g_{{}_{2}}}\}}
\Biggl.{\bf 1}_{\{i_{g_{{}_{3}}}=~i_{g_{{}_{4}}}\ne 0\}}
{\bf 1}_{\{j_{g_{{}_{3}}}=~j_{g_{{}_{4}}}\}}
\zeta_{j_{q_1}}^{(i_{q_1})}\Biggr).
$$

\vspace{7mm}
\noindent
The last equality obviously agrees with
(\ref{a5}).

Let us consider the generalization of Theorem 1 for the case
of an arbitrary complete orthonormal systems  
of functions in the space $L_2([t,T])$ 
and $\psi_1(\tau),\ldots,\psi_k(\tau)\in L_2([t, T]).$

\vspace{2mm}

{\bf Theorem~2}\ \cite{10a} (Sect.~1.11), \cite{11} (Sect.~15).
{\it Suppose that
$\psi_1(\tau),\ldots,\psi_k(\tau)\in L_2([t, T])$ and
$\{\phi_j(x)\}_{j=0}^{\infty}$ is an arbitrary complete orthonormal system  
of functions in the space $L_2([t,T]).$
Then the following expansion

\vspace{1mm}
$$
J[\psi^{(k)}]_{T,t}=
\hbox{\vtop{\offinterlineskip\halign{
\hfil#\hfil\cr
{\rm l.i.m.}\cr
$\stackrel{}{{}_{p_1,\ldots,p_k\to \infty}}$\cr
}} }
\sum\limits_{j_1=0}^{p_1}\ldots
\sum\limits_{j_k=0}^{p_k}
C_{j_k\ldots j_1}\Biggl(
\prod_{l=1}^k\zeta_{j_l}^{(i_l)}+\sum\limits_{r=1}^{[k/2]}
(-1)^r \times
\Biggr.
$$

\vspace{2mm}
\begin{equation}
\label{leto6000}
\times
\sum_{\stackrel{(\{\{g_1, g_2\}, \ldots, 
\{g_{2r-1}, g_{2r}\}\}, \{q_1, \ldots, q_{k-2r}\})}
{{}_{\{g_1, g_2, \ldots, 
g_{2r-1}, g_{2r}, q_1, \ldots, q_{k-2r}\}=\{1, 2, \ldots, k\}}}}
\prod\limits_{s=1}^r
{\bf 1}_{\{i_{g_{{}_{2s-1}}}=~i_{g_{{}_{2s}}}\ne 0\}}
\Biggl.{\bf 1}_{\{j_{g_{{}_{2s-1}}}=~j_{g_{{}_{2s}}}\}}
\prod_{l=1}^{k-2r}\zeta_{j_{q_l}}^{(i_{q_l})}\Biggr)
\end{equation}

\vspace{6mm}
\noindent
con\-verg\-ing in the mean-square sense is valid,
where $[x]$ is an integer part of a real number $x;$
another notations are the same as in Theorem~{\rm 1}.}

\vspace{2mm}

As noted above, in a number of works of the author 
\cite{5}-\cite{10axx1}, \cite{12}, \cite{15a}
Theorem 1 has been adapted for the iterated  Stratonovich stochastic integrals
(\ref{str}) of multiplicities 1 to 8.
Let us first present some old results as the following theorem.

\vspace{2mm}

{\bf Theorem 3} \cite{5}-\cite{10axx1}, \cite{12}, \cite{15a}. 
{\it Suppose that 
$\{\phi_j(x)\}_{j=0}^{\infty}$ is a complete orthonormal system of 
Legendre polynomials or trigonometric functions in the space $L_2([t, T]).$
At the same time $\psi_2(\tau)$ is a continuously differentiable 
function on $[t, T]$ and $\psi_1(\tau), \psi_3(\tau)$ are twice 
continuously differentiable functions on $[t, T]$. Then

\begin{equation}
\label{a}
J^{*}[\psi^{(2)}]_{T,t}=
\hbox{\vtop{\offinterlineskip\halign{
\hfil#\hfil\cr
{\rm l.i.m.}\cr
$\stackrel{}{{}_{p_1,p_2\to \infty}}$\cr
}} }\sum_{j_1=0}^{p_1}\sum_{j_2=0}^{p_2}
C_{j_2j_1}\zeta_{j_1}^{(i_1)}\zeta_{j_2}^{(i_2)}\ \ \ (i_1,i_2=1,\ldots,m),
\end{equation}

\vspace{1mm}
\begin{equation}
\label{feto19000ab}
J^{*}[\psi^{(3)}]_{T,t}=
\hbox{\vtop{\offinterlineskip\halign{
\hfil#\hfil\cr
{\rm l.i.m.}\cr
$\stackrel{}{{}_{p_1,p_2,p_3\to \infty}}$\cr
}} }\sum_{j_1=0}^{p_1}\sum_{j_2=0}^{p_2}\sum_{j_3=0}^{p_3}
C_{j_3 j_2 j_1}\zeta_{j_1}^{(i_1)}\zeta_{j_2}^{(i_2)}\zeta_{j_3}^{(i_3)}\ \ \
(i_1,i_2,i_3=0, 1,\ldots,m),
\end{equation}

\vspace{1mm}
\begin{equation}
\label{feto19000a}
J^{*}[\psi^{(3)}]_{T,t}=
\hbox{\vtop{\offinterlineskip\halign{
\hfil#\hfil\cr
{\rm l.i.m.}\cr
$\stackrel{}{{}_{p\to \infty}}$\cr
}} }
\sum\limits_{j_1,j_2,j_3=0}^{p}
C_{j_3 j_2 j_1}\zeta_{j_1}^{(i_1)}\zeta_{j_2}^{(i_2)}\zeta_{j_3}^{(i_3)}\ \ \
(i_1,i_2,i_3=1,\ldots,m),
\end{equation}

\vspace{1mm}
\begin{equation}
\label{uu}
J^{*}[\psi^{(4)}]_{T,t}=
\hbox{\vtop{\offinterlineskip\halign{
\hfil#\hfil\cr
{\rm l.i.m.}\cr
$\stackrel{}{{}_{p\to \infty}}$\cr
}} }
\sum\limits_{j_1,j_2,j_3,j_4=0}^{p}
C_{j_4 j_3 j_2 j_1}\zeta_{j_1}^{(i_1)}
\zeta_{j_2}^{(i_2)}\zeta_{j_3}^{(i_3)}\zeta_{j_4}^{(i_4)}\ \ \
(i_1,i_2,i_3,i_4=0, 1,\ldots,m),
\end{equation}

\vspace{5mm}
\noindent
where $J^{*}[\psi^{(k)}]_{T,t}$ is defined by {\rm (\ref{str})} and
$\psi_l(\tau)\equiv 1$ $(l=1,\ldots,4)$ in {\rm (\ref{feto19000ab})}, 
{\rm (\ref{uu});} another notations are the same as in Theorems {\rm 1, 2.}
}

\vspace{2mm}

Recently, a new approach to the expansion and mean-square 
approximation of iterated Stratonovich stochastic integrals has been obtained
\cite{10a} (Sect.~2.10--2.16), \cite{12} (Sect.~13--19), 
\cite{15a} (Sect.~5--11), \cite{arxiv-4} (Sect.~7--13).
Let us formulate four theorems that were obtained using this approach.

\vspace{2mm}

{\bf Theorem 4}\ \cite{10a}, \cite{12}, \cite{15a}, \cite{arxiv-4}.\
{\it Suppose 
that $\{\phi_j(x)\}_{j=0}^{\infty}$ is a complete orthonormal system of 
Legendre polynomials or trigonometric functions in the space $L_2([t, T]).$
Furthermore, let $\psi_1(\tau), \psi_2(\tau),$ $\psi_3(\tau)$ are continuously dif\-ferentiable 
nonrandom functions on $[t, T].$ 
Then, for the 
iterated Stra\-to\-no\-vich stochastic integral of third multiplicity

$$
J^{*}[\psi^{(3)}]_{T,t}={\int\limits_t^{*}}^T\psi_3(t_3)
{\int\limits_t^{*}}^{t_3}\psi_2(t_2)
{\int\limits_t^{*}}^{t_2}\psi_1(t_1)
d{\bf w}_{t_1}^{(i_1)}
d{\bf w}_{t_2}^{(i_2)}d{\bf w}_{t_3}^{(i_3)}\ \ \ (i_1,i_2,i_3=0,1,\ldots,m)
$$

\vspace{4mm}
\noindent
the following 
relations

\vspace{-1mm}
\begin{equation}
\label{fin1}
J^{*}[\psi^{(3)}]_{T,t}
=\hbox{\vtop{\offinterlineskip\halign{
\hfil#\hfil\cr
{\rm l.i.m.}\cr
$\stackrel{}{{}_{p\to \infty}}$\cr
}} }
\sum\limits_{j_1, j_2, j_3=0}^{p}
C_{j_3 j_2 j_1}\zeta_{j_1}^{(i_1)}\zeta_{j_2}^{(i_2)}\zeta_{j_3}^{(i_3)},
\end{equation}

\vspace{3mm}
\begin{equation}
\label{fin2}
{\sf M}\left\{\left(
J^{*}[\psi^{(3)}]_{T,t}-
\sum\limits_{j_1, j_2, j_3=0}^{p}
C_{j_3 j_2 j_1}\zeta_{j_1}^{(i_1)}\zeta_{j_2}^{(i_2)}\zeta_{j_3}^{(i_3)}\right)^2\right\}
\le \frac{C}{p}
\end{equation}

\vspace{5mm}
\noindent
are fulfilled, where $i_1, i_2, i_3=0,1,\ldots,m$ in {\rm (\ref{fin1})} and 
$i_1, i_2, i_3=1,\ldots,m$ in {\rm (\ref{fin2})},
constant $C$ is independent of $p,$

$$
C_{j_3 j_2 j_1}=\int\limits_t^T\psi_3(t_3)\phi_{j_3}(t_3)
\int\limits_t^{t_3}\psi_2(t_2)\phi_{j_2}(t_2)
\int\limits_t^{t_2}\psi_1(t_1)\phi_{j_1}(t_1)dt_1dt_2dt_3
$$

\vspace{4mm}
\noindent
and
$$
\zeta_{j}^{(i)}=
\int\limits_t^T \phi_{j}(\tau) d{\bf f}_{\tau}^{(i)}
$$ 

\vspace{2mm}
\noindent
are independent standard Gaussian random variables for various 
$i$ or $j$ {\rm (}in the case when $i\ne 0${\rm );} 
another notations are the same as in Theorems~{\rm 1, 2}.}

\vspace{2mm}

{\bf Theorem 5}\ \cite{10a}, \cite{12}, \cite{15a}, \cite{arxiv-4}.\ {\it Let
$\{\phi_j(x)\}_{j=0}^{\infty}$ be a complete orthonormal system of 
Legendre polynomials or trigonometric functions in the space $L_2([t, T]).$
Furthermore, let $\psi_1(\tau), \ldots,$ $\psi_4(\tau)$ be continuously dif\-ferentiable 
nonrandom functions on $[t, T].$ 
Then, for the 
iterated Stra\-to\-no\-vich stochastic integral of fourth multiplicity

\begin{equation}
\label{fin0}
J^{*}[\psi^{(4)}]_{T,t}={\int\limits_t^{*}}^T\psi_4(t_4)
{\int\limits_t^{*}}^{t_4}\psi_3(t_3)
{\int\limits_t^{*}}^{t_3}\psi_2(t_2)
{\int\limits_t^{*}}^{t_2}\psi_1(t_1)
d{\bf w}_{t_1}^{(i_1)}
d{\bf w}_{t_2}^{(i_2)}d{\bf w}_{t_3}^{(i_3)}d{\bf w}_{t_4}^{(i_4)}
\end{equation}

\vspace{4mm}
\noindent
the following 
relations

\begin{equation}
\label{fin3}
J^{*}[\psi^{(4)}]_{T,t}
=\hbox{\vtop{\offinterlineskip\halign{
\hfil#\hfil\cr
{\rm l.i.m.}\cr
$\stackrel{}{{}_{p\to \infty}}$\cr
}} }
\sum\limits_{j_1, j_2, j_3,j_4=0}^{p}
C_{j_4j_3 j_2 j_1}\zeta_{j_1}^{(i_1)}\zeta_{j_2}^{(i_2)}\zeta_{j_3}^{(i_3)}\zeta_{j_4}^{(i_4)},
\end{equation}

\vspace{3mm}

\begin{equation}
\label{fin4}
{\sf M}\left\{\left(
J^{*}[\psi^{(4)}]_{T,t}-
\sum\limits_{j_1, j_2, j_3, j_4=0}^{p}
C_{j_4 j_3 j_2 j_1}\zeta_{j_1}^{(i_1)}\zeta_{j_2}^{(i_2)}\zeta_{j_3}^{(i_3)}
\zeta_{j_4}^{(i_4)}
\right)^2\right\}
\le \frac{C}{p^{1-\varepsilon}}
\end{equation}

\vspace{5mm}
\noindent
are fulfilled, where $i_1, \ldots , i_4=0,1,\ldots,m$ in {\rm (\ref{fin0}),} {\rm (\ref{fin3})} 
and $i_1, \ldots, i_4=1,\ldots,m$ in {\rm (\ref{fin4}),}
constant $C$ does not depend on $p,$
$\varepsilon$ is an arbitrary
small positive real number 
for the case of complete orthonormal system of 
Legendre polynomials in the space $L_2([t, T])$
and $\varepsilon=0$ for the case of
complete orthonormal system of 
trigonometric functions in the space $L_2([t, T]),$

$$
C_{j_4 j_3 j_2 j_1}=
$$

$$
=
\int\limits_t^T\psi_4(t_4)\phi_{j_4}(t_4)
\int\limits_t^{t_4}\psi_3(t_3)\phi_{j_3}(t_3)
\int\limits_t^{t_3}\psi_2(t_2)\phi_{j_2}(t_2)
\int\limits_t^{t_2}\psi_1(t_1)\phi_{j_1}(t_1)dt_1dt_2dt_3dt_4;
$$

\vspace{4mm}
\noindent
another notations are the same as in Theorem~{\rm 4}.}

\vspace{2mm}

{\bf Theorem 6}\ \cite{10a}, \cite{12}, \cite{15a}, \cite{arxiv-4}.\
{\it Assume 
that $\{\phi_j(x)\}_{j=0}^{\infty}$ is a complete orthonormal system of 
Legendre polynomials or trigonometric functions in the space $L_2([t, T])$
and $\psi_1(\tau), \ldots,$ $\psi_5(\tau)$ are continuously dif\-ferentiable 
nonrandom functions on $[t, T].$ 
Then, for the 
iterated Stra\-to\-no\-vich stochastic integral of fifth multiplicity

\begin{equation}
\label{fin7}
J^{*}[\psi^{(5)}]_{T,t}={\int\limits_t^{*}}^T\psi_5(t_5)
\ldots
{\int\limits_t^{*}}^{t_2}\psi_1(t_1)
d{\bf w}_{t_1}^{(i_1)}
\ldots d{\bf w}_{t_5}^{(i_5)}
\end{equation}

\vspace{4mm}
\noindent
the following 
relations

\begin{equation}
\label{fin8}
J^{*}[\psi^{(5)}]_{T,t}
=\hbox{\vtop{\offinterlineskip\halign{
\hfil#\hfil\cr
{\rm l.i.m.}\cr
$\stackrel{}{{}_{p\to \infty}}$\cr
}} }
\sum\limits_{j_1,\ldots,j_5=0}^{p}
C_{j_5 \ldots j_1}\zeta_{j_1}^{(i_1)}\ldots \zeta_{j_5}^{(i_5)},
\end{equation}

\vspace{3mm}

\begin{equation}
\label{fin9}
{\sf M}\left\{\left(
J^{*}[\psi^{(5)}]_{T,t}-
\sum\limits_{j_1, \ldots, j_5=0}^{p}
C_{j_5 \ldots j_1}\zeta_{j_1}^{(i_1)}\ldots
\zeta_{j_5}^{(i_5)}
\right)^2\right\}
\le \frac{C}{p^{1-\varepsilon}}
\end{equation}

\vspace{5mm}
\noindent
are fulfilled, where $i_1, \ldots , i_5=0,1,\ldots,m$ in {\rm (\ref{fin7}),} {\rm (\ref{fin8})} 
and $i_1, \ldots, i_5=1,\ldots,m$ in {\rm (\ref{fin9}),}
constant $C$ is independent of $p,$
$\varepsilon$ is an arbitrary
small positive real number 
for the case of complete orthonormal system of 
Legendre polynomials in the space $L_2([t, T])$
and $\varepsilon=0$ for the case of
complete orthonormal system of 
trigonometric functions in the space $L_2([t, T]),$

$$
C_{j_5 \ldots j_1}=
\int\limits_t^T\psi_5(t_5)\phi_{j_5}(t_5)\ldots
\int\limits_t^{t_2}\psi_1(t_1)\phi_{j_1}(t_1)dt_1\ldots dt_5;
$$

\vspace{3mm}
\noindent
another notations are the same as in Theorems~{\rm 4, 5}.}

\vspace{2mm}

{\bf Theorem 7}\ \cite{10a}, \cite{12}, \cite{15a}, \cite{arxiv-4}.\
{\it Suppose that 
$\{\phi_j(x)\}_{j=0}^{\infty}$ is a complete orthonormal system of 
Legendre polynomials or trigonometric functions in the space $L_2([t, T]).$
Then, for the 
iterated Stratonovich stochastic integral of sixth multiplicity

\begin{equation}
\label{after10001qu1}
J_{T,t}^{*(i_1\ldots i_6)}={\int\limits_t^{*}}^T
\ldots
{\int\limits_t^{*}}^{t_2}
d{\bf w}_{t_1}^{(i_1)}
\ldots d{\bf w}_{t_6}^{(i_6)}
\end{equation}

\vspace{3mm}
\noindent
the following 
expansion 

\vspace{-1mm}
$$
J_{T,t}^{*(i_1\ldots i_6)}
=\hbox{\vtop{\offinterlineskip\halign{
\hfil#\hfil\cr
{\rm l.i.m.}\cr
$\stackrel{}{{}_{p\to \infty}}$\cr
}} }
\sum\limits_{j_1, \ldots, j_6=0}^{p}
C_{j_6 \ldots j_1}\zeta_{j_1}^{(i_1)}\ldots
\zeta_{j_6}^{(i_6)}
$$

\vspace{4mm}
\noindent
that converges in the mean-square sense is valid, where
$i_1, \ldots, i_6=0, 1,\ldots,m,$

$$
C_{j_6 \ldots j_1}=
\int\limits_t^T\phi_{j_6}(t_6)\ldots
\int\limits_t^{t_2}\phi_{j_1}(t_1)dt_1\ldots dt_6;
$$

\vspace{3mm}
\noindent
another notations are the same as in Theorems~{\rm 4--6}.}

\vspace{2mm}

The results of Theorems~3--7 were developed in 
\cite{10a} (Chapter~2), \cite{12}, \cite{15a}, \cite{arxiv-4}. 
In particular, analogues of Theorem~7 for iterated Stratonovich stochastic
integrals of multiplicities 7 and 8 were obtained in \cite{10a} (Sect.~2.36, 2.37).
In addition, the variants of Theorems 3--7
were obtained
for the case when $\{\phi_j(x)\}_{j=0}^{\infty}$ is an arbitrary complete orthonormal system
of functions in $L_2([t, T])$ \cite{10a} (Sect.~2.1.4, 2.23, 2.24, 2.31--2.34),
\cite{12}, \cite{15a}, \cite{arxiv-4}.

As we mentioned above,
Theorems 1 and 2 allow us to accurately calculate 
the mean-square 
approximation error for iterated Ito stochastic integrals
(see Theorem 8 below).

Assume that $J[\psi^{(k)}]_{T,t}^{p_1 \ldots p_k}$ is the approximation 
of (\ref{ito}), which is
the expression on the right-hand side of (\ref{leto6000}) before passing to the limit

\vspace{1mm}
$$
J[\psi^{(k)}]_{T,t}^{p_1 \ldots p_k}=
\sum\limits_{j_1=0}^{p_1}\ldots
\sum\limits_{j_k=0}^{p_k}
C_{j_k\ldots j_1}\Biggl(
\prod_{l=1}^k\zeta_{j_l}^{(i_l)}+\sum\limits_{r=1}^{[k/2]}
(-1)^r \times
\Biggr.
$$

\vspace{4mm}
$$
\times
\sum_{\stackrel{(\{\{g_1, g_2\}, \ldots, 
\{g_{2r-1}, g_{2r}\}\}, \{q_1, \ldots, q_{k-2r}\})}
{{}_{\{g_1, g_2, \ldots, 
g_{2r-1}, g_{2r}, q_1, \ldots, q_{k-2r}\}=\{1, 2, \ldots, k\}}}}
\prod\limits_{s=1}^r
{\bf 1}_{\{i_{g_{{}_{2s-1}}}=~i_{g_{{}_{2s}}}\ne 0\}}
\Biggl.{\bf 1}_{\{j_{g_{{}_{2s-1}}}=~j_{g_{{}_{2s}}}\}}
\prod_{l=1}^{k-2r}\zeta_{j_{q_l}}^{(i_{q_l})}\Biggr),
$$

\vspace{6mm}
\noindent
where $[x]$ is an integer part of a real number $x;$
another notations are the same as in Theorems~{\rm 1, 2}.

Let us denote

$$
E_k^{p_1,\ldots,p_k}\stackrel{{\rm def}}
{=}{\sf M}\left\{\left(J[\psi^{(k)}]_{T,t}-
J[\psi^{(k)}]_{T,t}^{p_1,\ldots,p_k}\right)^2\right\},
$$

\vspace{3mm}
$$
E_k^p\stackrel{{\rm def}}{=}E_k^{p_1,\ldots,p_k}\ \ \hbox{if}\ \ 
p_1=\ldots=p_k=p,
$$

\vspace{2mm}
$$
I_k\stackrel{{\rm def}}{=}\left\Vert K\right\Vert^2_{L_2([t,T]^k)}=\int\limits_{[t,T]^k}
K^2(t_1,\ldots,t_k)dt_1\ldots dt_k.
$$

\vspace{4mm}

In \cite{7}-\cite{11}, \cite{15b} it was shown that

\begin{equation}
\label{star00011}
E_k^{p_1,\ldots,p_k}\le k!\left(I_k-\sum_{j_1=0}^{p_1}\ldots
\sum_{j_k=0}^{p_k}C^2_{j_k\ldots j_1}\right)
\end{equation}

\vspace{4mm}
\noindent
if $i_1,\ldots,i_k=1,\ldots,m$ and $0<T-t<\infty$ or 
$i_1,\ldots,i_k=0, 1,\ldots,m$ and $0<T-t<1.$

Moreover, in \cite{10a} (Sect.~1.1.9, 1.11, 1.12), \cite{11} (Sect.~6, 15, 16)
the following estimate is obtained

\vspace{3mm}
$$
{\sf M}\left\{\left(J[\psi^{(k)}]_{T,t}-
J[\psi^{(k)}]_{T,t}^{p_1,\ldots,p_k}\right)^{2n}\right\}\le
$$

\vspace{2mm}
\begin{equation}
\label{99999}
\le
(k!)^{n} (2n-1)^{nk}
\left(I_k-\sum_{j_1=0}^{p_1}\ldots
\sum_{j_k=0}^{p_k}C^2_{j_k\ldots j_1}\right)^n,
\end{equation}

\vspace{5mm}
\noindent
where $n\in \mathbb{N}$.

The value $E_k^{p}$
can be calculated exactly.

\vspace{2mm}

{\bf Theorem 8} \cite{10a} (Sect.~1.12), \cite{15b} (Sect.~6).
{\it Suppose that $\{\phi_j(x)\}_{j=0}^{\infty}$ 
is an arbitrary complete orthonormal system  
of functions in the space $L_2([t,T])$ and
$\psi_1(\tau),\ldots,\psi_k(\tau)\in L_2([t, T]),$  $i_1,\ldots, i_k=1,\ldots,m$.
Then

$$
E_k^p=I_k- 
$$

\begin{equation}
\label{tttr11}
-\sum_{j_1,\ldots, j_k=0}^{p}
C_{j_k\ldots j_1}
{\sf M}\left\{J[\psi^{(k)}]_{T,t}
\sum\limits_{(j_1,\ldots,j_k)}
\int\limits_t^T \phi_{j_k}(t_k)
\ldots
\int\limits_t^{t_{2}}\phi_{j_{1}}(t_{1})
d{\bf f}_{t_1}^{(i_1)}\ldots
d{\bf f}_{t_k}^{(i_k)}\right\},
\end{equation}

\vspace{6mm}
\noindent
where $i_1,\ldots,i_k = 1,\ldots,m;$
expression 

\vspace{-1mm}
$$
\sum\limits_{(j_1,\ldots,j_k)}
$$ 

\vspace{3mm}
\noindent
means the sum with respect to all
possible permutations 
$(j_1,\ldots,j_k)$. At the same time if 
$j_r$ swapped with $j_q$ in the permutation $(j_1,\ldots,j_k),$
then $i_r$ swapped with $i_q$ in the permutation
$(i_1,\ldots,i_k);$
another notations are the same as in Theorems {\rm 1, 2.}
}

Note that 

$$
{\sf M}\left\{J[\psi^{(k)}]_{T,t}
\int\limits_t^T \phi_{j_k}(t_k)
\ldots
\int\limits_t^{t_{2}}\phi_{j_{1}}(t_{1})
d{\bf f}_{t_1}^{(i_1)}\ldots
d{\bf f}_{t_k}^{(i_k)}\right\}=C_{j_k\ldots j_1}.
$$

\vspace{5mm}

Then from Theorem 8 for pairwise different $i_1,\ldots,i_k$ 
and for $i_1=\ldots=i_k$
we obtain

$$
E_k^p= I_k- \sum_{j_1,\ldots,j_k=0}^{p}
C_{j_k\ldots j_1}^2,
$$

\vspace{2mm}
$$ 
E_k^p= I_k - \sum_{j_1,\ldots,j_k=0}^{p}
C_{j_k\ldots j_1}\Biggl(\sum\limits_{(j_1,\ldots,j_k)}
C_{j_k\ldots j_1}\Biggr).
$$

\vspace{6mm}

Consider some examples of the application of Theorem 8
$(i_1,\ldots,i_5=1,\ldots,m)$

\vspace{1mm}
$$
E_2^p
=I_2
-\sum_{j_1,j_2=0}^p
C_{j_2j_1}^2-
\sum_{j_1,j_2=0}^p
C_{j_2j_1}C_{j_1j_2}\ \ \ (i_1=i_2),
$$

\vspace{4mm}
$$
E_3^p=I_3
-\sum_{j_3,j_2,j_1=0}^p C_{j_3j_2j_1}^2-
\sum_{j_3,j_2,j_1=0}^p C_{j_3j_1j_2}C_{j_3j_2j_1}\ \ \ (i_1=i_2\ne i_3),
$$

\vspace{4mm}
$$
E_3^p=I_3-
\sum_{j_3,j_2,j_1=0}^p C_{j_3j_2j_1}^2-
\sum_{j_3,j_2,j_1=0}^p C_{j_2j_3j_1}C_{j_3j_2j_1}\ \ \ (i_1\ne i_2=i_3),
$$

\vspace{4mm}
$$
E_3^p=I_3
-\sum_{j_3,j_2,j_1=0}^p C_{j_3j_2j_1}^2-
\sum_{j_3,j_2,j_1=0}^p C_{j_3j_2j_1}C_{j_1j_2j_3}\ \ \ (i_1=i_3\ne i_2),
$$

\vspace{4mm}
\begin{equation}
\label{usl1}
E^p_4 = I_4 - \sum_{j_1,\ldots,j_4=0}^{p}
C_{j_4\ldots j_1}\Biggl(\sum\limits_{(j_1,j_2)}
C_{j_4\ldots j_1}\Biggr)\ \ \ (i_1=i_2\ne i_3, i_4;\ i_3\ne i_4),
\end{equation}

\vspace{4mm}
\begin{equation}
\label{usl2}
E^p_4 = I_4 - \sum_{j_1,\ldots,j_4=0}^{p}
C_{j_4\ldots j_1}\Biggl(\sum\limits_{(j_1,j_3)}
C_{j_4\ldots j_1}\Biggr)\ \ \ (i_1=i_3\ne i_2, i_4;\ i_2\ne i_4),
\end{equation}

\vspace{4mm}
\begin{equation}
\label{usl3}
E^p_4 = I_4 - \sum_{j_1,\ldots,j_4=0}^{p}
C_{j_4\ldots j_1}\Biggl(\sum\limits_{(j_1,j_4)}
C_{j_4\ldots j_1}\Biggr)\ \ \ (i_1=i_4\ne i_2, i_3;\ i_2\ne i_3),
\end{equation}

\vspace{4mm}
\begin{equation}
\label{usl4}
E^p_4 = I_4 - \sum_{j_1,\ldots,j_4=0}^{p}
C_{j_4\ldots j_1}\Biggl(\sum\limits_{(j_2,j_3)}
C_{j_4\ldots j_1}\Biggr)\ \ \ (i_2=i_3\ne i_1, i_4;\ i_1\ne i_4),
\end{equation}

\vspace{4mm}
\begin{equation}
\label{usl5}
E^p_4 = I_4 - \sum_{j_1,\ldots,j_4=0}^{p}
C_{j_4\ldots j_1}\Biggl(\sum\limits_{(j_2,j_4)}
C_{j_4\ldots j_1}\Biggr)\ \ \ (i_2=i_4\ne i_1, i_3;\ i_1\ne i_3),
\end{equation}

\vspace{4mm}
\begin{equation}
\label{usl6}
E^p_4 = I_4 - \sum_{j_1,\ldots,j_4=0}^{p}
C_{j_4\ldots j_1}\Biggl(\sum\limits_{(j_3,j_4)}
C_{j_4\ldots j_1}\Biggr)\ \ \ (i_3=i_4\ne i_1, i_2;\ i_1\ne i_2),
\end{equation}

\vspace{4mm}
\begin{equation}
\label{usl7}
E_4^p = I_4 -
\sum_{j_1,\ldots,j_4=0}^{p}
C_{j_4\ldots j_1}\Biggl(\sum\limits_{(j_1,j_2,j_3)}
C_{j_4\ldots j_1}\Biggr)\ \ \ (i_1=i_2=i_3\ne i_4),
\end{equation}

\vspace{4mm}
\begin{equation}
\label{usl8}
E_4^p = I_4 -
 \sum_{j_1,\ldots,j_4=0}^{p}
C_{j_4\ldots j_1}\Biggl(\sum\limits_{(j_2,j_3,j_4)}
C_{j_4\ldots j_1}\Biggr)\ \ \ (i_2=i_3=i_4\ne i_1),
\end{equation}

\vspace{4mm}
\begin{equation}
\label{usl9}
E_4^p = I_4 -
 \sum_{j_1,\ldots,j_4=0}^{p}
C_{j_4\ldots j_1}\Biggl(\sum\limits_{(j_1,j_2,j_4)}
C_{j_4\ldots j_1}\Biggr)\ \ \ (i_1=i_2=i_4\ne i_3),
\end{equation}

\vspace{4mm}
\begin{equation}
\label{usl10}
E_4^p = I_4 -
 \sum_{j_1,\ldots,j_4=0}^{p}
C_{j_4\ldots j_1}\Biggl(\sum\limits_{(j_1,j_3,j_4)}
C_{j_4\ldots j_1}\Biggr)\ \ \ (i_1=i_3=i_4\ne i_2),
\end{equation}

\vspace{4mm}
\begin{equation}
\label{usl11}
E^p_4 = I_4 - \sum_{j_1,\ldots,j_4=0}^{p}
C_{j_4\ldots j_1}\Biggl(\sum\limits_{(j_1,j_2)}\Biggl(
\sum\limits_{(j_3,j_4)}
C_{j_4\ldots j_1}\Biggr)\Biggr)\ \ \ (i_1=i_2\ne i_3=i_4),
\end{equation}

\vspace{4mm}
\begin{equation}
\label{usl12}
E^p_4 = I_4 - \sum_{j_1,\ldots,j_4=0}^{p}
C_{j_4\ldots j_1}\Biggl(\sum\limits_{(j_1,j_3)}\Biggl(
\sum\limits_{(j_2,j_4)}
C_{j_4\ldots j_1}\Biggr)\Biggr)\ \ \ (i_1=i_3\ne i_2=i_4),
\end{equation}

\vspace{4mm}
\begin{equation}
\label{usl13}
E^p_4 = I_4 - \sum_{j_1,\ldots,j_4=0}^{p}
C_{j_4\ldots j_1}\Biggl(\sum\limits_{(j_1,j_4)}\Biggl(
\sum\limits_{(j_2,j_3)}
C_{j_4\ldots j_1}\Biggr)\Biggr)\ \ \ (i_1=i_4\ne i_2=i_3),
\end{equation}

\vspace{5mm}
$$
E_5^p = I_5 - \sum_{j_1,\ldots,j_5=0}^{p}
C_{j_5\ldots j_1}\Biggl(\sum\limits_{(j_2,j_4)}\Biggl(
\sum\limits_{(j_3,j_5)}
C_{j_5\ldots j_1}\Biggr)\Biggr)\ \ \ (i_1\ne i_2=i_4\ne i_3=i_5\ne i_1),
$$

\vspace{5mm}
$$
E^p_5 = I_5 - \sum_{j_1,\ldots,j_5=0}^{p}
C_{j_5\ldots j_1}\Biggl(\sum\limits_{(j_4,j_5)}\Biggl(
\sum\limits_{(j_1,j_2,j_3)}
C_{j_5\ldots j_1}\Biggr)\Biggr)\ \ \ (i_1=i_2=i_3\ne i_4=i_5),
$$

\vspace{5mm}
$$
E^p_5 = I_5 - \sum_{j_1,\ldots,j_5=0}^{p}
C_{j_5\ldots j_1}\Biggl(\sum\limits_{(j_1,j_3,j_4,j_5)}
C_{j_5\ldots j_1}\Biggr)\ \ \ (i_1=i_3=i_4=i_5\ne i_2).
$$

\vspace{5mm}

\section{Expansions and Approximations of Specific Iterated 
Ito and Stratonovich Stochastic 
Integrals of Multiplicities 1 to 6 Using Legendre Polynomials}

\vspace{5mm}

In this section, we provide considerable practical material (based on 
Theorems 1--7) on expansions and approximations of 
iterated Ito and 
Stratonovich stochastic integrals of the following form

\vspace{-1mm}
\begin{equation}
\label{k1000}
I_{(l_1\ldots l_k)T,t}^{(i_1\ldots i_k)}
=\int\limits_t^T(t-t_k)^{l_k} \ldots \int\limits_t^{t_{2}}
(t-t_1)^{l_1} d{\bf f}_{t_1}^{(i_1)}\ldots
d{\bf f}_{t_k}^{(i_k)},
\end{equation}

\begin{equation}
\label{k1001}
I_{(l_1\ldots l_k)T,t}^{*(i_1\ldots i_k)}
=\int\limits_t^{*T}(t-t_k)^{l_k} \ldots \int\limits_t^{*t_{2}}
(t-t_1)^{l_1} d{\bf f}_{t_1}^{(i_1)}\ldots
d{\bf f}_{t_k}^{(i_k)},
\end{equation}

\vspace{4mm}
\noindent
where $i_1,\ldots, i_k=1,\dots,m,$\ \  $l_1,\ldots,l_k=0, 1,\ldots$

The complete orthonormal system of Legendre polynomials in the 
space $L_2([t,T])$ looks as follows

\begin{equation}
\label{4009}
\phi_j(x)=\sqrt{\frac{2j+1}{T-t}}P_j\left(\left(
x-\frac{T+t}{2}\right)\frac{2}{T-t}\right),\ \ \ j=0, 1, 2,\ldots,
\end{equation}

\vspace{3mm}
\noindent
where $P_j(x)$ is the Legendre polynomial. It is well known
that the polynomials $P_j(x)$ can be represented, for example, in the form

\vspace{-1mm}
$$
P_j(x)=\frac{1}{2^j j!} \frac{d^j}{dx^j}\left(x^2-1\right)^j.
$$

\vspace{3mm}

Consider some well known properties of the polynomials $P_j(x)$

\vspace{-1mm}
$$
P_j(1)=1,\ \ \ P_{j+1}(-1)=-P_j(-1),\ \ \ j=0, 1, 2,\ldots,
$$

\vspace{1mm}
$$
\frac{dP_{j+1}(x)}{dx}-\frac{dP_{j-1}(x)}{dx}=(2j+1)P_j(x),
$$

\vspace{1mm}
$$
xP_{j}(x)=\frac{(j+1)P_{j+1}(x)+jP_{j-1}(x)}{2j+1},\ \ \ j=1, 2,\ldots,
$$

\vspace{1mm}
$$
\int\limits_{-1}^1 x^k P_j(x)dx=0,\ \ \ k=0, 1, 2,\ldots,j-1,
$$

\vspace{1mm}
$$
\int\limits_{-1}^1 P_k(x)P_j(x)dx=\left\{
\begin{matrix}
0\ \  &\hbox{if}\ \  k\ne j\cr\cr
2/(2j+1)\ \  &\hbox{if}\ \  k=j
\end{matrix}\ \ ,\right.
$$

\vspace{2mm}
$$
P_n(x)P_m(x)=\sum_{k=0}^m K_{m,n,k} P_{n+m-2k}(x),
$$

\vspace{4mm}
\noindent
where

\vspace{-1mm}
$$
K_{m,n,k}=\frac{a_{m-k}a_k a_{n-k}}{a_{m+n-k}}\cdot
\frac{2n+2m-4k+1}{2n+2m-2k+1},\ \ \ a_k=\frac{(2k-1)!!}{k!},\ \ \ m\le n.
$$

\vspace{4mm}

Using the above properties,
system of 
functions (\ref{4009}) and Theorems 1--7, 
we obtain the following expansions of iterated Ito 
and Stratonovich stochastic integrals (\ref{k1000}) and
(\ref{k1001}) based on multiple Fourier--Legendre series

\begin{equation}
\label{4001}
I_{(0)T,t}^{(i_1)}=\sqrt{T-t}\zeta_0^{(i_1)},
\end{equation}

\vspace{1mm}

\begin{equation}
\label{4002}
I_{(1)T,t}^{(i_1)}=-\frac{(T-t)^{3/2}}{2}\left(\zeta_0^{(i_1)}+
\frac{1}{\sqrt{3}}\zeta_1^{(i_1)}\right),
\end{equation}

\vspace{3mm}

\begin{equation}
\label{4003}
I_{(2)T,t}^{(i_1)}=\frac{(T-t)^{5/2}}{3}\left(\zeta_0^{(i_1)}+
\frac{\sqrt{3}}{2}\zeta_1^{(i_1)}+
\frac{1}{2\sqrt{5}}\zeta_2^{(i_1)}\right),
\end{equation}

\vspace{3mm}

\begin{equation}
I_{(00)T,t}^{*(i_1 i_2)}=
\frac{T-t}{2}\left(\zeta_0^{(i_1)}\zeta_0^{(i_2)}+\sum_{i=1}^{\infty}
\frac{1}{\sqrt{4i^2-1}}\left(
\zeta_{i-1}^{(i_1)}\zeta_{i}^{(i_2)}-
\zeta_i^{(i_1)}\zeta_{i-1}^{(i_2)}\right)\right),
\label{4004}
\end{equation}

\vspace{3mm}

$$
I_{(00)T,t}^{(i_1 i_2)}=
\frac{T-t}{2}\left(\zeta_0^{(i_1)}\zeta_0^{(i_2)}+\sum_{i=1}^{\infty}
\frac{1}{\sqrt{4i^2-1}}\left(
\zeta_{i-1}^{(i_1)}\zeta_{i}^{(i_2)}-
\zeta_i^{(i_1)}\zeta_{i-1}^{(i_2)}\right) - {\bf 1}_{\{i_1=i_2\}}\right),
$$

\vspace{6mm}

$$
I_{(01)T,t}^{*(i_1 i_2)}=-\frac{T-t}{2}I_{(00)T,t}^{*(i_1 i_2)}
-\frac{(T-t)^2}{4}\Biggl(
\frac{1}{\sqrt{3}}\zeta_0^{(i_1)}\zeta_1^{(i_2)}+\Biggr.
$$

\vspace{1mm}
\begin{equation}
\label{ud111}
+\Biggl.\sum_{i=0}^{\infty}\Biggl(
\frac{(i+2)\zeta_i^{(i_1)}\zeta_{i+2}^{(i_2)}
-(i+1)\zeta_{i+2}^{(i_1)}\zeta_{i}^{(i_2)}}
{\sqrt{(2i+1)(2i+5)}(2i+3)}-
\frac{\zeta_i^{(i_1)}\zeta_{i}^{(i_2)}}{(2i-1)(2i+3)}\Biggr)\Biggr),
\end{equation}

\vspace{6mm}

$$
I_{(10)T,t}^{*(i_1 i_2)}=-\frac{T-t}{2}I_{(00)T,t}^{*(i_1 i_2)}
-\frac{(T-t)^2}{4}\Biggl(
\frac{1}{\sqrt{3}}\zeta_0^{(i_2)}\zeta_1^{(i_1)}+\Biggr.
$$

\vspace{1mm}
\begin{equation}
\label{4006}
+\Biggl.\sum_{i=0}^{\infty}\Biggl(
\frac{(i+1)\zeta_{i+2}^{(i_2)}\zeta_{i}^{(i_1)}
-(i+2)\zeta_{i}^{(i_2)}\zeta_{i+2}^{(i_1)}}
{\sqrt{(2i+1)(2i+5)}(2i+3)}+
\frac{\zeta_i^{(i_1)}\zeta_{i}^{(i_2)}}{(2i-1)(2i+3)}\Biggr)\Biggr),
\end{equation}

\vspace{6mm}
\noindent 
or
$$
I_{(01)T,t}^{*(i_1i_2)}
=\hbox{\vtop{\offinterlineskip\halign{
\hfil#\hfil\cr
{\rm l.i.m.}\cr
$\stackrel{}{{}_{p\to \infty}}$\cr
}} }
\sum_{j_1,j_2=0}^{p}
C_{j_2j_1}^{01}
\zeta_{j_1}^{(i_1)}\zeta_{j_2}^{(i_2)},
$$

\vspace{2mm}
$$
I_{(10)T,t}^{*(i_1i_2)}
=\hbox{\vtop{\offinterlineskip\halign{
\hfil#\hfil\cr
{\rm l.i.m.}\cr
$\stackrel{}{{}_{p\to \infty}}$\cr
}} }
\sum_{j_1,j_2=0}^{p}
C_{j_2j_1}^{10}
\zeta_{j_1}^{(i_1)}\zeta_{j_2}^{(i_2)},
$$

\vspace{3mm}
\noindent
where

\vspace{-2mm}
$$
C_{j_2j_1}^{01}
=\frac{\sqrt{(2j_1+1)(2j_2+1)}}{8}(T-t)^{2}\bar
C_{j_2j_1}^{01},
$$

\vspace{1mm}
$$
C_{j_2j_1}^{10}
=\frac{\sqrt{(2j_1+1)(2j_2+1)}}{8}(T-t)^{2}\bar
C_{j_2j_1}^{10},
$$

\vspace{1mm}
$$
\bar C_{j_2j_1}^{01}=-\int\limits_{-1}^{1}(1+y)P_{j_2}(y)
\int\limits_{-1}^{y}
P_{j_1}(x)dx dy,
$$

\vspace{1mm}
$$
\bar C_{j_2j_1}^{10}=-\int\limits_{-1}^{1}P_{j_2}(y)
\int\limits_{-1}^{y}
(1+x)P_{j_1}(x)dx dy;
$$

\vspace{6mm}

$$
I_{(10)T,t}^{(i_1 i_2)}=
I_{(10)T,t}^{*(i_1 i_2)}+
\frac{1}{4}{\bf 1}_{\{i_1=i_2\}}(T-t)^2\ \ \ {\rm w.\ p.\ 1},
$$

\vspace{3mm}
$$
I_{(01)T,t}^{(i_1 i_2)}=
I_{(01)T,t}^{*(i_1 i_2)}+
\frac{1}{4}{\bf 1}_{\{i_1=i_2\}}(T-t)^2\ \ \ {\rm w.\ p.\ 1},
$$

\vspace{8mm}

$$
I_{(01)T,t}^{(i_1 i_2)}=
-\frac{T-t}{2}
I_{(00)T,t}^{(i_1 i_2)}
-\frac{(T-t)^2}{4}\Biggl(
\frac{1}{\sqrt{3}}\zeta_0^{(i_1)}\zeta_1^{(i_2)}+\Biggr.
$$

\vspace{2mm}
$$
+\Biggl.\sum_{i=0}^{\infty}\Biggl(
\frac{(i+2)\zeta_i^{(i_1)}\zeta_{i+2}^{(i_2)}
-(i+1)\zeta_{i+2}^{(i_1)}\zeta_{i}^{(i_2)}}
{\sqrt{(2i+1)(2i+5)}(2i+3)}-
\frac{\zeta_i^{(i_1)}\zeta_{i}^{(i_2)}}{(2i-1)(2i+3)}\Biggr)\Biggr),
$$

\vspace{8mm}

$$
I_{(10)T,t}^{(i_1 i_2)}=
-\frac{T-t}{2}I_{(00)T,t}^{(i_1 i_2)}
-\frac{(T-t)^2}{4}\Biggl(
\frac{1}{\sqrt{3}}\zeta_0^{(i_2)}\zeta_1^{(i_1)}+\Biggr.
$$

\vspace{2mm}
$$
+\Biggl.\sum_{i=0}^{\infty}\Biggl(
\frac{(i+1)\zeta_{i+2}^{(i_2)}\zeta_{i}^{(i_1)}
-(i+2)\zeta_{i}^{(i_2)}\zeta_{i+2}^{(i_1)}}
{\sqrt{(2i+1)(2i+5)}(2i+3)}+
\frac{\zeta_i^{(i_1)}\zeta_{i}^{(i_2)}}{(2i-1)(2i+3)}\Biggr)\Biggr),
$$

\vspace{8mm}
\noindent
or
$$
I_{(01)T,t}^{(i_1 i_2)}=
\hbox{\vtop{\offinterlineskip\halign{
\hfil#\hfil\cr
{\rm l.i.m.}\cr
$\stackrel{}{{}_{p\to \infty}}$\cr
}} }
\sum_{j_1,j_2=0}^{p}
C_{j_2j_1}^{01}\Biggl(\zeta_{j_1}^{(i_1)}\zeta_{j_2}^{(i_2)}
-{\bf 1}_{\{i_1=i_2\}}
{\bf 1}_{\{j_1=j_2\}}\Biggr),
$$

\vspace{5mm}
$$
I_{(10)T,t}^{(i_1 i_2)}=
\hbox{\vtop{\offinterlineskip\halign{
\hfil#\hfil\cr
{\rm l.i.m.}\cr
$\stackrel{}{{}_{p\to \infty}}$\cr
}} }
\sum_{j_1,j_2=0}^{p}
C_{j_2j_1}^{10}\Biggl(\zeta_{j_1}^{(i_1)}\zeta_{j_2}^{(i_2)}
-{\bf 1}_{\{i_1=i_2\}}
{\bf 1}_{\{j_1=j_2\}}\Biggr),
$$

\vspace{6mm}

\begin{equation}
\label{good1}
I_{(000)T,t}^{*(i_1 i_2 i_3)}=
\hbox{\vtop{\offinterlineskip\halign{
\hfil#\hfil\cr
{\rm l.i.m.}\cr
$\stackrel{}{{}_{p\to \infty}}$\cr
}} }
\sum\limits_{j_1, j_2, j_3=0}^{p}
C_{j_3 j_2 j_1}\zeta_{j_1}^{(i_1)}\zeta_{j_2}^{(i_2)}\zeta_{j_3}^{(i_3)},
\end{equation}

\vspace{6mm}

$$
I_{(000)T,t}^{(i_1i_2i_3)}
=\hbox{\vtop{\offinterlineskip\halign{
\hfil#\hfil\cr
{\rm l.i.m.}\cr
$\stackrel{}{{}_{p\to \infty}}$\cr
}} }
\sum_{j_1,j_2,j_3=0}^{p}
C_{j_3j_2j_1}
\Biggl(
\zeta_{j_1}^{(i_1)}\zeta_{j_2}^{(i_2)}\zeta_{j_3}^{(i_3)}
-{\bf 1}_{\{i_1=i_2\}}
{\bf 1}_{\{j_1=j_2\}}
\zeta_{j_3}^{(i_3)}-
\Biggr.
$$
\begin{equation}
\label{zzz1}
\Biggl.
-{\bf 1}_{\{i_2=i_3\}}
{\bf 1}_{\{j_2=j_3\}}
\zeta_{j_1}^{(i_1)}-
{\bf 1}_{\{i_1=i_3\}}
{\bf 1}_{\{j_1=j_3\}}
\zeta_{j_2}^{(i_2)}\Biggr),
\end{equation}

\vspace{6mm}

$$
I_{(000)T,t}^{(i_1 i_1 i_1)}
=\frac{1}{6}(T-t)^{3/2}\biggl(
\left(\zeta_0^{(i_1)}\right)^3-3
\zeta_0^{(i_1)}\biggr)\ \ \ \hbox{w.\ p.\ 1},
$$

\vspace{3mm}

\begin{equation}
\label{dest4}
I_{(000)T,t}^{*(i_1 i_1 i_1)}
=\frac{1}{6}(T-t)^{3/2}
\left(\zeta_0^{(i_1)}\right)^3\ \ \ \hbox{w.\ p.\ 1},
\end{equation}

\vspace{3mm}
\noindent
where

\vspace{-2mm}
\begin{equation}
\label{zzz2}
C_{j_3j_2j_1}
=\frac{\sqrt{(2j_1+1)(2j_2+1)(2j_3+1)}}{8}(T-t)^{3/2}\bar
C_{j_3j_2j_1},
\end{equation}

\begin{equation}
\label{zzz3xx}
\bar C_{j_3j_2j_1}=\int\limits_{-1}^{1}P_{j_3}(z)
\int\limits_{-1}^{z}P_{j_2}(y)
\int\limits_{-1}^{y}
P_{j_1}(x)dx dy dz;
\end{equation}

\vspace{5mm}

$$
I_{(000)T,t}^{(i_1 i_2 i_3)}=I_{(000)T,t}^{*(i_1 i_2 i_3)}
+{\bf 1}_{\{i_1=i_2\ne 0\}}\frac{1}{2}I_{(1)T,t}^{(i_3)}-
$$

\vspace{1mm}
$$
-
{\bf 1}_{\{i_2=i_3\ne 0\}}\frac{1}{2}\left((T-t)
I_{(0)T,t}^{(i_1)}+I_{(1)T,t}^{(i_1)}\right)\ \ \ \hbox{w.\ p.\ 1},
$$

\vspace{6mm}

$$
I_{(02)T,t}^{*(i_1 i_2)}=-\frac{(T-t)^2}{4}I_{(00)T,t}^{*(i_1 i_2)}
-(T-t)I_{(01)T,t}^{*(i_1 i_2)}+
\frac{(T-t)^3}{8}\Biggl[
\frac{2}{3\sqrt{5}}\zeta_2^{(i_2)}\zeta_0^{(i_1)}+\Biggr.
$$

\vspace{1mm}
$$
+\frac{1}{3}\zeta_0^{(i_1)}\zeta_0^{(i_2)}+
\sum_{i=0}^{\infty}\Biggl(
\frac{(i+2)(i+3)\zeta_{i+3}^{(i_2)}\zeta_{i}^{(i_1)}
-(i+1)(i+2)\zeta_{i}^{(i_2)}\zeta_{i+3}^{(i_1)}}
{\sqrt{(2i+1)(2i+7)}(2i+3)(2i+5)}+
\Biggr.
$$

\vspace{1mm}
\begin{equation}
\label{leto1}
\Biggl.\Biggl.+
\frac{(i^2+i-3)\zeta_{i+1}^{(i_2)}\zeta_{i}^{(i_1)}
-(i^2+3i-1)\zeta_{i}^{(i_2)}\zeta_{i+1}^{(i_1)}}
{\sqrt{(2i+1)(2i+3)}(2i-1)(2i+5)}\Biggr)\Biggr],
\end{equation}

\vspace{6mm}

$$
I_{(20)T,t}^{*(i_1 i_2)}=-\frac{(T-t)^2}{4}I_{(00)T,t}^{*(i_1 i_2)}
-(T-t)I_{(10)T,t}^{*(i_1 i_2)}+
\frac{(T-t)^3}{8}\Biggl[
\frac{2}{3\sqrt{5}}\zeta_0^{(i_2)}\zeta_2^{(i_1)}+\Biggr.
$$

\vspace{1mm}
$$
+\frac{1}{3}\zeta_0^{(i_1)}\zeta_0^{(i_2)}+
\sum_{i=0}^{\infty}\Biggl(
\frac{(i+1)(i+2)\zeta_{i+3}^{(i_2)}\zeta_{i}^{(i_1)}
-(i+2)(i+3)\zeta_{i}^{(i_2)}\zeta_{i+3}^{(i_1)}}
{\sqrt{(2i+1)(2i+7)}(2i+3)(2i+5)}+
\Biggr.
$$

\vspace{1mm}
\begin{equation}
\label{leto2}
\Biggl.\Biggl.+
\frac{(i^2+3i-1)\zeta_{i+1}^{(i_2)}\zeta_{i}^{(i_1)}
-(i^2+i-3)\zeta_{i}^{(i_2)}\zeta_{i+1}^{(i_1)}}
{\sqrt{(2i+1)(2i+3)}(2i-1)(2i+5)}\Biggr)\Biggr],
\end{equation}

\vspace{6mm}

$$
I_{(11)T,t}^{*(i_1 i_2)}=-\frac{(T-t)^2}{4}I_{(00)T,t}^{*(i_1 i_2)}-
\frac{(T-t)}{2}\left(I_{(10)T,t}^{*(i_1 i_2)}+
I_{(01)T,t}^{*(i_1 i_2)}\right)+
$$

\vspace{1mm}
$$
+\frac{(T-t)^3}{8}\Biggl[
\frac{1}{3}\zeta_1^{(i_1)}\zeta_1^{(i_2)}+
\sum_{i=0}^{\infty}\Biggl(
\frac{(i+1)(i+3)\left(\zeta_{i+3}^{(i_2)}\zeta_{i}^{(i_1)}
-\zeta_{i}^{(i_2)}\zeta_{i+3}^{(i_1)}\right)}
{\sqrt{(2i+1)(2i+7)}(2i+3)(2i+5)}+
\Biggr.
$$

\vspace{1mm}
\begin{equation}
\label{leto3}
\Biggl.\Biggl.+
\frac{(i+1)^2\left(\zeta_{i+1}^{(i_2)}\zeta_{i}^{(i_1)}
-\zeta_{i}^{(i_2)}\zeta_{i+1}^{(i_1)}\right)}
{\sqrt{(2i+1)(2i+3)}(2i-1)(2i+5)}\Biggr)\Biggr],
\end{equation}

\vspace{8mm}
or
$$
I_{(02)T,t}^{*(i_1 i_2)}=
\hbox{\vtop{\offinterlineskip\halign{
\hfil#\hfil\cr
{\rm l.i.m.}\cr
$\stackrel{}{{}_{p\to \infty}}$\cr
}} }
\sum_{j_1,j_2=0}^{p}
C_{j_2j_1}^{02}\zeta_{j_1}^{(i_1)}\zeta_{j_2}^{(i_2)},
$$

\vspace{3mm}

$$
I_{(20)T,t}^{*(i_1 i_2)}=
\hbox{\vtop{\offinterlineskip\halign{
\hfil#\hfil\cr
{\rm l.i.m.}\cr
$\stackrel{}{{}_{p\to \infty}}$\cr
}} }
\sum_{j_1,j_2=0}^{p}
C_{j_2j_1}^{20}\zeta_{j_1}^{(i_1)}\zeta_{j_2}^{(i_2)},
$$

\vspace{3mm}

$$
I_{(11)T,t}^{*(i_1 i_2)}=
\hbox{\vtop{\offinterlineskip\halign{
\hfil#\hfil\cr
{\rm l.i.m.}\cr
$\stackrel{}{{}_{p\to \infty}}$\cr
}} }
\sum_{j_1,j_2=0}^{p}
C_{j_2j_1}^{11}\zeta_{j_1}^{(i_1)}\zeta_{j_2}^{(i_2)},
$$

\vspace{5mm}
\noindent
where

\vspace{-2mm}
$$
C_{j_2j_1}^{02}=
\frac{\sqrt{(2j_1+1)(2j_2+1)}}{16}(T-t)^{3}\bar
C_{j_2j_1}^{02},
$$

\vspace{3mm}
$$
C_{j_2j_1}^{20}=
\frac{\sqrt{(2j_1+1)(2j_2+1)}}{16}(T-t)^{3}\bar
C_{j_2j_1}^{20},
$$

\vspace{3mm}

$$
C_{j_2j_1}^{11}=
\frac{\sqrt{(2j_1+1)(2j_2+1)}}{16}(T-t)^{3}\bar
C_{j_2j_1}^{11}, 
$$

\vspace{3mm}
$$
\bar C_{j_2j_1}^{02}=
\int\limits_{-1}^{1}P_{j_2}(y)(y+1)^2
\int\limits_{-1}^{y}
P_{j_1}(x)dx dy,
$$

\vspace{3mm}
$$
\bar C_{j_2j_1}^{20}=
\int\limits_{-1}^{1}P_{j_2}(y)
\int\limits_{-1}^{y}
P_{j_1}(x)(x+1)^2 dx dy,
$$

\vspace{3mm}
$$
\bar C_{j_2j_1}^{11}=
\int\limits_{-1}^{1}P_{j_2}(y)(y+1)
\int\limits_{-1}^{y}
P_{j_1}(x)(x+1)dx dy;
$$

\vspace{8mm}

$$
I_{(11)T,t}^{*(i_1 i_1)}=\frac{1}{2}\left(I_{(1)T,t}^{(i_1)}
\right)^2\ \ \ \hbox{w.\ p.\ 1,}
$$

\begin{equation}
\label{seg1}
I_{(02)T,t}^{(i_1 i_2)}=
I_{(02)T,t}^{*(i_1 i_2)}-
\frac{1}{6}{\bf 1}_{\{i_1=i_2\}}(T-t)^3\ \ \ \hbox{w.\ p.\ 1},
\end{equation}

\begin{equation}
\label{seg2}
I_{(20)T,t}^{(i_1 i_2)}=
I_{(20)T,t}^{*(i_1 i_2)}-
\frac{1}{6}{\bf 1}_{\{i_1=i_2\}}(T-t)^3\ \ \ \hbox{w.\ p.\ 1},
\end{equation}

\begin{equation}
\label{uh11}
I_{(11)T,t}^{(i_1 i_2)}=
I_{(11)T,t}^{*(i_1 i_2)}-
\frac{1}{6}{\bf 1}_{\{i_1=i_2\}}(T-t)^3\ \ \ \hbox{w.\ p.\ 1},
\end{equation}

\vspace{6mm}

$$
I_{(02)T,t}^{(i_1 i_2)}
=-\frac{(T-t)^2}{4}I_{(00)T,t}^{(i_1 i_2)}
-(T-t) I_{(01)T,t}^{(i_1 i_2)}+
\frac{(T-t)^3}{8}\Biggl[
\frac{2}{3\sqrt{5}}\zeta_2^{(i_2)}\zeta_0^{(i_1)}+\Biggr.
$$

\vspace{1mm}
$$
+\frac{1}{3}\zeta_0^{(i_1)}\zeta_0^{(i_2)}+
\sum_{i=0}^{\infty}\Biggl(
\frac{(i+2)(i+3)\zeta_{i+3}^{(i_2)}\zeta_{i}^{(i_1)}
-(i+1)(i+2)\zeta_{i}^{(i_2)}\zeta_{i+3}^{(i_1)}}
{\sqrt{(2i+1)(2i+7)}(2i+3)(2i+5)}+
\Biggr.
$$

\vspace{1mm}
$$
\Biggl.\Biggl.+
\frac{(i^2+i-3)\zeta_{i+1}^{(i_2)}\zeta_{i}^{(i_1)}
-(i^2+3i-1)\zeta_{i}^{(i_2)}\zeta_{i+1}^{(i_1)}}
{\sqrt{(2i+1)(2i+3)}(2i-1)(2i+5)}\Biggr)\Biggr] - 
$$

\vspace{1mm}
\begin{equation}
\label{seak1}
-\frac{1}{24}{\bf 1}_{\{i_1=i_2\}}{(T-t)^3},
\end{equation}

\vspace{6mm}

$$
I_{(20)T,t}^{(i_1 i_2)}=-\frac{(T-t)^2}{4}
I_{(00)T,t}^{(i_1 i_2)}
-(T-t) I_{(10)T,t}^{(i_1 i_2)}+
\frac{(T-t)^3}{8}\Biggl[
\frac{2}{3\sqrt{5}}\zeta_0^{(i_2)}\zeta_2^{(i_1)}+\Biggr.
$$
$$
+\frac{1}{3}\zeta_0^{(i_1)}\zeta_0^{(i_2)}+
\sum_{i=0}^{\infty}\Biggl(
\frac{(i+1)(i+2)\zeta_{i+3}^{(i_2)}\zeta_{i}^{(i_1)}
-(i+2)(i+3)\zeta_{i}^{(i_2)}\zeta_{i+3}^{(i_1)}}
{\sqrt{(2i+1)(2i+7)}(2i+3)(2i+5)}+
\Biggr.
$$

\vspace{1mm}
$$
\Biggl.\Biggl.+
\frac{(i^2+3i-1)\zeta_{i+1}^{(i_2)}\zeta_{i}^{(i_1)}
-(i^2+i-3)\zeta_{i}^{(i_2)}\zeta_{i+1}^{(i_1)}}
{\sqrt{(2i+1)(2i+3)}(2i-1)(2i+5)}\Biggr)\Biggr] - 
$$

\vspace{1mm}
\begin{equation}
\label{seak2}
-
\frac{1}{24}{\bf 1}_{\{i_1=i_2\}}{(T-t)^3},
\end{equation}

\vspace{6mm}

$$
I_{(11)T,t}^{(i_1 i_2)}
=-\frac{(T-t)^2}{4}I_{(00)T,t}^{(i_1 i_2)}
-\frac{T-t}{2}\left(
I_{(10)T,t}^{(i_1 i_2)}+
I_{(01)T,t}^{(i_1 i_2)}\right)+
$$

\vspace{1mm}
$$
+
\frac{(T-t)^3}{8}\Biggl[
\frac{1}{3}\zeta_1^{(i_1)}\zeta_1^{(i_2)}+\Biggr.
\sum_{i=0}^{\infty}\Biggl(
\frac{(i+1)(i+3)\left(\zeta_{i+3}^{(i_2)}\zeta_{i}^{(i_1)}
-\zeta_{i}^{(i_2)}\zeta_{i+3}^{(i_1)}\right)}
{\sqrt{(2i+1)(2i+7)}(2i+3)(2i+5)}+
\Biggr.
$$

\vspace{1mm}
$$
\Biggl.\Biggl.
+\frac{(i+1)^2\left(\zeta_{i+1}^{(i_2)}\zeta_{i}^{(i_1)}
-\zeta_{i}^{(i_2)}\zeta_{i+1}^{(i_1)}\right)}
{\sqrt{(2i+1)(2i+3)}(2i-1)(2i+5)}\Biggr)\Biggr] - 
$$

\vspace{1mm}
\begin{equation}
\label{seak3}
-
\frac{1}{24}{\bf 1}_{\{i_1=i_2\}}{(T-t)^3},
\end{equation}

\vspace{8mm}
or
$$
I_{(02)T,t}^{(i_1 i_2)}=
\hbox{\vtop{\offinterlineskip\halign{
\hfil#\hfil\cr
{\rm l.i.m.}\cr
$\stackrel{}{{}_{p\to \infty}}$\cr
}} }
\sum_{j_1,j_2=0}^p
C_{j_2j_1}^{02}\Biggl(\zeta_{j_1}^{(i_1)}\zeta_{j_2}^{(i_2)}
-{\bf 1}_{\{i_1=i_2\}}
{\bf 1}_{\{j_1=j_2\}}\Biggr),
$$

\vspace{4mm}
$$
I_{(20)T,t}^{(i_1 i_2)}=
\hbox{\vtop{\offinterlineskip\halign{
\hfil#\hfil\cr
{\rm l.i.m.}\cr
$\stackrel{}{{}_{p\to \infty}}$\cr
}} }
\sum_{j_1,j_2=0}^{p}
C_{j_2j_1}^{20}\Biggl(\zeta_{j_1}^{(i_1)}\zeta_{j_2}^{(i_2)}
-{\bf 1}_{\{i_1=i_2\}}
{\bf 1}_{\{j_1=j_2\}}\Biggr),
$$

\vspace{4mm}
$$
I_{(11)T,t}^{(i_1 i_2)}=
\hbox{\vtop{\offinterlineskip\halign{
\hfil#\hfil\cr
{\rm l.i.m.}\cr
$\stackrel{}{{}_{p\to \infty}}$\cr
}} }
\sum_{j_1,j_2=0}^{p}
C_{j_2j_1}^{11}\Biggl(\zeta_{j_1}^{(i_1)}\zeta_{j_2}^{(i_2)}
-{\bf 1}_{\{i_1=i_2\}}
{\bf 1}_{\{j_1=j_2\}}\Biggr),
$$

\vspace{6mm}

\begin{equation}
\label{gg1}
I_{(3)T,t}^{(i_1)}=-\frac{(T-t)^{7/2}}{4}\left(\zeta_0^{(i_1)}+
\frac{3\sqrt{3}}{5}\zeta_1^{(i_1)}+
\frac{1}{\sqrt{5}}\zeta_2^{(i_1)}+
\frac{1}{5\sqrt{7}}\zeta_3^{(i_1)}\right),
\end{equation}

\vspace{5mm}

$$
I_{(0000)T,t}^{*(i_1 i_2 i_3 i_4)}=
\hbox{\vtop{\offinterlineskip\halign{
\hfil#\hfil\cr
{\rm l.i.m.}\cr
$\stackrel{}{{}_{p\to \infty}}$\cr
}} }
\sum\limits_{j_1, j_2, j_3, j_4=0}^{p}
C_{j_4 j_3 j_2 j_1}\zeta_{j_1}^{(i_1)}\zeta_{j_2}^{(i_2)}\zeta_{j_3}^{(i_3)}
\zeta_{j_4}^{(i_4)},
$$

\vspace{6mm}

$$
I_{(0000)T,t}^{(i_1 i_2 i_3 i_4)}
=
\hbox{\vtop{\offinterlineskip\halign{
\hfil#\hfil\cr
{\rm l.i.m.}\cr
$\stackrel{}{{}_{p\to \infty}}$\cr
}} }
\sum_{j_1,j_2,j_3,j_4=0}^{p}
C_{j_4 j_3 j_2 j_1}\Biggl(
\zeta_{j_1}^{(i_1)}\zeta_{j_2}^{(i_2)}\zeta_{j_3}^{(i_3)}\zeta_{j_4}^{(i_4)}
-\Biggr.
$$
$$
-
{\bf 1}_{\{i_1=i_2\}}
{\bf 1}_{\{j_1=j_2\}}
\zeta_{j_3}^{(i_3)}
\zeta_{j_4}^{(i_4)}
-
{\bf 1}_{\{i_1=i_3\}}
{\bf 1}_{\{j_1=j_3\}}
\zeta_{j_2}^{(i_2)}
\zeta_{j_4}^{(i_4)}-
$$
$$
-
{\bf 1}_{\{i_1=i_4\}}
{\bf 1}_{\{j_1=j_4\}}
\zeta_{j_2}^{(i_2)}
\zeta_{j_3}^{(i_3)}
-
{\bf 1}_{\{i_2=i_3\}}
{\bf 1}_{\{j_2=j_3\}}
\zeta_{j_1}^{(i_1)}
\zeta_{j_4}^{(i_4)}-
$$
$$
-
{\bf 1}_{\{i_2=i_4\}}
{\bf 1}_{\{j_2=j_4\}}
\zeta_{j_1}^{(i_1)}
\zeta_{j_3}^{(i_3)}
-
{\bf 1}_{\{i_3=i_4\}}
{\bf 1}_{\{j_3=j_4\}}
\zeta_{j_1}^{(i_1)}
\zeta_{j_2}^{(i_2)}+
$$
$$
+
{\bf 1}_{\{i_1=i_2\}}
{\bf 1}_{\{j_1=j_2\}}
{\bf 1}_{\{i_3=i_4\}}
{\bf 1}_{\{j_3=j_4\}}+
{\bf 1}_{\{i_1=i_3\}}
{\bf 1}_{\{j_1=j_3\}}
{\bf 1}_{\{i_2=i_4\}}
{\bf 1}_{\{j_2=j_4\}}+
$$
\begin{equation}
\label{zzz10}
+\Biggl.
{\bf 1}_{\{i_1=i_4\}}
{\bf 1}_{\{j_1=j_4\}}
{\bf 1}_{\{i_2=i_3\}}
{\bf 1}_{\{j_2=j_3\}}\Biggr),
\end{equation}

\vspace{6mm}

$$
I_{(0000)T,t}^{(i_1i_1i_1i_1)}=
\frac{1}{24}(T-t)^2
\left(\left(\zeta_0^{(i_1)}\right)^4-
6\left(\zeta_0^{(i_1)}\right)^2+3\right)\ \ \ \hbox{w.\ p.\ 1},
$$

\begin{equation}
\label{ud111ee}
I_{(0000)T,t}^{*(i_1i_1i_1i_1)}=
\frac{1}{24}(T-t)^2
\left(\zeta_0^{(i_1)}\right)^4\ \ \ \hbox{w.\ p.\ 1},
\end{equation}

\vspace{2mm}
\noindent
where

\begin{equation}
\label{zzz11}
C_{j_4j_3j_2j_1}=
\frac{\sqrt{(2j_1+1)(2j_2+1)(2j_3+1)(2j_4+1)}}{16}(T-t)^{2}\bar
C_{j_4j_3j_2j_1},
\end{equation}

\begin{equation}
\label{zzz12}
\bar C_{j_4j_3j_2j_1}=\int\limits_{-1}^{1}P_{j_4}(u)
\int\limits_{-1}^{u}P_{j_3}(z)
\int\limits_{-1}^{z}P_{j_2}(y)
\int\limits_{-1}^{y}
P_{j_1}(x)dx dy dz du;
\end{equation}

\vspace{4mm}

$$
I_{(001)T,t}^{*(i_1i_2i_3)}
=\hbox{\vtop{\offinterlineskip\halign{
\hfil#\hfil\cr
{\rm l.i.m.}\cr
$\stackrel{}{{}_{p\to \infty}}$\cr
}} }
\sum_{j_1,j_2,j_3=0}^{p}
C_{j_3 j_2 j_1}^{001}
\zeta_{j_1}^{(i_1)}\zeta_{j_2}^{(i_2)}\zeta_{j_3}^{(i_3)},
$$

\vspace{3mm}
$$
I_{(010)T,t}^{*(i_1i_2i_3)}
=\hbox{\vtop{\offinterlineskip\halign{
\hfil#\hfil\cr
{\rm l.i.m.}\cr
$\stackrel{}{{}_{p\to \infty}}$\cr
}} }
\sum_{j_1,j_2,j_3=0}^{p}
C_{j_3 j_2 j_1}^{010}
\zeta_{j_1}^{(i_1)}\zeta_{j_2}^{(i_2)}\zeta_{j_3}^{(i_3)},
$$

\vspace{3mm}
$$
I_{(100)T,t}^{*(i_1i_2i_3)}
=\hbox{\vtop{\offinterlineskip\halign{
\hfil#\hfil\cr
{\rm l.i.m.}\cr
$\stackrel{}{{}_{p\to \infty}}$\cr
}} }
\sum_{j_1,j_2,j_3=0}^{p}
C_{j_3 j_2 j_1}^{100}
\zeta_{j_1}^{(i_1)}\zeta_{j_2}^{(i_2)}\zeta_{j_3}^{(i_3)},
$$

\vspace{6mm}

$$
I_{(001)T,t}^{(i_1i_2i_3)}
=\hbox{\vtop{\offinterlineskip\halign{
\hfil#\hfil\cr
{\rm l.i.m.}\cr
$\stackrel{}{{}_{p\to \infty}}$\cr
}} }
\sum_{j_1,j_2,j_3=0}^{p}
C_{j_3j_2j_1}^{001}\Biggl(
\zeta_{j_1}^{(i_1)}\zeta_{j_2}^{(i_2)}\zeta_{j_3}^{(i_3)}
-{\bf 1}_{\{i_1=i_2\}}
{\bf 1}_{\{j_1=j_2\}}
\zeta_{j_3}^{(i_3)}-
\Biggr.
$$
\begin{equation}
\label{sss1}
\Biggl.
-{\bf 1}_{\{i_2=i_3\}}
{\bf 1}_{\{j_2=j_3\}}
\zeta_{j_1}^{(i_1)}-
{\bf 1}_{\{i_1=i_3\}}
{\bf 1}_{\{j_1=j_3\}}
\zeta_{j_2}^{(i_2)}\Biggr),
\end{equation}

\vspace{6mm}

$$
I_{(010)T,t}^{(i_1i_2i_3)}
=\hbox{\vtop{\offinterlineskip\halign{
\hfil#\hfil\cr
{\rm l.i.m.}\cr
$\stackrel{}{{}_{p\to \infty}}$\cr
}} }
\sum_{j_1,j_2,j_3=0}^{p}
C_{j_3j_2j_1}^{010}\Biggl(
\zeta_{j_1}^{(i_1)}\zeta_{j_2}^{(i_2)}\zeta_{j_3}^{(i_3)}
-{\bf 1}_{\{i_1=i_2\}}
{\bf 1}_{\{j_1=j_2\}}
\zeta_{j_3}^{(i_3)}-
\Biggr.
$$
\begin{equation}
\label{sss2}
\Biggl.
-{\bf 1}_{\{i_2=i_3\}}
{\bf 1}_{\{j_2=j_3\}}
\zeta_{j_1}^{(i_1)}-
{\bf 1}_{\{i_1=i_3\}}
{\bf 1}_{\{j_1=j_3\}}
\zeta_{j_2}^{(i_2)}\Biggr),
\end{equation}

\vspace{6mm}

$$
I_{(100)T,t}^{(i_1i_2i_3)}
=\hbox{\vtop{\offinterlineskip\halign{
\hfil#\hfil\cr
{\rm l.i.m.}\cr
$\stackrel{}{{}_{p\to \infty}}$\cr
}} }
\sum_{j_1,j_2,j_3=0}^{p}
C_{j_3j_2j_1}^{100}\Biggl(
\zeta_{j_1}^{(i_1)}\zeta_{j_2}^{(i_2)}\zeta_{j_3}^{(i_3)}
-{\bf 1}_{\{i_1=i_2\}}
{\bf 1}_{\{j_1=j_2\}}
\zeta_{j_3}^{(i_3)}-
\Biggr.
$$
\begin{equation}
\label{sss3}
\Biggl.
-{\bf 1}_{\{i_2=i_3\}}
{\bf 1}_{\{j_2=j_3\}}
\zeta_{j_1}^{(i_1)}-
{\bf 1}_{\{i_1=i_3\}}
{\bf 1}_{\{j_1=j_3\}}
\zeta_{j_2}^{(i_2)}\Biggr),
\end{equation}

\vspace{5mm}
\noindent
where

$$
C_{j_3j_2j_1}^{001}
=\frac{\sqrt{(2j_1+1)(2j_2+1)(2j_3+1)}}{16}(T-t)^{5/2}\bar
C_{j_3j_2j_1}^{001},
$$

$$
C_{j_3j_2j_1}^{010}
=\frac{\sqrt{(2j_1+1)(2j_2+1)(2j_3+1)}}{16}(T-t)^{5/2}\bar
C_{j_3j_2j_1}^{010},
$$

$$
C_{j_3j_2j_1}^{100}
=\frac{\sqrt{(2j_1+1)(2j_2+1)(2j_3+1)}}{16}(T-t)^{5/2}\bar
C_{j_3j_2j_1}^{100},
$$

\vspace{3mm}

$$
\bar C_{j_3j_2j_1}^{100}=-
\int\limits_{-1}^{1}P_{j_3}(z)
\int\limits_{-1}^{z}P_{j_2}(y)
\int\limits_{-1}^{y}
P_{j_1}(x)(x+1)dx dy dz,
$$

$$
\bar C_{j_3j_2j_1}^{010}=-
\int\limits_{-1}^{1}P_{j_3}(z)
\int\limits_{-1}^{z}P_{j_2}(y)(y+1)
\int\limits_{-1}^{y}
P_{j_1}(x)dx dy dz,
$$

$$
\bar C_{j_3j_2j_1}^{001}=-
\int\limits_{-1}^{1}P_{j_3}(z)(z+1)
\int\limits_{-1}^{z}P_{j_2}(y)
\int\limits_{-1}^{y}
P_{j_1}(x)dx dy dz;
$$

\vspace{5mm}

$$
I_{(lll)T,t}^{(i_1i_1i_1)}=
\frac{1}{6}\left(\left(I_{(l)T,t}^{(i_1)}\right)^3-
3I_{(l)T,t}^{(i_1)}\Delta_{l(T,t)}\right)\ \ \ \hbox{w.\ p.\ 1},
$$

\vspace{3mm}
$$
I_{(lll)T,t}^{*(i_1i_1i_1)}=
\frac{1}{6}\left(I_{(l)T,t}^{(i_1)}\right)^3\ \ \ \hbox{w.\ p.\ 1},
$$

\vspace{3mm}

$$
I_{(llll)T,t}^{(i_1i_1i_1i_1)}=
\frac{1}{24}\left(\left(I_{(l)T,t}^{(i_1)}\right)^4-
6\left(I_{(l)T,t}^{(i_1)}\right)^2\Delta_{(l)T,t}+3
\left(\Delta_{(l)T,t}\right)^2\right)\ \ \ \hbox{w.\ p.\ 1},
$$

\vspace{3mm}
$$
I_{(llll)T,t}^{*(i_1i_1i_1i_1)}=
\frac{1}{24}\left(I_{(l)T,t}^{(i_1)}\right)^4\ \ \ \hbox{w.\ p.\ 1},
$$

\vspace{5mm}
\noindent
where
$$
I_{(l)T,t}^{(i_1)}=\sum_{j=0}^l C_j^l \zeta_j^{(i_1)}\ \ \ \hbox{w.\ p.\ 1},
$$       

\vspace{3mm}
$$
\Delta_{l(T,t)}=\int\limits_t^T(t-s)^{2l}ds,\ \ \
C_j^l=\int\limits_t^T(t-s)^l\phi_j(s)ds;
$$ 

\vspace{6mm}

$$
I_{(00000)T,t}^{*(i_1 i_2 i_3 i_4 i_5)}=
\hbox{\vtop{\offinterlineskip\halign{
\hfil#\hfil\cr
{\rm l.i.m.}\cr
$\stackrel{}{{}_{p\to \infty}}$\cr
}} }
\sum\limits_{j_1, j_2, j_3, j_4, j_5=0}^{p}
C_{j_5j_4 j_3 j_2 j_1}
\zeta_{j_1}^{(i_1)}\zeta_{j_2}^{(i_2)}\zeta_{j_3}^{(i_3)}
\zeta_{j_4}^{(i_4)}\zeta_{j_5}^{(i_5)},
$$

\vspace{7mm}

$$
I_{(00000)T,t}^{(i_1 i_2 i_3 i_4 i_5)}
=
\hbox{\vtop{\offinterlineskip\halign{
\hfil#\hfil\cr
{\rm l.i.m.}\cr
$\stackrel{}{{}_{p\to \infty}}$\cr
}} }
\sum_{j_1,j_2,j_3,j_4,j_5=0}^p
C_{j_5 j_4 j_3 j_2 j_1}\Biggl(
\prod_{l=1}^5\zeta_{j_l}^{(i_l)}
-\Biggr.
$$
$$
-
{\bf 1}_{\{i_1=i_2\}}
{\bf 1}_{\{j_1=j_2\}}
\zeta_{j_3}^{(i_3)}
\zeta_{j_4}^{(i_4)}
\zeta_{j_5}^{(i_5)}-
{\bf 1}_{\{i_1=i_3\}}
{\bf 1}_{\{j_1=j_3\}}
\zeta_{j_2}^{(i_2)}
\zeta_{j_4}^{(i_4)}
\zeta_{j_5}^{(i_5)}-
$$
$$
-
{\bf 1}_{\{i_1=i_4\}}
{\bf 1}_{\{j_1=j_4\}}
\zeta_{j_2}^{(i_2)}
\zeta_{j_3}^{(i_3)}
\zeta_{j_5}^{(i_5)}-
{\bf 1}_{\{i_1=i_5\}}
{\bf 1}_{\{j_1=j_5\}}
\zeta_{j_2}^{(i_2)}
\zeta_{j_3}^{(i_3)}
\zeta_{j_4}^{(i_4)}-
$$
$$
-
{\bf 1}_{\{i_2=i_3\}}
{\bf 1}_{\{j_2=j_3\}}
\zeta_{j_1}^{(i_1)}
\zeta_{j_4}^{(i_4)}
\zeta_{j_5}^{(i_5)}-
{\bf 1}_{\{i_2=i_4\}}
{\bf 1}_{\{j_2=j_4\}}
\zeta_{j_1}^{(i_1)}
\zeta_{j_3}^{(i_3)}
\zeta_{j_5}^{(i_5)}-
$$
$$
-
{\bf 1}_{\{i_2=i_5\}}
{\bf 1}_{\{j_2=j_5\}}
\zeta_{j_1}^{(i_1)}
\zeta_{j_3}^{(i_3)}
\zeta_{j_4}^{(i_4)}
-{\bf 1}_{\{i_3=i_4\}}
{\bf 1}_{\{j_3=j_4\}}
\zeta_{j_1}^{(i_1)}
\zeta_{j_2}^{(i_2)}
\zeta_{j_5}^{(i_5)}-
$$
$$
-
{\bf 1}_{\{i_3=i_5\}}
{\bf 1}_{\{j_3=j_5\}}
\zeta_{j_1}^{(i_1)}
\zeta_{j_2}^{(i_2)}
\zeta_{j_4}^{(i_4)}
-{\bf 1}_{\{i_4=i_5\}}
{\bf 1}_{\{j_4=j_5\}}
\zeta_{j_1}^{(i_1)}
\zeta_{j_2}^{(i_2)}
\zeta_{j_3}^{(i_3)}+
$$
$$
+
{\bf 1}_{\{i_1=i_2\}}
{\bf 1}_{\{j_1=j_2\}}
{\bf 1}_{\{i_3=i_4\}}
{\bf 1}_{\{j_3=j_4\}}\zeta_{j_5}^{(i_5)}+
{\bf 1}_{\{i_1=i_2\}}
{\bf 1}_{\{j_1=j_2\}}
{\bf 1}_{\{i_3=i_5\}}
{\bf 1}_{\{j_3=j_5\}}\zeta_{j_4}^{(i_4)}+
$$
$$
+
{\bf 1}_{\{i_1=i_2\}}
{\bf 1}_{\{j_1=j_2\}}
{\bf 1}_{\{i_4=i_5\}}
{\bf 1}_{\{j_4=j_5\}}\zeta_{j_3}^{(i_3)}+
{\bf 1}_{\{i_1=i_3\}}
{\bf 1}_{\{j_1=j_3\}}
{\bf 1}_{\{i_2=i_4\}}
{\bf 1}_{\{j_2=j_4\}}\zeta_{j_5}^{(i_5)}+
$$
$$
+
{\bf 1}_{\{i_1=i_3\}}
{\bf 1}_{\{j_1=j_3\}}
{\bf 1}_{\{i_2=i_5\}}
{\bf 1}_{\{j_2=j_5\}}\zeta_{j_4}^{(i_4)}+
{\bf 1}_{\{i_1=i_3\}}
{\bf 1}_{\{j_1=j_3\}}
{\bf 1}_{\{i_4=i_5\}}
{\bf 1}_{\{j_4=j_5\}}\zeta_{j_2}^{(i_2)}+
$$
$$
+
{\bf 1}_{\{i_1=i_4\}}
{\bf 1}_{\{j_1=j_4\}}
{\bf 1}_{\{i_2=i_3\}}
{\bf 1}_{\{j_2=j_3\}}\zeta_{j_5}^{(i_5)}+
{\bf 1}_{\{i_1=i_4\}}
{\bf 1}_{\{j_1=j_4\}}
{\bf 1}_{\{i_2=i_5\}}
{\bf 1}_{\{j_2=j_5\}}\zeta_{j_3}^{(i_3)}+
$$
$$
+
{\bf 1}_{\{i_1=i_4\}}
{\bf 1}_{\{j_1=j_4\}}
{\bf 1}_{\{i_3=i_5\}}
{\bf 1}_{\{j_3=j_5\}}\zeta_{j_2}^{(i_2)}+
{\bf 1}_{\{i_1=i_5\}}
{\bf 1}_{\{j_1=j_5\}}
{\bf 1}_{\{i_2=i_3\}}
{\bf 1}_{\{j_2=j_3\}}\zeta_{j_4}^{(i_4)}+
$$
$$
+
{\bf 1}_{\{i_1=i_5\}}
{\bf 1}_{\{j_1=j_5\}}
{\bf 1}_{\{i_2=i_4\}}
{\bf 1}_{\{j_2=j_4\}}\zeta_{j_3}^{(i_3)}+
{\bf 1}_{\{i_1=i_5\}}
{\bf 1}_{\{j_1=j_5\}}
{\bf 1}_{\{i_3=i_4\}}
{\bf 1}_{\{j_3=j_4\}}\zeta_{j_2}^{(i_2)}+
$$
$$
+
{\bf 1}_{\{i_2=i_3\}}
{\bf 1}_{\{j_2=j_3\}}
{\bf 1}_{\{i_4=i_5\}}
{\bf 1}_{\{j_4=j_5\}}\zeta_{j_1}^{(i_1)}+
{\bf 1}_{\{i_2=i_4\}}
{\bf 1}_{\{j_2=j_4\}}
{\bf 1}_{\{i_3=i_5\}}
{\bf 1}_{\{j_3=j_5\}}\zeta_{j_1}^{(i_1)}+
$$
\begin{equation}
\label{sss4}
+\Biggl.
{\bf 1}_{\{i_2=i_5\}}
{\bf 1}_{\{j_2=j_5\}}
{\bf 1}_{\{i_3=i_4\}}
{\bf 1}_{\{j_3=j_4\}}\zeta_{j_1}^{(i_1)}\Biggr),
\end{equation}

\vspace{7mm}
$$         
I_{(00000)T,t}^{(i_1i_1i_1i_1i_1)}=
\frac{1}{120}(T-t)^{5/2}
\left(\left(\zeta_0^{(i_1)}\right)^5-
10\left(\zeta_0^{(i_1)}\right)^3+15\zeta_0^{(i_1)}\right)\ \ \ 
\hbox{w.\ p.\ 1},
$$

\vspace{2mm}
$$
I_{(00000)T,t}^{*(i_1i_1i_1i_1i_1)}=
\frac{1}{120}(T-t)^{5/2}\left(\zeta_0^{(i_1)}\right)^5\ \ \ \hbox{w.\ p.\ 1},
$$

\vspace{2mm}
\noindent
where 

$$
C_{j_5j_4 j_3 j_2 j_1}=
\frac{\sqrt{(2j_1+1)(2j_2+1)(2j_3+1)(2j_4+1)(2j_5+1)}}{32}(T-t)^{5/2}\bar
C_{j_5j_4 j_3 j_2 j_1},
$$

$$
\bar C_{j_5j_4 j_3 j_2 j_1}=
\int\limits_{-1}^{1}P_{j_5}(v)
\int\limits_{-1}^{v}P_{j_4}(u)
\int\limits_{-1}^{u}P_{j_3}(z)
\int\limits_{-1}^{z}P_{j_2}(y)
\int\limits_{-1}^{y}
P_{j_1}(x)dx dy dz du dv;
$$

\vspace{4mm}

$$
I_{(0001)T,t}^{*(i_1i_2i_3)}
=\hbox{\vtop{\offinterlineskip\halign{
\hfil#\hfil\cr
{\rm l.i.m.}\cr
$\stackrel{}{{}_{p\to \infty}}$\cr
}} }
\sum_{j_1,j_2,j_3,j_4=0}^{p}
C_{j_4j_3 j_2 j_1}^{0001}
\zeta_{j_1}^{(i_1)}\zeta_{j_2}^{(i_2)}\zeta_{j_3}^{(i_3)}\zeta_{j_4}^{(i_4)},
$$

\vspace{2mm}
$$
I_{(0010)T,t}^{*(i_1i_2i_3)}
=\hbox{\vtop{\offinterlineskip\halign{
\hfil#\hfil\cr
{\rm l.i.m.}\cr
$\stackrel{}{{}_{p\to \infty}}$\cr
}} }
\sum_{j_1,j_2,j_3,j_4=0}^{p}
C_{j_4j_3 j_2 j_1}^{0010}
\zeta_{j_1}^{(i_1)}\zeta_{j_2}^{(i_2)}\zeta_{j_3}^{(i_3)}\zeta_{j_4}^{(i_4)},
$$

\vspace{2mm}

$$
I_{(0100)T,t}^{*(i_1i_2i_3)}
=\hbox{\vtop{\offinterlineskip\halign{
\hfil#\hfil\cr
{\rm l.i.m.}\cr
$\stackrel{}{{}_{p\to \infty}}$\cr
}} }
\sum_{j_1,j_2,j_3,j_4=0}^{p}
C_{j_4j_3 j_2 j_1}^{0100}
\zeta_{j_1}^{(i_1)}\zeta_{j_2}^{(i_2)}\zeta_{j_3}^{(i_3)}\zeta_{j_4}^{(i_4)},
$$

\vspace{2mm}

$$
I_{(1000)T,t}^{*(i_1i_2i_3)}
=\hbox{\vtop{\offinterlineskip\halign{
\hfil#\hfil\cr
{\rm l.i.m.}\cr
$\stackrel{}{{}_{p\to \infty}}$\cr
}} }
\sum_{j_1,j_2,j_3,j_4=0}^{p}
C_{j_4j_3 j_2 j_1}^{1000}
\zeta_{j_1}^{(i_1)}\zeta_{j_2}^{(i_2)}\zeta_{j_3}^{(i_3)}\zeta_{j_4}^{(i_4)},
$$

\vspace{6mm}

$$
I_{(0001)T,t}^{(i_1 i_2 i_3 i_4)}
=\hbox{\vtop{\offinterlineskip\halign{
\hfil#\hfil\cr
{\rm l.i.m.}\cr
$\stackrel{}{{}_{p\to \infty}}$\cr
}} }
\sum_{j_1,j_2,j_3,j_4=0}^{p}
C_{j_4 j_3 j_2 j_1}^{0001}\Biggl(
\zeta_{j_1}^{(i_1)}\zeta_{j_2}^{(i_2)}\zeta_{j_3}^{(i_3)}\zeta_{j_4}^{(i_4)}
-\Biggr.
$$
$$
-
{\bf 1}_{\{i_1=i_2\}}
{\bf 1}_{\{j_1=j_2\}}
\zeta_{j_3}^{(i_3)}
\zeta_{j_4}^{(i_4)}
-
{\bf 1}_{\{i_1=i_3\}}
{\bf 1}_{\{j_1=j_3\}}
\zeta_{j_2}^{(i_2)}
\zeta_{j_4}^{(i_4)}-
$$
$$
-
{\bf 1}_{\{i_1=i_4\}}
{\bf 1}_{\{j_1=j_4\}}
\zeta_{j_2}^{(i_2)}
\zeta_{j_3}^{(i_3)}
-
{\bf 1}_{\{i_2=i_3\}}
{\bf 1}_{\{j_2=j_3\}}
\zeta_{j_1}^{(i_1)}
\zeta_{j_4}^{(i_4)}-
$$
$$
-
{\bf 1}_{\{i_2=i_4\}}
{\bf 1}_{\{j_2=j_4\}}
\zeta_{j_1}^{(i_1)}
\zeta_{j_3}^{(i_3)}
-
{\bf 1}_{\{i_3=i_4\}}
{\bf 1}_{\{j_3=j_4\}}
\zeta_{j_1}^{(i_1)}
\zeta_{j_2}^{(i_2)}+
$$
$$
+
{\bf 1}_{\{i_1=i_2\}}
{\bf 1}_{\{j_1=j_2\}}
{\bf 1}_{\{i_3=i_4\}}
{\bf 1}_{\{j_3=j_4\}}+
{\bf 1}_{\{i_1=i_3\}}
{\bf 1}_{\{j_1=j_3\}}
{\bf 1}_{\{i_2=i_4\}}
{\bf 1}_{\{j_2=j_4\}}+
$$
$$
+\Biggl.
{\bf 1}_{\{i_1=i_4\}}
{\bf 1}_{\{j_1=j_4\}}
{\bf 1}_{\{i_2=i_3\}}
{\bf 1}_{\{j_2=j_3\}}\Biggr),
$$

\vspace{6mm}

$$
I_{(0010)T,t}^{(i_1 i_2 i_3 i_4)}
=\hbox{\vtop{\offinterlineskip\halign{
\hfil#\hfil\cr
{\rm l.i.m.}\cr
$\stackrel{}{{}_{p\to \infty}}$\cr
}} }
\sum_{j_1,j_2,j_3,j_4=0}^{p}
C_{j_4 j_3 j_2 j_1}^{0010}\Biggl(
\zeta_{j_1}^{(i_1)}\zeta_{j_2}^{(i_2)}\zeta_{j_3}^{(i_3)}\zeta_{j_4}^{(i_4)}
-\Biggr.
$$
$$
-
{\bf 1}_{\{i_1=i_2\}}
{\bf 1}_{\{j_1=j_2\}}
\zeta_{j_3}^{(i_3)}
\zeta_{j_4}^{(i_4)}
-
{\bf 1}_{\{i_1=i_3\}}
{\bf 1}_{\{j_1=j_3\}}
\zeta_{j_2}^{(i_2)}
\zeta_{j_4}^{(i_4)}-
$$
$$
-
{\bf 1}_{\{i_1=i_4\}}
{\bf 1}_{\{j_1=j_4\}}
\zeta_{j_2}^{(i_2)}
\zeta_{j_3}^{(i_3)}
-
{\bf 1}_{\{i_2=i_3\}}
{\bf 1}_{\{j_2=j_3\}}
\zeta_{j_1}^{(i_1)}
\zeta_{j_4}^{(i_4)}-
$$
$$
-
{\bf 1}_{\{i_2=i_4\}}
{\bf 1}_{\{j_2=j_4\}}
\zeta_{j_1}^{(i_1)}
\zeta_{j_3}^{(i_3)}
-
{\bf 1}_{\{i_3=i_4\}}
{\bf 1}_{\{j_3=j_4\}}
\zeta_{j_1}^{(i_1)}
\zeta_{j_2}^{(i_2)}+
$$
$$
+
{\bf 1}_{\{i_1=i_2\}}
{\bf 1}_{\{j_1=j_2\}}
{\bf 1}_{\{i_3=i_4\}}
{\bf 1}_{\{j_3=j_4\}}+
{\bf 1}_{\{i_1=i_3\}}
{\bf 1}_{\{j_1=j_3\}}
{\bf 1}_{\{i_2=i_4\}}
{\bf 1}_{\{j_2=j_4\}}+
$$
$$
+\Biggl.
{\bf 1}_{\{i_1=i_4\}}
{\bf 1}_{\{j_1=j_4\}}
{\bf 1}_{\{i_2=i_3\}}
{\bf 1}_{\{j_2=j_3\}}\Biggr),
$$

\vspace{6mm}

$$
I_{(0100)T,t}^{(i_1 i_2 i_3 i_4)}
=\hbox{\vtop{\offinterlineskip\halign{
\hfil#\hfil\cr
{\rm l.i.m.}\cr
$\stackrel{}{{}_{p\to \infty}}$\cr
}} }
\sum_{j_1,j_2,j_3,j_4=0}^{p}
C_{j_4 j_3 j_2 j_1}^{0100}\Biggl(
\zeta_{j_1}^{(i_1)}\zeta_{j_2}^{(i_2)}\zeta_{j_3}^{(i_3)}\zeta_{j_4}^{(i_4)}
-\Biggr.
$$
$$
-
{\bf 1}_{\{i_1=i_2\}}
{\bf 1}_{\{j_1=j_2\}}
\zeta_{j_3}^{(i_3)}
\zeta_{j_4}^{(i_4)}
-
{\bf 1}_{\{i_1=i_3\}}
{\bf 1}_{\{j_1=j_3\}}
\zeta_{j_2}^{(i_2)}
\zeta_{j_4}^{(i_4)}-
$$
$$
-
{\bf 1}_{\{i_1=i_4\}}
{\bf 1}_{\{j_1=j_4\}}
\zeta_{j_2}^{(i_2)}
\zeta_{j_3}^{(i_3)}
-
{\bf 1}_{\{i_2=i_3\}}
{\bf 1}_{\{j_2=j_3\}}
\zeta_{j_1}^{(i_1)}
\zeta_{j_4}^{(i_4)}-
$$
$$
-
{\bf 1}_{\{i_2=i_4\}}
{\bf 1}_{\{j_2=j_4\}}
\zeta_{j_1}^{(i_1)}
\zeta_{j_3}^{(i_3)}
-
{\bf 1}_{\{i_3=i_4\}}
{\bf 1}_{\{j_3=j_4\}}
\zeta_{j_1}^{(i_1)}
\zeta_{j_2}^{(i_2)}+
$$
$$
+
{\bf 1}_{\{i_1=i_2\}}
{\bf 1}_{\{j_1=j_2\}}
{\bf 1}_{\{i_3=i_4\}}
{\bf 1}_{\{j_3=j_4\}}+
{\bf 1}_{\{i_1=i_3\}}
{\bf 1}_{\{j_1=j_3\}}
{\bf 1}_{\{i_2=i_4\}}
{\bf 1}_{\{j_2=j_4\}}+
$$
$$
+\Biggl.
{\bf 1}_{\{i_1=i_4\}}
{\bf 1}_{\{j_1=j_4\}}
{\bf 1}_{\{i_2=i_3\}}
{\bf 1}_{\{j_2=j_3\}}\Biggr),
$$

\vspace{6mm}

$$
I_{(1000)T,t}^{(i_1 i_2 i_3 i_4)}
=\hbox{\vtop{\offinterlineskip\halign{
\hfil#\hfil\cr
{\rm l.i.m.}\cr
$\stackrel{}{{}_{p\to \infty}}$\cr
}} }
\sum_{j_1,j_2,j_3,j_4=0}^{p}
C_{j_4 j_3 j_2 j_1}^{1000}\Biggl(
\zeta_{j_1}^{(i_1)}\zeta_{j_2}^{(i_2)}\zeta_{j_3}^{(i_3)}\zeta_{j_4}^{(i_4)}
-\Biggr.
$$
$$
-
{\bf 1}_{\{i_1=i_2\}}
{\bf 1}_{\{j_1=j_2\}}
\zeta_{j_3}^{(i_3)}
\zeta_{j_4}^{(i_4)}
-
{\bf 1}_{\{i_1=i_3\}}
{\bf 1}_{\{j_1=j_3\}}
\zeta_{j_2}^{(i_2)}
\zeta_{j_4}^{(i_4)}-
$$
$$
-
{\bf 1}_{\{i_1=i_4\}}
{\bf 1}_{\{j_1=j_4\}}
\zeta_{j_2}^{(i_2)}
\zeta_{j_3}^{(i_3)}
-
{\bf 1}_{\{i_2=i_3\}}
{\bf 1}_{\{j_2=j_3\}}
\zeta_{j_1}^{(i_1)}
\zeta_{j_4}^{(i_4)}-
$$
$$
-
{\bf 1}_{\{i_2=i_4\}}
{\bf 1}_{\{j_2=j_4\}}
\zeta_{j_1}^{(i_1)}
\zeta_{j_3}^{(i_3)}
-
{\bf 1}_{\{i_3=i_4\}}
{\bf 1}_{\{j_3=j_4\}}
\zeta_{j_1}^{(i_1)}
\zeta_{j_2}^{(i_2)}+
$$
$$
+
{\bf 1}_{\{i_1=i_2\}}
{\bf 1}_{\{j_1=j_2\}}
{\bf 1}_{\{i_3=i_4\}}
{\bf 1}_{\{j_3=j_4\}}+
{\bf 1}_{\{i_1=i_3\}}
{\bf 1}_{\{j_1=j_3\}}
{\bf 1}_{\{i_2=i_4\}}
{\bf 1}_{\{j_2=j_4\}}+
$$
$$
+\Biggl.
{\bf 1}_{\{i_1=i_4\}}
{\bf 1}_{\{j_1=j_4\}}
{\bf 1}_{\{i_2=i_3\}}
{\bf 1}_{\{j_2=j_3\}}\Biggr),
$$

\vspace{3mm}
\noindent
where

$$
C_{j_4j_3j_2j_1}^{0001}
=\frac{\sqrt{(2j_1+1)(2j_2+1)(2j_3+1)(2j_4+1)}}{32}(T-t)^{3}\bar
C_{j_4j_3j_2j_1}^{0001},
$$

\vspace{2mm}
$$
C_{j_3j_2j_1}^{0010}
=\frac{\sqrt{(2j_1+1)(2j_2+1)(2j_3+1)(2j_4+1)}}{32}(T-t)^{3}\bar
C_{j_4j_3j_2j_1}^{0010},
$$

\vspace{2mm}
$$
C_{j_4j_3j_2j_1}^{0100}=
\frac{\sqrt{(2j_1+1)(2j_2+1)(2j_3+1)(2j_4+1)}}{32}(T-t)^{3}\bar
C_{j_3j_2j_1}^{0100},
$$

\vspace{2mm}
$$
C_{j_4j_3j_2j_1}^{1000}
=\frac{\sqrt{(2j_1+1)(2j_2+1)(2j_3+1)(2j_4+1)}}{32}(T-t)^{3}\bar
C_{j_4j_3j_2j_1}^{1000},
$$

\vspace{2mm}
$$
\bar C_{j_4j_3j_2j_1}^{1000}=-
\int\limits_{-1}^{1}P_{j_4}(u)
\int\limits_{-1}^{u}P_{j_3}(z)
\int\limits_{-1}^{z}P_{j_2}(y)
\int\limits_{-1}^{y}
P_{j_1}(x)(x+1)dx dy dz du,
$$

\vspace{2mm}
$$
\bar C_{j_4j_3j_2j_1}^{0100}=-
\int\limits_{-1}^{1}P_{j_4}(u)
\int\limits_{-1}^{u}P_{j_3}(z)
\int\limits_{-1}^{z}P_{j_2}(y)(y+1)
\int\limits_{-1}^{y}
P_{j_1}(x)dx dy dz du,
$$

\vspace{2mm}
$$
\bar C_{j_4j_3j_2j_1}^{0010}=-
\int\limits_{-1}^{1}P_{j_4}(u)
\int\limits_{-1}^{u}P_{j_3}(z)(z+1)
\int\limits_{-1}^{z}P_{j_2}(y)
\int\limits_{-1}^{y}
P_{j_1}(x)dx dy dz du,
$$

\vspace{2mm}
$$
\bar C_{j_4j_3j_2j_1}^{0001}=-
\int\limits_{-1}^{1}P_{j_4}(u)(u+1)
\int\limits_{-1}^{u}P_{j_3}(z)
\int\limits_{-1}^{z}P_{j_2}(y)
\int\limits_{-1}^{y}
P_{j_1}(x)dx dy dz du;
$$

\vspace{5mm}

$$
I_{(000000)T,t}^{*(i_1 i_2 i_3 i_4 i_5 i_6)}=
\hbox{\vtop{\offinterlineskip\halign{
\hfil#\hfil\cr
{\rm l.i.m.}\cr
$\stackrel{}{{}_{p\to \infty}}$\cr
}} }
\sum\limits_{j_1, j_2, j_3, j_4, j_5, j_6=0}^{p}
C_{j_6j_5j_4 j_3 j_2 j_1}
\zeta_{j_1}^{(i_1)}\zeta_{j_2}^{(i_2)}\zeta_{j_3}^{(i_3)}
\zeta_{j_4}^{(i_4)}\zeta_{j_5}^{(i_5)}\zeta_{j_6}^{(i_6)},
$$

\vspace{7mm}

$$
I_{(000000)T,t}^{(i_1 i_2 i_3 i_4 i_5 i_6)}
=\hbox{\vtop{\offinterlineskip\halign{
\hfil#\hfil\cr
{\rm l.i.m.}\cr
$\stackrel{}{{}_{p\to \infty}}$\cr
}} }\sum_{j_1,j_2,j_3,j_4,j_5,j_6=0}^{p}
C_{j_6 j_5 j_4 j_3 j_2 j_1}\Biggl(
\prod_{l=1}^6
\zeta_{j_l}^{(i_l)}
-\Biggr.
$$
$$
-
{\bf 1}_{\{j_1=j_6\}}
{\bf 1}_{\{i_1=i_6\}}
\zeta_{j_2}^{(i_2)}
\zeta_{j_3}^{(i_3)}
\zeta_{j_4}^{(i_4)}
\zeta_{j_5}^{(i_5)}-
{\bf 1}_{\{j_2=j_6\}}
{\bf 1}_{\{i_2=i_6\}}
\zeta_{j_1}^{(i_1)}
\zeta_{j_3}^{(i_3)}
\zeta_{j_4}^{(i_4)}
\zeta_{j_5}^{(i_5)}-
$$
$$
-
{\bf 1}_{\{j_3=j_6\}}
{\bf 1}_{\{i_3=i_6\}}
\zeta_{j_1}^{(i_1)}
\zeta_{j_2}^{(i_2)}
\zeta_{j_4}^{(i_4)}
\zeta_{j_5}^{(i_5)}-
{\bf 1}_{\{j_4=j_6\}}
{\bf 1}_{\{i_4=i_6\}}
\zeta_{j_1}^{(i_1)}
\zeta_{j_2}^{(i_2)}
\zeta_{j_3}^{(i_3)}
\zeta_{j_5}^{(i_5)}-
$$
$$
-
{\bf 1}_{\{j_5=j_6\}}
{\bf 1}_{\{i_5=i_6\}}
\zeta_{j_1}^{(i_1)}
\zeta_{j_2}^{(i_2)}
\zeta_{j_3}^{(i_3)}
\zeta_{j_4}^{(i_4)}-
{\bf 1}_{\{j_1=j_2\}}
{\bf 1}_{\{i_1=i_2\}}
\zeta_{j_3}^{(i_3)}
\zeta_{j_4}^{(i_4)}
\zeta_{j_5}^{(i_5)}
\zeta_{j_6}^{(i_6)}-
$$
$$
-
{\bf 1}_{\{j_1=j_3\}}
{\bf 1}_{\{i_1=i_3\}}
\zeta_{j_2}^{(i_2)}
\zeta_{j_4}^{(i_4)}
\zeta_{j_5}^{(i_5)}
\zeta_{j_6}^{(i_6)}-
{\bf 1}_{\{j_1=j_4\}}
{\bf 1}_{\{i_1=i_4\}}
\zeta_{j_2}^{(i_2)}
\zeta_{j_3}^{(i_3)}
\zeta_{j_5}^{(i_5)}
\zeta_{j_6}^{(i_6)}-
$$
$$
-
{\bf 1}_{\{j_1=j_5\}}
{\bf 1}_{\{i_1=i_5\}}
\zeta_{j_2}^{(i_2)}
\zeta_{j_3}^{(i_3)}
\zeta_{j_4}^{(i_4)}
\zeta_{j_6}^{(i_6)}-
{\bf 1}_{\{j_2=j_3\}}
{\bf 1}_{\{i_2=i_3\}}
\zeta_{j_1}^{(i_1)}
\zeta_{j_4}^{(i_4)}
\zeta_{j_5}^{(i_5)}
\zeta_{j_6}^{(i_6)}-
$$
$$
-
{\bf 1}_{\{j_2=j_4\}}
{\bf 1}_{\{i_2=i_4\}}
\zeta_{j_1}^{(i_1)}
\zeta_{j_3}^{(i_3)}
\zeta_{j_5}^{(i_5)}
\zeta_{j_6}^{(i_6)}-
{\bf 1}_{\{j_2=j_5\}}
{\bf 1}_{\{i_2=i_5\}}
\zeta_{j_1}^{(i_1)}
\zeta_{j_3}^{(i_3)}
\zeta_{j_4}^{(i_4)}
\zeta_{j_6}^{(i_6)}-
$$
$$
-
{\bf 1}_{\{j_3=j_4\}}
{\bf 1}_{\{i_3=i_4\}}
\zeta_{j_1}^{(i_1)}
\zeta_{j_2}^{(i_2)}
\zeta_{j_5}^{(i_5)}
\zeta_{j_6}^{(i_6)}-
{\bf 1}_{\{j_3=j_5\}}
{\bf 1}_{\{i_3=i_5\}}
\zeta_{j_1}^{(i_1)}
\zeta_{j_2}^{(i_2)}
\zeta_{j_4}^{(i_4)}
\zeta_{j_6}^{(i_6)}-
$$
$$
-
{\bf 1}_{\{j_4=j_5\}}
{\bf 1}_{\{i_4=i_5\}}
\zeta_{j_1}^{(i_1)}
\zeta_{j_2}^{(i_2)}
\zeta_{j_3}^{(i_3)}
\zeta_{j_6}^{(i_6)}+
$$
$$
+
{\bf 1}_{\{j_1=j_2\}}
{\bf 1}_{\{i_1=i_2\}}
{\bf 1}_{\{j_3=j_4\}}
{\bf 1}_{\{i_3=i_4\}}
\zeta_{j_5}^{(i_5)}
\zeta_{j_6}^{(i_6)}
+
{\bf 1}_{\{j_1=j_2\}}
{\bf 1}_{\{i_1=i_2\}}
{\bf 1}_{\{j_3=j_5\}}
{\bf 1}_{\{i_3=i_5\}}
\zeta_{j_4}^{(i_4)}
\zeta_{j_6}^{(i_6)}+
$$
$$
+
{\bf 1}_{\{j_1=j_2\}}
{\bf 1}_{\{i_1=i_2\}}
{\bf 1}_{\{j_4=j_5\}}
{\bf 1}_{\{i_4=i_5\}}
\zeta_{j_3}^{(i_3)}
\zeta_{j_6}^{(i_6)}
+
{\bf 1}_{\{j_1=j_3\}}
{\bf 1}_{\{i_1=i_3\}}
{\bf 1}_{\{j_2=j_4\}}
{\bf 1}_{\{i_2=i_4\}}
\zeta_{j_5}^{(i_5)}
\zeta_{j_6}^{(i_6)}+
$$
$$
+
{\bf 1}_{\{j_1=j_3\}}
{\bf 1}_{\{i_1=i_3\}}
{\bf 1}_{\{j_2=j_5\}}
{\bf 1}_{\{i_2=i_5\}}
\zeta_{j_4}^{(i_4)}
\zeta_{j_6}^{(i_6)}
+
{\bf 1}_{\{j_1=j_3\}}
{\bf 1}_{\{i_1=i_3\}}
{\bf 1}_{\{j_4=j_5\}}
{\bf 1}_{\{i_4=i_5\}}
\zeta_{j_2}^{(i_2)}
\zeta_{j_6}^{(i_6)}+
$$
$$
+
{\bf 1}_{\{j_1=j_4\}}
{\bf 1}_{\{i_1=i_4\}}
{\bf 1}_{\{j_2=j_3\}}
{\bf 1}_{\{i_2=i_3\}}
\zeta_{j_5}^{(i_5)}
\zeta_{j_6}^{(i_6)}
+
{\bf 1}_{\{j_1=j_4\}}
{\bf 1}_{\{i_1=i_4\}}
{\bf 1}_{\{j_2=j_5\}}
{\bf 1}_{\{i_2=i_5\}}
\zeta_{j_3}^{(i_3)}
\zeta_{j_6}^{(i_6)}+
$$
$$
+
{\bf 1}_{\{j_1=j_4\}}
{\bf 1}_{\{i_1=i_4\}}
{\bf 1}_{\{j_3=j_5\}}
{\bf 1}_{\{i_3=i_5\}}
\zeta_{j_2}^{(i_2)}
\zeta_{j_6}^{(i_6)}
+
{\bf 1}_{\{j_1=j_5\}}
{\bf 1}_{\{i_1=i_5\}}
{\bf 1}_{\{j_2=j_3\}}
{\bf 1}_{\{i_2=i_3\}}
\zeta_{j_4}^{(i_4)}
\zeta_{j_6}^{(i_6)}+
$$
$$
+
{\bf 1}_{\{j_1=j_5\}}
{\bf 1}_{\{i_1=i_5\}}
{\bf 1}_{\{j_2=j_4\}}
{\bf 1}_{\{i_2=i_4\}}
\zeta_{j_3}^{(i_3)}
\zeta_{j_6}^{(i_6)}
+
{\bf 1}_{\{j_1=j_5\}}
{\bf 1}_{\{i_1=i_5\}}
{\bf 1}_{\{j_3=j_4\}}
{\bf 1}_{\{i_3=i_4\}}
\zeta_{j_2}^{(i_2)}
\zeta_{j_6}^{(i_6)}+
$$
$$
+
{\bf 1}_{\{j_2=j_3\}}
{\bf 1}_{\{i_2=i_3\}}
{\bf 1}_{\{j_4=j_5\}}
{\bf 1}_{\{i_4=i_5\}}
\zeta_{j_1}^{(i_1)}
\zeta_{j_6}^{(i_6)}
+
{\bf 1}_{\{j_2=j_4\}}
{\bf 1}_{\{i_2=i_4\}}
{\bf 1}_{\{j_3=j_5\}}
{\bf 1}_{\{i_3=i_5\}}
\zeta_{j_1}^{(i_1)}
\zeta_{j_6}^{(i_6)}+
$$
$$
+
{\bf 1}_{\{j_2=j_5\}}
{\bf 1}_{\{i_2=i_5\}}
{\bf 1}_{\{j_3=j_4\}}
{\bf 1}_{\{i_3=i_4\}}
\zeta_{j_1}^{(i_1)}
\zeta_{j_6}^{(i_6)}
+
{\bf 1}_{\{j_6=j_1\}}
{\bf 1}_{\{i_6=i_1\}}
{\bf 1}_{\{j_3=j_4\}}
{\bf 1}_{\{i_3=i_4\}}
\zeta_{j_2}^{(i_2)}
\zeta_{j_5}^{(i_5)}+
$$
$$
+
{\bf 1}_{\{j_6=j_1\}}
{\bf 1}_{\{i_6=i_1\}}
{\bf 1}_{\{j_3=j_5\}}
{\bf 1}_{\{i_3=i_5\}}
\zeta_{j_2}^{(i_2)}
\zeta_{j_4}^{(i_4)}
+
{\bf 1}_{\{j_6=j_1\}}
{\bf 1}_{\{i_6=i_1\}}
{\bf 1}_{\{j_2=j_5\}}
{\bf 1}_{\{i_2=i_5\}}
\zeta_{j_3}^{(i_3)}
\zeta_{j_4}^{(i_4)}+
$$
$$
+
{\bf 1}_{\{j_6=j_1\}}
{\bf 1}_{\{i_6=i_1\}}
{\bf 1}_{\{j_2=j_4\}}
{\bf 1}_{\{i_2=i_4\}}
\zeta_{j_3}^{(i_3)}
\zeta_{j_5}^{(i_5)}
+
{\bf 1}_{\{j_6=j_1\}}
{\bf 1}_{\{i_6=i_1\}}
{\bf 1}_{\{j_4=j_5\}}
{\bf 1}_{\{i_4=i_5\}}
\zeta_{j_2}^{(i_2)}
\zeta_{j_3}^{(i_3)}+
$$
$$
+
{\bf 1}_{\{j_6=j_1\}}
{\bf 1}_{\{i_6=i_1\}}
{\bf 1}_{\{j_2=j_3\}}
{\bf 1}_{\{i_2=i_3\}}
\zeta_{j_4}^{(i_4)}
\zeta_{j_5}^{(i_5)}
+
{\bf 1}_{\{j_6=j_2\}}
{\bf 1}_{\{i_6=i_2\}}
{\bf 1}_{\{j_3=j_5\}}
{\bf 1}_{\{i_3=i_5\}}
\zeta_{j_1}^{(i_1)}
\zeta_{j_4}^{(i_4)}+
$$
$$
+
{\bf 1}_{\{j_6=j_2\}}
{\bf 1}_{\{i_6=i_2\}}
{\bf 1}_{\{j_4=j_5\}}
{\bf 1}_{\{i_4=i_5\}}
\zeta_{j_1}^{(i_1)}
\zeta_{j_3}^{(i_3)}
+
{\bf 1}_{\{j_6=j_2\}}
{\bf 1}_{\{i_6=i_2\}}
{\bf 1}_{\{j_3=j_4\}}
{\bf 1}_{\{i_3=i_4\}}
\zeta_{j_1}^{(i_1)}
\zeta_{j_5}^{(i_5)}+
$$
$$
+
{\bf 1}_{\{j_6=j_2\}}
{\bf 1}_{\{i_6=i_2\}}
{\bf 1}_{\{j_1=j_5\}}
{\bf 1}_{\{i_1=i_5\}}
\zeta_{j_3}^{(i_3)}
\zeta_{j_4}^{(i_4)}
+
{\bf 1}_{\{j_6=j_2\}}
{\bf 1}_{\{i_6=i_2\}}
{\bf 1}_{\{j_1=j_4\}}
{\bf 1}_{\{i_1=i_4\}}
\zeta_{j_3}^{(i_3)}
\zeta_{j_5}^{(i_5)}+
$$
$$
+
{\bf 1}_{\{j_6=j_2\}}
{\bf 1}_{\{i_6=i_2\}}
{\bf 1}_{\{j_1=j_3\}}
{\bf 1}_{\{i_1=i_3\}}
\zeta_{j_4}^{(i_4)}
\zeta_{j_5}^{(i_5)}
+
{\bf 1}_{\{j_6=j_3\}}
{\bf 1}_{\{i_6=i_3\}}
{\bf 1}_{\{j_2=j_5\}}
{\bf 1}_{\{i_2=i_5\}}
\zeta_{j_1}^{(i_1)}
\zeta_{j_4}^{(i_4)}+
$$
$$
+
{\bf 1}_{\{j_6=j_3\}}
{\bf 1}_{\{i_6=i_3\}}
{\bf 1}_{\{j_4=j_5\}}
{\bf 1}_{\{i_4=i_5\}}
\zeta_{j_1}^{(i_1)}
\zeta_{j_2}^{(i_2)}
+
{\bf 1}_{\{j_6=j_3\}}
{\bf 1}_{\{i_6=i_3\}}
{\bf 1}_{\{j_2=j_4\}}
{\bf 1}_{\{i_2=i_4\}}
\zeta_{j_1}^{(i_1)}
\zeta_{j_5}^{(i_5)}+
$$
$$
+
{\bf 1}_{\{j_6=j_3\}}
{\bf 1}_{\{i_6=i_3\}}
{\bf 1}_{\{j_1=j_5\}}
{\bf 1}_{\{i_1=i_5\}}
\zeta_{j_2}^{(i_2)}
\zeta_{j_4}^{(i_4)}
+
{\bf 1}_{\{j_6=j_3\}}
{\bf 1}_{\{i_6=i_3\}}
{\bf 1}_{\{j_1=j_4\}}
{\bf 1}_{\{i_1=i_4\}}
\zeta_{j_2}^{(i_2)}
\zeta_{j_5}^{(i_5)}+
$$
$$
+
{\bf 1}_{\{j_6=j_3\}}
{\bf 1}_{\{i_6=i_3\}}
{\bf 1}_{\{j_1=j_2\}}
{\bf 1}_{\{i_1=i_2\}}
\zeta_{j_4}^{(i_4)}
\zeta_{j_5}^{(i_5)}
+
{\bf 1}_{\{j_6=j_4\}}
{\bf 1}_{\{i_6=i_4\}}
{\bf 1}_{\{j_3=j_5\}}
{\bf 1}_{\{i_3=i_5\}}
\zeta_{j_1}^{(i_1)}
\zeta_{j_2}^{(i_2)}+
$$
$$
+
{\bf 1}_{\{j_6=j_4\}}
{\bf 1}_{\{i_6=i_4\}}
{\bf 1}_{\{j_2=j_5\}}
{\bf 1}_{\{i_2=i_5\}}
\zeta_{j_1}^{(i_1)}
\zeta_{j_3}^{(i_3)}
+
{\bf 1}_{\{j_6=j_4\}}
{\bf 1}_{\{i_6=i_4\}}
{\bf 1}_{\{j_2=j_3\}}
{\bf 1}_{\{i_2=i_3\}}
\zeta_{j_1}^{(i_1)}
\zeta_{j_5}^{(i_5)}+
$$
$$
+
{\bf 1}_{\{j_6=j_4\}}
{\bf 1}_{\{i_6=i_4\}}
{\bf 1}_{\{j_1=j_5\}}
{\bf 1}_{\{i_1=i_5\}}
\zeta_{j_2}^{(i_2)}
\zeta_{j_3}^{(i_3)}
+
{\bf 1}_{\{j_6=j_4\}}
{\bf 1}_{\{i_6=i_4\}}
{\bf 1}_{\{j_1=j_3\}}
{\bf 1}_{\{i_1=i_3\}}
\zeta_{j_2}^{(i_2)}
\zeta_{j_5}^{(i_5)}+
$$
$$
+
{\bf 1}_{\{j_6=j_4\}}
{\bf 1}_{\{i_6=i_4\}}
{\bf 1}_{\{j_1=j_2\}}
{\bf 1}_{\{i_1=i_2\}}
\zeta_{j_3}^{(i_3)}
\zeta_{j_5}^{(i_5)}
+
{\bf 1}_{\{j_6=j_5\}}
{\bf 1}_{\{i_6=i_5\}}
{\bf 1}_{\{j_3=j_4\}}
{\bf 1}_{\{i_3=i_4\}}
\zeta_{j_1}^{(i_1)}
\zeta_{j_2}^{(i_2)}+
$$
$$
+
{\bf 1}_{\{j_6=j_5\}}
{\bf 1}_{\{i_6=i_5\}}
{\bf 1}_{\{j_2=j_4\}}
{\bf 1}_{\{i_2=i_4\}}
\zeta_{j_1}^{(i_1)}
\zeta_{j_3}^{(i_3)}
+
{\bf 1}_{\{j_6=j_5\}}
{\bf 1}_{\{i_6=i_5\}}
{\bf 1}_{\{j_2=j_3\}}
{\bf 1}_{\{i_2=i_3\}}
\zeta_{j_1}^{(i_1)}
\zeta_{j_4}^{(i_4)}+
$$
$$
+
{\bf 1}_{\{j_6=j_5\}}
{\bf 1}_{\{i_6=i_5\}}
{\bf 1}_{\{j_1=j_4\}}
{\bf 1}_{\{i_1=i_4\}}
\zeta_{j_2}^{(i_2)}
\zeta_{j_3}^{(i_3)}
+
{\bf 1}_{\{j_6=j_5\}}
{\bf 1}_{\{i_6=i_5\}}
{\bf 1}_{\{j_1=j_3\}}
{\bf 1}_{\{i_1=i_3\}}
\zeta_{j_2}^{(i_2)}
\zeta_{j_4}^{(i_4)}+
$$
$$
+
{\bf 1}_{\{j_6=j_5\}}
{\bf 1}_{\{i_6=i_5\}}
{\bf 1}_{\{j_1=j_2\}}
{\bf 1}_{\{i_1=i_2\}}
\zeta_{j_3}^{(i_3)}
\zeta_{j_4}^{(i_4)}-
$$
$$
-
{\bf 1}_{\{j_6=j_1\}}
{\bf 1}_{\{i_6=i_1\}}
{\bf 1}_{\{j_2=j_5\}}
{\bf 1}_{\{i_2=i_5\}}
{\bf 1}_{\{j_3=j_4\}}
{\bf 1}_{\{i_3=i_4\}}-
$$
$$
-
{\bf 1}_{\{j_6=j_1\}}
{\bf 1}_{\{i_6=i_1\}}
{\bf 1}_{\{j_2=j_4\}}
{\bf 1}_{\{i_2=i_4\}}
{\bf 1}_{\{j_3=j_5\}}
{\bf 1}_{\{i_3=i_5\}}-
$$
$$
-
{\bf 1}_{\{j_6=j_1\}}
{\bf 1}_{\{i_6=i_1\}}
{\bf 1}_{\{j_2=j_3\}}
{\bf 1}_{\{i_2=i_3\}}
{\bf 1}_{\{j_4=j_5\}}
{\bf 1}_{\{i_4=i_5\}}-
$$
$$
-               
{\bf 1}_{\{j_6=j_2\}}
{\bf 1}_{\{i_6=i_2\}}
{\bf 1}_{\{j_1=j_5\}}
{\bf 1}_{\{i_1=i_5\}}
{\bf 1}_{\{j_3=j_4\}}
{\bf 1}_{\{i_3=i_4\}}-
$$
$$
-
{\bf 1}_{\{j_6=j_2\}}
{\bf 1}_{\{i_6=i_2\}}
{\bf 1}_{\{j_1=j_4\}}
{\bf 1}_{\{i_1=i_4\}}
{\bf 1}_{\{j_3=j_5\}}
{\bf 1}_{\{i_3=i_5\}}-
$$
$$
-
{\bf 1}_{\{j_6=j_2\}}
{\bf 1}_{\{i_6=i_2\}}
{\bf 1}_{\{j_1=j_3\}}
{\bf 1}_{\{i_1=i_3\}}
{\bf 1}_{\{j_4=j_5\}}
{\bf 1}_{\{i_4=i_5\}}-
$$
$$
-
{\bf 1}_{\{j_6=j_3\}}
{\bf 1}_{\{i_6=i_3\}}
{\bf 1}_{\{j_1=j_5\}}
{\bf 1}_{\{i_1=i_5\}}
{\bf 1}_{\{j_2=j_4\}}
{\bf 1}_{\{i_2=i_4\}}-
$$
$$
-
{\bf 1}_{\{j_6=j_3\}}
{\bf 1}_{\{i_6=i_3\}}
{\bf 1}_{\{j_1=j_4\}}
{\bf 1}_{\{i_1=i_4\}}
{\bf 1}_{\{j_2=j_5\}}
{\bf 1}_{\{i_2=i_5\}}-
$$
$$
-
{\bf 1}_{\{j_3=j_6\}}
{\bf 1}_{\{i_3=i_6\}}
{\bf 1}_{\{j_1=j_2\}}
{\bf 1}_{\{i_1=i_2\}}
{\bf 1}_{\{j_4=j_5\}}
{\bf 1}_{\{i_4=i_5\}}-
$$
$$
-
{\bf 1}_{\{j_6=j_4\}}
{\bf 1}_{\{i_6=i_4\}}
{\bf 1}_{\{j_1=j_5\}}
{\bf 1}_{\{i_1=i_5\}}
{\bf 1}_{\{j_2=j_3\}}
{\bf 1}_{\{i_2=i_3\}}-
$$
$$
-
{\bf 1}_{\{j_6=j_4\}}
{\bf 1}_{\{i_6=i_4\}}
{\bf 1}_{\{j_1=j_3\}}
{\bf 1}_{\{i_1=i_3\}}
{\bf 1}_{\{j_2=j_5\}}
{\bf 1}_{\{i_2=i_5\}}-
$$
$$
-
{\bf 1}_{\{j_6=j_4\}}
{\bf 1}_{\{i_6=i_4\}}
{\bf 1}_{\{j_1=j_2\}}
{\bf 1}_{\{i_1=i_2\}}
{\bf 1}_{\{j_3=j_5\}}
{\bf 1}_{\{i_3=i_5\}}-
$$
$$
-
{\bf 1}_{\{j_6=j_5\}}
{\bf 1}_{\{i_6=i_5\}}
{\bf 1}_{\{j_1=j_4\}}
{\bf 1}_{\{i_1=i_4\}}
{\bf 1}_{\{j_2=j_3\}}
{\bf 1}_{\{i_2=i_3\}}-
$$
$$
-
{\bf 1}_{\{j_6=j_5\}}
{\bf 1}_{\{i_6=i_5\}}
{\bf 1}_{\{j_1=j_2\}}
{\bf 1}_{\{i_1=i_2\}}
{\bf 1}_{\{j_3=j_4\}}
{\bf 1}_{\{i_3=i_4\}}-
$$
$$
\Biggl.-
{\bf 1}_{\{j_6=j_5\}}
{\bf 1}_{\{i_6=i_5\}}
{\bf 1}_{\{j_1=j_3\}}
{\bf 1}_{\{i_1=i_3\}}
{\bf 1}_{\{j_2=j_4\}}
{\bf 1}_{\{i_2=i_4\}}\Biggr),
$$

\vspace{6mm}

$$         
I_{(000000)T,t}^{(i_1i_1i_1i_1i_1i_1)}=
\frac{1}{720}(T-t)^{3}
\left(\left(\zeta_0^{(i_1)}\right)^6-
15\left(\zeta_0^{(i_1)}\right)^4+45\left(\zeta_0^{(i_1)}\right)^2-
15\right)\ \ \ 
\hbox{w.\ p.\ 1},
$$

\vspace{3mm}
$$
I_{(000000)T,t}^{*(i_1i_1i_1i_1i_1i_1)}=
\frac{1}{720}(T-t)^{3}\left(\zeta_0^{(i_1)}\right)^6\ \ \ \hbox{w.\ p.\ 1},
$$

\vspace{3mm}
\noindent
where

$$
C_{j_6j_5j_4 j_3 j_2 j_1}
=\frac{\sqrt{(2j_1+1)(2j_2+1)(2j_3+1)
(2j_4+1)(2j_5+1)(2j_6+1)}}{64}(T-t)^{3}\bar
C_{j_6j_5j_4 j_3 j_2 j_1},
$$

\vspace{3mm}

$$
\bar C_{j_6j_5j_4 j_3 j_2 j_1}=
\int\limits_{-1}^{1}P_{j_6}(w)
\int\limits_{-1}^{w}P_{j_5}(v)
\int\limits_{-1}^{v}P_{j_4}(u)
\int\limits_{-1}^{u}P_{j_3}(z)
\int\limits_{-1}^{z}P_{j_2}(y)
\int\limits_{-1}^{y}
P_{j_1}(x)dx dy dz du dv dw.
$$

\vspace{5mm}

It should be noted that instead of the expansion (\ref{good1})
we may to consider the following expansion, which is derived by direct 
calculation

\vspace{1mm}
$$
I_{(000)T,t}^{*(i_1 i_2 i_3)}=-\frac{1}{T-t}\left(
I_{(0)T,t}^{(i_3)}I_{(10)T,t}^{*(i_2 i_1)}+
I_{(0)T,t}^{(i_1)}I_{(10)T,t}^{*(i_2 i_3)}\right)+
\frac{1}{2}I_{(0)T,t}^{(i_3)}\left(
I_{(00)T,t}^{*(i_1 i_2)}-I_{(00)T,t}^{*(i_2 i_1)}\right)-
$$

\vspace{1mm}
\begin{equation}
\label{4004.ii}
-(T-t)^{3/2}\left(\frac{1}{6}\zeta_0^{(i_1)}\zeta_0^{(i_3)}
\left(\zeta_0^{(i_2)}+\sqrt{3}\zeta_1^{(i_2)}-\frac{1}{\sqrt{5}}
\zeta_2^{(i_2)}\right)
+\frac{1}{4}D^{(i_1 i_2 i_3)}_{T,t}\right),
\end{equation}

\vspace{5mm}
\noindent
where

$$
D^{(i_1 i_2 i_3)}_{T,t}=
\sum\limits_{\stackrel{i=1,\ j=0,\ k=i}{{}_{2i\ge k+i-j\ge -2;\
k+i-j\ - {\rm even}}}}^{\infty}
N_{ijk}K_{i+1,k+1,\frac{k+i-j}{2}+1}\zeta_i^{(i_1)}\zeta_j^{(i_2)}
\zeta_k^{(i_3)}+
$$
$$
+\sum\limits_{\stackrel{i=1,\ j=0,\ 1\le k\le i-1}{{}_{2k\ge k+i-j\ge -2;\
k+i-j\ - {\rm even}}}}^{\infty}
N_{ijk}K_{k+1,i+1,\frac{k+i-j}{2}+1}\zeta_i^{(i_1)}\zeta_j^{(i_2)}
\zeta_k^{(i_3)}-
$$
$$
-\sum\limits_{\stackrel{i=1,\ j=0,\ k=i+2}{{}_{2i+2\ge k+i-j\ge 0;\ k+i-j\ 
- {\rm even}}}}^{\infty}
N_{ijk}K_{i+1,k-1,\frac{k+i-j}{2}}\zeta_i^{(i_1)}\zeta_j^{(i_2)}
\zeta_k^{(i_3)}-
$$
$$
-\sum\limits_{\stackrel{i=1,\ j=0,\ 1\le k\le i+1}
{{}_{2k-2\ge k+i-j\ge 0;\ k+i-j\ - {\rm even}}}}^{\infty}
N_{ijk}K_{k-1,i+1,\frac{k+i-j}{2}}\zeta_i^{(i_1)}\zeta_j^{(i_2)}
\zeta_k^{(i_3)}-
$$
$$
-\sum\limits_{\stackrel{i=1,\ j=0,\ k=i-2, k\ge 1}
{{}_{2i-2\ge k+i-j\ge 0;\ k+i-j\ - {\rm even}}}}^{\infty}
N_{ijk}K_{i-1,k+1,\frac{k+i-j}{2}}\zeta_i^{(i_1)}\zeta_j^{(i_2)}
\zeta_k^{(i_3)}-
$$
$$
-\sum\limits_{\stackrel{i=1,\ j=0,\ 1\le k \le i-3} 
{{}_{2k+2\ge k+i-j\ge 0;\ k+i-j\ - {\rm even}}}}^{\infty}
N_{ijk}K_{k+1,i-1,\frac{k+i-j}{2}}\zeta_i^{(i_1)}\zeta_j^{(i_2)}
\zeta_k^{(i_3)}+
$$
$$
+\sum\limits_{\stackrel{i=1,\ j=0,\ k=i}
{{}_{2i\ge k+i-j\ge 2;\ k+i-j\ - {\rm even}}}}^{\infty}
N_{ijk}K_{i-1,k-1,\frac{k+i-j}{2}-1}\zeta_i^{(i_1)}\zeta_j^{(i_2)}
\zeta_k^{(i_3)}+
$$
$$
+\sum\limits_{\stackrel{i=1,\ j=0\ 1\le k\le i-1}
{{}_{2k\ge k+i-j\ge 2;\ k+i-j\ - {\rm even}}}}^{\infty}
N_{ijk}K_{k-1,i-1,\frac{k+i-j}{2}-1}\zeta_i^{(i_1)}\zeta_j^{(i_2)}
\zeta_k^{(i_3)},
$$

\vspace{3mm}
\noindent
where

$$
N_{ijk}=\sqrt{\frac{1}{(2k+1)(2j+1)(2i+1)}}\ ,
$$

\vspace{2mm}
$$
K_{m,n,k}=\frac{a_{m-k}a_k a_{n-k}}{a_{m+n-k}}\cdot
\frac{2n+2m-4k+1}{2n+2m-2k+1},\ \ \ a_k=\frac{(2k-1)!!}{k!},\ \ \ m\le n.
$$

\vspace{5mm}

However, as we will see further, the expansion (\ref{zzz1})
is more convenient for practical implementation then
(\ref{4004.ii}).

Also note the following relation between iterated
Ito and Stratonovich stochastic integrals

$$
I_{(0000)T,t}^{(i_1 i_2 i_3 i_4)}=
I_{(0000)T,t}^{*(i_1 i_2 i_3 i_4)}+
\frac{1}{2}{\bf 1}_{\{i_1=i_2\}}I_{(10)T,t}^{*(i_3 i_4)}-
\frac{1}{2}{\bf 1}_{\{i_2=i_3\}}\biggl(
I_{(10)T,t}^{*(i_1 i_4)}-
I_{(01)T,t}^{*(i_1 i_4)}\biggr)-
$$

\vspace{1mm}
$$
-
\frac{1}{2}{\bf 1}_{\{i_3=i_4\}}\biggl(
(T-t) I_{(00)T,t}^{*(i_1 i_2)}+
I_{(01)T,t}^{*(i_1 i_2)}\biggr)
+\frac{1}{8}(T-t)^2{\bf 1}_{\{i_1=i_2\}}{\bf 1}_{\{i_3=i_4\}}\ \ \
\hbox{w.\ p.\ 1}.
$$

\vspace{5mm}

Let 
$$
I_{(l_1\ldots l_k)T,t}^{(i_1\ldots i_k)q},\ \ \
I_{(l_1\ldots l_k)T,t}^{*(i_1\ldots i_k)q}
$$

\vspace{3mm}
\noindent
be approximations of the iterated Ito and Stratonovich
stochastic integrals 

$$
I_{(l_1\ldots l_k)T,t}^{(i_1\ldots i_k)},\ \ \
I_{(l_1\ldots l_k)T,t}^{*(i_1\ldots i_k)}
$$ 

\vspace{3mm}
\noindent
defined by 
(\ref{k1000}), (\ref{k1001}), i.e. we replace
$\infty$ with $q$ in the expansions of these stochastic integrals.
For example, $I_{(00)T,t}^{*(i_1 i_2)q}$ be the
approximation
of the iterated 
Stratonovich stochastic integral $I_{(00)T,t}^{*(i_1 i_2)}$ obtained from 
(\ref{4004}) by replacing $\infty$ with $q$, etc.

It is easy to prove that

\begin{equation}
\label{fff09}
{\sf M}\biggl\{\left(I_{(00)T,t}^{*(i_1 i_2)}-
I_{(00)T,t}^{*(i_1 i_2)q}
\right)^2\biggr\}
=\frac{(T-t)^2}{2}\left(\frac{1}{2}-\sum_{i=1}^q
\frac{1}{4i^2-1}\right)\ \ \ (i_1\ne i_2).
\end{equation}

\vspace{4mm}

Moreover, using Theorem 8, we obtain for $i_1\ne i_2$

$$
{\sf M}\biggl\{\left(I_{(10)T,t}^{*(i_1 i_2)}-I_{(10)T,t}^{*(i_1 i_2)q}
\right)^2\biggr\}=
{\sf M}\biggl\{\left(I_{(01)T,t}^{*(i_1 i_2)}-
I_{(01)T,t}^{*(i_1 i_2)q}\right)^2\biggr\}=
$$

\vspace{2mm}
$$
=\frac{(T-t)^4}{16}\left(\frac{5}{9}-
2\sum_{i=2}^q\frac{1}{4i^2-1}-
\sum_{i=1}^q
\frac{1}{(2i-1)^2(2i+3)^2}
-\right.
$$

\vspace{2mm}
\begin{equation}
\label{fff09xxx}
\left.-\sum_{i=0}^q\frac{(i+2)^2+(i+1)^2}{(2i+1)(2i+5)(2i+3)^2}
\right).
\end{equation}

\vspace{6mm}

For the case $i_1=i_2$,
using Theorem 8, we have

$$
{\sf M}\biggl\{\left(I_{(10)T,t}^{(i_1 i_1)}-
I_{(10)T,t}^{(i_1 i_1)q}
\right)^2\biggr\}=
{\sf M}\biggl\{\left(I_{(01)T,t}^{(i_1 i_1)}-
I_{(01)T,t}^{(i_1 i_1)q}\right)^2\biggr\}=
$$

\vspace{2mm}
\begin{equation}
\label{fff09xxxx}
=\frac{(T-t)^4}{16}\left(\frac{1}{9}-
\sum_{i=0}^{q}
\frac{1}{(2i+1)(2i+5)(2i+3)^2}
-2\sum_{i=1}^{q}
\frac{1}{(2i-1)^2(2i+3)^2}\right).
\end{equation}

\vspace{5mm}

In Tables 1--3 we have calculations according to the formulas 
(\ref{fff09})--(\ref{fff09xxxx}) for various values of $q.$ 
In the given tables $\varepsilon$
means the right-hand sides of these formulas.

Let us consider (\ref{ud111}), (\ref{4006}) for $i_1=i_2$

$$
I_{(01)T,t}^{*(i_1 i_1)}
=-\frac{(T-t)^2}{4}\Biggl(
\left(\zeta_0^{(i_1)}\right)^2+
\frac{1}{\sqrt{3}}\zeta_0^{(i_1)}\zeta_1^{(i_1)}+\Biggr.
$$

\vspace{1mm}
\begin{equation}
\label{leto1000}
\Biggl.
+\sum_{i=0}^{\infty}\Biggl(\frac{1}{\sqrt{(2i+1)(2i+5)}(2i+3)}
\zeta_i^{(i_1)}\zeta_{i+2}^{(i_1)}-
\frac{1}{(2i-1)(2i+3)}\left(\zeta_i^{(i_1)}\right)^2
\Biggr)\Biggr),
\end{equation}

\vspace{5mm}

$$
I_{(10)T,t}^{*(i_1 i_1)}
=-\frac{(T-t)^2}{4}\Biggl(
\left(\zeta_0^{(i_1)}\right)^2+
\frac{1}{\sqrt{3}}\zeta_0^{(i_1)}\zeta_1^{(i_1)}+\Biggr.
$$

\vspace{1mm}
\begin{equation}
\label{leto1001}
\Biggl.
+\sum_{i=0}^{\infty}\Biggl(-\frac{1}{\sqrt{(2i+1)(2i+5)}(2i+3)}
\zeta_i^{(i_1)}\zeta_{i+2}^{(i_1)}+
\frac{1}{(2i-1)(2i+3)}\left(\zeta_i^{(i_1)}\right)^2
\Biggr)\Biggr).
\end{equation}

\vspace{5mm}

From (\ref{leto1000}), (\ref{leto1001}), 
considering (\ref{4001}) and (\ref{4002}), we obtain

\begin{equation}
\label{leto1002}
I_{(10)T,t}^{*(i_1 i_1)}+I_{(01)T,t}^{*(i_1 i_1)}=
-\frac{(T-t)^2}{2}\left(\left(\zeta_0^{(i_1)}\right)^2+
\frac{1}{\sqrt{3}}\zeta_0^{(i_1)}\zeta_1^{(i_1)}\right)=
I_{(0)T,t}^{(i_1)}I_{(1)T,t}^{(i_1)}\ \ \ \hbox{w.\ p.\ 1.}
\end{equation}

\vspace{2mm}

\begin{table}
\centering
\caption{Confirmation of the formula (\ref{fff09})}      
\label{tab:1}      

\begin{tabular}{p{2.3cm}p{1.7cm}p{1.7cm}p{1.7cm}p{2.3cm}p{2.3cm}p{2.3cm}}

\hline\noalign{\smallskip}

$2\varepsilon/(T-t)^2$&0.1667&0.0238&0.0025&$2.4988\cdot 10^{-4}$&$2.4999\cdot 10^{-5}$\\

\noalign{\smallskip}\hline\noalign{\smallskip}

$q$&1&10&100&1000&10000\\

\noalign{\smallskip}\hline\noalign{\smallskip}
\end{tabular}
\vspace{3mm}
\end{table}

\begin{table}
\centering
\caption{Confirmation of the formula (\ref{fff09xxx})}
\label{tab:2}      

\begin{tabular}{p{2.3cm}p{1.7cm}p{1.7cm}p{1.7cm}p{2.3cm}p{2.3cm}p{2.3cm}}

\hline\noalign{\smallskip}

$16\varepsilon/(T-t)^4$&0.3797&0.0581&0.0062&$6.2450\cdot 10^{-4}$&$6.2495\cdot 10^{-5}$\\

\noalign{\smallskip}\hline\noalign{\smallskip}

$q$&1&10&100&1000&10000\\

\noalign{\smallskip}\hline\noalign{\smallskip}
\end{tabular}
\vspace{3mm}
\end{table}

\begin{table}
\centering
\caption{Confirmation of the formula (\ref{fff09xxxx})}
\label{tab:3}      

\begin{tabular}{p{2.1cm}p{1.2cm}p{2.1cm}p{2.1cm}p{2.3cm}p{2.3cm}p{2.3cm}}

\hline\noalign{\smallskip}

$16\varepsilon/(T-t)^4$&0.0070&$4.3551\cdot 10^{-5}$&$6.0076\cdot 10^{-8}$&$6.2251\cdot 10^{-11}$&$6.3178\cdot 10^{-14}$\\

\noalign{\smallskip}\hline\noalign{\smallskip}

$q$&1&10&100&1000&10000\\

\noalign{\smallskip}\hline\noalign{\smallskip}
\end{tabular}
\vspace{3mm}
\end{table}

Obtaining (\ref{leto1002}), we supposed that the formulas
(\ref{ud111}), (\ref{4006})
are valid w.~p.~1. The complete proof of this fact will 
be given in Sect.~5, 6.

Applying
the Ito formula 
and standard relations between 
iterated Ito and 
Stratonovich stochastic integrals, 
it is easy to get the equality (\ref{leto1002}).

Furthermore, using the Ito formula, we obtain

\begin{equation}
\label{leto1010}
I_{(11)T,t}^{*(i_1 i_1)}=\frac{\left(I_{(1)T,t}^{(i_1)}
\right)^2}{2}\ \ \ \hbox{w.\ p.\ 1.}
\end{equation}

\vspace{3mm}

In addition, applying the Ito formula, we have

\vspace{-1mm}
\begin{equation}
\label{tqtq}
I_{(20)T,t}^{(i_1 i_1)}+I_{(02)T,t}^{(i_1 i_1)}=
I_{(0)T,t}^{(i_1)}I_{(2)T,t}^{(i_1)}-
\frac{(T-t)^3}{3}\ \ \ \hbox{w.\ p.\ 1.}
\end{equation}

\vspace{4mm}

From (\ref{tqtq}), considering the formulas (\ref{seg1}), (\ref{seg2}), 
we get

\vspace{-1mm}
\begin{equation}
\label{leto1012}
I_{(20)T,t}^{*(i_1 i_1)}+I_{(02)T,t}^{*(i_1 i_1)}=
I_{(0)T,t}^{(i_1)}I_{(2)T,t}^{(i_1)}\ \ \ \hbox{w.\ p.\ 1.}
\end{equation}

\vspace{4mm}

Let us check whether the formulas (\ref{leto1010}), (\ref{leto1012}) 
follow
from (\ref{leto1})--(\ref{leto3}), if we suppose $i_1=i_2$ in the last ones.
From (\ref{leto1})--(\ref{leto3}) for
$i_1=i_2$ we obtain

$$
I_{(20)T,t}^{*(i_1 i_1)}+I_{(02)T,t}^{*(i_1 i_1)}=
-\frac{(T-t)^2}{2}I_{(00)T,t}^{*(i_1 i_1)}
-(T-t)\left(I_{(10)T,t}^{*(i_1 i_1)}+I_{(01)T,t}^{*(i_1 i_1)}\right)+
$$

\begin{equation}
\label{leto1020}
+\frac{(T-t)^3}{4}
\left(\frac{1}{3}\left(\zeta_0^{(i_1)}\right)^2+
\frac{2}{3\sqrt{5}}\zeta_2^{(i_1)}\zeta_0^{(i_1)}\right),
\end{equation}

\vspace{3mm}

\begin{equation}
\label{leto1021}
I_{(11)T,t}^{*(i_1 i_1)}=
-\frac{(T-t)^2}{4}I_{(00)T,t}^{*(i_1 i_1)}
-\frac{T-t}{2}\left(I_{(10)T,t}^{*(i_1 i_1)}+I_{(01)T,t}^{*(i_1 i_1)}\right)
+\frac{(T-t)^3}{24}
\left(\zeta_1^{(i_1)}\right)^2.
\end{equation}

\vspace{5mm}

It is easy to see that from (\ref{leto1020}) and (\ref{leto1021}), 
considering (\ref{leto1002}) and (\ref{4001})--(\ref{4004}), 
we actually obtain the equalities  
(\ref{leto1010}) and (\ref{leto1012}), and it indirectly confirm 
the correctness of 
the formulas (\ref{leto1})--(\ref{leto3}).

Obtaining (\ref{leto1010}), (\ref{leto1012}), we supposed that the formulas
(\ref{leto1})--(\ref{leto3})
are valid w.~p.~1. The complete proof of this fact will 
be given in Sect.~5, 6.

On the basis of 
the presented 
expansions of 
iterated stochastic integrals we 
can see that increasing of multiplicities of these integrals 
or degree indexes of their weight functions 
leads
to noticeable complication of formulas 
for mentioned expansions. 

However, increasing of mentioned parameters leads to increasing 
of orders of smallness with respect to $T-t$ in the mean-square sense 
for iterated stochastic integrals that leads to a sharp decrease  
of member 
quantities
in expansions of iterated stochastic 
integrals, which are required for achieving the acceptable accuracy
of approximation. In this context, let us consider the approach 
to the approximation of iterated stochastic integrals, which 
provides a possibility to obtain the mean-square approximations of 
the required accuracy without using the 
complex expansions like (\ref{4004.ii}).

Let us consider the following approximation of iterated Ito stochastic integral 
$I_{(000)T,t}^{(i_1i_2i_3)}$ using (\ref{zzz1})

$$
I_{(000)T,t}^{(i_1i_2i_3)q_1}
=\sum_{j_1,j_2,j_3=0}^{q_1}
C_{j_3j_2j_1}\Biggl(
\zeta_{j_1}^{(i_1)}\zeta_{j_2}^{(i_2)}\zeta_{j_3}^{(i_3)}
-{\bf 1}_{\{i_1=i_2\}}
{\bf 1}_{\{j_1=j_2\}}
\zeta_{j_3}^{(i_3)}-
\Biggr.
$$
\begin{equation}
\label{sad001}
\Biggl.
-{\bf 1}_{\{i_2=i_3\}}
{\bf 1}_{\{j_2=j_3\}}
\zeta_{j_1}^{(i_1)}-
{\bf 1}_{\{i_1=i_3\}}
{\bf 1}_{\{j_1=j_3\}}
\zeta_{j_2}^{(i_2)}\Biggr),
\end{equation}

\vspace{3mm}
\noindent
where $C_{j_3j_2j_1}$ is defined by (\ref{zzz2}), (\ref{zzz3xx}).

In particular, from (\ref{sad001}) for 
$i_1\ne i_2$, 
$i_2\ne i_3$, $i_1\ne i_3$
we obtain

\begin{equation}
\label{38}
I_{(000)T,t}^{(i_1i_2i_3)q_1}=
\sum_{j_1,j_2,j_3=0}^{q_1}
C_{j_3j_2j_1}
\zeta_{j_1}^{(i_1)}\zeta_{j_2}^{(i_2)}\zeta_{j_3}^{(i_3)}.
\end{equation}

\vspace{3mm}

Furthermore, using Theorem 8 for $k=3$,
we get

$$
{\sf M}\left\{\left(
I_{(000)T,t}^{(i_1i_2 i_3)}-
I_{(000)T,t}^{(i_1i_2 i_3)q_1}\right)^2\right\}=
$$

\begin{equation}
\label{39}
=
\frac{(T-t)^{3}}{6}
-\sum_{j_1,j_2,j_3=0}^{q_1}
C_{j_3j_2j_1}^2\ \ \ (i_1\ne i_2, i_1\ne i_3, i_2\ne i_3),
\end{equation}

\vspace{3mm}

$$
{\sf M}\left\{\left(
I_{(000)T,t}^{(i_1i_2 i_3)}-
I_{(000)T,t}^{(i_1i_2 i_3)q_1}\right)^2\right\}=
$$

\begin{equation}
\label{39a}
=
\frac{(T-t)^{3}}{6}-\sum_{j_1,j_2,j_3=0}^{q_1}
C_{j_3j_2j_1}^2
-\sum_{j_1,j_2,j_3=0}^{q_1}
C_{j_2j_3j_1}C_{j_3j_2j_1}\ \ \ (i_1\ne i_2=i_3),
\end{equation}

\vspace{3mm}

$$
{\sf M}\left\{\left(
I_{(000)T,t}^{(i_1i_2 i_3)}-
I_{(000)T,t}^{(i_1i_2 i_3)q_1}\right)^2\right\}=
$$

\begin{equation}
\label{39b}
=
\frac{(T-t)^{3}}{6}-\sum_{j_1,j_2,j_3=0}^{q_1}
C_{j_3j_2j_1}^2
-\sum_{j_1,j_2,j_3=0}^{q_1}
C_{j_3j_2j_1}C_{j_1j_2j_3}\ \ \ (i_1=i_3\ne i_2),
\end{equation}

\vspace{3mm}

$$
{\sf M}\left\{\left(
I_{(000)T,t}^{(i_1i_2 i_3)}-
I_{(000)T,t}^{(i_1i_2 i_3)q_1}\right)^2\right\}=
$$

\begin{equation}
\label{39c}
=
\frac{(T-t)^{3}}{6}-\sum_{j_1,j_2,j_3=0}^{q_1}
C_{j_3j_2j_1}^2
-\sum_{j_1,j_2,j_3=0}^{q_1}
C_{j_3j_1j_2}C_{j_3j_2j_1}\ \ \ (i_1=i_2\ne i_3).
\end{equation}

\vspace{5mm}

From the other hand, from (\ref{star00011}) for $k=3$ we obtain

\begin{equation}
\label{leto1041}
{\sf M}\left\{\left(
I_{(000)T,t}^{(i_1i_2 i_3)}-
I_{(000)T,t}^{(i_1i_2 i_3)q_1}\right)^2\right\}\le
6\left(\frac{(T-t)^{3}}{6}-\sum_{j_1,j_2,j_3=0}^{q_1}
C_{j_3j_2j_1}^2\right),
\end{equation}

\vspace{4mm}
\noindent
where $i_1, i_2, i_3=1,\ldots,m$.

We may act similarly with more complicated 
iterated stochastic integrals. For example, for  
approximation of the stochastic integral
$I_{(0000)T,t}^{(i_1 i_2 i_3 i_4)}$ 
we can write (see (\ref{zzz10}))

\vspace{1mm}
$$
I_{(0000)T,t}^{(i_1 i_2 i_3 i_4)q_2}=
\sum_{j_1,j_2,j_3,j_4=0}^{q_2}
C_{j_4 j_3 j_2 j_1}\Biggl(
\zeta_{j_1}^{(i_1)}\zeta_{j_2}^{(i_2)}\zeta_{j_3}^{(i_3)}\zeta_{j_4}^{(i_4)}
-\Biggr.
$$
$$
-
{\bf 1}_{\{i_1=i_2\}}
{\bf 1}_{\{j_1=j_2\}}
\zeta_{j_3}^{(i_3)}
\zeta_{j_4}^{(i_4)}
-
{\bf 1}_{\{i_1=i_3\}}
{\bf 1}_{\{j_1=j_3\}}
\zeta_{j_2}^{(i_2)}
\zeta_{j_4}^{(i_4)}-
$$
$$
-
{\bf 1}_{\{i_1=i_4\}}
{\bf 1}_{\{j_1=j_4\}}
\zeta_{j_2}^{(i_2)}
\zeta_{j_3}^{(i_3)}
-
{\bf 1}_{\{i_2=i_3\}}
{\bf 1}_{\{j_2=j_3\}}
\zeta_{j_1}^{(i_1)}
\zeta_{j_4}^{(i_4)}-
$$
$$
-
{\bf 1}_{\{i_2=i_4\}}
{\bf 1}_{\{j_2=j_4\}}
\zeta_{j_1}^{(i_1)}
\zeta_{j_3}^{(i_3)}
-
{\bf 1}_{\{i_3=i_4\}}
{\bf 1}_{\{j_3=j_4\}}
\zeta_{j_1}^{(i_1)}
\zeta_{j_2}^{(i_2)}+
$$
$$
+
{\bf 1}_{\{i_1=i_2\}}
{\bf 1}_{\{j_1=j_2\}}
{\bf 1}_{\{i_3=i_4\}}
{\bf 1}_{\{j_3=j_4\}}
+
{\bf 1}_{\{i_1=i_3\}}
{\bf 1}_{\{j_1=j_3\}}
{\bf 1}_{\{i_2=i_4\}}
{\bf 1}_{\{j_2=j_4\}}+
$$
\begin{equation}
\label{res100}
+\Biggl.
{\bf 1}_{\{i_1=i_4\}}
{\bf 1}_{\{j_1=j_4\}}
{\bf 1}_{\{i_2=i_3\}}
{\bf 1}_{\{j_2=j_3\}}\Biggr),
\end{equation}

\vspace{4mm}
\noindent
where $C_{j_4 j_3 j_2 j_1}$ is defined by (\ref{zzz11}), (\ref{zzz12}).

Moreover, according to (\ref{star00011}) for $k=4$, we get

$$
{\sf M}\left\{\left(
I_{(0000)T,t}^{(i_1i_2 i_3 i_4)}-
I_{(0000)T,t}^{(i_1i_2 i_3 i_4)q_2}\right)^2\right\}\le
24\left(\frac{(T-t)^{4}}{24}-\sum_{j_1,j_2,j_3,j_4=0}^{q_2}
C_{j_4j_3j_2j_1}^2\right),
$$

\vspace{3mm}
\noindent
where
$i_1, i_2, i_3, i_4=1,\ldots,m$.

For pairwise different $i_1, i_2, i_3, i_4=1,\ldots,m$ from 
Theorem 8 we obtain

\begin{equation}
\label{r7}
{\sf M}\left\{\left(
I_{(0000)T,t}^{(i_1i_2 i_3i_4)}-
I_{(0000)T,t}^{(i_1i_2 i_3i_4)q_2}\right)^2\right\}=
\frac{(T-t)^{4}}{24}
-\sum_{j_1,j_2,j_3,j_4=0}^{q_2}
C_{j_4j_3j_2j_1}^2.
\end{equation}

\vspace{5mm}

Using Theorem 8, we can calculate exactly the left-hand
side of (\ref{r7})
for any possible combinations
of $i_1, i_2, i_3, i_4$. These relations were obtained in 
\cite{10a}-\cite{10axx1}, \cite{15b}. For example,

\vspace{1mm}
$$
{\sf M}\left\{\left(
I_{(0000)T,t}^{(i_1i_2 i_3 i_4)}-
I_{(0000)T,t}^{(i_1i_2 i_3 i_4)q_2}\right)^2\right\}=
$$

\vspace{1mm}
$$
= \frac{(T-t)^{4}}{24} - \sum_{j_1,j_2,j_3,j_4=0}^{q_2}
C_{j_4j_3j_2j_1}\Biggl(\sum\limits_{(j_1,j_2)}\Biggl(
\sum\limits_{(j_3,j_4)}
C_{j_4j_3j_2j_1}\Biggr)\Biggr)\ \ \ (i_1=i_2\ne i_3=i_4),
$$

\vspace{6mm}

$$
{\sf M}\left\{\left(
I_{(0000)T,t}^{(i_1i_2 i_3 i_4)}-
I_{(0000)T,t}^{(i_1i_2 i_3 i_4)q_2}\right)^2\right\}=
$$

\vspace{1mm}
$$
=\frac{(T-t)^{4}}{24} -
\sum_{j_1,\ldots,j_4=0}^{q_2}
C_{j_4\ldots j_1}\Biggl(\sum\limits_{(j_1,j_2,j_3)}
C_{j_4\ldots j_1}\Biggr)\ \ \ (i_1=i_2=i_3\ne i_4),
$$

\vspace{5mm}
\noindent
where
$$
\sum\limits_{(j_1,j_2)}
$$

\vspace{3mm}
\noindent
means the sum with respect to permutations $(j_1,j_2)$.

Assume that $q_1=6$. In Tables 4--10 we have the exact 
values of coefficients 
$\bar C_{j_3j_2j_1},$ $j_1,j_2,j_3=0, 1,\ldots,6.$
Note that in \cite{Kuz-Kuz}, \cite{Mikh-1}
the database with 270,000 exactly
calculated Fourier--Legendre coefficients was described.

Calculating the value  
(\ref{39}) for $q_1=6,$ 
$i_1\ne i_2,$ $i_1\ne i_3,$ $i_3\ne i_2$,
we obtain the following  
approximate equality

$$
{\sf M}\left\{\left(
I_{(000)T,t}^{(i_1i_2 i_3)}-
I_{(000)T,t}^{(i_1i_2 i_3)q_1}\right)^2\right\}\approx
0.01956(T-t)^3.
$$

\vspace{4mm}

Let us choose, for example, $q_2=2.$ In Tables 11--19  
we have the exact values of coefficients 
$\bar C_{j_4j_3j_2j_1}$
($j_1,j_2,j_3,j_4=0, 1, 2$).  
In the case of pairwise different
$i_1, i_2, i_3, i_4$ we have from (\ref{r7}) the following  
approximate equality

\begin{equation}
\label{46000}
{\sf M}\left\{\left(
I_{(0000)T,t}^{(i_1i_2i_3 i_4)}-
I_{(0000)T,t}^{(i_1i_2i_3 i_4)q_2}\right)^2\right\}\approx
0.0236084(T-t)^4.
\end{equation}

\vspace{4mm}

Let us consider the following four approximations of iterated Ito 
stochastic integrals (see (\ref{sss1})--(\ref{sss4}))

\vspace{2mm}
$$
I_{(001)T,t}^{(i_1i_2i_3)q_3}
=
\sum_{j_1,j_2,j_3=0}^{q_3}
C_{j_3j_2j_1}^{001}\Biggl(
\zeta_{j_1}^{(i_1)}\zeta_{j_2}^{(i_2)}\zeta_{j_3}^{(i_3)}
-{\bf 1}_{\{i_1=i_2\}}
{\bf 1}_{\{j_1=j_2\}}
\zeta_{j_3}^{(i_3)}-
\Biggr.
$$
\begin{equation}
\label{r9}
\Biggl.
-{\bf 1}_{\{i_2=i_3\}}
{\bf 1}_{\{j_2=j_3\}}
\zeta_{j_1}^{(i_1)}-
{\bf 1}_{\{i_1=i_3\}}
{\bf 1}_{\{j_1=j_3\}}
\zeta_{j_2}^{(i_2)}\Biggr),
\end{equation}

\vspace{5mm}

$$
I_{(010)T,t}^{(i_1i_2i_3)q_3}
=
\sum_{j_1,j_2,j_3=0}^{q_3}
C_{j_3j_2j_1}^{010}\Biggl(
\zeta_{j_1}^{(i_1)}\zeta_{j_2}^{(i_2)}\zeta_{j_3}^{(i_3)}
-{\bf 1}_{\{i_1=i_2\}}
{\bf 1}_{\{j_1=j_2\}}
\zeta_{j_3}^{(i_3)}-
\Biggr.
$$
\begin{equation}
\label{r10}
\Biggl.
-{\bf 1}_{\{i_2=i_3\}}
{\bf 1}_{\{j_2=j_3\}}
\zeta_{j_1}^{(i_1)}-
{\bf 1}_{\{i_1=i_3\}}
{\bf 1}_{\{j_1=j_3\}}
\zeta_{j_2}^{(i_2)}\Biggr),
\end{equation}

\vspace{5mm}

$$
I_{(100)T,t}^{(i_1i_2i_3)q_3}
=
\sum_{j_1,j_2,j_3=0}^{q_3}
C_{j_3j_2j_1}^{100}\Biggl(
\zeta_{j_1}^{(i_1)}\zeta_{j_2}^{(i_2)}\zeta_{j_3}^{(i_3)}
-{\bf 1}_{\{i_1=i_2\}}
{\bf 1}_{\{j_1=j_2\}}
\zeta_{j_3}^{(i_3)}-
\Biggr.
$$
\begin{equation}
\label{r10a}
\Biggl.
-{\bf 1}_{\{i_2=i_3\}}
{\bf 1}_{\{j_2=j_3\}}
\zeta_{j_1}^{(i_1)}-
{\bf 1}_{\{i_1=i_3\}}
{\bf 1}_{\{j_1=j_3\}}
\zeta_{j_2}^{(i_2)}\Biggr),
\end{equation}

\vspace{8mm}

$$
I_{(00000)T,t}^{(i_1i_2i_3i_4 i_5)q_4}=
\sum_{j_1,j_2,j_3,j_4,j_5=0}^{q_4}
C_{j_5 j_4 j_3 j_2 j_1}\Biggl(
\prod_{l=1}^5\zeta_{j_l}^{(i_l)}
-\Biggr.
$$

$$
-
{\bf 1}_{\{i_1=i_2\}}
{\bf 1}_{\{j_1=j_2\}}
\zeta_{j_3}^{(i_3)}
\zeta_{j_4}^{(i_4)}
\zeta_{j_5}^{(i_5)}-
{\bf 1}_{\{i_1=i_3\}}
{\bf 1}_{\{j_1=j_3\}}
\zeta_{j_2}^{(i_2)}
\zeta_{j_4}^{(i_4)}
\zeta_{j_5}^{(i_5)}-
$$
$$
-
{\bf 1}_{\{i_1=i_4\}}
{\bf 1}_{\{j_1=j_4\}}
\zeta_{j_2}^{(i_2)}
\zeta_{j_3}^{(i_3)}
\zeta_{j_5}^{(i_5)}-
{\bf 1}_{\{i_1=i_5\}}
{\bf 1}_{\{j_1=j_5\}}
\zeta_{j_2}^{(i_2)}
\zeta_{j_3}^{(i_3)}
\zeta_{j_4}^{(i_4)}-
$$
$$
-
{\bf 1}_{\{i_2=i_3\}}
{\bf 1}_{\{j_2=j_3\}}
\zeta_{j_1}^{(i_1)}
\zeta_{j_4}^{(i_4)}
\zeta_{j_5}^{(i_5)}-
{\bf 1}_{\{i_2=i_4\}}
{\bf 1}_{\{j_2=j_4\}}
\zeta_{j_1}^{(i_1)}
\zeta_{j_3}^{(i_3)}
\zeta_{j_5}^{(i_5)}-
$$
$$
-
{\bf 1}_{\{i_2=i_5\}}
{\bf 1}_{\{j_2=j_5\}}
\zeta_{j_1}^{(i_1)}
\zeta_{j_3}^{(i_3)}
\zeta_{j_4}^{(i_4)}
-{\bf 1}_{\{i_3=i_4\}}
{\bf 1}_{\{j_3=j_4\}}
\zeta_{j_1}^{(i_1)}
\zeta_{j_2}^{(i_2)}
\zeta_{j_5}^{(i_5)}-
$$
$$
-
{\bf 1}_{\{i_3=i_5\}}
{\bf 1}_{\{j_3=j_5\}}
\zeta_{j_1}^{(i_1)}
\zeta_{j_2}^{(i_2)}
\zeta_{j_4}^{(i_4)}
-{\bf 1}_{\{i_4=i_5\}}
{\bf 1}_{\{j_4=j_5\}}
\zeta_{j_1}^{(i_1)}
\zeta_{j_2}^{(i_2)}
\zeta_{j_3}^{(i_3)}+
$$
$$
+
{\bf 1}_{\{i_1=i_2\}}
{\bf 1}_{\{j_1=j_2\}}
{\bf 1}_{\{i_3=i_4\}}
{\bf 1}_{\{j_3=j_4\}}\zeta_{j_5}^{(i_5)}+
{\bf 1}_{\{i_1=i_2\}}
{\bf 1}_{\{j_1=j_2\}}
{\bf 1}_{\{i_3=i_5\}}
{\bf 1}_{\{j_3=j_5\}}\zeta_{j_4}^{(i_4)}+
$$
$$
+
{\bf 1}_{\{i_1=i_2\}}
{\bf 1}_{\{j_1=j_2\}}
{\bf 1}_{\{i_4=i_5\}}
{\bf 1}_{\{j_4=j_5\}}\zeta_{j_3}^{(i_3)}+
{\bf 1}_{\{i_1=i_3\}}
{\bf 1}_{\{j_1=j_3\}}
{\bf 1}_{\{i_2=i_4\}}
{\bf 1}_{\{j_2=j_4\}}\zeta_{j_5}^{(i_5)}+
$$
$$
+
{\bf 1}_{\{i_1=i_3\}}
{\bf 1}_{\{j_1=j_3\}}
{\bf 1}_{\{i_2=i_5\}}
{\bf 1}_{\{j_2=j_5\}}\zeta_{j_4}^{(i_4)}+
{\bf 1}_{\{i_1=i_3\}}
{\bf 1}_{\{j_1=j_3\}}
{\bf 1}_{\{i_4=i_5\}}
{\bf 1}_{\{j_4=j_5\}}\zeta_{j_2}^{(i_2)}+
$$
$$
+
{\bf 1}_{\{i_1=i_4\}}
{\bf 1}_{\{j_1=j_4\}}
{\bf 1}_{\{i_2=i_3\}}
{\bf 1}_{\{j_2=j_3\}}\zeta_{j_5}^{(i_5)}+
{\bf 1}_{\{i_1=i_4\}}
{\bf 1}_{\{j_1=j_4\}}
{\bf 1}_{\{i_2=i_5\}}
{\bf 1}_{\{j_2=j_5\}}\zeta_{j_3}^{(i_3)}+
$$
$$
+
{\bf 1}_{\{i_1=i_4\}}
{\bf 1}_{\{j_1=j_4\}}
{\bf 1}_{\{i_3=i_5\}}
{\bf 1}_{\{j_3=j_5\}}\zeta_{j_2}^{(i_2)}+
{\bf 1}_{\{i_1=i_5\}}
{\bf 1}_{\{j_1=j_5\}}
{\bf 1}_{\{i_2=i_3\}}
{\bf 1}_{\{j_2=j_3\}}\zeta_{j_4}^{(i_4)}+
$$
$$
+
{\bf 1}_{\{i_1=i_5\}}
{\bf 1}_{\{j_1=j_5\}}
{\bf 1}_{\{i_2=i_4\}}
{\bf 1}_{\{j_2=j_4\}}\zeta_{j_3}^{(i_3)}+
{\bf 1}_{\{i_1=i_5\}}
{\bf 1}_{\{j_1=j_5\}}
{\bf 1}_{\{i_3=i_4\}}
{\bf 1}_{\{j_3=j_4\}}\zeta_{j_2}^{(i_2)}+
$$
$$
+
{\bf 1}_{\{i_2=i_3\}}
{\bf 1}_{\{j_2=j_3\}}
{\bf 1}_{\{i_4=i_5\}}
{\bf 1}_{\{j_4=j_5\}}\zeta_{j_1}^{(i_1)}+
{\bf 1}_{\{i_2=i_4\}}
{\bf 1}_{\{j_2=j_4\}}
{\bf 1}_{\{i_3=i_5\}}
{\bf 1}_{\{j_3=j_5\}}\zeta_{j_1}^{(i_1)}+
$$
\begin{equation}
\label{r11}
+\Biggl.
{\bf 1}_{\{i_2=i_5\}}
{\bf 1}_{\{j_2=j_5\}}
{\bf 1}_{\{i_3=i_4\}}
{\bf 1}_{\{j_3=j_4\}}\zeta_{j_1}^{(i_1)}\Biggr).
\end{equation}

\begin{table}
\centering
\caption{Coefficients $\bar C_{0j_2j_1}$}
\label{tab:4}      

\begin{tabular}{p{1.3cm}p{1.3cm}p{1.3cm}p{1.3cm}p{1.3cm}p{1.3cm}p{1.3cm}p{1.3cm}}

\hline\noalign{\smallskip}

&$j_1=0$&$j_1=1$&$j_1=2$&$j_1=3$&$j_1=4$&$j_1=5$&$j_1=6$\\

\noalign{\smallskip}\hline\noalign{\smallskip}

$j_2=0$&$\frac{4}{3}$&$\frac{-2}{3}$&$\frac{2}{15}$&$0$&0&0&0\\

\noalign{\smallskip}

$j_2=1$&$0$&$\frac{2}{15}$&$\frac{-2}{15}$&$\frac{4}{105}$&0&0&0\\

\noalign{\smallskip}

$j_2=2$&$\frac{-4}{15}$&$\frac{2}{15}$&$\frac{2}{105}$&$\frac{-2}{35}$&
$\frac{2}{105}$&0&0\\

\noalign{\smallskip}

$j_2=3$&$0$&$\frac{-2}{35}$&$\frac{2}{35}$&$\frac{2}{315}$&
$\frac{-2}{63}$&$\frac{8}{693}$&0\\

\noalign{\smallskip}

$j_2=4$&$0$&$0$&$\frac{-8}{315}$&$\frac{2}{63}$&
$\frac{2}{693}$&$\frac{-2}{99}$&$\frac{10}{1287}$\\

\noalign{\smallskip}

$j_2=5$&$0$&$0$&$0$&$\frac{-10}{693}$&
$\frac{2}{99}$&$\frac{2}{1287}$&$\frac{-2}{143}$\\

\noalign{\smallskip}

$j_2=6$&$0$&$0$&$0$&$0$&$\frac{-4}{429}$&$\frac{2}{143}$&$\frac{2}{2145}$\\

\noalign{\smallskip}\hline\noalign{\smallskip}
\end{tabular}
\end{table}

\begin{table}
\centering
\caption{Coefficients $\bar C_{1j_2j_1}$}
\label{tab:5}      

\begin{tabular}{p{1.3cm}p{1.3cm}p{1.3cm}p{1.3cm}p{1.3cm}p{1.3cm}p{1.3cm}p{1.3cm}}

\hline\noalign{\smallskip}

&$j_1=0$&$j_1=1$&$j_1=2$&$j_1=3$&$j_1=4$&$j_1=5$&$j_1=6$\\

\noalign{\smallskip}\hline\noalign{\smallskip}

$j_2=0$&$\frac{2}{3}$&$\frac{-4}{15}$&$0$&$\frac{2}{105}$&0&0&0\\

\noalign{\smallskip}

$j_2=1$&$\frac{2}{15}$&$0$&$\frac{-4}{105}$&$0$&$\frac{2}{315}$&0&0\\

\noalign{\smallskip}

$j_2=2$&$\frac{-2}{15}$&$\frac{8}{105}$&$0$&$\frac{-2}{105}$&
0&$\frac{4}{1155}$&0\\

\noalign{\smallskip}

$j_2=3$&$\frac{-2}{35}$&0&$\frac{8}{315}$&0&
$\frac{-38}{3465}$&0&$\frac{20}{9009}$\\

\noalign{\smallskip}

$j_2=4$&$0$&$\frac{-4}{315}$&0&$\frac{46}{3465}$&
0&$\frac{-64}{9009}$&0\\

\noalign{\smallskip}

$j_2=5$&$0$&$0$&$\frac{-4}{693}$&0&
$\frac{74}{9009}$&0&$\frac{-32}{6435}$\\

\noalign{\smallskip}

$j_2=6$&$0$&$0$&$0$&$\frac{-10}{3003}$&$0$&$\frac{4}{715}$&$0$\\

\noalign{\smallskip}\hline\noalign{\smallskip}
\end{tabular}
\end{table}

\begin{table}
\centering
\caption{Coefficients $\bar C_{2j_2j_1}$}
\label{tab:6}      

\begin{tabular}{p{1.3cm}p{1.3cm}p{1.3cm}p{1.3cm}p{1.3cm}p{1.3cm}p{1.3cm}p{1.3cm}}

\hline\noalign{\smallskip}

&$j_1=0$&$j_1=1$&$j_1=2$&$j_1=3$&$j_1=4$&$j_1=5$&$j_1=6$\\

\noalign{\smallskip}\hline\noalign{\smallskip}

$j_2=0$&$\frac{2}{15}$&0&$\frac{-4}{105}$&0&$\frac{2}{315}$&0&0\\

\noalign{\smallskip}

$j_2=1$&$\frac{2}{15}$&$\frac{-4}{105}$&0&$\frac{-2}{315}$&0&$\frac{8}{3465}$&0\\

\noalign{\smallskip}

$j_2=2$&$\frac{2}{105}$&0&$0$&0&
$\frac{-2}{495}$&0&$\frac{4}{3003}$\\

\noalign{\smallskip}

$j_2=3$&$\frac{-2}{35}$&$\frac{8}{315}$&0&$\frac{-2}{3465}$&
0&$\frac{-116}{45045}$&0\\

\noalign{\smallskip}

$j_2=4$&$\frac{-8}{315}$&0&$\frac{4}{495}$&0&
$\frac{-2}{6435}$&0&$\frac{-16}{9009}$\\

\noalign{\smallskip}

$j_2=5$&$0$&$\frac{-4}{693}$&0&$\frac{38}{9009}$&
0&$\frac{-8}{45045}$&0\\

\noalign{\smallskip}

$j_2=6$&$0$&$0$&$\frac{-8}{3003}$&$0$&$\frac{118}{45045}$&$0$&$\frac{-4}{36465}$\\

\noalign{\smallskip}\hline\noalign{\smallskip}
\end{tabular}
\end{table}

\begin{table}
\centering
\caption{Coefficients $\bar C_{3j_2j_1}$}
\label{tab:7}      

\begin{tabular}{p{1.3cm}p{1.3cm}p{1.3cm}p{1.3cm}p{1.3cm}p{1.3cm}p{1.3cm}p{1.3cm}}

\hline\noalign{\smallskip}

&$j_1=0$&$j_1=1$&$j_1=2$&$j_1=3$&$j_1=4$&$j_1=5$&$j_1=6$\\

\noalign{\smallskip}\hline\noalign{\smallskip}

$j_2=0$&$0$&$\frac{2}{105}$&$0$&$\frac{-4}{315}$&$0$&$\frac{2}{693}$&0\\
\noalign{\smallskip}
$j_2=1$&$\frac{4}{105}$&0&$\frac{-2}{315}$&0&$\frac{-8}{3465}$&0&$\frac{10}{9009}$\\
\noalign{\smallskip}
$j_2=2$&$\frac{2}{35}$&$\frac{-2}{105}$&$0$&$\frac{4}{3465}$&
$0$&$\frac{-74}{45045}$&0\\
\noalign{\smallskip}
$j_2=3$&$\frac{2}{315}$&$0$&$\frac{-2}{3465}$&0&
$\frac{16}{45045}$&0&$\frac{-10}{9009}$\\
\noalign{\smallskip}
$j_2=4$&$\frac{-2}{63}$&$\frac{46}{3465}$&0&$\frac{-32}{45045}$&
0&$\frac{2}{9009}$&0\\
\noalign{\smallskip}
$j_2=5$&$\frac{-10}{693}$&0&$\frac{38}{9009}$&0&
$\frac{-4}{9009}$&0&$\frac{122}{765765}$\\
\noalign{\smallskip}
$j_2=6$&$0$&$\frac{-10}{3003}$&$0$&$\frac{20}{9009}$&$0$&$\frac{-226}{765765}$&$0$\\

\noalign{\smallskip}\hline\noalign{\smallskip}
\end{tabular}
\end{table}

\begin{table}
\centering
\caption{Coefficients $\bar C_{4j_2j_1}$}
\label{tab:8}      

\begin{tabular}{p{1.3cm}p{1.3cm}p{1.3cm}p{1.3cm}p{1.3cm}p{1.3cm}p{1.3cm}p{1.3cm}}

\hline\noalign{\smallskip}

&$j_1=0$&$j_1=1$&$j_1=2$&$j_1=3$&$j_1=4$&$j_1=5$&$j_1=6$\\

\noalign{\smallskip}\hline\noalign{\smallskip}

$j_2=0$&$0$&0&$\frac{2}{315}$&0&$\frac{-4}{693}$&0&$\frac{2}{1287}$\\
\noalign{\smallskip}
$j_2=1$&0&$\frac{2}{315}$&0&$\frac{-8}{3465}$&0&$\frac{-10}{9009}$&0\\
\noalign{\smallskip}
$j_2=2$&$\frac{2}{105}$&0&$\frac{-2}{495}$&0&
$\frac{4}{6435}$&0&$\frac{-38}{45045}$\\
\noalign{\smallskip}
$j_2=3$&$\frac{2}{63}$&$\frac{-38}{3465}$&0&$\frac{16}{45045}$&
0&$\frac{2}{9009}$&0\\
\noalign{\smallskip}
$j_2=4$&$\frac{2}{693}$&0&$\frac{-2}{6435}$&0&
0&0&$\frac{2}{13923}$\\
\noalign{\smallskip}
$j_2=5$&$\frac{-2}{99}$&$\frac{74}{9009}$&0&$\frac{-4}{9009}$&
0&$\frac{-2}{153153}$&0\\
\noalign{\smallskip}
$j_2=6$&$\frac{-4}{429}$&$0$&$\frac{118}{45045}$&$0$&$\frac{-4}{13923}$&$0$&
$\frac{-2}{188955}$\\
\noalign{\smallskip}\hline\noalign{\smallskip}
\end{tabular}
\end{table}

\begin{table}
\centering
\caption{Coefficients $\bar C_{5j_2j_1}$}
\label{tab:9}      

\begin{tabular}{p{1.3cm}p{1.3cm}p{1.3cm}p{1.3cm}p{1.3cm}p{1.3cm}p{1.3cm}p{1.3cm}}

\hline\noalign{\smallskip}

&$j_1=0$&$j_1=1$&$j_1=2$&$j_1=3$&$j_1=4$&$j_1=5$&$j_1=6$\\

\noalign{\smallskip}\hline\noalign{\smallskip}

$j_2=0$&$0$&0&0&$\frac{2}{693}$&0&$\frac{-4}{1287}$&0\\
\noalign{\smallskip}

$j_2=1$&0&0&$\frac{8}{3465}$&0&$\frac{-10}{9009}$&0&$\frac{-4}{6435}$\\
\noalign{\smallskip}
$j_2=2$&0&$\frac{4}{1155}$&0&$\frac{-74}{45045}$&
0&$\frac{16}{45045}$&0\\
\noalign{\smallskip}
$j_2=3$&$\frac{8}{693}$&0&$\frac{-116}{45045}$&0&
$\frac{2}{9009}$&0&$\frac{8}{58905}$\\
\noalign{\smallskip}
$j_2=4$&$\frac{2}{99}$&$\frac{-64}{9009}$&0&$\frac{2}{9009}$&
0&$\frac{4}{153153}$&0\\
\noalign{\smallskip}
$j_2=5$&$\frac{2}{1287}$&$0$&$\frac{-8}{45045}$&0&
$\frac{-2}{153153}$&0&$\frac{4}{415701}$\\
\noalign{\smallskip}
$j_2=6$&$\frac{-2}{143}$&$\frac{4}{715}$&$0$&$\frac{-226}{765765}$&
$0$&$\frac{-8}{415701}$&$0$\\
\noalign{\smallskip}\hline\noalign{\smallskip}
\end{tabular}
\end{table}

\begin{table}
\centering
\caption{Coefficients $\bar C_{6j_2j_1}$}
\label{tab:10}      

\begin{tabular}{p{1.3cm}p{1.3cm}p{1.3cm}p{1.3cm}p{1.3cm}p{1.3cm}p{1.3cm}p{1.3cm}}

\hline\noalign{\smallskip}

&$j_1=0$&$j_1=1$&$j_1=2$&$j_1=3$&$j_1=4$&$j_1=5$&$j_1=6$\\

\noalign{\smallskip}\hline\noalign{\smallskip}

$j_2=0$&$0$&0&0&$0$&$\frac{2}{1287}$&$0$&$\frac{-4}{2145}$\\
\noalign{\smallskip}
$j_2=1$&0&0&$0$&$\frac{10}{9009}$&$0$&$\frac{-4}{6435}$&$0$\\
\noalign{\smallskip}
$j_2=2$&0&$0$&$\frac{4}{3003}$&$0$&
$\frac{-38}{45045}$&$0$&$\frac{8}{36465}$\\
\noalign{\smallskip}
$j_2=3$&$0$&$\frac{20}{9009}$&0&$\frac{-10}{9009}$&0&
$\frac{8}{58905}$&0\\
\noalign{\smallskip}
$j_2=4$&$\frac{10}{1287}$&$0$&$\frac{-16}{9009}$&0&$\frac{2}{13923}$&
0&$\frac{4}{188955}$\\
\noalign{\smallskip}
$j_2=5$&$\frac{2}{143}$&$\frac{-32}{6435}$&0&
$\frac{122}{765765}$&0&$\frac{4}{415701}$&0\\
\noalign{\smallskip}
$j_2=6$&$\frac{2}{2145}$&$0$&$\frac{-4}{36465}$&$0$&$\frac{-2}{188955}$&
$0$&$0$\\
\noalign{\smallskip}\hline\noalign{\smallskip}
\end{tabular}
\end{table}

Assume that 
$q_3=2,$ $q_4=1.$  In 
Tables 20--36 we have the exact values of Fo\-u\-ri\-er--Le\-gen\-dre
coefficients 
$\bar C_{j_3j_2j_1}^{001},$
$\bar C_{j_3j_2j_1}^{010},$
$\bar C_{j_3j_2j_1}^{100}$
($j_1,j_2,j_3=0, 1, 2),$
$\bar C_{j_5j_4j_3j_2j_1}$
($j_1,\ldots,j_5=0, 1$).

In the case of pairwise different 
$i_1, \ldots, i_5$ from 
Tables  
20--36 we have

$$
{\sf M}\left\{\left(
I_{(100)T,t}^{(i_1i_2 i_3)}-
I_{(100)T,t}^{(i_1i_2 i_3)q_3}\right)^2\right\}=
$$

$$
=
\frac{(T-t)^{5}}{60}-\sum_{j_1,j_2,j_3=0}^{2}
\left(C_{j_3j_2j_1}^{100}\right)^2\approx 0.00815429(T-t)^5,
$$

\vspace{4mm}

$$
{\sf M}\left\{\left(
I_{(010)T,t}^{(i_1i_2 i_3)}-
I_{(010)T,t}^{(i_1i_2 i_3)q_3}\right)^2\right\}=
$$

$$
=
\frac{(T-t)^{5}}{20}-\sum_{j_1,j_2,j_3=0}^{2}
\left(C_{j_3j_2j_1}^{010}\right)^2\approx 0.0173903(T-t)^5,
$$

\vspace{4mm}

$$
{\sf M}\left\{\left(
I_{(001)T,t}^{(i_1i_2 i_3)}-
I_{(001)T,t}^{(i_1i_2 i_3)q_3}\right)^2\right\}=
$$

$$
=
\frac{(T-t)^5}{10}-\sum_{j_1,j_2,j_3=0}^{2}
\left(C_{j_3j_2j_1}^{001}\right)^2
\approx 0.0252801(T-t)^5,
$$

\vspace{4mm}

$$
{\sf M}\left\{\left(
I_{(00000)T,t}^{(i_1i_2i_3i_4 i_5)}-
I_{(00000)T,t}^{(i_1i_2i_3i_4 i_5)q_4}\right)^2\right\}=
$$

$$
=
\frac{(T-t)^5}{120}-\sum_{j_1,j_2,j_3,j_4,j_5=0}^{1}
C_{j_5j_4j_3j_2j_1}^2\approx 0.00759105(T-t)^5.
$$

\vspace{4mm}

\begin{table}
\centering
\caption{Coefficients $\bar C_{00j_2j_1}$}
\label{tab:11}      

\begin{tabular}{p{1.3cm}p{1.3cm}p{1.3cm}p{1.3cm}}

\hline\noalign{\smallskip}

&$j_1=0$&$j_1=1$&$j_1=2$\\

\noalign{\smallskip}\hline\noalign{\smallskip}

$j_2=0$&$\frac{2}{3}$&$\frac{-2}{5}$&$\frac{2}{15}$\\
\noalign{\smallskip}
$j_2=1$&$\frac{-2}{15}$&$\frac{2}{15}$&$\frac{-2}{21}$\\
\noalign{\smallskip}
$j_2=2$&$\frac{-2}{15}$&$\frac{2}{35}$&$\frac{2}{105}$\\
\noalign{\smallskip}\hline\noalign{\smallskip}
\end{tabular}
\end{table}

\begin{table}
\centering
\caption{Coefficients $\bar C_{10j_2j_1}$}
\label{tab:12}      
\begin{tabular}{p{1.3cm}p{1.3cm}p{1.3cm}p{1.3cm}}
\hline\noalign{\smallskip}
&$j_1=0$&$j_1=1$&$j_1=2$\\
\noalign{\smallskip}\hline\noalign{\smallskip}
$j_2=0$&$\frac{2}{5}$&$\frac{-2}{9}$&$\frac{2}{35}$\\
\noalign{\smallskip}
$j_2=1$&$\frac{-2}{45}$&$\frac{2}{35}$&$\frac{-2}{45}$\\
\noalign{\smallskip}
$j_2=2$&$\frac{-2}{21}$&$\frac{2}{45}$&$\frac{2}{315}$\\
\noalign{\smallskip}\hline\noalign{\smallskip}
\end{tabular}
\end{table}

\begin{table}
\centering
\caption{Coefficients $\bar C_{02j_2j_1}$}
\label{tab:13}      
\begin{tabular}{p{1.3cm}p{1.3cm}p{1.3cm}p{1.3cm}}
\hline\noalign{\smallskip}
&$j_1=0$&$j_1=1$&$j_1=2$\\
\noalign{\smallskip}\hline\noalign{\smallskip}
$j_2=0$&$\frac{-2}{15}$&$\frac{2}{21}$&$\frac{-4}{105}$\\
\noalign{\smallskip}
$j_2=1$&$\frac{2}{35}$&$\frac{-4}{105}$&$\frac{2}{105}$\\
\noalign{\smallskip}
$j_2=2$&$\frac{4}{105}$&$\frac{-2}{105}$&$0$\\
\noalign{\smallskip}\hline\noalign{\smallskip}
\end{tabular}
\end{table}

\begin{table}
\centering
\caption{Coefficients $\bar C_{01j_2j_1}$}
\label{tab:14}      
\begin{tabular}{p{1.3cm}p{1.3cm}p{1.3cm}p{1.3cm}}
\hline\noalign{\smallskip}
&$j_1=0$&$j_1=1$&$j_1=2$\\
\noalign{\smallskip}\hline\noalign{\smallskip}
$j_2=0$&$\frac{2}{15}$&$\frac{-2}{45}$&$\frac{-2}{105}$\\
\noalign{\smallskip}
$j_2=1$&$\frac{2}{45}$&$\frac{-2}{105}$&$\frac{2}{315}$\\
\noalign{\smallskip}
$j_2=2$&$\frac{-2}{35}$&$\frac{2}{63}$&$\frac{-2}{315}$\\
\noalign{\smallskip}\hline\noalign{\smallskip}
\end{tabular}
\vspace{3mm}
\end{table}

\begin{table}
\centering
\caption{Coefficients $\bar C_{11j_2j_1}$}
\label{tab:15}      
\begin{tabular}{p{1.3cm}p{1.3cm}p{1.3cm}p{1.3cm}}
\hline\noalign{\smallskip}
&$j_1=0$&$j_1=1$&$j_1=2$\\
\noalign{\smallskip}\hline\noalign{\smallskip}
$j_2=0$&$\frac{2}{15}$&$\frac{-2}{35}$&$0$\\
\noalign{\smallskip}
$j_2=1$&$\frac{2}{105}$&$0$&$\frac{-2}{315}$\\
\noalign{\smallskip}
$j_2=2$&$\frac{-4}{105}$&$\frac{2}{105}$&$0$\\
\noalign{\smallskip}\hline\noalign{\smallskip}
\end{tabular}
\vspace{3mm}
\end{table}

\begin{table}
\centering
\caption{Coefficients $\bar C_{20j_2j_1}$}
\label{tab:16}      
\begin{tabular}{p{1.3cm}p{1.3cm}p{1.3cm}p{1.3cm}}
\hline\noalign{\smallskip}
&$j_1=0$&$j_1=1$&$j_1=2$\\
\noalign{\smallskip}\hline\noalign{\smallskip}
$j_2=0$&$\frac{2}{15}$&$\frac{-2}{35}$&$0$\\
\noalign{\smallskip}
$j_2=1$&$\frac{2}{105}$&$0$&$\frac{-2}{315}$\\
\noalign{\smallskip}
$j_2=2$&$\frac{-4}{105}$&$\frac{2}{105}$&$0$\\
\noalign{\smallskip}\hline\noalign{\smallskip}
\end{tabular}
\vspace{3mm}
\end{table}

\begin{table}
\centering
\caption{Coefficients $\bar C_{21j_2j_1}$}
\label{tab:17}      
\begin{tabular}{p{1.3cm}p{1.3cm}p{1.3cm}p{1.3cm}}
\hline\noalign{\smallskip}
&$j_1=0$&$j_1=1$&$j_1=2$\\
\noalign{\smallskip}\hline\noalign{\smallskip}
$j_2=0$&$\frac{2}{21}$&$\frac{-2}{45}$&$\frac{2}{315}$\\
\noalign{\smallskip}
$j_2=1$&$\frac{2}{315}$&$\frac{2}{315}$&$\frac{-2}{225}$\\
\noalign{\smallskip}
$j_2=2$&$\frac{-2}{105}$&$\frac{2}{225}$&$\frac{2}{1155}$\\
\noalign{\smallskip}\hline\noalign{\smallskip}
\end{tabular}
\vspace{3mm}
\end{table}

\begin{table}
\centering
\caption{Coefficients $\bar C_{12j_2j_1}$}
\label{tab:18}      
\begin{tabular}{p{1.3cm}p{1.3cm}p{1.3cm}p{1.3cm}}
\hline\noalign{\smallskip}
&$j_1=0$&$j_1=1$&$j_1=2$\\
\noalign{\smallskip}\hline\noalign{\smallskip}
$j_2=0$&$\frac{-2}{35}$&$\frac{2}{45}$&$\frac{-2}{105}$\\
\noalign{\smallskip}
$j_2=1$&$\frac{2}{63}$&$\frac{-2}{105}$&$\frac{2}{225}$\\
\noalign{\smallskip}
$j_2=2$&$\frac{2}{105}$&$\frac{-2}{225}$&$\frac{-2}{3465}$\\
\noalign{\smallskip}\hline\noalign{\smallskip}
\end{tabular}
\vspace{3mm}
\end{table}

\begin{table}
\centering
\caption{Coefficients $\bar C_{22j_2j_1}$}
\label{tab:19}      
\begin{tabular}{p{1.3cm}p{1.3cm}p{1.3cm}p{1.3cm}}
\hline\noalign{\smallskip}
&$j_1=0$&$j_1=1$&$j_1=2$\\
\noalign{\smallskip}\hline\noalign{\smallskip}
$j_2=0$&$\frac{2}{105}$&$\frac{-2}{315}$&$0$\\
\noalign{\smallskip}
$j_2=1$&$\frac{2}{315}$&$0$&$\frac{-2}{1155}$\\
\noalign{\smallskip}
$j_2=2$&$0$&$\frac{2}{3465}$&$0$\\
\noalign{\smallskip}\hline\noalign{\smallskip}
\end{tabular}
\end{table}

\begin{table}
\centering
\caption{Coefficients $\bar C_{0j_2j_1}^{001}$}
\label{tab:20}      
\begin{tabular}{p{1.3cm}p{1.3cm}p{1.3cm}p{1.3cm}}
\hline\noalign{\smallskip}
&$j_1=0$&$j_1=1$&$j_1=2$\\
\noalign{\smallskip}\hline\noalign{\smallskip}
$j_2=0$&$-2$&$\frac{14}{15}$&$\frac{-2}{15}$\\
\noalign{\smallskip}
$j_2=1$&$\frac{-2}{15}$&$\frac{-2}{15}$&$\frac{6}{35}$\\
\noalign{\smallskip}
$j_2=2$&$\frac{2}{5}$&$\frac{-22}{105}$&$\frac{-2}{105}$\\
\noalign{\smallskip}\hline\noalign{\smallskip}
\end{tabular}
\end{table}

\begin{table}
\centering
\caption{Coefficients $\bar C_{1j_2j_1}^{001}$}
\label{tab:21}      
\begin{tabular}{p{1.3cm}p{1.3cm}p{1.3cm}p{1.3cm}}
\hline\noalign{\smallskip}
&$j_1=0$&$j_1=1$&$j_1=2$\\
\noalign{\smallskip}\hline\noalign{\smallskip}
$j_2=0$&$\frac{-6}{5}$&$\frac{22}{45}$&$\frac{-2}{105}$\\
\noalign{\smallskip}
$j_2=1$&$\frac{-2}{9}$&$\frac{-2}{105}$&$\frac{26}{315}$\\
\noalign{\smallskip}
$j_2=2$&$\frac{22}{105}$&$\frac{-38}{315}$&$\frac{-2}{315}$\\
\noalign{\smallskip}\hline\noalign{\smallskip}
\end{tabular}
\end{table}

\begin{table}
\centering
\caption{Coefficients $\bar C_{2j_2j_1}^{001}$}
\label{tab:22}      
\begin{tabular}{p{1.3cm}p{1.3cm}p{1.3cm}p{1.3cm}}
\hline\noalign{\smallskip}
&$j_1=0$&$j_1=1$&$j_1=2$\\
\noalign{\smallskip}\hline\noalign{\smallskip}
$j_2=0$&$\frac{-2}{5}$&$\frac{2}{21}$&$\frac{4}{105}$\\
\noalign{\smallskip}
$j_2=1$&$\frac{-22}{105}$&$\frac{4}{105}$&$\frac{2}{105}$\\
\noalign{\smallskip}
$j_2=2$&$0$&$\frac{-2}{105}$&$0$\\
\noalign{\smallskip}\hline\noalign{\smallskip}
\end{tabular}
\end{table}

\begin{table}
\centering
\caption{Coefficients $\bar C_{0j_2j_1}^{100}$}
\label{tab:23}      
\begin{tabular}{p{1.3cm}p{1.3cm}p{1.3cm}p{1.3cm}}
\hline\noalign{\smallskip}
&$j_1=0$&$j_1=1$&$j_1=2$\\
\noalign{\smallskip}\hline\noalign{\smallskip}
$j_2=0$&$\frac{-2}{3}$&$\frac{2}{15}$&$\frac{2}{15}$\\
\noalign{\smallskip}
$j_2=1$&$\frac{-2}{15}$&$\frac{-2}{45}$&$\frac{2}{35}$\\
\noalign{\smallskip}
$j_2=2$&$\frac{2}{15}$&$\frac{-2}{35}$&$\frac{-4}{105}$\\
\noalign{\smallskip}\hline\noalign{\smallskip}
\end{tabular}
\end{table}

Note that from (\ref{star00011}) for $k=5$ we obtain

$$
{\sf M}\left\{\left(
I_{(00000)T,t}^{(i_1i_2 i_3 i_4 i_5)}-
I_{(00000)T,t}^{(i_1i_2 i_3 i_4 i_5)q_4}\right)^2\right\}\le
120\left(\frac{(T-t)^{5}}{120}-\sum_{j_1,j_2,j_3,j_4,j_5=0}^{q_4}
C_{j_5j_4j_3j_2j_1}^2\right),
$$

\vspace{3mm}
\noindent
where $i_1, \ldots, i_5=1,\ldots,m$.

Moreover, from the inequality (\ref{star00011}) we get the 
following useful estimates

\vspace{1mm}

$$
{\sf M}\left\{\left(
I_{(01)T,t}^{(i_1i_2)}-
I_{(01)T,t}^{(i_1i_2)q}\right)^2\right\}\le
2\Biggl(\frac{(T-t)^{4}}{4}-\sum_{j_1,j_2=0}^{q}
\left(C_{j_2j_1}^{01}\right)^2\Biggr),
$$

\vspace{3mm}
$$
{\sf M}\left\{\left(
I_{(10)T,t}^{(i_1i_2)}-
I_{(10)T,t}^{(i_1i_2)q}\right)^2\right\}\le
2\Biggl(\frac{(T-t)^{4}}{12}-\sum_{j_1,j_2=0}^{q}
\left(C_{j_2j_1}^{10}\right)^2\Biggr),
$$

\vspace{3mm}
$$
{\sf M}\left\{\left(
I_{(100)T,t}^{(i_1i_2 i_3)}-
I_{(100)T,t}^{(i_1i_2 i_3)q}\right)^2\right\}\le
6\Biggl(\frac{(T-t)^{5}}{60}-\sum_{j_1,j_2,j_3=0}^{q}
\left(C_{j_3j_2j_1}^{100}\right)^2\Biggr),
$$

\vspace{3mm}
$$
{\sf M}\left\{\left(
I_{(010)T,t}^{(i_1i_2 i_3)}-
I_{(010)T,t}^{(i_1i_2 i_3)q}\right)^2\right\}\le
6\Biggl(\frac{(T-t)^{5}}{20}-\sum_{j_1,j_2,j_3=0}^{q}
\left(C_{j_3j_2j_1}^{010}\right)^2\Biggr),
$$

\vspace{3mm}

$$
{\sf M}\left\{\left(
I_{(001)T,t}^{(i_1i_2 i_3)}-
I_{(001)T,t}^{(i_1i_2 i_3)q}\right)^2\right\}\le
6\Biggl(\frac{(T-t)^5}{10}-\sum_{j_1,j_2,j_3=0}^{q}
\left(C_{j_3j_2j_1}^{001}\right)^2\Biggr),
$$

\vspace{3mm}

$$
{\sf M}\left\{\left(
I_{(20))T,t}^{(i_1i_2)}-
I_{(20)T,t}^{(i_1i_2)q}\right)^2\right\}\le
2\Biggl(\frac{(T-t)^6}{30}-\sum_{j_2,j_1=0}^{q}
\left(C_{j_2j_1}^{20}\right)^2\Biggr),
$$

\vspace{3mm}

$$
{\sf M}\left\{\left(
I_{(11)T,t}^{(i_1i_2)}-
I_{(11)T,t}^{(i_1i_2)q}\right)^2\right\}\le
2\Biggl(\frac{(T-t)^6}{18}-\sum_{j_2,j_1=0}^{q}
\left(C_{j_2j_1}^{11}\right)^2\Biggr),
$$

\vspace{3mm}

$$
{\sf M}\left\{\left(
I_{(02)T,t}^{(i_1i_2)}-
I_{(02)T,t}^{(i_1i_2)q}\right)^2\right\}\le
2\Biggl(\frac{(T-t)^6}{6}-\sum_{j_2,j_1=0}^{q}
\left(C_{j_2j_1}^{02}\right)^2\Biggr),
$$

\vspace{3mm}

$$
{\sf M}\left\{\left(
I_{(1000)T,t}^{(i_1i_2 i_3i_4)}-
I_{(1000)T,t}^{(i_1i_2 i_3i_4)q}\right)^2\right\}\le
24\Biggl(\frac{(T-t)^{6}}{360}-\sum_{j_1,j_2,j_3, j_4=0}^{q}
\left(C_{j_4j_3j_2j_1}^{1000}\right)^2\Biggr),
$$

\vspace{3mm}

$$
{\sf M}\left\{\left(
I_{(0100)T,t}^{(i_1i_2 i_3i_4)}-
I_{(0100)T,t}^{(i_1i_2 i_3i_4)q}\right)^2\right\}\le
24\Biggl(\frac{(T-t)^{6}}{120}-\sum_{j_1,j_2,j_3, j_4=0}^{q}
\left(C_{j_4j_3j_2j_1}^{0100}\right)^2\Biggr),
$$

\vspace{3mm}

$$
{\sf M}\left\{\left(
I_{(0010)T,t}^{(i_1i_2 i_3i_4)}-
I_{(0010)T,t}^{(i_1i_2 i_3 i_4)q}\right)^2\right\}\le
24\Biggl(\frac{(T-t)^6}{60}-\sum_{j_1,j_2,j_3, j_4=0}^{q}
\left(C_{j_4j_3j_2j_1}^{0010}\right)^2\Biggr),
$$

\vspace{3mm}

$$
{\sf M}\left\{\left(
I_{(0001)T,t}^{(i_1i_2 i_3 i_4)}-
I_{(0001)T,t}^{(i_1i_2 i_3 i_4)q}\right)^2\right\}\le
24\Biggl(\frac{(T-t)^6}{36}-\sum_{j_1,j_2,j_3, j_4=0}^{q}
\left(C_{j_4j_3j_2j_1}^{0001}\right)^2\Biggr),
$$

\vspace{3mm}

$$
{\sf M}\left\{\left(
I_{(000000)T,t}^{(i_1 i_2 i_3 i_4 i_5 i_6)}-
I_{(000000)T,t}^{(i_1 i_2 i_3 i_4 i_5 i_6)q}\right)^2\right\}\le
720\left(\frac{(T-t)^{6}}{720}-\sum_{j_1,j_2,j_3,j_4,j_5,j_6=0}^{q}
C_{j_6 j_5 j_4 j_3 j_2 j_1}^2\right).
$$

\vspace{5mm}

\begin{table}
\centering
\caption{Coefficients $\bar C_{1j_2j_1}^{100}$}
\label{tab:24}      
\begin{tabular}{p{1.3cm}p{1.3cm}p{1.3cm}p{1.3cm}}
\hline\noalign{\smallskip}
&$j_1=0$&$j_1=1$&$j_1=2$\\
\noalign{\smallskip}\hline\noalign{\smallskip}
$j_2=0$&$\frac{-2}{5}$&$\frac{2}{45}$&$\frac{2}{21}$\\
\noalign{\smallskip}
$j_2=1$&$\frac{-2}{15}$&$\frac{-2}{105}$&$\frac{4}{105}$\\
\noalign{\smallskip}
$j_2=2$&$\frac{2}{35}$&$\frac{-2}{63}$&$\frac{-2}{105}$\\
\noalign{\smallskip}\hline\noalign{\smallskip}
\end{tabular}
\end{table}

\begin{table}
\centering
\caption{Coefficients $\bar C_{2j_2j_1}^{100}$}
\label{tab:25}      
\begin{tabular}{p{1.3cm}p{1.3cm}p{1.3cm}p{1.3cm}}
\hline\noalign{\smallskip}
&$j_1=0$&$j_1=1$&$j_1=2$\\
\noalign{\smallskip}\hline\noalign{\smallskip}
$j_2=0$&$\frac{-2}{15}$&$\frac{-2}{105}$&$\frac{4}{105}$\\
\noalign{\smallskip}
$j_2=1$&$\frac{-2}{21}$&$\frac{-2}{315}$&$\frac{2}{105}$\\
\noalign{\smallskip}
$j_2=2$&$\frac{-2}{105}$&$\frac{-2}{315}$&$0$\\
\noalign{\smallskip}\hline\noalign{\smallskip}
\end{tabular}
\end{table}

\begin{table}
\centering
\caption{Coefficients $\bar C_{0j_2j_1}^{010}$}
\label{tab:26}      
\begin{tabular}{p{1.3cm}p{1.3cm}p{1.3cm}p{1.3cm}}
\hline\noalign{\smallskip}
&$j_1=0$&$j_1=1$&$j_1=2$\\
\noalign{\smallskip}\hline\noalign{\smallskip}
$j_2=0$&$\frac{-4}{3}$&$\frac{8}{15}$&$0$\\
\noalign{\smallskip}
$j_2=1$&$\frac{-4}{15}$&$0$&$\frac{8}{105}$\\
\noalign{\smallskip}
$j_2=2$&$\frac{4}{15}$&$\frac{-16}{105}$&$0$\\
\noalign{\smallskip}\hline\noalign{\smallskip}
\end{tabular}
\end{table}

\begin{table}
\centering
\caption{Coefficients $\bar C_{1j_2j_1}^{010}$}
\label{tab:27}      
\begin{tabular}{p{1.3cm}p{1.3cm}p{1.3cm}p{1.3cm}}
\hline\noalign{\smallskip}
&$j_1=0$&$j_1=1$&$j_1=2$\\
\noalign{\smallskip}\hline\noalign{\smallskip}
$j_2=0$&$\frac{-4}{5}$&$\frac{4}{15}$&$\frac{4}{105}$\\
\noalign{\smallskip}
$j_2=1$&$\frac{-4}{15}$&$\frac{4}{105}$&$\frac{4}{105}$\\
\noalign{\smallskip}
$j_2=2$&$\frac{4}{35}$&$\frac{-8}{105}$&$0$\\
\noalign{\smallskip}\hline\noalign{\smallskip}
\end{tabular}
\end{table}

\begin{table}
\centering
\caption{Coefficients $\bar C_{2j_2j_1}^{010}$}
\label{tab:28}      
\begin{tabular}{p{1.3cm}p{1.3cm}p{1.3cm}p{1.3cm}}
\hline\noalign{\smallskip}
&$j_1=0$&$j_1=1$&$j_1=2$\\
\noalign{\smallskip}\hline\noalign{\smallskip}
$j_2=0$&$\frac{-4}{15}$&$\frac{4}{105}$&$\frac{4}{105}$\\
\noalign{\smallskip}
$j_2=1$&$\frac{-4}{21}$&$\frac{4}{105}$&$\frac{4}{315}$\\
\noalign{\smallskip}
$j_2=2$&$\frac{-4}{105}$&$0$&$0$\\
\noalign{\smallskip}\hline\noalign{\smallskip}
\end{tabular}
\end{table}

\begin{table}
\centering
\caption{Coefficients $\bar C_{000j_2j_1}$}
\label{tab:29}      
\begin{tabular}{p{1.3cm}p{1.3cm}p{1.3cm}}
\hline\noalign{\smallskip}
&$j_1=0$&$j_1=1$\\
\noalign{\smallskip}\hline\noalign{\smallskip}
$j_2=0$&$\frac{4}{15}$&$\frac{-8}{45}$\\
\noalign{\smallskip}
$j_2=1$&$\frac{-4}{45}$&$\frac{8}{105}$\\
\noalign{\smallskip}\hline\noalign{\smallskip}
\end{tabular}
\end{table}

\begin{table}
\centering
\caption{Coefficients $\bar C_{010j_2j_1}$}
\label{tab:30}      
\begin{tabular}{p{1.3cm}p{1.3cm}p{1.3cm}}
\hline\noalign{\smallskip}
&$j_1=0$&$j_1=1$\\
\noalign{\smallskip}\hline\noalign{\smallskip}
$j_2=0$&$\frac{4}{45}$&$\frac{-16}{315}$\\
\noalign{\smallskip}
$j_2=1$&$\frac{-4}{315}$&$\frac{4}{315}$\\
\noalign{\smallskip}\hline\noalign{\smallskip}
\end{tabular}
\end{table}

\begin{table}
\centering
\caption{Coefficients $\bar C_{110j_2j_1}$}
\label{tab:31}      
\begin{tabular}{p{1.3cm}p{1.3cm}p{1.3cm}}
\hline\noalign{\smallskip}
&$j_1=0$&$j_1=1$\\
\noalign{\smallskip}\hline\noalign{\smallskip}
$j_2=0$&$\frac{8}{105}$&$\frac{-2}{45}$\\
\noalign{\smallskip}
$j_2=1$&$\frac{-4}{315}$&$\frac{4}{315}$\\
\noalign{\smallskip}\hline\noalign{\smallskip}
\end{tabular}
\end{table}

In addition, from Theorem 8 for $k=2$ we have

$$
{\sf M}\biggl\{\left(I_{(10)T,t}^{(i_1 i_2)}-
I_{(10)T,t}^{(i_1 i_2)q}
\right)^2\biggr\}
=\frac{(T-t)^4}{12}
-\sum_{j_1,j_2=0}^q
\left(C_{j_2j_1}^{10}\right)^2-
\sum_{j_1,j_2=0}^q C_{j_2j_1}^{10}C_{j_1j_2}^{10}\ \ \ (i_1=i_2),
$$

\vspace{3mm}

$$
{\sf M}\biggl\{\left(I_{(10)T,t}^{(i_1 i_2)}-
I_{(10)T,t}^{(i_1 i_2)q}
\right)^2\biggr\}=\frac{(T-t)^4}{12}
-\sum_{j_1,j_2=0}^q
\left(C_{j_2j_1}^{10}\right)^2\ \ \ (i_1\ne i_2),
$$

\vspace{3mm}

$$
{\sf M}\biggl\{\left(I_{(01)T,t}^{(i_1 i_2)}-
I_{(01)T,t}^{(i_1 i_2)q}
\right)^2\biggr\}
=\frac{(T-t)^4}{4}
-\sum_{j_1,j_2=0}^q
\left(C_{j_2j_1}^{01}\right)^2-
\sum_{j_1,j_2=0}^q C_{j_2j_1}^{01}C_{j_1j_2}^{01}\ \ \ (i_1=i_2),
$$

\vspace{3mm}

$$
{\sf M}\biggl\{\left(I_{(01)T,t}^{(i_1 i_2)}-
I_{(01)T,t}^{(i_1 i_2)q}
\right)^2\biggr\}=\frac{(T-t)^4}{4}
-\sum_{j_1,j_2=0}^q
\left(C_{j_2j_1}^{01}\right)^2\ \ \ (i_1\ne i_2),
$$

$$
{\sf M}\biggl\{\left(I_{(20)T,t}^{(i_1 i_2)}-
I_{(20)T,t}^{(i_1 i_2)q}
\right)^2\biggr\}
=\frac{(T-t)^6}{30}
-\sum_{j_1,j_2=0}^q
\left(C_{j_2j_1}^{20}\right)^2-
\sum_{j_1,j_2=0}^q C_{j_2j_1}^{20}C_{j_1j_2}^{20}\ \ \ (i_1=i_2),
$$

\vspace{3mm}

$$
{\sf M}\biggl\{\left(I_{(20)T,t}^{(i_1 i_2)}-
I_{(20)T,t}^{(i_1 i_2)q}
\right)^2\biggr\}=\frac{(T-t)^6}{30}
-\sum_{j_1,j_2=0}^q
\left(C_{j_2j_1}^{20}\right)^2\ \ \ (i_1\ne i_2),
$$

\vspace{3mm}

$$
{\sf M}\biggl\{\left(I_{(11)T,t}^{(i_1 i_2)}-
I_{(11)T,t}^{(i_1 i_2)q}
\right)^2\biggr\}=
\frac{(T-t)^6}{18}
-\sum_{j_1,j_2=0}^q
\left(C_{j_2j_1}^{11}\right)^2-
\sum_{j_1,j_2=0}^q C_{j_2j_1}^{11}C_{j_1j_2}^{11}\ \ \ (i_1=i_2),
$$

\vspace{3mm}

$$
{\sf M}\biggl\{\left(I_{(11)T,t}^{(i_1 i_2)}-
I_{(11)T,t}^{(i_1 i_2)q}
\right)^2\biggr\}=\frac{(T-t)^6}{18}
-\sum_{j_1,j_2=0}^q
\left(C_{j_2j_1}^{11}\right)^2\ \ \ (i_1\ne i_2),
$$

\vspace{3mm}

$$
{\sf M}\biggl\{\left(I_{(02)}^{(i_1 i_2)}-
I_{(02)T,t}^{(i_1 i_2)q}
\right)^2\biggr\}=
\frac{(T-t)^6}{6}
-\sum_{j_1,j_2=0}^q
\left(C_{j_2j_1}^{02}\right)^2-
\sum_{j_1,j_2=0}^q C_{j_2j_1}^{02}C_{j_1j_2}^{02}\ \ \ (i_1=i_2),
$$

\vspace{3mm}

$$
{\sf M}\biggl\{\left(I_{(02)T,t}^{(i_1 i_2)}-
I_{(02)T,t}^{(i_1 i_2)q}
\right)^2\biggr\}=\frac{(T-t)^6}{6}
-\sum_{j_1,j_2=0}^q
\left(C_{j_2j_1}^{02}\right)^2\ \ \ (i_1\ne i_2).
$$

\vspace{3mm}

\begin{table}
\centering
\caption{Coefficients $\bar C_{011j_2j_1}$}
\label{tab:32}      
\begin{tabular}{p{1.3cm}p{1.3cm}p{1.3cm}}
\hline\noalign{\smallskip}
&$j_1=0$&$j_1=1$\\
\noalign{\smallskip}\hline\noalign{\smallskip}
$j_2=0$&$\frac{8}{315}$&$\frac{-4}{315}$\\
\noalign{\smallskip}
$j_2=1$&$0$&$\frac{2}{945}$\\
\noalign{\smallskip}\hline\noalign{\smallskip}
\end{tabular}
\end{table}

\begin{table}
\centering
\caption{Coefficients $\bar C_{001j_2j_1}$}
\label{tab:33}      
\begin{tabular}{p{1.3cm}p{1.3cm}p{1.3cm}}
\hline\noalign{\smallskip}
&$j_1=0$&$j_1=1$\\
\noalign{\smallskip}\hline\noalign{\smallskip}
$j_2=0$&$0$&$\frac{4}{315}$\\
\noalign{\smallskip}
$j_2=1$&$\frac{8}{315}$&$\frac{-2}{105}$\\
\noalign{\smallskip}\hline\noalign{\smallskip}
\end{tabular}
\end{table}

\begin{table}
\centering
\caption{Coefficients $\bar C_{100j_2j_1}$}
\label{tab:34}      
\begin{tabular}{p{1.3cm}p{1.3cm}p{1.3cm}}
\hline\noalign{\smallskip}
&$j_1=0$&$j_1=1$\\
\noalign{\smallskip}\hline\noalign{\smallskip}
$j_2=0$&$\frac{8}{45}$&$\frac{-4}{35}$\\
\noalign{\smallskip}
$j_2=1$&$\frac{-16}{315}$&$\frac{2}{45}$\\
\noalign{\smallskip}\hline\noalign{\smallskip}
\end{tabular}
\end{table}

\begin{table}
\centering
\caption{Coefficients $\bar C_{101j_2j_1}$}
\label{tab:35}      
\begin{tabular}{p{1.3cm}p{1.3cm}p{1.3cm}}
\hline\noalign{\smallskip}
&$j_1=0$&$j_1=1$\\
\noalign{\smallskip}\hline\noalign{\smallskip}
$j_2=0$&$\frac{4}{315}$&$0$\\
\noalign{\smallskip}
$j_2=1$&$\frac{4}{315}$&$\frac{-8}{945}$\\
\noalign{\smallskip}\hline\noalign{\smallskip}
\end{tabular}
\end{table}

\begin{table}
\centering
\caption{Coefficients $\bar C_{111j_2j_1}$}
\label{tab:36}      
\begin{tabular}{p{1.3cm}p{1.3cm}p{1.3cm}}
\hline\noalign{\smallskip}
&$j_1=0$&$j_1=1$\\
\noalign{\smallskip}\hline\noalign{\smallskip}
$j_2=0$&$\frac{2}{105}$&$\frac{-8}{945}$\\
\noalign{\smallskip}
$j_2=1$&$\frac{2}{945}$&$0$\\
\noalign{\smallskip}\hline\noalign{\smallskip}
\end{tabular}
\end{table}

Clearly, expansions for iterated Stratonovich stochastic integrals
(see above) are simpler 
than expansions for
iterated Ito stochastic integrals (see Theorems 1, 2, and
(\ref{a1})--(\ref{a6})). However, the calculation of the mean-square
approximation error for iterated Stratonovich
stochastic integrals turns out to be much more difficult than for 
iterated Ito stochastic integrals.
Below we consider how we can estimate or calculate exactly
(for some particular cases)
the mean-square
approximation error for iterated Stratonovich
stochastic integrals.

As we mentioned above, on the basis of 
the presented 
approximations of 
iterated Stratonovich stochastic integrals we 
can see that increasing of multiplicities of these integrals 
leads to increasing 
of orders of smallness with respect to $T-t$
in the mean-square sense 
for iterated Stratonovich stochastic integrals
($T-t\ll 1$ 
since the length of integration interval $[t, T]$ 
for iterated Stratonovich 
stochastic integrals 
plays the role of integration step for the numerical 
methods for Ito SDEs, i.e. $T-t$ is already fairly small).
This leads to a sharp decrease  
of member 
quantities
in the appro\-xi\-ma\-ti\-ons of iterated Stratonovich stochastic 
integrals,
which are required for achieving the acceptable accuracy
of approximation.

From (\ref{fff09}) $(i_1\ne i_2)$ we obtain

\vspace{2mm}
$$
{\sf M}\left\{\left(I_{(00)T,t}^{*(i_1 i_2)}-
I_{(00)T,t}^{*(i_1 i_2)q}
\right)^2\right\}=\frac{(T-t)^2}{2}
\sum\limits_{i=q+1}^{\infty}\frac{1}{4i^2-1}\le 
$$

\vspace{2mm}
\begin{equation}
\label{teacxx}
\le \frac{(T-t)^2}{2}\int\limits_{q}^{\infty}
\frac{1}{4x^2-1}dx
=-\frac{(T-t)^2}{8}{\rm ln}\left|
1-\frac{2}{2q+1}\right|\le 
C_1\frac{(T-t)^2}{q},
\end{equation}

\vspace{5mm}
\noindent
where $C_1$ is a constant.

It is easy to notice that for a sufficiently
small $T-t$ (recall that $T-t\ll 1$ since it is a step of integration
for numerical schemes for Ito SDEs) there 
exists a constant $C_2$ such that

\begin{equation}
\label{teac3xx}
{\sf M}\left\{\left(
I_{(l_1\ldots l_k)T,t}^{*(i_1\ldots i_k)}-
I_{(l_1\ldots l_k)T,t}^{*(i_1\ldots i_k)q}\right)^2\right\}
\le C_2 {\sf M}\left\{\left(I_{(00)T,t}^{*(i_1 i_2)}-
I_{(00)T,t}^{*(i_1 i_2)q}\right)^2\right\},
\end{equation}

\vspace{4mm}
\noindent
where 
$I_{(l_1\ldots l_k)T,t}^{*(i_1\ldots i_k)q}$
is an approximation of the iterated Stratonovich stochastic integral 
$I_{(l_1\ldots l_k)T,t}^{*(i_1\ldots i_k)}.$

From (\ref{teacxx}) and (\ref{teac3xx}) we finally obtain

\vspace{1mm}
\begin{equation}
\label{teac4}
{\sf M}\left\{\left(
I_{(l_1\ldots l_k)T,t}^{*(i_1\ldots i_k)}-
I_{(l_1\ldots l_k)T,t}^{*(i_1\ldots i_k)q}\right)^2\right\}
\le C \frac{(T-t)^2}{q},
\end{equation}

\vspace{4mm}
\noindent
where constant $C$ is independent of $T-t$.

The same idea can be found in \cite{1995} in the framework of 
the method of approximation of iterated Stratonovich stochastic 
integrals based
on the trigonometric expansion of the
Brownian bridge process \cite{1988}. 
Note that, in contrast to the estimate (\ref{teac4}), 
the constant $C$ in Theorems 4--6 does not depend on $p.$

We can get more information about the numbers $q$ (these
numbers are different for different iterated Stratonovich
stochastic integrals)
using the another approach.
Since for pairwise different $i_1,\ldots,i_k=1,\ldots,m$

\vspace{-3mm}
$$
J^{*}[\psi^{(k)}]_{T,t}=J[\psi^{(k)}]_{T,t}\ \ \ \hbox{w.\ p.\ 1,}
$$

\vspace{4mm}
\noindent
where $J[\psi^{(k)}]_{T,t},$ $J^{*}[\psi^{(k)}]_{T,t}$
are defined by (\ref{ito}) and (\ref{str}) correspondingly,
then 
for pairwise different 
$i_1,\ldots,i_6=1,\ldots,m$ from Theorem 8 we obtain

\vspace{1mm}
$$
{\sf M}\left\{\left(
I_{(01)T,t}^{*(i_1i_2)}-
I_{(01)T,t}^{*(i_1i_2)q}\right)^2\right\}=
\frac{(T-t)^{4}}{4}-\sum_{j_1,j_2=0}^{q}
\left(C_{j_2j_1}^{01}\right)^2,
$$

\vspace{3mm}
$$
{\sf M}\left\{\left(
I_{(10)T,t}^{*(i_1i_2)}-
I_{(10)T,t}^{*(i_1i_2)q}\right)^2\right\}=
\frac{(T-t)^{4}}{12}-\sum_{j_1,j_2=0}^{q}
\left(C_{j_2j_1}^{10}\right)^2,
$$

\vspace{3mm}
$$
{\sf M}\left\{\left(
I_{(000)T,t}^{*(i_1i_2 i_3)}-
I_{(000)T,t}^{*(i_1i_2 i_3)q}\right)^2\right\}=
\frac{(T-t)^{3}}{6}-\sum_{j_3,j_2,j_1=0}^{q}
C_{j_3j_2j_1}^2,
$$

\vspace{3mm}

$$
{\sf M}\left\{\left(
I_{(0000)T,t}^{*(i_1i_2 i_3 i_4)}-
I_{(0000)T,t}^{*(i_1i_2 i_3 i_4)q}\right)^2\right\}=
\frac{(T-t)^{4}}{24}-\sum_{j_1,j_2,j_3,j_4=0}^{q}
C_{j_4j_3j_2j_1}^2,
$$

\vspace{3mm}

$$
{\sf M}\left\{\left(
I_{(100)T,t}^{*(i_1i_2 i_3)}-
I_{(100)T,t}^{*(i_1i_2 i_3)q}\right)^2\right\}=
\frac{(T-t)^{5}}{60}-\sum_{j_1,j_2,j_3=0}^{q}
\left(C_{j_3j_2j_1}^{100}\right)^2,
$$

\vspace{3mm}

$$
{\sf M}\left\{\left(
I_{(010))T,t}^{*(i_1i_2 i_3)}-
I_{(010)T,t}^{*(i_1i_2 i_3)q}\right)^2\right\}=
\frac{(T-t)^{5}}{20}-\sum_{j_1,j_2,j_3=0}^{q}
\left(C_{j_3j_2j_1}^{010}\right)^2,
$$

\vspace{3mm}

$$
{\sf M}\left\{\left(
I_{(001)T,t}^{*(i_1i_2 i_3)}-
I_{(001)T,t}^{*(i_1i_2 i_3)q}\right)^2\right\}=
\frac{(T-t)^5}{10}-\sum_{j_1,j_2,j_3=0}^{q}
\left(C_{j_3j_2j_1}^{001}\right)^2,
$$

\vspace{3mm}

$$
{\sf M}\left\{\left(
I_{(00000)T,t}^{*(i_1 i_2 i_3 i_4 i_5)}-
I_{(00000)T,t}^{*(i_1 i_2 i_3 i_4 i_5)q}\right)^2\right\}=
\frac{(T-t)^{5}}{120}-\sum_{j_1,j_2,j_3,j_4,j_5=0}^{q}
C_{j_5 i_4 i_3 i_2 j_1}^2,
$$

\vspace{3mm}

$$
{\sf M}\left\{\left(
I_{(20)T,t}^{*(i_1i_2)}-
I_{(20)T,t}^{*(i_1i_2)q}\right)^2\right\}=
\frac{(T-t)^6}{30}-\sum_{j_2,j_1=0}^{q}
\left(C_{j_2j_1}^{20}\right)^2,
$$

\vspace{3mm}

$$
{\sf M}\left\{\left(
I_{(11)T,t}^{*(i_1i_2)}-
I_{(11)T,t}^{*(i_1i_2)q}\right)^2\right\}=
\frac{(T-t)^6}{18}-\sum_{j_2,j_1=0}^{q}
\left(C_{j_2j_1}^{11}\right)^2,
$$

\vspace{3mm}

$$
{\sf M}\left\{\left(
I_{(02)T,t}^{*(i_1i_2)}-
I_{(02)T,t}^{*(i_1i_2)q}\right)^2\right\}=
\frac{(T-t)^6}{6}-\sum_{j_2,j_1=0}^{q}
\left(C_{j_2j_1}^{02}\right)^2,
$$

\vspace{3mm}

$$
{\sf M}\left\{\left(
I_{(1000)T,t}^{*(i_1i_2 i_3i_4)}-
I_{(1000)T,t}^{*(i_1i_2 i_3i_4)q}\right)^2\right\}=
\frac{(T-t)^{6}}{360}-\sum_{j_1,j_2,j_3, j_4=0}^{q}
\left(C_{j_4j_3j_2j_1}^{1000}\right)^2,
$$

\vspace{3mm}

$$
{\sf M}\left\{\left(
I_{(0100)T,t}^{*(i_1i_2 i_3i_4)}-
I_{(0100)T,t}^{*(i_1i_2 i_3i_4)q}\right)^2\right\}=
\frac{(T-t)^{6}}{120}-\sum_{j_1,j_2,j_3, j_4=0}^{q}
\left(C_{j_4j_3j_2j_1}^{0100}\right)^2,
$$

\vspace{3mm}

$$
{\sf M}\left\{\left(
I_{(0010)T,t}^{*(i_1i_2 i_3i_4)}-
I_{(0010)T,t}^{*(i_1i_2 i_3 i_4)q}\right)^2\right\}=
\frac{(T-t)^6}{60}-\sum_{j_1,j_2,j_3, j_4=0}^{q}
\left(C_{j_4j_3j_2j_1}^{0010}\right)^2,
$$

\vspace{3mm}

$$
{\sf M}\left\{\left(
I_{(0001)T,t}^{*(i_1i_2 i_3 i_4)}-
I_{(0001)T,t}^{*(i_1i_2 i_3 i_4)q}\right)^2\right\}=
\frac{(T-t)^6}{36}-\sum_{j_1,j_2,j_3, j_4=0}^{q}
\left(C_{j_4j_3j_2j_1}^{0001}\right)^2,
$$
\vspace{2mm}

$$
{\sf M}\left\{\left(
I_{(000000)T,t}^{*(i_1 i_2 i_3 i_4 i_5 i_6)}-
I_{(000000)T,t}^{*(i_1 i_2 i_3 i_4 i_5 i_6)q}\right)^2\right\}=
\frac{(T-t)^{6}}{720}-\sum_{j_1,j_2,j_3,j_4,j_5,j_6=0}^{q}
C_{j_6 j_5 j_4 j_3 j_2 j_1}^2.
$$

\vspace{5mm}

\section{Legendre Polynomials of Trigonometry?}

\vspace{5mm}

This section is devoted to the comparative analysis of the efficiency of 
application of Legendre polynomials and trigonometric functions to the 
mean-square approximation of iterated Ito and 
Stratonovich stochastic 
integrals.

Using Theorems 1, 2, 8 and the complete orthonormal
system of trigonometric functions in the space $L_2([t, T])$, we obtain
for $i_1\ne i_2$ ($i_1, i_2=1,\ldots,m)$

\begin{equation}
\label{43}
I_{(00)T,t}^{(i_1 i_2)}=\frac{1}{2}(T-t)\Biggl(
\zeta_{0}^{(i_1)}\zeta_{0}^{(i_2)}
+\frac{1}{\pi}
\sum_{r=1}^{\infty}\frac{1}{r}\left(
\zeta_{2r}^{(i_1)}\zeta_{2r-1}^{(i_2)}-
\zeta_{2r-1}^{(i_1)}\zeta_{2r}^{(i_2)}
\right.\biggr.
+\left.\sqrt{2}\left(\zeta_{2r-1}^{(i_1)}\zeta_{0}^{(i_2)}-
\zeta_{0}^{(i_1)}\zeta_{2r-1}^{(i_2)}\right)\right)
\Biggr),
\end{equation}

\vspace{2mm}
\begin{equation}
\label{801}
{\sf M}\left\{\left(I_{(00)T,t}^{(i_1 i_2)}-
I_{(00)T,t}^{(i_1 i_2)q}
\right)^2\right\}
=\frac{3(T-t)^{2}}{2\pi^2}\Biggl(\frac{\pi^2}{6}-
\sum_{r=1}^q \frac{1}{r^2}\Biggr),
\end{equation}

\vspace{1mm}
\begin{equation}
\label{1000x}
I_{(00)T,t}^{(i_1 i_2)q}=\frac{1}{2}(T-t)\Biggl(
\zeta_{0}^{(i_1)}\zeta_{0}^{(i_2)}
+\frac{1}{\pi}
\sum_{r=1}^{q}\frac{1}{r}\left(
\zeta_{2r}^{(i_1)}\zeta_{2r-1}^{(i_2)}-
\zeta_{2r-1}^{(i_1)}\zeta_{2r}^{(i_2)}
\right.\biggr.
+\left.\sqrt{2}\left(\zeta_{2r-1}^{(i_1)}\zeta_{0}^{(i_2)}-
\zeta_{0}^{(i_1)}\zeta_{2r-1}^{(i_2)}\right)\right)
\Biggr),
\end{equation}

\vspace{4mm}
\noindent
where 
$$
\zeta_{j}^{(i)}=
\int\limits_t^T \phi_{j}(s) d{\bf f}_s^{(i)}
$$

\vspace{2mm}
\noindent
are independent standard Gaussian random variables
for various
$i$ or $j$,

\vspace{1mm}
$$
\phi_j(s)=\frac{1}{\sqrt{T-t}}
\begin{cases}
1\ &\hbox{for}\ j=0\cr\cr
\sqrt{2}{\rm sin}(2\pi r(s-t)/(T-t))\ &\hbox{for}\ j=2r-1\cr\cr
\sqrt{2}{\rm cos}(2\pi r(s-t)/(T-t))\ &\hbox{for}\ j=2r
\end{cases},
$$

\vspace{5mm}
\noindent
where $r=1, 2,\ldots;$
another notations are the same as in Theorems 1, 2.

The expansion (\ref{43})
was first derived by Milstein G.N. in \cite{1988} on the base of the Karhunen--Loeve
expansion of the Brownian bridge process.

However, this approach has an obvious drawback. Indeed, we have too complex 
formulas (in comparison with (\ref{4002}), (\ref{4003})) 
for the following stochastic integrals with Gaussian distribution

\begin{equation}
\label{mut1}
I_{(1)T,t}^{(i_1)}=-
\frac{{(T-t)}^{3/2}}{2}
\Biggl(\zeta_0^{(i_1)}-\frac{\sqrt{2}}{\pi}\sum_{r=1}^{\infty}
\frac{1}{r}\zeta_{2r-1}^{(i_1)}\Biggr),
\end{equation}

\vspace{1mm}
\begin{equation}
\label{mut2}
I_{(2)T,t}^{(i_1)}=
(T-t)^{5/2}\Biggl(
\frac{1}{3}\zeta_0^{(i_1)}+\frac{1}{\sqrt{2}\pi^2}
\sum_{r=1}^{\infty}\frac{1}{r^2}\zeta_{2r}^{(i_1)}
-\frac{1}{\sqrt{2}\pi}\sum_{r=1}^{\infty}
\frac{1}{r}\zeta_{2r-1}^{(i_1)}
\Biggr),
\end{equation}

\vspace{4mm}
\noindent
where
$i_1=1,\ldots,m.$

In \cite{1988} Milstein G.N. proposed the following mean-square 
approximations on the base of the expansions
(\ref{43}), (\ref{mut1})

\begin{equation}
\label{444}
I_{(1)T,t}^{(i_1)q}=-\frac{{(T-t)}^{3/2}}{2}
\Biggl(\zeta_0^{(i_1)}-\frac{\sqrt{2}}{\pi}\biggl(\sum_{r=1}^{q}
\frac{1}{r}
\zeta_{2r-1}^{(i_1)}+\sqrt{\alpha_q}\xi_q^{(i_1)}\biggr)
\Biggr),
\end{equation}

\vspace{2mm}
$$
I_{(00)T,t}^{(i_1 i_2)q}=\frac{1}{2}(T-t)\Biggl(
\zeta_{0}^{(i_1)}\zeta_{0}^{(i_2)}
+\frac{1}{\pi}
\sum_{r=1}^{q}\frac{1}{r}\left(
\zeta_{2r}^{(i_1)}\zeta_{2r-1}^{(i_2)}-
\zeta_{2r-1}^{(i_1)}\zeta_{2r}^{(i_2)}+
\right.\Biggr.
$$

\vspace{1mm}
\begin{equation}
\label{555}
+\Biggl.\left.\sqrt{2}\left(\zeta_{2r-1}^{(i_1)}\zeta_{0}^{(i_2)}-
\zeta_{0}^{(i_1)}\zeta_{2r-1}^{(i_2)}\right)\right)
+\frac{\sqrt{2}}{\pi}\sqrt{\alpha_q}\left(
\xi_q^{(i_1)}\zeta_0^{(i_2)}-\zeta_0^{(i_1)}\xi_q^{(i_2)}\right)\Biggr),
\end{equation}

\vspace{4mm}
\noindent
where
\begin{equation}
\label{333}
\xi_q^{(i)}=\frac{1}{\sqrt{\alpha_q}}\sum_{r=q+1}^{\infty}
\frac{1}{r}~\zeta_{2r-1}^{(i)},\ \ \
\alpha_q=\frac{\pi^2}{6}-\sum_{r=1}^q\frac{1}{r^2},
\end{equation}

\vspace{3mm}
\noindent
where
$\zeta_0^{(i)},$ $\zeta_{2r}^{(i)},$
$\zeta_{2r-1}^{(i)},$ $\xi_q^{(i)}$ ($r=1,\ldots,q;$
$i=1,\ldots,m$)
are independent standard Gaussian random variables.

Obviously, for the approximations (\ref{444}) and (\ref{555}) we obtain

$$
{\sf M}\left\{\left(I_{(1)T,t}^{(i_1)}-
I_{(1)T,t}^{(i_1)q}
\right)^2\right\}=0,
$$

\vspace{1mm}
\begin{equation}
\label{8010}
{\sf M}\left\{\left(I_{(00)T,t}^{(i_1 i_2)}-
I_{(00)T,t}^{(i_1 i_2)q}
\right)^2\right\}
=\frac{(T-t)^{2}}{2\pi^2}\Biggl(\frac{\pi^2}{6}-
\sum_{r=1}^q \frac{1}{r^2}\Biggr).
\end{equation}

\vspace{4mm}

This idea has been developed in \cite{1995}.
For example, the approximation $I_{(2)T,t}^{(i_1)q},$ 
which corresponds to (\ref{444}), (\ref{555}),
has the form \cite{1995}

\begin{equation}
\label{1970}
I_{(2)T,t}^{(i_1)q}=
(T-t)^{5/2}\Biggl(
\frac{1}{3}\zeta_0^{(i_1)}+\frac{1}{\sqrt{2}\pi^2}
\Biggl(\sum_{r=1}^{q}\frac{1}{r^2}\zeta_{2r}^{(i_1)}+\sqrt{\beta_q}
\mu_q^{(i_1)}\Biggr)\Biggr.
\Biggl.-\frac{1}{\sqrt{2}\pi}\Biggl(\sum_{r=1}^q
\frac{1}{r}\zeta_{2r-1}^{(i_1)}+\sqrt{\alpha_q}\xi_q^{(i_1)}\Biggr)\Biggr),
\end{equation}

\vspace{3mm}

$$
{\sf M}\left\{\left(I_{(2)T,t}^{(i_1)}-
I_{(2)T,t}^{(i_1)q}
\right)^2\right\}=0,
$$

\vspace{5mm}
\noindent
where
$\xi_q^{(i)},$ $\alpha_q$ have the form (\ref{333}) and

\begin{equation}
\label{333a}
\mu_q^{(i)}=\frac{1}{\sqrt{\beta_q}}\sum_{r=q+1}^{\infty}
\frac{1}{r^2}~\zeta_{2r}^{(i)},\ \ \
\beta_q=\frac{\pi^4}{90}-\sum_{r=1}^q\frac{1}{r^4},
\end{equation}

\vspace{4mm}
\noindent
where
$\zeta_0^{(i)},$ $\zeta_{2r}^{(i)},$
$\zeta_{2r-1}^{(i)},$ $\xi_q^{(i)},$ $\mu_q^{(i)}$ ($r=1,\ldots,q;$
$i=1,\ldots,m$) are independent
standard Gaussian random variables.

Nevetheless, the expansions (\ref{444}), (\ref{1970}) are too complex for
the numerical modeling of two Gaussian random variables
$I_{(1)T,t}^{(i_1)},$ $I_{(2)T,t}^{(i_1)}.$

Further, we will see that the using of random 
variables $\xi_q^{(i)}$ and 
$\mu_q^{(i)}$ will drastically
complicate the approximation of the stochastic 
integral $I_{(000)T,t}^{(i_1 i_2 i_3)};$
$i_1,i_2,i_3=1,\ldots,m.$ 
This is due to the fact that for this approach 
the number $q$ is fixed for all stochastic integrals 
included into the considered collection
\cite{1995}. However, it is clear that due 
to the smallness of $T-t$, the number $q$ for $I_{(000)T,t}^{(i_1 i_2 i_3)}$
could be taken significantly 
less than in the formula (\ref{555}) (see for comparison 
the case of Legendre polynomials). 
This feature is also valid for the formulas (\ref{444}), (\ref{1970}).

To obtain the expansion for (\ref{str}) on the base of the  
approach from \cite{1988}
the truncated 
trigonometric expansions of components of the multidimensional Wiener  
process ${\bf f}_s$ must be
iteratively substituted in the single integrals, and the integrals
must be calculated, starting from the innermost integral.
This is a complicated procedure that obviously does not lead to a general
expansion of (\ref{str}) valid for an arbitrary multiplicity $k.$
For this reason, only expansions of simplest single, double, and triple
integrals (\ref{str}) were obtained (see \cite{1995}, \cite{1988},
\cite{KPS}-\cite{Zapad-9}).

At that, in \cite{1988} the case $\psi_1(s), \psi_2(s)\equiv 1$ and
$i_1, i_2=0, 1,\ldots,m$ is considered. In 
\cite{1995}, \cite{KPS}-\cite{Zapad-9} the attempt to consider the case 
$\psi_1(s), \psi_2(s), \psi_3(s)\equiv 1$ and 
$i_1, i_2, i_3=0, 1,\ldots,m$ is realized.

Note that the mean-square
convergence of $J_{(111)T,t}^{*(i_1 i_2 i_3)q}$ to
$J_{(111)T,t}^{*(i_1 i_2 i_3)}$ if $q\to\infty$
was not proved rigorously in
\cite{1995}
(Sect.~5.8, pp.~202--204), \cite{KPS} (pp.~82-84),
\cite{Zapad-2} (pp.~438-439),  
\cite{Zapad-9} (pp.~263-264) 
within the 
frames of the Milstein approach \cite{1988}
together with the Wong--Zakai approximation \cite{W-Z-1}-\cite{Watanabe}
(see discussion in Sect.~8 for detail).

Consider the approximation
$I_{(00)T,t}^{(i_1 i_2)q}$ 
of the iterated stochastic 
integral $I_{(00)T,t}^{(i_1 i_2)}$ obtained from 
(\ref{4004}) by replacing $\infty$ with $q$

\begin{equation}
\label{401}
I_{(00)T,t}^{(i_1 i_2)q}=
\frac{T-t}{2}\Biggl(\zeta_0^{(i_1)}\zeta_0^{(i_2)}+\sum_{i=1}^{q}
\frac{1}{\sqrt{4i^2-1}}\left(
\zeta_{i-1}^{(i_1)}\zeta_{i}^{(i_2)}-
\zeta_i^{(i_1)}\zeta_{i-1}^{(i_2)}\right)\Biggr)\ \ \ (i_1\ne i_2).
\end{equation}

\vspace{3mm}

Let us compare computational costs for 
the approximations (\ref{555}), (\ref{401}).
It is not difficult to show that \cite{1}-\cite{10axx1}

\vspace{-3mm}
\begin{equation}
\label{400}
{\sf M}\left\{\left(I_{(00)T,t}^{(i_1 i_2)}-
I_{(00)T,t}^{(i_1 i_2)q}
\right)^2\right\}
=\frac{(T-t)^2}{2}\Biggl(\frac{1}{2}-\sum_{i=1}^q
\frac{1}{4i^2-1}\Biggr).
\end{equation}

\vspace{3mm}

Let us compare (\ref{401}) with (\ref{555}) and (\ref{400}) with (\ref{8010}).
Consider minimal natural numbers $q_{\rm trig}$ and 
$q_{\rm pol},$ which satisfy to (see Table 37)

\vspace{1mm}
$$
\frac{(T-t)^2}{2}\Biggl(\frac{1}{2}-\sum_{i=1}^{q_{\rm pol}}
\frac{1}{4i^2-1}\Biggr)\le (T-t)^3,\ \ \
\frac{(T-t)^{2}}{2\pi^2}\Biggl(\frac{\pi^2}{6}-
\sum_{r=1}^{q_{\rm trig}}\frac{1}{r^2}\Biggr)\le (T-t)^3.
$$

\vspace{4mm}

Thus, we have

\vspace{-1mm}
$$
\frac{q_{\rm pol}}{q_{\rm trig}}\ \ \approx\ \ 1.67,\ \ 
2.22,\ \ 2.43,\ \ 2.36,\ \ 2.41,\ \ 
2.43,\ \ 2.45,\ \ 2.45.
$$

\vspace{3mm}

From the other hand, the  
formula (\ref{555}) includes $(4q+4)m$ independent
standard Gaussian random variables. At the same time the folmula
(\ref{401}) includes only $(2q+2)m$ independent
standard Gaussian random variables. Moreover, the formula
(\ref{401}) is simpler than the formula (\ref{555}).
Thus, in this case we can talk about approximately equal computational costs
for the formulas (\ref{555}) and (\ref{401}).

\begin{table}
\centering
\caption{Numbers $q_{\rm trig},$ $q_{\rm pol}$}
\label{tab:37}      
\begin{tabular}{p{1.1cm}p{1.1cm}p{1.1cm}p{1.1cm}p{1.1cm}p{1.1cm}p{1.1cm}p{1.1cm}p{1.1cm}}
\hline\noalign{\smallskip}
$T-t$&$2^{-5}$&$2^{-6}$&$2^{-7}$&$2^{-8}$&$2^{-9}$&$2^{-10}$&$2^{-11}$&$2^{-12}$\\
\noalign{\smallskip}\hline\noalign{\smallskip}
$q_{\rm trig}$&3&4&7&14&27&53&105&209\\
\noalign{\smallskip}\hline\noalign{\smallskip}
$q_{\rm trig}^{*}$&6&11&20&40&79&157&312&624\\
\noalign{\smallskip}\hline\noalign{\smallskip}
$q_{\rm pol}$&5&9&17&33&65&129&257&513\\
\noalign{\smallskip}\hline\noalign{\smallskip}
\end{tabular}
\end{table}

There is one important feature. 
As we mentioned above, further we will see that the using of random 
variables $\xi_q^{(i)}$ and 
$\mu_q^{(i)}$ will drastically
complicate the approximation of 
the stochastic integral $I_{(000)T,t}^{(i_1 i_2 i_3)};$
$i_1,i_2,i_3=1,\ldots,m.$ 
This is due to the fact that
the number $q$ is fixed for all stochastic integrals, which 
included into the considered collection (the case of trigonometric functions). 
However, it is clear that due 
to the smallness of $T-t$, the number $q$ for $I_{(000)T,t}^{(i_1 i_2 i_3)}$
could be chosen significantly 
less than in the formula (\ref{555}). 
This feature is also valid for the formulas (\ref{444}), (\ref{1970}).
However, for the case of Legendre polynomials we can choose different 
numbers $q$ for different stochastic integrals (see Sect.~3).

From the other hand, if we will not use the random 
variables $\xi_q^{(i)}$ and 
$\mu_q^{(i)},$ then the mean-square error of approximation of the
stochastic integral $I_{(00)T,t}^{(i_1 i_2)}$ will be three times larger
(see (\ref{801})). Moreover, in this case the stochastic integrals
$I_{(1)T,t}^{(i_1)}$, $I_{(2)T,t}^{(i_1)}$ (with Gaussian distribution)
will be approximated worse. 

Consider minimal natural numbers $q_{\rm trig}^{*}$, 
which satisfy the condition (see Table 37)

$$
\frac{3(T-t)^{2}}{2\pi^2}\Biggl(\frac{\pi^2}{6}-
\sum_{r=1}^{q_{\rm trig}^{*}}\frac{1}{r^2}\Biggr)\le (T-t)^3.
$$

\vspace{4mm}

In this situation we can talk about
the advantage of Ledendre polynomials ($q_{\rm trig}^{*} > q_{pol}$ and
(\ref{1000x}) is more complex than (\ref{401})).

Using Theorems 1, 2 for the system of trigonometric functions, we have 
$(i_1\ne i_2,\ i_1\ne i_3,\ i_2\ne i_3)$ \cite{3}-\cite{10axx1}
(also see \cite{1}, \cite{1a})

\vspace{2mm}
$$
I_{(000)T,t}^{(i_1 i_2 i_3)q}
=(T-t)^{3/2}\Biggl(\frac{1}{6}
\zeta_{0}^{(i_1)}\zeta_{0}^{(i_2)}\zeta_{0}^{(i_3)}+\Biggr.
\frac{\sqrt{\alpha_q}}{2\sqrt{2}\pi}\left(
\xi_q^{(i_1)}\zeta_0^{(i_2)}\zeta_0^{(i_3)}-\xi_q^{(i_3)}\zeta_0^{(i_2)}
\zeta_0^{(i_1)}\right)+
$$
$$
+\frac{1}{2\sqrt{2}\pi^2}\sqrt{\beta_q}\left(
\mu_q^{(i_1)}\zeta_0^{(i_2)}\zeta_0^{(i_3)}-2\mu_q^{(i_2)}\zeta_0^{(i_1)}
\zeta_0^{(i_3)}+\mu_q^{(i_3)}\zeta_0^{(i_1)}\zeta_0^{(i_2)}\right)+
$$
$$
+
\frac{1}{2\sqrt{2}}\sum_{r=1}^{q}
\Biggl(\frac{1}{\pi r}\left(
\zeta_{2r-1}^{(i_1)}
\zeta_{0}^{(i_2)}\zeta_{0}^{(i_3)}-
\zeta_{2r-1}^{(i_3)}
\zeta_{0}^{(i_2)}\zeta_{0}^{(i_1)}\right)+\Biggr.
$$
$$
\Biggl.+
\frac{1}{\pi^2 r^2}\left(
\zeta_{2r}^{(i_1)}
\zeta_{0}^{(i_2)}\zeta_{0}^{(i_3)}-
2\zeta_{2r}^{(i_2)}
\zeta_{0}^{(i_3)}\zeta_{0}^{(i_1)}+
\zeta_{2r}^{(i_3)}
\zeta_{0}^{(i_2)}\zeta_{0}^{(i_1)}\right)\Biggr)+
$$
$$
+
\sum_{r=1}^{q}
\Biggl(\frac{1}{4\pi r}\left(
\zeta_{2r}^{(i_1)}
\zeta_{2r-1}^{(i_2)}\zeta_{0}^{(i_3)}-
\zeta_{2r-1}^{(i_1)}
\zeta_{2r}^{(i_2)}\zeta_{0}^{(i_3)}-
\zeta_{2r-1}^{(i_2)}
\zeta_{2r}^{(i_3)}\zeta_{0}^{(i_1)}+
\zeta_{2r-1}^{(i_3)}
\zeta_{2r}^{(i_2)}\zeta_{0}^{(i_1)}\right)+\Biggr.
$$
$$
+
\frac{1}{8\pi^2 r^2}\left(
3\zeta_{2r-1}^{(i_1)}
\zeta_{2r-1}^{(i_2)}\zeta_{0}^{(i_3)}+
\zeta_{2r}^{(i_1)}
\zeta_{2r}^{(i_2)}\zeta_{0}^{(i_3)}-
6\zeta_{2r-1}^{(i_1)}
\zeta_{2r-1}^{(i_3)}\zeta_{0}^{(i_2)}+\right.
$$
\begin{equation}
\label{x5}
\Biggl.\left.
+
3\zeta_{2r-1}^{(i_2)}
\zeta_{2r-1}^{(i_3)}\zeta_{0}^{(i_1)}-
2\zeta_{2r}^{(i_1)}
\zeta_{2r}^{(i_3)}\zeta_{0}^{(i_2)}+
\zeta_{2r}^{(i_3)}
\zeta_{2r}^{(i_2)}\zeta_{0}^{(i_1)}\right)\Biggr)
\Biggl.+D_{T,t}^{(i_1i_2i_3)q}\Biggr),
\end{equation}

\vspace{5mm}
\noindent
where

\vspace{-2mm}
$$
D_{T,t}^{(i_1i_2i_3)q}=
\frac{1}{2\pi^2}\sum_{\stackrel{r,l=1}{{}_{r\ne l}}}^{q}
\Biggl(\frac{1}{r^2-l^2}\biggl(
\zeta_{2r}^{(i_1)}
\zeta_{2l}^{(i_2)}\zeta_{0}^{(i_3)}-
\zeta_{2r}^{(i_2)}
\zeta_{0}^{(i_1)}\zeta_{2l}^{(i_3)}+\biggr.\Biggr.
$$
$$
\Biggl.+\biggl.
\frac{r}{l}
\zeta_{2r-1}^{(i_1)}
\zeta_{2l-1}^{(i_2)}\zeta_{0}^{(i_3)}-\frac{l}{r}
\zeta_{0}^{(i_1)}
\zeta_{2r-1}^{(i_2)}\zeta_{2l-1}^{(i_3)}\biggr)-
\frac{1}{rl}\zeta_{2r-1}^{(i_1)}
\zeta_{0}^{(i_2)}\zeta_{2l-1}^{(i_3)}\Biggr)+
$$
$$
+
\frac{1}{4\sqrt{2}\pi^2}\Biggl(
\sum_{r,m=1}^{q}\Biggl(\frac{2}{rm}
\left(-\zeta_{2r-1}^{(i_1)}
\zeta_{2m-1}^{(i_2)}\zeta_{2m}^{(i_3)}+
\zeta_{2r-1}^{(i_1)}
\zeta_{2r}^{(i_2)}\zeta_{2m-1}^{(i_3)}+
\right.\Biggr.\Biggr.
$$
$$
\left.+
\zeta_{2r-1}^{(i_1)}
\zeta_{2m}^{(i_2)}\zeta_{2m-1}^{(i_3)}-
\zeta_{2r}^{(i_1)}
\zeta_{2r-1}^{(i_2)}\zeta_{2m-1}^{(i_3)}\right)+
$$
$$
+\frac{1}{m(r+m)}
\left(-\zeta_{2(m+r)}^{(i_1)}
\zeta_{2r}^{(i_2)}\zeta_{2m}^{(i_3)}-
\zeta_{2(m+r)-1}^{(i_1)}
\zeta_{2r-1}^{(i_2)}\zeta_{2m}^{(i_3)}-
\right.
$$
$$
\Biggl.\left.
-\zeta_{2(m+r)-1}^{(i_1)}
\zeta_{2r}^{(i_2)}\zeta_{2m-1}^{(i_3)}+
\zeta_{2(m+r)}^{(i_1)}
\zeta_{2r-1}^{(i_2)}\zeta_{2m-1}^{(i_3)}\right)\Biggr)+
$$
$$
+
\sum_{m=1}^{q}\sum_{l=m+1}^{q}\Biggl(\frac{1}{m(l-m)}
\left(\zeta_{2(l-m)}^{(i_1)}
\zeta_{2l}^{(i_2)}\zeta_{2m}^{(i_3)}+
\zeta_{2(l-m)-1}^{(i_1)}
\zeta_{2l-1}^{(i_2)}\zeta_{2m}^{(i_3)}-
\right.\Biggr.
$$
$$
\left.
-\zeta_{2(l-m)-1}^{(i_1)}
\zeta_{2l}^{(i_2)}\zeta_{2m-1}^{(i_3)}+
\zeta_{2(l-m)}^{(i_1)}
\zeta_{2l-1}^{(i_2)}\zeta_{2m-1}^{(i_3)}\right)+
$$
$$
+
\frac{1}{l(l-m)}
\left(-\zeta_{2(l-m)}^{(i_1)}
\zeta_{2m}^{(i_2)}\zeta_{2l}^{(i_3)}+
\zeta_{2(l-m)-1}^{(i_1)}
\zeta_{2m-1}^{(i_2)}\zeta_{2l}^{(i_3)}-
\right.
$$
$$
\Biggl.
\Biggl.
\Biggl.
\left.
-\zeta_{2(l-m)-1}^{(i_1)}
\zeta_{2m}^{(i_2)}\zeta_{2l-1}^{(i_3)}-
\zeta_{2(l-m)}^{(i_1)}
\zeta_{2m-1}^{(i_2)}\zeta_{2l-1}^{(i_3)}\right)\Biggr)\Biggr),
$$

\vspace{6mm}
\noindent
where
$\zeta_0^{(i)},$ $\zeta_{2r}^{(i)},$
$\zeta_{2r-1}^{(i)},$ $\xi_q^{(i)},$ $\mu_q^{(i)}$ ($r=1,\ldots,q;$
$i=1,\ldots,m$) are independent
standard Gaussian random variables (see (\ref{333}), (\ref{333a})).

The mean-square error of approximation (\ref{x5})
$(i_1\ne i_2,\ i_1\ne i_3,\ i_2\ne i_3)$ has the following form 
\cite{3}-\cite{10axx1}
(also see \cite{1}, \cite{1a})

\vspace{1mm}
$$
{\sf M}\left\{\left(I_{(000)T,t}^{(i_1 i_2 i_3)}-
I_{(000)T,t}^{(i_1 i_2 i_3)q}\right)^2\right\}=
(T-t)^3\Biggl(\frac{4}{45}-\frac{1}{4\pi^2}\sum_{r=1}^q\frac{1}{r^2}-
\Biggl.
$$

\vspace{1mm}
\begin{equation}
\label{x7}
\Biggl.-\frac{55}{32\pi^4}\sum_{r=1}^q\frac{1}{r^4}-
\frac{1}{4\pi^4}\sum_{\stackrel{r,l=1}{{}_{r\ne l}}}^q
\frac{5l^4+4r^4-3r^2l^2}{r^2 l^2 \left(r^2-l^2\right)^2}\Biggr).
\end{equation}

\vspace{3mm}

\begin{table}                
\centering
\caption{Confirmation of the formula (\ref{x7})}
\label{tab:38}      
\begin{tabular}{p{2.1cm}p{1.7cm}p{1.7cm}p{2.1cm}p{2.3cm}p{2.3cm}p{2.3cm}}
\hline\noalign{\smallskip}
$\varepsilon/(T-t)^3$&0.0459&0.0072&$7.5722\cdot 10^{-4}$
&$7.5973\cdot 10^{-5}$&
$7.5990\cdot 10^{-6}$\\
\noalign{\smallskip}\hline\noalign{\smallskip}
$q$&1&10&100&1000&10000\\
\noalign{\smallskip}\hline\noalign{\smallskip}
\end{tabular}
\end{table}

In Table 38 we can see the numerical confirmation of 
the formula (\ref{x7}) ($\varepsilon$ is the right-hand side of (\ref{x7})).

As we mentioned above, the Milstein expansion \cite{1988} (i.e. expansion
based on the Karhunen--Loeve expansion of the Brownian bridge
process)
for iterated stochastic integrals leads to iterated application
of the operation of 
limit transition. The analogue of (\ref{x5}) for iterated
Stratonovich stochastic integrals has been derived in 
\cite{1995}, \cite{KPS}-\cite{Zapad-9} on the base of the Milstein expansion
together with the Wong--Zakai approximation \cite{W-Z-1}-\cite{Watanabe}
(without rigorous proof).
It means that the authors in 
\cite{1995}
(Sect.~5.8, pp.~202--204), \cite{KPS} (pp.~82-84),
\cite{Zapad-2} (pp.~438-439),  
\cite{Zapad-9} (pp.~263-264) 
formally
could not use the double sum with the upper limit $q$ in the analogue of
(\ref{x5}). From the other hand the correctness of (\ref{x5}) follows 
directly from Theorems 1, 2. Note that (\ref{x5}) has been
obtained  reasonably
for the first time in \cite{3}. The version of (\ref{x5}) without
using the random variables $\xi_q^{(i)}$ and $\mu_q^{(i)}$ 
can be found in \cite{1} (1997).

The mean-square error (\ref{x7}) has been obtained for the first time
in \cite{3} on the base of the simplified variant of 
Theorem 8 (the case of pairwise different $i_1,\ldots,i_k$).

As we noted above, the number $q$ must be the same
in (\ref{444}), (\ref{555}), (\ref{x5}). This is the main 
drawback of this approach,
because really the number $q$ in (\ref{x5})
can be chosen essentially smaller than in (\ref{555}).

Note that in (\ref{x5}) we can replace 
$I_{(000)T,t}^{(i_1 i_2 i_3)q}$ on $I_{(000)T,t}^{*(i_1 i_2 i_3)q}$
and (\ref{x5}) then will be valid for any $i_1, i_2, i_3 = 1,\ldots,m$
(see Theorem 3).

Let us compare 
the efficiency of application of Legendre polynomials 
and trigonometric functions for approximation of 
the iterated stochastic integrals
$I_{(00)T,t}^{(i_1 i_2)},$ $I_{(000)T,t}^{(i_1 i_2 i_3)}$ $(i_1\ne i_2,\ 
i_1\ne i_3,\ i_2\ne i_3).$

Consider the following conditions $(i_1\ne i_2,\ 
i_1\ne i_3,\ i_2\ne i_3)$

\begin{equation}
\label{z1}
\frac{(T-t)^2}{2}\Biggl(\frac{1}{2}-\sum_{i=1}^{q}
\frac{1}{4i^2-1}\Biggr)\le (T-t)^4,
\end{equation}

\vspace{1mm}
\begin{equation}
\label{z2}
(T-t)^{3}\Biggl(\frac{1}{6}-\sum_{j_1,j_2,j_3=0}^{q_1}
\frac{\left(C_{j_3j_2j_1}\right)^2}{(T-t)^3}\Biggr)\le (T-t)^4,
\end{equation}

\vspace{1mm}
\begin{equation}
\label{z3}
\frac{(T-t)^{2}}{2\pi^2}\Biggl(\frac{\pi^2}{6}-
\sum_{r=1}^{p}\frac{1}{r^2}\Biggr)\le (T-t)^4,
\end{equation}

\vspace{1mm}
\begin{equation}
\label{z4}
(T-t)^3\Biggl(\frac{4}{45}-\frac{1}{4\pi^2}\sum_{r=1}^{p_1}\frac{1}{r^2}
\Biggl.
\Biggl.-\frac{55}{32\pi^4}\sum_{r=1}^{p_1}\frac{1}{r^4}-
\frac{1}{4\pi^4}\sum_{\stackrel{r,l=1}{{}_{r\ne l}}}^{p_1}
\frac{5l^4+4r^4-3r^2l^2}{r^2 l^2 \left(r^2-l^2\right)^2}\Biggr)\le (T-t)^4,
\end{equation}

\vspace{4mm}
\noindent
where
$C_{j_3j_2j_1}$ is defined by
(\ref{zzz2}).

In Tables 39 and 40 we can see minimal numbers 
$q,$ $q_1,$ $p,$ $p_1,$ which satisfy the conditions
(\ref{z1})--(\ref{z4}). As we mentioned above, 
the numbers $q,$ $q_1$ are different. At that $q_1\ll q$ 
(the case of Legendre polynomials). Moreover,
we cannot take different numbers 
$p,$ $p_1$ for the case of trigonometric functions. Thus, we must 
choose $q=p$ in (\ref{444}), (\ref{555}), (\ref{x5}). This leads
to huge computational costs (see very complex formula (\ref{x5})).
From the other hand, we can choose different numbers $q$
in (\ref{444}), (\ref{555}), (\ref{x5}). At that we must exclude
random variables $\xi_q^{(i)},$ $\mu_q^{(i)}$ from 
(\ref{444}), (\ref{555}), (\ref{x5}). 
At this situation we have

\begin{equation}
\label{zzz3}
\frac{3(T-t)^{2}}{2\pi^2}\Biggl(\frac{\pi^2}{6}-
\sum_{r=1}^{p^{*}}\frac{1}{r^2}\Biggr)\le (T-t)^4,
\end{equation}

\vspace{1mm}
\begin{equation}
\label{zzz4}
(T-t)^3\Biggl(\frac{5}{36}-\frac{1}{2\pi^2}\sum_{r=1}^{p_1^{*}}\frac{1}{r^2}
\Biggl.
\Biggl.-\frac{79}{32\pi^4}\sum_{r=1}^{p_1^{*}}\frac{1}{r^4}-
\frac{1}{4\pi^4}\sum_{\stackrel{r,l=1}{{}_{r\ne l}}}^{p_1^{*}}
\frac{5l^4+4r^4-3r^2l^2}{r^2 l^2 \left(r^2-l^2\right)^2}\Biggr)\le (T-t)^4,
\end{equation}

\vspace{4mm}
\noindent
where the left-hand sides of (\ref{zzz3}), (\ref{zzz4}) correspond
to (\ref{555}), (\ref{x5}) but without $\xi_q^{(i)},$ $\mu_q^{(i)}$.
In Table 40 we can see minimal numbers 
$p^{*},$ $p_1^{*}$, which satisfy the conditions
(\ref{zzz3}), (\ref{zzz4}). Moreover,

\begin{equation}
\label{1001x}
{\sf M}\left\{\left(I_{(1)T,t}^{(i_1)}-
I_{(1)T,t}^{(i_1)q}
\right)^2\right\}
=\frac{(T-t)^{3}}{2\pi^2}\Biggl(\frac{\pi^2}{6}-
\sum_{r=1}^q \frac{1}{r^2}\Biggr),
\end{equation}

\vspace{4mm}
\noindent
where
$$
I_{(1)T,t}^{(i_1)q}=-
\frac{{(T-t)}^{3/2}}{2}
\Biggl(\zeta_0^{(i_1)}-\frac{\sqrt{2}}{\pi}\sum_{r=1}^{q}
\frac{1}{r}\zeta_{2r-1}^{(i_1)}\Biggr).
$$

\vspace{3mm}

It is not difficult to see that numbers $q_{\rm trig}$ in Table 37 
correspond to minimal numbers $q_{\rm trig}$, which satisfy 
the condition

$$
\frac{(T-t)^{3}}{2\pi^2}\Biggl(\frac{\pi^2}{6}-
\sum_{r=1}^{q_{\rm trig}} \frac{1}{r^2}\Biggr)\le (T-t)^4.
$$

\vspace{3mm}

From the other hand, the right-hand side of (\ref{4002}) 
includes only 2 random variables.
In this situation we again can talk about
the advantage of Ledendre polynomials.

\begin{table}
\centering
\caption{Numbers $q,$ $q_1$}
\label{tab:39}      
\begin{tabular}{p{1.1cm}p{2.1cm}p{2.1cm}p{2.1cm}p{2.1cm}p{2.1cm}p{2.1cm}}
\hline\noalign{\smallskip}
$T-t$&$0.08222$&$0.05020$&$0.02310$&$0.01956$\\
\noalign{\smallskip}\hline\noalign{\smallskip}
$q$&19&51&235&328\\
\noalign{\smallskip}\hline\noalign{\smallskip}
$q_1$&1&2&5&6\\
\noalign{\smallskip}\hline\noalign{\smallskip}
\end{tabular}
\end{table}

\begin{table}
\centering
\caption{Numbers $p,$ $p_1,$ $p^{*},$ $p_1^{*}$}
\label{tab:40}      
\begin{tabular}{p{1.1cm}p{2.1cm}p{2.1cm}p{2.1cm}p{2.1cm}p{2.1cm}p{2.1cm}}
\hline\noalign{\smallskip}
$T-t$&$0.08222$&$0.05020$&$0.02310$&$0.01956$\\
\noalign{\smallskip}\hline\noalign{\smallskip}
$p$&8&21&96&133\\
\noalign{\smallskip}\hline\noalign{\smallskip}
$p_1$&1&1&3&4\\
\noalign{\smallskip}\hline\noalign{\smallskip}
$p^{*}$&23&61&286&398\\
\noalign{\smallskip}\hline\noalign{\smallskip}
$p_1^{*}$&1&2&4&5\\
\noalign{\smallskip}\hline\noalign{\smallskip}
\end{tabular}
\end{table}

\begin{table}
\centering
\caption{Confirmation of the formula (\ref{zzz4})}
\label{tab:41}      
\begin{tabular}{p{2.1cm}p{1.7cm}p{1.7cm}p{2.1cm}p{2.3cm}p{2.3cm}p{2.3cm}}
\hline\noalign{\smallskip}
$\varepsilon/(T-t)^3$&$0.0629$&$0.0097$&$0.0010$&$1.0129\cdot 10^{-4}$&
$1.0132\cdot 10^{-5}$\\
\noalign{\smallskip}\hline\noalign{\smallskip}
$q$&1&10&100&1000&10000\\
\noalign{\smallskip}\hline\noalign{\smallskip}
\end{tabular}
\end{table}

In Table 41 we can see the numerical confirmation of 
(\ref{zzz4}) ($\varepsilon$ is the left-hand side of (\ref{zzz4})).

Let us compare computational costs for the approximation
$I_{(10)T,t}^{*(i_1 i_2)q}$ 
obtained from
(\ref{4006}) by replacing $\infty$ with $q$
(the case of Legendre olynomials) and for the approximation
$I_{(10)T,t}^{*(i_1 i_2)q}$ 
obtained
by Theorem 3 (the case of trigonometric functions)

$$
I_{(10)T,t}^{*(i_1 i_2)q}=-(T-t)^{2}\Biggl(\frac{1}{6}
\zeta_{0}^{(i_1)}\zeta_{0}^{(i_2)}-\frac{1}{2\sqrt{2}\pi}
\sqrt{\alpha_q}\xi_q^{(i_2)}\zeta_0^{(i_1)}+\Biggr.
$$
$$
+\frac{1}{2\sqrt{2}\pi^2}\sqrt{\beta_q}\Biggl(
\mu_q^{(i_2)}\zeta_0^{(i_1)}-2\mu_q^{(i_1)}\zeta_0^{(i_2)}\Biggr)+
$$
$$
+\frac{1}{2\sqrt{2}}\sum_{r=1}^{q}
\Biggl(-\frac{1}{\pi r}
\zeta_{2r-1}^{(i_2)}
\zeta_{0}^{(i_1)}+
\frac{1}{\pi^2 r^2}\left(
\zeta_{2r}^{(i_2)}
\zeta_{0}^{(i_1)}-
2\zeta_{2r}^{(i_1)}
\zeta_{0}^{(i_2)}\right)\Biggr)-
$$
$$
-
\frac{1}{2\pi^2}\sum_{\stackrel{r,l=1}{{}_{r\ne l}}}^{q}
\frac{1}{r^2-l^2}\Biggl(
\zeta_{2r}^{(i_1)}
\zeta_{2l}^{(i_2)}+
\frac{l}{r}
\zeta_{2r-1}^{(i_1)}
\zeta_{2l-1}^{(i_2)}
\Biggr)+
$$
$$
+
\sum_{r=1}^{q}
\Biggl(\frac{1}{4\pi r}\left(
\zeta_{2r}^{(i_1)}
\zeta_{2r-1}^{(i_2)}-
\zeta_{2r-1}^{(i_1)}
\zeta_{2r}^{(i_2)}\right)+
$$
\begin{equation}
\label{t1}
+
\Biggl.\Biggl.
\frac{1}{8\pi^2 r^2}\left(
3\zeta_{2r-1}^{(i_1)}
\zeta_{2r-1}^{(i_2)}+
\zeta_{2r}^{(i_2)}
\zeta_{2r}^{(i_1)}\right)\Biggr)\Biggr).
\end{equation}

\vspace{5mm}

For the formula (\ref{t1}) ($i_1\ne i_2$) from Theorem 8 we obtain
\cite{3}-\cite{10axx1}

\vspace{0.5mm}
$$
{\sf M}\left\{\left(I_{(01)T,t}^{(i_1 i_2)}-
I_{(01)T,t}^{(i_1 i_2)q}\right)^2\right\}=
\frac{(T-t)^4}{4}\Biggl(\frac{1}{9}-
\frac{1}{2\pi^2}\sum_{r=1}^q \frac{1}{r^2}-\Biggr.
$$

\vspace{1mm}
\begin{equation}
\label{t2}
\Biggl.-\frac{5}{8\pi^4}\sum_{r=1}^q \frac{1}{r^4}-
\frac{1}{\pi^4}\sum_{\stackrel{k,l=1}{{}_{k\ne l}}}^q
\frac{k^2+l^2}{l^2\left(l^2-k^2\right)^2}\Biggr).
\end{equation}

\vspace{3mm}

In Table 42 we can see the numerical confirmation of 
(\ref{t2}) ($\varepsilon$ is the right-hand side of (\ref{t2})).

\begin{table}
\centering
\caption{Confirmation of the formulas (\ref{t2})}
\label{tab:42}      
\begin{tabular}{p{2.1cm}p{1.7cm}p{1.7cm}p{2.1cm}p{2.3cm}p{2.3cm}p{2.3cm}}
\hline\noalign{\smallskip}
$4\varepsilon/(T-t)^4$&0.0540&0.0082&$8.4261\cdot 10^{-4}$
&$8.4429\cdot 10^{-5}$&
$8.4435\cdot 10^{-6}$\\
\noalign{\smallskip}\hline\noalign{\smallskip}
$q$&1&10&100&1000&10000\\
\noalign{\smallskip}\hline\noalign{\smallskip}
\end{tabular}
\end{table}

Let us compare the complexity of approximation based on the formula
(\ref{4006}) with the complexity of approximation (\ref{t1}).
The formula (\ref{t1}) includes the double sum

$$
\frac{1}{2\pi^2}\sum_{\stackrel{r,l=1}{{}_{r\ne l}}}^{q}
\frac{1}{r^2-l^2}\Biggl(
\zeta_{2r}^{(i_1)}
\zeta_{2l}^{(i_2)}+
\frac{l}{r}
\zeta_{2r-1}^{(i_1)}
\zeta_{2l-1}^{(i_2)}
\Biggr).
$$

\vspace{3mm}

Thus, the approximation (\ref{t1}) is more complex than the
approximation based on the formula (\ref{4006})
even if we take identical numbers $q$ in these approximations.
As we noted above, the number $q$ in (\ref{t1}) must be equal to the
number $q$ from the 
formula (\ref{555}), so it is much larger than the number $q$ 
from the approximation based on the formula (\ref{4006}). 
As s result, we have an obvious advantage of the Legendre polynomials
in computational costs in the considered case. 
As we mentioned above,
if we will not use the random 
variables $\xi_q^{(i)}$ and 
$\mu_q^{(i)},$ then the number $q$ in (\ref{t1}) can be chosen smaller, but
the mean-square error of approximation of the
stochastic integral $I_{(00)T,t}^{(i_1 i_2)}$ will be three times larger
(see (\ref{801})). Moreover, in this case the stochastic integrals
$I_{(1)T,t}^{(i_1)}$, $I_{(2)T,t}^{(i_1)}$ (with Gaussian distribution)
will be approximated worse. In this situation we can again talk about
the advantage of Ledendre polynomials.

Summing up the results of this section, we obtain to the 
following conclusions.

\vspace{2mm}

(I) We can talk about the approximately equal computational costs
for the formulas (\ref{555}) and (\ref{401}). This means that
computational costs for the implementation of Milstein scheme (explicit
one-step strong numerical method with the order of accuracy 1.0 
for Ito SDEs) 
for the case of Legendre polynomials and for the case 
of trigonometric functions are approximately the same.

\vspace{2mm}

(II) If we will not use the random 
variables $\xi_q^{(i)}$ (see (\ref{555})), then 
the mean-square error of approximation of the
stochastic integral $I_{(00)T,t}^{(i_1 i_2)}$ will be three times larger
(see (\ref{801})). 
In this situation we can talk about
the advantage of Ledendre polynomials within the frames 
of the Milstein scheme for Ito SDEs.
Moreover, in this case the stochastic integrals 
$I_{(1)T,t}^{(i_1)}$, $I_{(2)T,t}^{(i_1)}$ (with Gaussian distribution)
will be approximated worse.

\vspace{2mm}

(III) If we talk about an explicit one-step 
strong Taylor--Ito scheme 
of the order of accuracy $\gamma=1.5$ for Ito SDEs, then 
the numbers $q,$ $q_1$ (see (\ref{sad001}), (\ref{401})) are different. 
At that $q_1\ll q$
(the case of Legendre polynomials). 
The number $q$ must be the same in (\ref{444}), (\ref{555}), (\ref{x5}) 
(the case of trigonometric functions).
This leads to 
huge computational costs (see very complex formula (\ref{x5})).
From the other hand, we can choose different numbers $q$
in (\ref{444}), (\ref{555}), (\ref{x5}). At that we must exclude
the random variables $\xi_q^{(i)},$ $\mu_q^{(i)}$ from 
(\ref{444}), (\ref{555}), (\ref{x5}). This leads to another 
problems which we discussed above (see Conclusion (II)).

\vspace{2mm}

(IV) In addition, 
the author of this article supposes that the effect described in Conclusion 
(III) will be more impressive when 
analyzing more complex families of iterated Ito and Stratonovich 
stochastic
integrals (when $\gamma=$ $2.0,$ $2.5,$ $3.0,$ $\ldots $). 
This supposition is based on the fact that the polynomial 
system of functions has the significant advantage (in comparison with 
the trigonometric system of functions) 
for approximation of iterated stochastic 
integrals for which not all weight functions are equal to 1.

\vspace{5mm}

\section{Convergence With Probability 1 of Expansions 
of Iterated Stochastic Integrals of Multiplicities 1 and 2}

\vspace{5mm}

Let us address now to the convergence with probability 1. 
Note that proving Theorem 1 \cite{10a} (Theorem 1.1, Sect.~1.1.3) or Theorem 2
\cite{10a} (Theorem 1.16, Sect.~1.11) 
we obtained the following representation

\vspace{-2mm}
$$
J[\psi^{(k)}]_{T,t}=
\sum\limits_{j_1=0}^{p_1}\ldots
\sum\limits_{j_k=0}^{p_k}
C_{j_k\ldots j_1}\Biggl(
\prod_{l=1}^k\zeta_{j_l}^{(i_l)}+\sum\limits_{r=1}^{[k/2]}
(-1)^r \times
\Biggr.
$$

\vspace{1mm}
$$
\times
\sum_{\stackrel{(\{\{g_1, g_2\}, \ldots, 
\{g_{2r-1}, g_{2r}\}\}, \{q_1, \ldots, q_{k-2r}\})}
{{}_{\{g_1, g_2, \ldots, 
g_{2r-1}, g_{2r}, q_1, \ldots, q_{k-2r}\}=\{1, 2, \ldots, k\}}}}
\prod\limits_{s=1}^r
{\bf 1}_{\{i_{g_{{}_{2s-1}}}=~i_{g_{{}_{2s}}}\ne 0\}}
\Biggl.{\bf 1}_{\{j_{g_{{}_{2s-1}}}=~j_{g_{{}_{2s}}}\}}
\prod_{l=1}^{k-2r}\zeta_{j_{q_l}}^{(i_{q_l})}\Biggr)+R_{T,t}^{p_1,\ldots,p_k}
$$

\vspace{3mm}
\noindent
\ \ \ \hbox{w.\ p.\ 1,} where

\vspace{-2mm}
\begin{equation}
\label{y007}
R_{T,t}^{p_1,\ldots,p_k}=\sum_{(t_1,\ldots,t_k)}
\int\limits_{t}^{T}
\ldots
\int\limits_{t}^{t_2}
R_{p_1,\ldots,p_k}(t_1,\ldots,t_k)
d{\bf w}_{t_1}^{(i_1)}
\ldots
d{\bf w}_{t_k}^{(i_k)},
\end{equation}

\vspace{1mm}
$$
R_{p_1,\ldots,p_k}(t_1,\ldots,t_k)=K(t_1,\ldots,t_k)-
\sum_{j_1=0}^{p_1}\ldots
\sum_{j_k=0}^{p_k}
C_{j_k\ldots j_1}
\prod_{l=1}^k\phi_{j_l}(t_l),
$$

\vspace{2mm}
\noindent
where permutations $(t_1,\ldots,t_k)$ when summing 
in (\ref{y007}) are performed only 
in the values $d{\bf w}_{t_1}^{(i_1)}$ $\ldots$
$d{\bf w}_{t_k}^{(i_k)}$. At the same time the indexes near 
upper limits of integration in the iterated  stochastic integrals 
are changed correspondently and if $t_r$ swapped with $t_q$ in the  
permutation $(t_1,\ldots,t_k)$, then $i_r$ swapped with $i_q$ in the 
permutation $(i_1,\ldots,i_k)$. Another notations are the same as
in Theorems 1, 2.

Let us consider in detail the following
expansion of iterated  Ito stochastic integral

\begin{equation}
I_{(00)T,t}^{(i_1 i_2)}=
\frac{T-t}{2}\Biggl(\zeta_0^{(i_1)}\zeta_0^{(i_2)}+\sum_{i=1}^{\infty}
\frac{1}{\sqrt{4i^2-1}}\left(
\zeta_{i-1}^{(i_1)}\zeta_{i}^{(i_2)}-
\zeta_i^{(i_1)}\zeta_{i-1}^{(i_2)}\right)-{\bf 1}_{\{i_1=i_2\}}\Biggr).
\label{hqye}
\end{equation}

\vspace{3mm}

If $i_1=i_2$, then from (\ref{hqye}) we obtain the following equality

$$
I_{(00)T,t}^{(i_1 i_1)}=
\frac{1}{2}(T-t)\left(\left(\zeta_0^{(i_1)}\right)^2-1\right),
$$

\vspace{3mm}
\noindent
which is correct w.~p.~1 
and can be obtained using the Ito formula.

Let us consider the case $i_1\ne i_2$. In this case

$$
I_{(00)T,t}^{*(i_1 i_2)}=I_{(00)T,t}^{(i_1 i_2)}\ \ \ \hbox{w.\ p.\ 1.}
$$

\vspace{3mm}

First, note the well-known fact.

{\bf Lemma 1.}\ {\it If for the sequence of random variables
$\xi_p$ and for some
$\alpha>0$ the number series

$$
\sum\limits_{p=1}^{\infty}{\sf M}\left\{\left|\xi_p\right|^{\alpha}\right\}
$$

\vspace{3mm}
\noindent
converges, then the sequence $\xi_p$ converges to zero w.~p.~{\rm 1}.}

In our specific case $(i_1\ne i_2)$

$$
I_{(00)T,t}^{(i_1 i_2)}=I_{(00)T,t}^{(i_1 i_2)p}+\xi_p,\ \ \
\xi_p
=\frac{T-t}{2}\sum_{i=p+1}^{\infty}
\frac{1}{\sqrt{4i^2-1}}\left(
\zeta_{i-1}^{(i_1)}\zeta_{i}^{(i_2)}-
\zeta_i^{(i_1)}\zeta_{i-1}^{(i_2)}\right),
$$

\vspace{2mm}
\begin{equation}
\label{90}
I_{(00)T,t}^{(i_1 i_2)p}=
\frac{T-t}{2}\Biggl(\zeta_0^{(i_1)}\zeta_0^{(i_2)}+\sum_{i=1}^{p}
\frac{1}{\sqrt{4i^2-1}}\left(
\zeta_{i-1}^{(i_1)}\zeta_{i}^{(i_2)}-
\zeta_i^{(i_1)}\zeta_{i-1}^{(i_2)}\right)\Biggr).
\end{equation}

\vspace{5mm}

Let 
$$
R_{T,t}^{p_1,p_2}\stackrel{\sf def}{=}R_{T,t}^p,\ \ \ 
R_{p_1p_2}(t_1,t_2)\stackrel{\sf def}{=}R_p(t_1,t_2)\ \ \
\hbox{for}\ \ \ p_1=p_2=p.
$$ 

\vspace{2mm}

Then

\vspace{-1mm}
$$
\xi_p=R_{T,t}^p=\int\limits_t^T\int\limits_t^{t_2}
R_p(t_1,t_2)d{\bf f}_{t_1}^{(i_1)}d{\bf f}_{t_2}^{(i_2)}+
\int\limits_t^T\int\limits_{t}^{t_1}
R_p(t_1,t_2)d{\bf f}_{t_2}^{(i_2)}d{\bf f}_{t_1}^{(i_1)},
$$

\vspace{2mm}
\begin{equation}
\label{2017rock}
{\sf M}\{|\xi_p|^2\}=
\int\limits_t^T\int\limits_t^{t_2}
\left(R_p(t_1,t_2)\right)^2dt_1 dt_2 +
\int\limits_t^T\int\limits_{t}^{t_1}
\left(R_p(t_1,t_2)\right)^2 dt_2 dt_1
=\int\limits_{[t, T]^2}\left(R_p(t_1,t_2)\right)^2dt_1 dt_2,
\end{equation}

\vspace{2mm}
\begin{equation}
\label{2017rock1}
{\sf M}\left\{|\xi_p|^2\right\}=
\frac{(T-t)^2}{2}\sum_{i=p+1}^{\infty}
\frac{1}{4i^2-1},
\end{equation}

\vspace{2mm}
$$
R_p(t_1,t_2)=K(t_1,t_2)-\sum\limits_{j_1,j_2=0}^p
C_{j_2 j_1}\phi_{j_1}(t_1)\phi_{j_2}(t_2),
$$

\vspace{2mm}
\begin{equation}
\label{yuye}
\sum\limits_{i=p+1}^{\infty}\frac{1}{4i^2-1}\le \int\limits_{p}^{\infty}
\frac{1}{4x^2-1}dx
=-\frac{1}{4}{\rm ln}\biggl|
1-\frac{2}{2p+1}\biggr|\le \frac{C}{p},
\end{equation}

\vspace{5mm}
\noindent
where constant $C$ does not depend on $p.$

Therefore, taking $\alpha=2$ in Lemma 1, we 
cannot
prove the convergence of 
$\xi_p$ to zero w.~p.~1, 
since the series 

\vspace{-2mm}
$$
\sum\limits_{p=1}^{\infty}{\sf M}\left\{|\xi_p|^{2}\right\}
$$

\vspace{3mm}
\noindent
will be 
majorized by the 
divergent
Dirichlet series
with the index 1. Let us take $\alpha=4$ and estimate the value
${\sf M}\left\{|\xi_p|^4\right\}$.

According to (\ref{99999}), we can write

\begin{equation}
\label{e1}
{\sf M}\left\{(R_{T,t}^{p_1,\ldots,p_k})^{2n}\right\}\le
C_{n,k}\left(
\int\limits_{[t, T]^k}
R^2_{p_1\ldots p_k}(t_1,\dots,t_k)dt_1\ldots dt_k
\right)^n,
\end{equation}

\vspace{3mm}
\noindent
where $(k!)^{n} (2n-1)^{nk}.$

From (\ref{e1}) for $k=2$, $n=2$ and
(\ref{2017rock})--(\ref{yuye}) we obtain

\begin{equation}
\label{91}
{\sf M}\left\{|\xi_p|^4\right\}\le K
\left(\int\limits_{[t,T]^2}
R_p^2(t_1,t_2)dt_1dt_2\right)^2\le \frac{K_1}{p^2}
\end{equation}

\vspace{3mm}
\noindent
and
\begin{equation}
\label{hhq}
\sum_{p=1}^{\infty}
{\sf M}\left\{|\xi_p|^4\right\}\le {K_1}
\sum_{p=1}^{\infty}\frac{1}{p^2}<\infty,
\end{equation}

\vspace{3mm}
\noindent
where constants $K,$ $K_1$ do not depend on $p.$

Since the series in (\ref{hhq}) converges, then according to Lemma 1
we obtain that
$\xi_p \to 0$ when $p\to \infty$\ w.~p.~1. Then  

\vspace{-2mm}
$$
I_{(00)T,t}^{(i_1 i_2)p}\to
I_{(00)T,t}^{(i_1 i_2)}\ \ \hbox{when}\ \ p\to \infty\ \ \hbox{w.\ p.\ 1.}
$$

\vspace{3mm}

Let us consider the stochastic integrals 
$I_{(01)T,t}^{*(i_1 i_2)},$
$I_{(10)T,t}^{*(i_1 i_2)}$
whose expansions 
look as (\ref{ud111}), (\ref{4006}).

Consider the case $i_1\ne i_2$. In this case

$$
I_{(01)T,t}^{*(i_1 i_2)}=I_{(01)T,t}^{(i_1 i_2)},\ \ \
I_{(10)T,t}^{*(i_1 i_2)}=I_{(10)T,t}^{(i_1 i_2)},\ \ \
I_{(00)T,t}^{*(i_1 i_2)}=I_{(00)T,t}^{(i_1 i_2)}\ \ \ \hbox{w.\ p.\ 1,}
$$

\vspace{2mm}
\noindent
and

\vspace{-1mm}
$$
I_{(01)T,t}^{(i_1 i_2)}=-\frac{T-t}{2}I_{(00)T,t}^{(i_1 i_2)p}
-\frac{(T-t)^2}{4}\Biggl(
\frac{1}{\sqrt{3}}\zeta_0^{(i_1)}\zeta_1^{(i_2)}+\Biggr.
$$

\vspace{2mm}
$$
+\Biggl.\sum_{i=0}^{p}\Biggl(
\frac{(i+2)\zeta_i^{(i_1)}\zeta_{i+2}^{(i_2)}
-(i+1)\zeta_{i+2}^{(i_1)}\zeta_{i}^{(i_2)}}
{\sqrt{(2i+1)(2i+5)}(2i+3)}-
\frac{\zeta_i^{(i_1)}\zeta_{i}^{(i_2)}}{(2i-1)(2i+3)}\Biggr)\Biggr)
+ \xi_p^{(01)},
$$

\vspace{6mm}

$$
I_{(10)T,t}^{(i_1 i_2)}=-\frac{T-t}{2}I_{(00)T,t}^{(i_1 i_2)p}
-\frac{(T-t)^2}{4}\Biggl(
\frac{1}{\sqrt{3}}\zeta_0^{(i_2)}\zeta_1^{(i_1)}+\Biggr.
$$

\vspace{2mm}
$$
+\Biggl.\sum_{i=0}^{p}\Biggl(
\frac{(i+1)\zeta_{i+2}^{(i_2)}\zeta_{i}^{(i_1)}
-(i+2)\zeta_{i}^{(i_2)}\zeta_{i+2}^{(i_1)}}
{\sqrt{(2i+1)(2i+5)}(2i+3)}+
\frac{\zeta_i^{(i_1)}\zeta_{i}^{(i_2)}}{(2i-1)(2i+3)}\Biggr)\Biggr)
+\xi_p^{(10)},
$$

\vspace{6mm}
\noindent
where

\vspace{-1mm}
$$
\xi_p^{(01)}=
-\frac{(T-t)^2}{4}
\Biggl(\sum_{i=p+1}^{\infty}
\frac{1}{\sqrt{4i^2-1}}\left(
\zeta_{i-1}^{(i_1)}\zeta_{i}^{(i_2)}-
\zeta_i^{(i_1)}\zeta_{i-1}^{(i_2)}\right)\Biggr.+
$$

\vspace{2mm}
$$
+
\Biggl.\sum_{i=p+1}^{\infty}\Biggl(
\frac{(i+2)\zeta_i^{(i_1)}\zeta_{i+2}^{(i_2)}
-(i+1)\zeta_{i+2}^{(i_1)}\zeta_{i}^{(i_2)}}
{\sqrt{(2i+1)(2i+5)}(2i+3)}-
\frac{\zeta_i^{(i_1)}\zeta_{i}^{(i_2)}}{(2i-1)(2i+3)}\Biggr)\Biggr),
$$

\vspace{6mm}

$$
\xi_p^{(10)}=
-\frac{(T-t)^2}{4}
\Biggl(\sum_{i=p+1}^{\infty}
\frac{1}{\sqrt{4i^2-1}}\left(
\zeta_{i-1}^{(i_1)}\zeta_{i}^{(i_2)}-
\zeta_i^{(i_1)}\zeta_{i-1}^{(i_2)}\right)\Biggr.+
$$                                     

\vspace{2mm}
$$
+
\Biggl.\sum_{i=p+1}^{\infty}\Biggl(
\frac{(i+1)\zeta_i^{(i_1)}\zeta_{i+2}^{(i_2)}
-(i+2 )\zeta_{i+2}^{(i_1)}\zeta_{i}^{(i_2)}}
{\sqrt{(2i+1)(2i+5)}(2i+3)}-
\frac{\zeta_i^{(i_1)}\zeta_{i}^{(i_2)}}{(2i-1)(2i+3)}\Biggr)\Biggr).
$$

\vspace{6mm}

Then
$$
{\sf M}\left\{\left|\xi_p^{(01)}\right|^2\right\}=
\int\limits_{[t, T]^2}\left(R_p^{(01)}(t_1,t_2)\right)^2dt_1 dt_2=
\frac{(T-t)^4}{16}\times
$$

\begin{equation}
\label{2017mas1}
\times
\sum\limits_{i=p+1}^{\infty}
\Biggl(\frac{2}{4i^2-1}+\frac{(i+2)^2+(i+1)^2}{(2i+1)(2i+5)(2i+3)^2}
+\frac{1}{(2i-1)^2(2i+3)^2}\Biggr) 
\le 
C\sum\limits_{i=p+1}^{\infty}\frac{1}{i^2}\le \frac{K}{p},
\end{equation}

\vspace{5mm}
\noindent
where constants $C,$ $K$ do not depend on $p$.

Analogously, we obtain

\begin{equation}
\label{2017mas2}
{\sf M}\left\{\left|\xi_p^{(10)}\right|^2\right\}=
\int\limits_{[t, T]^2}\biggl(R_p^{(10)}(t_1,t_2)\biggr)^2dt_1 dt_2
\le \frac{K}{p},
\end{equation}

\vspace{3mm}
\noindent
where constant $K$ does not depend on $p$.

According (\ref{e1}) when $k=2$, $n=2$ and
(\ref{2017mas1}), (\ref{2017mas2}), we obtain

$$
{\sf M}\left\{\left|\xi_p^{(01)}\right|^4\right\}\le K
\left(\int\limits_{[t,T]^2}
\left(R_p^{(01)}(t_1,t_2)\right)^2 dt_1dt_2\right)^2\le \frac{K_1}{p^2},
$$

\vspace{2mm}
$$
{\sf M}\left\{\left|\xi_p^{(10)}\right|^4\right\}\le K
\left(\int\limits_{[t,T]^2}
\left(R_p^{(10)}(t_1,t_2)\right)^2 dt_1dt_2\right)^2\le \frac{K_1}{p^2},
$$

\vspace{3mm}
\noindent
and
\begin{equation}
\label{hhq11}
\sum_{p=1}^{\infty}
{\sf M}\left\{\left|\xi_p^{(01)}\right|^4\right\}\le {K_1}
\sum_{p=1}^{\infty}\frac{1}{p^2}<\infty,\ \ \
\sum_{p=1}^{\infty}
{\sf M}\left\{\left|\xi_p^{(10)}\right|^4\right\}\le {K_1}
\sum_{p=1}^{\infty}\frac{1}{p^2}<\infty,
\end{equation}

\vspace{3mm}
\noindent
where constant $K_1$ does not depend on $p.$

From (\ref{hhq11}) and Lemma 1 we obtain that
$\xi_p^{(01)},\ \xi_p^{(10)}\to 0$ when $p\to \infty$\ w.~p.~1. 
Then  

\vspace{-1mm}
$$
I_{(01)T,t}^{(i_1 i_2)p}\to
I_{(01)T,t}^{(i_1 i_2)},\ \ \
I_{(10)T,t}^{(i_1 i_2)p}\to
I_{(10)T,t}^{(i_1 i_2)}\ \ \ \hbox{when}\ \ \  p\to \infty\ \ \ 
\hbox{w.\ p.\ 1,}
$$

\vspace{3mm}
\noindent
where $i_1\ne i_2.$

Let us consider the case $i_1=i_2$

\vspace{2mm}
$$
I_{(01)T,t}^{(i_1 i_1)}=
-\frac{(T-t)^2}{4}\Biggl(
\left(\zeta_0^{(i_1)}\right)^2-1+
\frac{1}{\sqrt{3}}\zeta_0^{(i_1)}\zeta_1^{(i_1)}+\Biggr.
$$

\vspace{2mm}
$$
\Biggl.
+\sum_{i=0}^{p}\Biggl(\frac{1}{\sqrt{(2i+1)(2i+5)}(2i+3)}
\zeta_i^{(i_1)}\zeta_{i+2}^{(i_1)}-
\frac{1}{(2i-1)(2i+3)}\left(\zeta_i^{(i_1)}\right)^2
\Biggr)\Biggr)+
\mu_p^{(01)},
$$

\vspace{8mm}

$$
I_{(10)T,t}^{(i_1 i_1)}=
-\frac{(T-t)^2}{4}\Biggl(
\left(\zeta_0^{(i_1)}\right)^2-1+
\frac{1}{\sqrt{3}}\zeta_0^{(i_1)}\zeta_1^{(i_1)}+\Biggr.
$$

\vspace{2mm}
$$
\Biggl.
+\sum_{i=0}^{p}\Biggl(-\frac{1}{\sqrt{(2i+1)(2i+5)}(2i+3)}
\zeta_i^{(i_1)}\zeta_{i+2}^{(i_1)}+
\frac{1}{(2i-1)(2i+3)}\left(\zeta_i^{(i_1)}\right)^2
\Biggr)\Biggr)+
\mu_p^{(10)},
$$

\vspace{5mm}
\noindent
where

$$
\mu_p^{(01)}=-\frac{(T-t)^2}{4}
\sum_{i=p+1}^{\infty}
\Biggl(\frac{1}{\sqrt{(2i+1)(2i+5)}(2i+3)}
\zeta_i^{(i_1)}\zeta_{i+2}^{(i_1)}-
\frac{1}{(2i-1)(2i+3)}\left(\zeta_i^{(i_1)}\right)^2
\Biggr),
$$

\vspace{2mm}
$$
\mu_p^{(10)}=-\frac{(T-t)^2}{4}
\sum_{i=p+1}^{\infty}\Biggl(-\frac{1}{\sqrt{(2i+1)(2i+5)}(2i+3)}
\zeta_i^{(i_1)}\zeta_{i+2}^{(i_1)}+
\frac{1}{(2i-1)(2i+3)}\left(\zeta_i^{(i_1)}\right)^2
\Biggr).
$$

\vspace{5mm}

Then

$$
{\sf M}\biggl\{\left(\mu_p^{(01)}\right)^2\biggr\}=
{\sf M}\biggl\{\left(\mu_p^{(10)}\right)^2\biggr\}=\frac{(T-t)^4}{16}\times
$$

$$
\times\left(
\sum_{i=p+1}^{\infty}
\frac{1}{(2i+1)(2i+5)(2i+3)^2}+
\sum_{i=p+1}^{\infty}
\frac{2}{(2i-1)^2(2i+3)^2} + 
\Biggl(\sum_{i=p+1}^{\infty}
\frac{1}{(2i-1)(2i+3)}\Biggr)^2\right)
\le \frac{K}{p^2}
$$

\vspace{4mm}
\noindent
and

\vspace{-2mm}
\begin{equation}
\label{hhq112}
\sum_{p=1}^{\infty}
{\sf M}\left\{\left|\mu_p^{(01)}\right|^2\right\}\le {K}
\sum_{p=1}^{\infty}\frac{1}{p^2}<\infty,\ \ \
\sum_{p=1}^{\infty}
{\sf M}\left\{\left|\mu_p^{(10)}\right|^2\right\}\le {K}
\sum_{p=1}^{\infty}\frac{1}{p^2}<\infty,
\end{equation}

\vspace{5mm}
\noindent
where constant $K$ does not depend on $p.$

According Lemma 1 and (\ref{hhq112}), we obtain that
$\mu_p^{(01)},\ \mu_p^{(10)}\to 0$ when $p\to \infty$\ w.~p.~1. 
Then  

$$
I_{(01)T,t}^{(i_1 i_1)p}\to
I_{(01)T,t}^{(i_1 i_1)},\ \ \
I_{(10)T,t}^{(i_1 i_1)p}\to
I_{(10)T,t}^{(i_1 i_1)}\ \ \ \hbox{when}\ \ \  
p\to \infty\ \ \ \hbox{w.\ p.\ 1.}
$$

\vspace{3mm}

Analogously, we obtain

$$
I_{(02)T,t}^{(i_1 i_2)p}\to
I_{(02)T,t}^{(i_1 i_2)},\ \ \
I_{(11)T,t}^{(i_1 i_2)p}\to
I_{(11)T,t}^{(i_1 i_2)},\ \ \
I_{(20)T,t}^{(i_1 i_2)p}\to
I_{(20)T,t}^{(i_1 i_2)}\ \ \ \hbox{when}\ \ \ 
p\to \infty\ \ \ \hbox{w.\ p.\ 1,}
$$

\vspace{3mm}
\noindent
where $i_1, i_2=1,\ldots,m.$
This result based on the following 
truncated expansions of the stochastic
integrals
$I_{(02)T,t}^{(i_1 i_2)}$,
$I_{(20)T,t}^{(i_1 i_2)}$,
$I_{(11)T,t}^{(i_1 i_2)}$ (see (\ref{seak1})--(\ref{seak3}))

\vspace{4mm}

$$
I_{(02)T,t}^{(i_1 i_2)p}
=-\frac{(T-t)^2}{4}I_{(00)T,t}^{(i_1 i_2)p}
-(T-t) I_{(01)T,t}^{(i_1 i_2)p}+
\frac{(T-t)^3}{8}\Biggl[
\frac{2}{3\sqrt{5}}\zeta_2^{(i_2)}\zeta_0^{(i_1)}+\Biggr.
$$

\vspace{1mm}
$$
+\frac{1}{3}\zeta_0^{(i_1)}\zeta_0^{(i_2)}+
\sum_{i=0}^{p}\Biggl(
\frac{(i+2)(i+3)\zeta_{i+3}^{(i_2)}\zeta_{i}^{(i_1)}
-(i+1)(i+2)\zeta_{i}^{(i_2)}\zeta_{i+3}^{(i_1)}}
{\sqrt{(2i+1)(2i+7)}(2i+3)(2i+5)}+
\Biggr.
$$

\vspace{1mm}
$$
\Biggl.\Biggl.+
\frac{(i^2+i-3)\zeta_{i+1}^{(i_2)}\zeta_{i}^{(i_1)}
-(i^2+3i-1)\zeta_{i}^{(i_2)}\zeta_{i+1}^{(i_1)}}
{\sqrt{(2i+1)(2i+3)}(2i-1)(2i+5)}\Biggr)\Biggr] - 
$$

\vspace{1mm}
$$
-\frac{1}{24}{\bf 1}_{\{i_1=i_2\}}{(T-t)^3},
$$

\vspace{4mm}

$$
I_{(20)T,t}^{(i_1 i_2)p}=-\frac{(T-t)^2}{4}
I_{(00)T,t}^{(i_1 i_2)p}
-(T-t) I_{(10)T,t}^{(i_1 i_2)p}+
\frac{(T-t)^3}{8}\Biggl[
\frac{2}{3\sqrt{5}}\zeta_0^{(i_2)}\zeta_2^{(i_1)}+\Biggr.
$$
$$
+\frac{1}{3}\zeta_0^{(i_1)}\zeta_0^{(i_2)}+
\sum_{i=0}^{p}\Biggl(
\frac{(i+1)(i+2)\zeta_{i+3}^{(i_2)}\zeta_{i}^{(i_1)}
-(i+2)(i+3)\zeta_{i}^{(i_2)}\zeta_{i+3}^{(i_1)}}
{\sqrt{(2i+1)(2i+7)}(2i+3)(2i+5)}+
\Biggr.
$$

\vspace{1mm}
$$
\Biggl.\Biggl.+
\frac{(i^2+3i-1)\zeta_{i+1}^{(i_2)}\zeta_{i}^{(i_1)}
-(i^2+i-3)\zeta_{i}^{(i_2)}\zeta_{i+1}^{(i_1)}}
{\sqrt{(2i+1)(2i+3)}(2i-1)(2i+5)}\Biggr)\Biggr] - 
$$

\vspace{1mm}
$$
-
\frac{1}{24}{\bf 1}_{\{i_1=i_2\}}{(T-t)^3},
$$

\vspace{6mm}

$$
I_{(11)T,t}^{(i_1 i_2)p}
=-\frac{(T-t)^2}{4}I_{(00)T,t}^{(i_1 i_2)p}
-\frac{T-t}{2}\left(
I_{(10)T,t}^{(i_1 i_2)p}+
I_{(01)T,t}^{(i_1 i_2)p}\right)+
$$

\vspace{1mm}
$$
+
\frac{(T-t)^3}{8}\Biggl[
\frac{1}{3}\zeta_1^{(i_1)}\zeta_1^{(i_2)}+\Biggr.
\sum_{i=0}^{p}\Biggl(
\frac{(i+1)(i+3)\left(\zeta_{i+3}^{(i_2)}\zeta_{i}^{(i_1)}
-\zeta_{i}^{(i_2)}\zeta_{i+3}^{(i_1)}\right)}
{\sqrt{(2i+1)(2i+7)}(2i+3)(2i+5)}+
\Biggr.
$$

\vspace{1mm}
$$
\Biggl.\Biggl.
+\frac{(i+1)^2\left(\zeta_{i+1}^{(i_2)}\zeta_{i}^{(i_1)}
-\zeta_{i}^{(i_2)}\zeta_{i+1}^{(i_1)}\right)}
{\sqrt{(2i+1)(2i+3)}(2i-1)(2i+5)}\Biggr)\Biggr] - 
$$

\vspace{1mm}
$$
-
\frac{1}{24}{\bf 1}_{\{i_1=i_2\}}{(T-t)^3}.
$$

\vspace{5mm}

The expansions (\ref{4001})--(\ref{4003}), (\ref{gg1})
for the stochastic integrals
$I_{(0)T,t}^{(i_1)},$
$I_{(1)T,t}^{(i_1)},$
$I_{(2)T,t}^{(i_1)},$
$I_{(3)T,t}^{(i_1)}$
are initially 
correct w.~p.~1  (they include 1, 2, 3, and 4 members of 
expansion, correspondently). 

Apparently, using the proposed scheme 
we can prove convergence w.~p.~1 for other iterated
stochastic integrals. In the next section, we consider
the more general and effective approach.

\vspace{5mm}

\section{Convergence With Probability 1
of Expansion of Iterated
Ito Stochastic Integrals in Theorem 1 for the Case of Multiplicity
$k$ $(k\in\mathbb{N})$}

\vspace{5mm}

This section is written on the base of 
\cite{10a} (Sect.~1.7.2), \cite{9999}, \cite{new-new-2}.
Remind that in a lot of author's publications \cite{3}-\cite{arxiv-4}
the convergence in Theorem 1 has been considered in different
probabilistic
senses. For example, the mean-square convergence \cite{3} (2006) (also see
\cite{3a}-\cite{arxiv-4}) and convergence in the mean of degree
$2n$ $(n\in\mathbb{N})$ 
\cite{10a} (Sect.~1.1.9, 1.11, 1.12), \cite{11} (Sect.~6, 15, 16)
have been proved. On the examples
of specific iterated Ito stochastic integrals of mutiplicities 1 and 2
the convergence with probability 1 has been considered in the previous 
section
(also see \cite{4} (2007), \cite{4a}-\cite{7}, \cite{9a}-\cite{11},
\cite{13}). However, these examples are narrow particular 
cases of the iterated Ito stochastic integrals (\ref{ito}).

In this section, we formulate and prove the theorem 
\cite{10a} (Sect.~1.7.2), \cite{9999}, \cite{new-new-2} on 
convergence with probability 1 of the expansions 
of iterated Ito stochastic integrals from Theorem 1.

Let us remind the well-known fact from the mathematical analysis,
which is connected to existence
of iterated limits.

\vspace{2mm}

{\bf Proposition 1.}\ {\it Let $\bigl\{x_{n,m}\bigr\}_{n,m=1}^{\infty}$
be a double sequence and let there exists the limit

\vspace{-1mm}
$$
\lim\limits_{n,m\to\infty}x_{n,m}=a<\infty.
$$

\vspace{2mm}

Moreover, let there exist the limits

\vspace{-1mm}
$$
\lim\limits_{n\to\infty}x_{n,m}<\infty\ \ \ \hbox{for any}\ \ \ m,\ \ \ \
\lim\limits_{m\to\infty}x_{n,m}<\infty\ \ \ \hbox{for any}\ \ \ n.
$$

\vspace{2mm}

Then there exist the iterated limits

\vspace{-1mm}
$$
\lim\limits_{n\to\infty}\lim\limits_{m\to\infty}x_{n,m},\ \ \ 
\lim\limits_{m\to\infty}\lim\limits_{n\to\infty}x_{n,m}
$$

\vspace{2mm}
and moreover,

\vspace{-1mm}
$$
\lim\limits_{n\to\infty}\lim\limits_{m\to\infty}x_{n,m}=
\lim\limits_{m\to\infty}\lim\limits_{n\to\infty}x_{n,m}=a.
$$
}

\vspace{3mm}

{\bf Theorem 9} \cite{10a} (Sect.~1.7.2), \cite{9999}, \cite{new-new-2}.\
{\it Let 
$\psi_l(\tau)$ $(l=1,\ldots, k)$ are 
continuously differentiable nonrandom functions on the interval
$[t, T]$ and $\{\phi_j(x)\}_{j=0}^{\infty}$ is a complete
orthonormal system of Legendre polynomials or 
trigonometric functions in the space $L_2([t, T]).$
Then 

\vspace{2mm}
$$
J[\psi^{(k)}]_{T,t}^{p,\ldots,p}\ \to \ J[\psi^{(k)}]_{T,t}\ \ \ 
\hbox{if}\ \ \ p\to \infty
$$

\vspace{5mm}
\noindent
w.\ p.\ {\rm 1,} where $J[\psi^{(k)}]_{T,t}^{p,\ldots,p}$
is the expression on the right-hand side of {\rm (\ref{tyyy})}
before passing to the limit 
$\hbox{\vtop{\offinterlineskip\halign{
\hfil#\hfil\cr
{\rm l.i.m.}\cr
$\stackrel{}{{}_{p_1,\ldots,p_k\to \infty}}$\cr
}} }$ for the case $p_1=\ldots=p_k=p,$ i.e. {\rm (}see Theorem {\rm 1)}

\vspace{2mm}
$$
J[\psi^{(k)}]_{T,t}^{p,\ldots,p}=
\sum_{j_1=0}^{p}\ldots\sum_{j_k=0}^{p}
C_{j_k\ldots j_1}\Biggl(
\prod_{l=1}^k\zeta_{j_l}^{(i_l)}\ -
\Biggr.
$$

\vspace{3mm}
$$
-\ \Biggl.
\hbox{\vtop{\offinterlineskip\halign{
\hfil#\hfil\cr
{\rm l.i.m.}\cr
$\stackrel{}{{}_{N\to \infty}}$\cr
}} }\sum_{(l_1,\ldots,l_k)\in {\rm G}_k}
\phi_{j_{1}}(\tau_{l_1})
\Delta{\bf w}_{\tau_{l_1}}^{(i_1)}\ldots
\phi_{j_{k}}(\tau_{l_k})
\Delta{\bf w}_{\tau_{l_k}}^{(i_k)}\Biggr),
$$

\vspace{6mm}
\noindent
where $i_1,\ldots,i_k=1,\ldots,m$.}

\vspace{2mm}

{\bf Proof.} Let us consider the Parseval equality

\begin{equation}
\label{par1}
\int\limits_{[t,T]^k}K^2(t_1,\ldots,t_k)dt_1\ldots dt_k=
\lim\limits_{p_1,\ldots,p_k\to\infty}
\sum_{j_1=0}^{p_1}\ldots \sum_{j_k=0}^{p_k}
C_{j_k\ldots j_1}^2,
\end{equation}

\vspace{3mm}
\noindent
where
\begin{equation}
\label{pppx}
K(t_1,\ldots,t_k)=
\begin{cases}
\psi_1(t_1)\ldots \psi_k(t_k),\ &t_1<\ldots<t_k\\
~\\
~\\
0,\ &\hbox{\rm otherwise}
\end{cases}\ \ \ \ 
=\ \ \ \ 
\prod\limits_{l=1}^k
\psi_l(t_l)\ \prod\limits_{l=1}^{k-1}{\bf 1}_{\{t_l<t_{l+1}\}},\ 
\end{equation}

\vspace{3mm}
\noindent
where $t_1,\ldots,t_k\in [t, T]$ for $k\ge 2$ and 
$K(t_1)\equiv\psi_1(t_1)$ for $t_1\in[t, T],$ 
${\bf 1}_A$ denotes the indicator of the set $A$,
\begin{equation}
\label{ppppax}
C_{j_k\ldots j_1}=\int\limits_{[t,T]^k}
K(t_1,\ldots,t_k)\prod_{l=1}^{k}\phi_{j_l}(t_l)dt_1\ldots dt_k
\end{equation}

\vspace{2mm}
\noindent
is the Fourier coefficient.

Using (\ref{pppx}), we obtain
$$
C_{j_k\ldots j_1}=
\int\limits_t^T
\phi_{j_k}(t_k)\psi_k(t_k)\ldots \int\limits_t^{t_2}
\phi_{j_1}(t_1)\psi_1(t_1)dt_1\ldots dt_k.
$$

\vspace{2mm}

Further, we denote

\vspace{-1mm}
$$
\lim\limits_{p_1,\ldots,p_k\to\infty}
\sum_{j_1=0}^{p_1}\ldots \sum_{j_k=0}^{p_k}
C_{j_k\ldots j_1}^2\stackrel{\sf def}{=}
\sum_{j_1,\ldots,j_k=0}^{\infty}
C_{j_k\ldots j_1}^2.
$$

\vspace{4mm}

If $p_1=\ldots=p_k=p,$ then we also write

\vspace{1mm}
$$
\lim\limits_{p\to\infty}
\sum_{j_1=0}^{p}\ldots \sum_{j_k=0}^{p}
C_{j_k\ldots j_1}^2\stackrel{\sf def}{=}
\sum_{j_1,\ldots,j_k=0}^{\infty}
C_{j_k\ldots j_1}^2.
$$

\vspace{4mm}

From the other hand, for iterated limits we write

\vspace{1mm}
$$
\lim\limits_{p_1\to\infty}\ldots \lim\limits_{p_k\to\infty}
\sum_{j_1=0}^{p_1}\ldots \sum_{j_k=0}^{p_k}
C_{j_k\ldots j_1}^2\stackrel{\sf def}{=}
\sum_{j_1=0}^{\infty}\ldots
\sum_{j_k=0}^{\infty}
C_{j_k\ldots j_1}^2,
$$

\vspace{3mm}
$$
\lim\limits_{p_1\to\infty}\lim\limits_{p_2,\ldots,p_k\to\infty}
\sum_{j_1=0}^{p_1}\ldots \sum_{j_k=0}^{p_k}
C_{j_k\ldots j_1}^2\stackrel{\sf def}{=}
\sum_{j_1=0}^{\infty}
\sum_{j_2,\ldots,j_k=0}^{\infty}
C_{j_k\ldots j_1}^2
$$

\vspace{2mm}
\noindent
and so on.

Let us consider the following lemma.

\vspace{2mm}

{\bf Lemma 2.}\ {\it The following equalities are fulfilled

\vspace{1mm}
$$
\sum_{j_1,\ldots,j_k=0}^{\infty}
C_{j_k\ldots j_1}^2=
\sum_{j_1=0}^{\infty}\ldots
\sum_{j_k=0}^{\infty}
C_{j_k\ldots j_1}^2=
$$

\vspace{2mm}
\begin{equation}
\label{lem1}
=\sum_{j_k=0}^{\infty}\ldots
\sum_{j_1=0}^{\infty}
C_{j_k\ldots j_1}^2=
\sum_{j_{q_1}=0}^{\infty}\ldots
\sum_{j_{q_k}=0}^{\infty}
C_{j_k\ldots j_1}^2
\end{equation}

\vspace{5mm}
\noindent
for any permutation $(q_1,\ldots,q_k)$ such that
$\{q_1,\ldots,q_k\}=\{1,\ldots,k\}.$}

\vspace{2mm}

{\bf Proof.} Let us consider the value

\vspace{-3mm}
\begin{equation}
\label{21}
\sum_{j_{q_l}=0}^{p}\ldots
\sum_{j_{q_k}=0}^{p}
C_{j_k\ldots j_1}^2
\end{equation}

\vspace{3mm}
\noindent
for any permutation $(q_l,\ldots,q_k)$, where $l=1,2,\ldots,k$,
$\{q_1,\ldots,q_k\}=\{1,\ldots,k\}.$

Obviously, (\ref{21}) 
is the non-decreasing sequence with respect to $p$.
Moreover,

$$
\sum_{j_{q_l}=0}^{p}\ldots
\sum_{j_{q_k}=0}^{p}
C_{j_k\ldots j_1}^2\le 
\sum_{j_{q_1}=0}^{p}\sum_{j_{q_2}=0}^{p}\ldots
\sum_{j_{q_k}=0}^{p}
C_{j_k\ldots j_1}^2\le 
$$

\vspace{4mm}
$$
\le
\sum_{j_1,\ldots,j_k=0}^{\infty}
C_{j_k\ldots j_1}^2<\infty.
$$

\vspace{4mm}

Then the following limit

\vspace{1mm}
$$
\lim\limits_{p\to\infty}\sum\limits_{j_{q_l}=0}^p \ldots 
\sum\limits_{j_{q_k}=0}^{p}
C_{j_k\ldots j_1}^2=
\sum_{j_{q_l},\ldots,j_{q_k}=0}^{\infty}
C_{j_k\ldots j_1}^2
$$

\vspace{3mm}
\noindent
exists.

Let $p_l,\ldots,p_k$ simultaneously tend to infinity.
Then $g, r\to \infty$, where $g=\min\{p_l,\ldots,p_k\}$ and
$r=\max\{p_l,\ldots,p_k\}$. Moreover,

\vspace{2mm}
$$
\sum_{j_{q_l}=0}^{g}\ldots
\sum_{j_{q_k}=0}^{g}
C_{j_k\ldots j_1}^2\le 
\sum_{j_{q_l}=0}^{p_l}\ldots
\sum_{j_{q_k}=0}^{p_k}
C_{j_k\ldots j_1}^2\le
\sum_{j_{q_l}=0}^{r}\ldots
\sum_{j_{q_k}=0}^{r}
C_{j_k\ldots j_1}^2.
$$

\vspace{5mm}

This means that the existence of the limit 

\begin{equation}
\label{1c1c}
\lim\limits_{p\to\infty}\sum_{j_{q_l}=0}^{p}\ldots
\sum_{j_{q_k}=0}^{p}
C_{j_k\ldots j_1}^2
\end{equation}

\vspace{3mm}
\noindent
implies the existence of the limit 
\begin{equation}
\label{1d1d}
\lim\limits_{p_l,\ldots,p_k\to\infty}\sum_{j_{q_l}=0}^{p_l}\ldots
\sum_{j_{q_k}=0}^{p_k}
C_{j_k\ldots j_1}^2
\end{equation}

\vspace{3mm}
\noindent
and equality of the limits (\ref{1c1c}) and (\ref{1d1d}).
 
Taking into account the above reasoning, we have 

$$
\lim\limits_{p,q\to\infty}\sum_{j_{q_l}=0}^{q}\sum_{j_{q_{l+1}}=0}^{p}\ldots
\sum_{j_{q_k}=0}^{p}
C_{j_k\ldots j_1}^2=
\lim\limits_{p\to\infty}\sum_{j_{q_l}=0}^{p}\ldots
\sum_{j_{q_k}=0}^{p}
C_{j_k\ldots j_1}^2=
$$

\vspace{1mm}
\begin{equation}
\label{1h1h}
=\lim\limits_{p_l,\ldots,p_k\to\infty}\sum_{j_{q_l}=0}^{p_l}\ldots
\sum_{j_{q_k}=0}^{p_k}
C_{j_k\ldots j_1}^2.
\end{equation}

\vspace{4mm}

Since the limit
$$
\sum_{j_1,\ldots,j_k=0}^{\infty}
C_{j_k\ldots j_1}^2
$$

\vspace{3mm}
\noindent
exists (see the Parseval equality (\ref{par1})), then from Proposition 1
we have

\vspace{1mm}
$$
\sum_{j_{q_1}=0}^{\infty}\sum_{j_{q_2},\ldots,j_{q_k}=0}^{\infty}
C_{j_k\ldots j_1}^2=
\lim\limits_{q\to\infty}
\lim\limits_{p\to\infty}
\sum_{j_{q_1}=0}^{q}\sum_{j_{q_2}=0}^p \ldots \sum_{j_{q_k}=0}^{p}
C_{j_k\ldots j_1}^2=
$$

\vspace{3mm}
\begin{equation}
\label{1b1b}
=\lim\limits_{q,p\to\infty}
\sum_{j_{q_1}=0}^{q}\sum_{j_{q_2}=0}^p \ldots \sum_{j_{q_k}=0}^{p}
C_{j_k\ldots j_1}^2=
\sum_{j_1,\ldots,j_k=0}^{\infty}
C_{j_k\ldots j_1}^2.
\end{equation}

\vspace{5mm}

Using (\ref{1h1h}) and Proposition 1, we get

$$
\sum_{j_{q_2}=0}^{\infty}\sum_{j_{q_3},\ldots,j_{q_k}=0}^{\infty}
C_{j_k\ldots j_1}^2=
\lim\limits_{q\to\infty}
\lim\limits_{p\to\infty}
\sum_{j_{q_2}=0}^{q}\sum_{j_{q_3}=0}^p \ldots \sum_{j_{q_k}=0}^{p}
C_{j_k\ldots j_1}^2=
$$

\vspace{3mm}
\begin{equation}
\label{1a1a}
=\lim\limits_{q,p\to\infty}
\sum_{j_{q_2}=0}^{q}\sum_{j_{q_3}=0}^p \ldots \sum_{j_{q_k}=0}^{p}
C_{j_k\ldots j_1}^2=
\sum_{j_{q_2},\ldots,j_{q_k}=0}^{\infty}
C_{j_k\ldots j_1}^2.
\end{equation}

\vspace{5mm}

Combining (\ref{1a1a}) and (\ref{1b1b}), we obtain

$$
\sum_{j_{q_1}=0}^{\infty}\sum_{j_{q_2}=0}^{\infty}
\sum_{j_{q_3},\ldots,j_{q_k}=0}^{\infty}
C_{j_k\ldots j_1}^2=
\sum_{j_{1},\ldots,j_{k}=0}^{\infty}
C_{j_k\ldots j_1}^2.
$$

\vspace{3mm}

Repeating the previous steps, we complete the proof of Lemma
2.

Further, let us show that for $s=1,\ldots,k$

\vspace{2mm}
$$
\sum_{j_1=0}^{\infty}\ldots
\sum_{j_{s-1}=0}^{\infty}
\sum_{j_s=p+1}^{\infty}\sum_{j_{s+1}=0}^{\infty}\ldots \sum_{j_k=0}^{\infty}
C_{j_k\ldots j_1}^2=
$$

\vspace{3mm}
\begin{equation}
\label{d11}
=
\sum_{j_s=p+1}^{\infty}\sum_{j_{s-1}=0}^{\infty}\ldots
\sum_{j_{1}=0}^{\infty}
\sum_{j_{s+1}=0}^{\infty}\ldots \sum_{j_k=0}^{\infty}
C_{j_k\ldots j_1}^2.
\end{equation}

\vspace{5mm}

Using the arguments which we used when proving Lemma 2, we obtain

\vspace{2mm}
$$
\lim\limits_{n\to\infty}
\sum_{j_1=0}^{n}\ldots
\sum_{j_{s-1}=0}^{n}
\sum_{j_s=0}^{p}\sum_{j_{s+1}=0}^{n}\ldots \sum_{j_k=0}^{n}
C_{j_k\ldots j_1}^2=
$$

\vspace{3mm}
\begin{equation}
\label{ura0}
=\sum_{j_s=0}^{p}\ \sum_{j_{1},\ldots, j_{s-1}, j_{s+1},\ldots,j_k=0}^{\infty}
C_{j_k\ldots j_1}^2
=\sum_{j_s=0}^{p}\sum_{j_{q_1}=0}^{\infty}\ldots
\sum_{j_{q_{k-1}}=0}^{\infty}
C_{j_k\ldots j_1}^2
\end{equation}

\vspace{5mm}
\noindent
for any permutation $(q_1,\ldots,q_{k-1})$ such that
$\{q_1,\ldots,q_{k-1}\}=\{1,\ldots,s-1,s+1,\ldots,k\}$,
where $p$ is a fixed natural number.

Obviously, we have

\vspace{2mm}
$$
\sum_{j_s=0}^{p}\sum_{j_{q_1}=0}^{\infty}\ldots
\sum_{j_{q_{k-1}}=0}^{\infty}
C_{j_k\ldots j_1}^2=
\sum_{j_{q_1}=0}^{\infty}\ldots \sum_{j_s=0}^{p} \ldots
\sum_{j_{q_{k-1}}=0}^{\infty}C_{j_k\ldots j_1}^2= \ldots =
$$

\vspace{3mm}
\begin{equation}
\label{ura1}
=
\sum_{j_{q_1}=0}^{\infty}\ldots 
\sum_{j_{q_{k-1}}=0}^{\infty}
\sum_{j_s=0}^{p}
C_{j_k\ldots j_1}^2.
\end{equation}

\vspace{5mm}

Using (\ref{ura0}), (\ref{ura1}), and Lemma 2, we get

\vspace{2mm}
$$
\sum_{j_1=0}^{\infty}\ldots
\sum_{j_{s-1}=0}^{\infty}
\sum_{j_s=p+1}^{\infty}\sum_{j_{s+1}=0}^{\infty}\ldots \sum_{j_k=0}^{\infty}
C_{j_k\ldots j_1}^2=
\sum_{j_1=0}^{\infty}\ldots
\sum_{j_{s-1}=0}^{\infty}
\sum_{j_s=0}^{\infty}\sum_{j_{s+1}=0}^{\infty}\ldots \sum_{j_k=0}^{\infty}
C_{j_k\ldots j_1}^2-
$$

\vspace{3mm}
$$
-\sum_{j_1=0}^{\infty}\ldots
\sum_{j_{s-1}=0}^{\infty}
\sum_{j_s=0}^{p}\sum_{j_{s+1}=0}^{\infty}\ldots \sum_{j_k=0}^{\infty}
C_{j_k\ldots j_1}^2=
$$

\vspace{3mm}
$$
=
\sum_{j_s=0}^{\infty}
\sum_{j_{s-1}=0}^{\infty}\ldots
\sum_{j_1=0}^{\infty}\sum_{j_{s+1}=0}^{\infty}\ldots \sum_{j_k=0}^{\infty}
C_{j_k\ldots j_1}^2-
\sum_{j_s=0}^{p}
\sum_{j_{s-1}=0}^{\infty}\ldots
\sum_{j_1=0}^{\infty}\sum_{j_{s+1}=0}^{\infty}\ldots \sum_{j_k=0}^{\infty}
C_{j_k\ldots j_1}^2=
$$

\vspace{3mm}
$$
=
\sum_{j_s=p+1}^{\infty}
\sum_{j_{s-1}=0}^{\infty}\ldots
\sum_{j_1=0}^{\infty}\sum_{j_{s+1}=0}^{\infty}\ldots \sum_{j_k=0}^{\infty}
C_{j_k\ldots j_1}^2.
$$

\vspace{5mm}

The equality (\ref{d11}) is proved.

Using the Parseval equality and Lemma 2, we obtain

\vspace{2mm}
$$
\int\limits_{[t,T]^k}K^2(t_1,\ldots,t_k)dt_1\ldots dt_k-
\sum_{j_1=0}^{p}\ldots \sum_{j_k=0}^{p}
C_{j_k\ldots j_1}^2=
$$

\vspace{3mm}

$$
=\sum_{j_1,\ldots,j_k=0}^{\infty}
C_{j_k\ldots j_1}^2-
\sum_{j_1=0}^{p}\ldots \sum_{j_k=0}^{p}
C_{j_k\ldots j_1}^2=
$$

\vspace{5mm}

$$
=
\sum_{j_1=0}^{\infty}\ldots \sum_{j_k=0}^{\infty}
C_{j_k\ldots j_1}^2-
\sum_{j_1=0}^{p}\ldots \sum_{j_k=0}^{p}
C_{j_k\ldots j_1}^2=
$$

\vspace{5mm}

$$
=\sum_{j_1=0}^{p}\sum_{j_2=0}^{\infty}\ldots \sum_{j_k=0}^{\infty}
C_{j_k\ldots j_1}^2+
\sum_{j_1=p+1}^{\infty}\sum_{j_2=0}^{\infty}\ldots \sum_{j_k=0}^{\infty}
C_{j_k\ldots j_1}^2-
\sum_{j_1=0}^{p}\ldots \sum_{j_k=0}^{p}
C_{j_k\ldots j_1}^2=
$$

\vspace{6mm}

$$
=\sum_{j_1=0}^{p}\sum_{j_2=0}^{p}\sum_{j_3=0}^{\infty}
\ldots \sum_{j_k=0}^{\infty}
C_{j_k\ldots j_1}^2+
\sum_{j_1=0}^{p}\sum_{j_2=p+1}^{\infty}
\sum_{j_3=0}^{\infty}
\ldots \sum_{j_k=0}^{\infty}+
$$

\vspace{4mm}

$$
+\sum_{j_1=p+1}^{\infty}\sum_{j_2=0}^{\infty}\ldots \sum_{j_k=0}^{\infty}
C_{j_k\ldots j_1}^2-
\sum_{j_1=0}^{p}\ldots \sum_{j_k=0}^{p}
C_{j_k\ldots j_1}^2=\ldots =
$$

\vspace{6mm}

$$
=\sum_{j_1=p+1}^{\infty}\sum_{j_2=0}^{\infty}\ldots \sum_{j_k=0}^{\infty}
C_{j_k\ldots j_1}^2+
\sum_{j_1=0}^p
\sum_{j_2=p+1}^{\infty}\sum_{j_2=0}^{\infty}\ldots \sum_{j_k=0}^{\infty}
C_{j_k\ldots j_1}^2+
$$

\vspace{6mm}

$$
+\sum_{j_1=0}^p\sum_{j_2=0}^p
\sum_{j_3=p+1}^{\infty}\sum_{j_4=0}^{\infty}\ldots \sum_{j_k=0}^{\infty}
C_{j_k\ldots j_1}^2+ \ldots +
\sum_{j_1=0}^p\ldots \sum_{j_{k-1}=0}^p
\sum_{j_k=p+1}^{\infty}C_{j_k\ldots j_1}^2\le
$$

\vspace{6mm}

$$
\le\sum_{j_1=p+1}^{\infty}\sum_{j_2=0}^{\infty}\ldots \sum_{j_k=0}^{\infty}
C_{j_k\ldots j_1}^2+
\sum_{j_1=0}^{\infty}
\sum_{j_2=p+1}^{\infty}\sum_{j_2=0}^{\infty}\ldots \sum_{j_k=0}^{\infty}
C_{j_k\ldots j_1}^2+
$$

\vspace{4mm}

$$
+\sum_{j_1=0}^{\infty}\sum_{j_2=0}^{\infty}
\sum_{j_3=p+1}^{\infty}\sum_{j_4=0}^{\infty}\ldots \sum_{j_k=0}^{\infty}
C_{j_k\ldots j_1}^2+ \ldots +
\sum_{j_1=0}^{\infty}\ldots \sum_{j_{k-1}=0}^{\infty}
\sum_{j_k=p+1}^{\infty}C_{j_k\ldots j_1}^2=
$$

\vspace{4mm}

\begin{equation}
\label{aaap}
=\sum\limits_{s=1}^k \left(\sum_{j_1=0}^{\infty}\ldots
\sum_{j_{s-1}=0}^{\infty}
\sum_{j_s=p+1}^{\infty}\sum_{j_{s+1}=0}^{\infty}\ldots \sum_{j_k=0}^{\infty}
C_{j_k\ldots j_1}^2\right).
\end{equation}

\vspace{6mm}

Note that deriving (\ref{aaap}), we used the following

\vspace{1mm}
$$
\sum_{j_1=0}^{p}\ldots
\sum_{j_{s-1}=0}^{p}
\sum_{j_s=p+1}^{\infty}\sum_{j_{s+1}=0}^{\infty}\ldots \sum_{j_k=0}^{\infty}
C_{j_k\ldots j_1}^2\le
$$

\vspace{3mm}
$$
\le
\sum_{j_1=0}^{m_1}\ldots
\sum_{j_{s-1}=0}^{m_{s-1}}
\sum_{j_s=p+1}^{\infty}\sum_{j_{s+1}=0}^{\infty}\ldots \sum_{j_k=0}^{\infty}
C_{j_k\ldots j_1}^2\le
$$

\vspace{3mm}
$$
\le
\lim\limits_{m_{s-1}\to\infty}
\sum_{j_1=0}^{m_1}\ldots
\sum_{j_{s-1}=0}^{m_{s-1}}
\sum_{j_s=p+1}^{\infty}\sum_{j_{s+1}=0}^{\infty}\ldots \sum_{j_k=0}^{\infty}
C_{j_k\ldots j_1}^2=
$$

\vspace{3mm}
$$
=
\sum_{j_1=0}^{m_1}\ldots
\sum_{j_{s-2}=0}^{m_{s-2}}\sum_{j_{s-1}=0}^{\infty}
\sum_{j_s=p+1}^{\infty}\sum_{j_{s+1}=0}^{\infty}\ldots \sum_{j_k=0}^{\infty}
C_{j_k\ldots j_1}^2\le
$$

\vspace{1mm}
$$
\le\ldots\le
$$

\vspace{-1mm}
$$
\le\sum_{j_1=0}^{\infty}\ldots
\sum_{j_{s-1}=0}^{\infty}
\sum_{j_s=p+1}^{\infty}\sum_{j_{s+1}=0}^{\infty}\ldots \sum_{j_k=0}^{\infty}
C_{j_k\ldots j_1}^2,
$$

\vspace{5mm}
\noindent
where $m_1,\ldots,m_{s-1}>p.$

Denote
$$
C_{j_s\ldots j_1}(\tau)=
\int\limits_t^{\tau}
\phi_{j_s}(t_s)\psi_s(t_s)\ldots \int\limits_t^{t_2}
\phi_{j_1}(t_1)\psi_1(t_1)dt_1\ldots dt_s,
$$

\vspace{2mm}
\noindent
where
$s=1,\ldots,k-1.$

Let us remind the Dini Theorem, which we will use further.

\vspace{2mm}

{\bf Theorem (Dini).} {\it 
Let the functional sequence $u_n(x)$ 
be non-decreasing at each point of the interval $[a, b]$.
In addition, all the functions $u_n(x)$
of this sequence and the limit function $u(x)$ are continuous on the interval
$[a, b].$ Then the convergence $u_n(x)$ to 
$u(x)$ is uniform on the interval $[a,b].$}

\vspace{2mm}

For $s<k$ due to the Parseval equality, Dini Theorem
and (\ref{d11}) we obtain

\vspace{2mm}
$$
\sum_{j_1=0}^{\infty}\ldots
\sum_{j_{s-1}=0}^{\infty}
\sum_{j_s=p+1}^{\infty}\sum_{j_{s+1}=0}^{\infty}\ldots \sum_{j_k=0}^{\infty}
C_{j_k\ldots j_1}^2=
$$

\vspace{5mm}

$$
\stackrel{\hbox{(\ref{d11})}}{=}\ \ \ 
\sum_{j_s=p+1}^{\infty}
\sum_{j_{s-1}=0}^{\infty}\ldots
\sum_{j_{1}=0}^{\infty}
\sum_{j_{s+1}=0}^{\infty}\ldots \sum_{j_k=0}^{\infty}
C_{j_k\ldots j_1}^2=
$$

\vspace{5mm}

$$
\stackrel{\hbox{(Parseval Eq.)}}{=}\ \ \
\sum_{j_s=p+1}^{\infty}
\sum_{j_{s-1}=0}^{\infty}\ldots
\sum_{j_{1}=0}^{\infty}
\sum_{j_{s+1}=0}^{\infty}\ldots 
\sum_{j_{k-1}=0}^{\infty}
\int\limits_t^T \psi_k^2(t_k) \left(C_{j_{k-1}\ldots j_1}(t_k)\right)^2 dt_k=
$$

\vspace{5mm}

$$
\stackrel{\hbox{(Dini Th.)}}{=}\ \ \
\sum_{j_s=p+1}^{\infty}
\sum_{j_{s-1}=0}^{\infty}\ldots
\sum_{j_{1}=0}^{\infty}
\sum_{j_{s+1}=0}^{\infty}\ldots 
\sum_{j_{k-2}=0}^{\infty}
\int\limits_t^T \psi_k^2(t_k) 
\sum_{j_{k-1}=0}^{\infty}\left(C_{j_{k-1}\ldots j_1}(t_k)\right)^2 dt_k=
$$

\vspace{5mm}

$$
\stackrel{\hbox{(Parseval Eq.)}}{=}\ \ \
\sum_{j_s=p+1}^{\infty}
\sum_{j_{s-1}=0}^{\infty}\ldots
\sum_{j_{1}=0}^{\infty}
\sum_{j_{s+1}=0}^{\infty}\ldots 
\sum_{j_{k-2}=0}^{\infty}
\int\limits_t^T \psi_k^2(t_k) \int\limits_t^{t_k} \psi_{k-1}^2(t_{k-1}) 
\left(C_{j_{k-2}\ldots j_1}(t_{k-1})\right)^2\times
$$

\vspace{2mm}
$$
\times dt_{k-1}dt_k\le
$$

\vspace{5mm}
$$
\le C\sum_{j_s=p+1}^{\infty}
\sum_{j_{s-1}=0}^{\infty}\ldots
\sum_{j_{1}=0}^{\infty}
\sum_{j_{s+1}=0}^{\infty}\ldots 
\sum_{j_{k-2}=0}^{\infty}
\int\limits_t^T 
\left(C_{j_{k-2}\ldots j_1}(\tau)\right)^2 d\tau =
$$

\vspace{5mm}

$$
\stackrel{\hbox{(Dini Th.)}}{=}\ \ \
C\sum_{j_s=p+1}^{\infty}
\sum_{j_{s-1}=0}^{\infty}\ldots
\sum_{j_{1}=0}^{\infty}
\sum_{j_{s+1}=0}^{\infty}\ldots 
\sum_{j_{k-3}=0}^{\infty}
\int\limits_t^T 
\sum_{j_{k-2}=0}^{\infty}
\left(C_{j_{k-2}\ldots j_1}(\tau)\right)^2 d\tau =
$$

\vspace{5mm}

$$
\stackrel{\hbox{(Parseval Eq.)}}{=}\ \ \ C\sum_{j_s=p+1}^{\infty}
\sum_{j_{s-1}=0}^{\infty}\ldots
\sum_{j_{1}=0}^{\infty}
\sum_{j_{s+1}=0}^{\infty}\ldots 
\sum_{j_{k-3}=0}^{\infty}
\int\limits_t^T \int\limits_t^{\tau}
\psi_{k-2}^2(\theta)
\left(C_{j_{k-3}\ldots j_1}(\theta)\right)^2
d\theta d\tau\le 
$$

\vspace{5mm}

$$
\le K
\sum_{j_s=p+1}^{\infty}
\sum_{j_{s-1}=0}^{\infty}\ldots
\sum_{j_{1}=0}^{\infty}
\sum_{j_{s+1}=0}^{\infty}\ldots 
\sum_{j_{k-3}=0}^{\infty}
\int\limits_t^T
\left(C_{j_{k-3}\ldots j_1}(\tau)\right)^2
d\tau\le 
$$

\vspace{4mm}
$$
\le \ldots \le
$$

$$
\le C_k
\sum_{j_s=p+1}^{\infty}
\sum_{j_{s-1}=0}^{\infty}\ldots
\sum_{j_{1}=0}^{\infty}
\int\limits_t^T 
\left(C_{j_{s}\ldots j_1}(\tau)\right)^2 d\tau=
$$

\vspace{5mm}

\begin{equation}
\label{d14}
\stackrel{\hbox{(Dini Th.)}}{=}\ \ \ C_k
\sum_{j_s=p+1}^{\infty}
\sum_{j_{s-1}=0}^{\infty}\ldots
\sum_{j_{2}=0}^{\infty}
\int\limits_t^T  \sum_{j_{1}=0}^{\infty}
\left(C_{j_{s}\ldots j_1}(\tau)\right)^2 d\tau,
\end{equation}

\vspace{6mm}
\noindent
where constants $C,$ $K$ depend on $T-t$ and
constant $C_k$ depends on $k$ and $T-t.$

Let us explane more precisely how we obtain (\ref{d14}).
For any function $g(s)\in L_2([t,T])$ we have the following
Parseval equality

$$
\sum\limits_{j=0}^{\infty}\left(\int\limits_t^{\tau}
\phi_j(s)g(s)ds\right)^2=
\sum\limits_{j=0}^{\infty}\left(\int\limits_t^T
{\bf 1}_{\{s<\tau\}}\phi_j(s)g(s)ds\right)^2=
$$

\begin{equation}
\label{d15}
=\int\limits_t^T
\left({\bf 1}_{\{s<\tau\}}\right)^2 g^2(s)ds=
\int\limits_t^{\tau}
g^2(s)ds.
\end{equation}

\vspace{3mm}

The equality (\ref{d15}) has been applied repeatedly when we obtaining
(\ref{d14}).

Using the replacement of the integrating order in Riemann integrals, we have

\vspace{1mm}
$$
C_{j_s\ldots j_1}(\tau)=
\int\limits_t^{\tau}
\phi_{j_s}(t_s)\psi_s(t_s)\ldots \int\limits_t^{t_2}
\phi_{j_1}(t_1)\psi_1(t_1)dt_1\ldots dt_s=
$$

\vspace{2mm}
$$
=\int\limits_t^{\tau}
\phi_{j_1}(t_1)\psi_1(t_1)\int\limits_{t_1}^{\tau}
\phi_{j_2}(t_2)\psi_2(t_2)
\ldots
\int\limits_{t_{s-1}}^{\tau}
\phi_{j_s}(t_s)\psi_s(t_s)dt_s\ldots dt_2dt_1
\stackrel{\sf def}{=}
$$

\vspace{4mm}
$$
\stackrel{\sf def}{=}
{\tilde C}_{j_s\ldots j_1}(\tau).
$$

\vspace{7mm}

For $l=1,\ldots,s$ we will use the following notation

\vspace{2mm}
$$
{\tilde C}_{j_s\ldots j_l}(\tau,\theta)=
\int\limits_{\theta}^{\tau}
\phi_{j_l}(t_l)\psi_l(t_l)\int\limits_{t_l}^{\tau}
\phi_{j_{l+1}}(t_{l+1})\psi_{l+1}(t_{l+1})
\ldots
\int\limits_{t_{s-1}}^{\tau}
\phi_{j_s}(t_s)\psi_s(t_s)dt_s\ldots dt_{l+1}dt_l.
$$

\vspace{5mm}

Using the Parseval equality and Dini Theorem, from (\ref{d14}) we obtain

\vspace{2mm}
$$
\sum_{j_1=0}^{\infty}\ldots
\sum_{j_{s-1}=0}^{\infty}
\sum_{j_s=p+1}^{\infty}\sum_{j_{s+1}=0}^{\infty}\ldots \sum_{j_k=0}^{\infty}
C_{j_k\ldots j_1}^2\le
$$

\vspace{5mm}

$$
\le
C_k
\sum_{j_s=p+1}^{\infty}
\sum_{j_{s-1}=0}^{\infty}\ldots
\sum_{j_{2}=0}^{\infty}
\int\limits_t^T  \sum_{j_{1}=0}^{\infty}
\left(C_{j_{s}\ldots j_1}(\tau)\right)^2 d\tau=
$$

\vspace{5mm}

$$
=C_k
\sum_{j_s=p+1}^{\infty}
\sum_{j_{s-1}=0}^{\infty}\ldots
\sum_{j_{2}=0}^{\infty}
\int\limits_t^T  \sum_{j_{1}=0}^{\infty}
\left({\tilde C}_{j_{s}\ldots j_1}(\tau)\right)^2 d\tau=
$$

\vspace{5mm}

\begin{equation}
\label{molod1}
\stackrel{\hbox{(Parseval Eq.)}}{=}\ \ \ C_k
\sum_{j_s=p+1}^{\infty}
\sum_{j_{s-1}=0}^{\infty}\ldots
\sum_{j_{2}=0}^{\infty}
\int\limits_t^T\int\limits_t^{\tau}\psi_1^2(t_1)  
\left({\tilde C}_{j_{s}\ldots j_2}(\tau,t_1)\right)^2 dt_1d\tau=
\end{equation}

\vspace{5mm}

\begin{equation}
\label{molod2}
\stackrel{\hbox{(Dini Th.)}}{=}\ \ \ C_k
\sum_{j_s=p+1}^{\infty}
\sum_{j_{s-1}=0}^{\infty}\ldots
\sum_{j_{3}=0}^{\infty}
\int\limits_t^T\int\limits_t^{\tau}\psi_1^2(t_1)  
\sum_{j_{2}=0}^{\infty}
\left({\tilde C}_{j_{s}\ldots j_2}(\tau,t_1)\right)^2 dt_1d\tau=
\end{equation}

\vspace{5mm}

$$
\stackrel{\hbox{(Parseval Eq.)}}{=}\ \ \ C_k
\sum_{j_s=p+1}^{\infty}
\sum_{j_{s-1}=0}^{\infty}\ldots
\sum_{j_{3}=0}^{\infty}
\int\limits_t^T\int\limits_t^{\tau}\psi_1^2(t_1)  
\int\limits_{t_1}^{\tau}\psi_2^2(t_2)  
\left({\tilde C}_{j_{s}\ldots j_3}(\tau,t_2)\right)^2 dt_2dt_1d\tau \le
$$

\vspace{5mm}

$$
\le C_k
\sum_{j_s=p+1}^{\infty}
\sum_{j_{s-1}=0}^{\infty}\ldots
\sum_{j_{3}=0}^{\infty}
\int\limits_t^T\int\limits_t^{\tau}\psi_1^2(t_1)  
\int\limits_{t}^{\tau}\psi_2^2(t_2)  
\left({\tilde C}_{j_{s}\ldots j_3}(\tau,t_2)\right)^2 dt_2dt_1d\tau\le
$$

\vspace{5mm}

$$
\le C^{'}_k
\sum_{j_s=p+1}^{\infty}
\sum_{j_{s-1}=0}^{\infty}\ldots
\sum_{j_{3}=0}^{\infty}
\int\limits_t^T
\int\limits_{t}^{\tau}\psi_2^2(t_2)  
\left({\tilde C}_{j_{s}\ldots j_3}(\tau,t_2)\right)^2 dt_2d\tau
\le 
$$

\vspace{4mm}

$$
\le  \ldots \le
$$

\vspace{1mm}

$$
\le C^{''}_k
\sum_{j_s=p+1}^{\infty}
\int\limits_t^T\int\limits_t^{\tau}
\psi_{s-1}^2(t_{s-1})
\left({\tilde C}_{j_{s}}(\tau,t_{s-1})\right)^2 dt_{s-1} d\tau\le
$$

\vspace{4mm}

\begin{equation}
\label{la}
\le {\tilde C}_k
\sum_{j_s=p+1}^{\infty}
\int\limits_t^T\int\limits_t^{\tau}
\left(~\int\limits_{u}^{\tau}\phi_{j_s}(\theta)
\psi_s(\theta)d\theta\right)^2 du d\tau,
\end{equation}

\vspace{6mm}
\noindent
where constants $C^{'}_k,$ $C^{''}_k,$ $\tilde C_k$
depend on $k$ and $T-t.$

Let us explane more precisely how we obtain (\ref{la}).
For any function $g(s)\in L_2([t,T])$ we have the following
Parseval equality

$$
\sum\limits_{j=0}^{\infty}\left(\int\limits_{\theta}^{\tau}
\phi_j(s)g(s)ds\right)^2=
\sum\limits_{j=0}^{\infty}\left(\int\limits_t^T
{\bf 1}_{\{\theta<s<\tau\}}\phi_j(s)g(s)ds\right)^2=
$$

\begin{equation}
\label{d22}
=\int\limits_t^T
\left({\bf 1}_{\{\theta<s<\tau\}}\right)^2 g^2(s)ds=
\int\limits_{\theta}^{\tau}
g^2(s)ds.
\end{equation}

\vspace{3mm}

The equality (\ref{d22}) has been applied repeatedly when we obtaining
(\ref{la}).

Let us explane more precisely the passing from (\ref{molod1})
to (\ref{molod2}) (the same steps have been used when we 
deriving (\ref{la})).

We have

$$
\int\limits_t^T\int\limits_t^{\tau}\psi_1^2(t_1)  
\sum_{j_{2}=0}^{\infty}
\left({\tilde C}_{j_{s}\ldots j_2}(\tau,t_1)\right)^2 dt_1d\tau -
\sum_{j_{2}=0}^{n}\int\limits_t^T\int\limits_t^{\tau}\psi_1^2(t_1)  
\left({\tilde C}_{j_{s}\ldots j_2}(\tau,t_1)\right)^2 dt_1d\tau =
$$

\vspace{4mm}

$$
=\int\limits_t^T\int\limits_t^{\tau}\psi_1^2(t_1)  
\sum_{j_{2}=n+1}^{\infty}
\left({\tilde C}_{j_{s}\ldots j_2}(\tau,t_1)\right)^2 dt_1d\tau =
$$

\vspace{4mm}

\begin{equation}
\label{molod3}
=\lim\limits_{N\to\infty}
\sum\limits_{j=0}^{N-1}\int\limits_t^{\tau_j}\psi_1^2(t_1)  
\sum_{j_{2}=n+1}^{\infty}
\left({\tilde C}_{j_{s}\ldots j_2}(\tau_j,t_1)\right)^2 dt_1 \Delta\tau_j,
\end{equation}

\vspace{4mm}
\noindent
where $\{\tau_j\}_{j=0}^{N}$ is the partition of the 
interval $[t, T],$ which satisfies the condition (\ref{1111}).

Since the non-decreasing functional sequence $u_n(\tau_j,t_1)$ and its
limit function $u(\tau_j,t_1)$ are continuous on the
interval $[t,\tau_j]\subseteq [t, T]$ with respect to $t_1$,
where

$$
u_n(\tau_j,t_1)=
\sum_{j_{2}=0}^{n}
\left({\tilde C}_{j_{s}\ldots j_2}(\tau_j,t_1)\right)^2,
$$

$$
u(\tau_j,t_1)=
\sum_{j_{2}=0}^{\infty}
\left({\tilde C}_{j_{s}\ldots j_2}(\tau_j,t_1)\right)^2=
\int\limits_{t_1}^{\tau_j}
\psi_2^2(t_2)
\left({\tilde C}_{j_{s}\ldots j_3}(\tau_j,t_2)\right)^2 dt_2,
$$

\vspace{4mm}

\noindent 
then by Dini Theorem we have the uniform convergence
of $u_n(\tau_j,t_1)$ to $u(\tau_j,t_1)$ at the interval $[t,\tau_j]\subseteq
[t, T]$
with respect to $t_1.$ As a result, we obtain

\begin{equation}
\label{molod4}
\sum_{j_{2}=n+1}^{\infty}
\left({\tilde C}_{j_{s}\ldots j_2}(\tau_j,t_1)\right)^2<\varepsilon,\ \ \ 
t_1\in [t,\tau_j]
\end{equation}

\vspace{4mm}
\noindent
for $n>N(\varepsilon)$ ($N(\varepsilon)$ exists
for any $\varepsilon>0$ and it does not depend on $t_1$).

From (\ref{molod3}) and (\ref{molod4}) we obtain

$$
\lim\limits_{N\to\infty}
\sum\limits_{j=0}^{N-1}\int\limits_t^{\tau_j}\psi_1^2(t_1)  
\sum_{j_{2}=n+1}^{\infty}
\left({\tilde C}_{j_{s}\ldots j_2}(\tau_j,t_1)\right)^2 dt_1 \Delta\tau_j
\le
\varepsilon 
\lim\limits_{N\to\infty}
\sum\limits_{j=0}^{N-1}\int\limits_t^{\tau_j}\psi_1^2(t_1)  
dt_1 \Delta\tau_j= 
$$

\vspace{3mm}
\begin{equation}
\label{molod6}
=\varepsilon \int\limits_t^T
\int\limits_t^{\tau}\psi_1^2(t_1)  
dt_1 d\tau.
\end{equation}

\vspace{4mm}

Using (\ref{molod6}), we get

$$
\lim\limits_{n\to\infty}\int\limits_t^T\int\limits_t^{\tau}\psi_1^2(t_1)  
\sum_{j_{2}=n+1}^{\infty}
\left({\tilde C}_{j_{s}\ldots j_2}(\tau,t_1)\right)^2 dt_1d\tau = 0.
$$

\vspace{4mm}

This fact completes the proof of passing 
from (\ref{molod1})
to (\ref{molod2}).

Let us estimate the integral 
\begin{equation}
\label{st1}
\int\limits_{u}^{\tau}\phi_{j_s}(\theta)
\psi_s(\theta)d\theta
\end{equation}

\vspace{2mm}
\noindent
from (\ref{la}) for the case when $\{\phi_j(s)\}_{j=0}^{\infty}$
is a complete orthonormal system of Legendre polynomials or
trigonometric functions in the space $L_2([t,T])$.

Note that the estimates for the integral

\vspace{-2mm}
\begin{equation}
\label{st2}
\int\limits_{t}^{\tau}\phi_{j}(\theta)\psi(\theta)d\theta,\ \ \ j\ge p+1,
\end{equation}

\vspace{2mm}
\noindent
where $\psi(\theta)$ is a continuously
differentiable function on the interval $[t, T]$,
have been obtained in \cite{12}, \cite{15c}. The same estimates 
can also be found in early publications \cite{5a}-\cite{7}, 
\cite{9a}, \cite{10} and in the monographs \cite{10a}-\cite{10axx1}.

Let us estimate the integral (\ref{st1}) using the approach from
\cite{12}, \cite{15c}.

First, consider the case of Legendre polynomials.
Then $\phi_j(s)$ looks as follows

\vspace{1mm}
\begin{equation}
\label{ogo7}
\phi_j(\theta)=\sqrt{\frac{2j+1}{T-t}}P_j\left(\left(
\theta-\frac{T+t}{2}\right)\frac{2}{T-t}\right),\ \ \ j\ge 0,
\end{equation}

\vspace{5mm}
\noindent
where $P_j(x)$ $(j=0, 1, 2\ldots)$ is the
Legendre polynomial.

Further, we have 

$$
\int\limits_v^x\phi_{j}(\theta)\psi(\theta)d\theta=
\frac{\sqrt{T-t}\sqrt{2j+1}}{2}
\int\limits_{z(v)}^{z(x)}P_{j}(y)
\psi(u(y))dy=
$$

$$
=\frac{\sqrt{T-t}}{2\sqrt{2j+1}}\Biggl((P_{j+1}(z(x))-
P_{j-1}(z(x)))\psi(x)-
(P_{j+1}(z(v))-
P_{j-1}(z(v)))\psi(v)-
\Biggr.
$$

\begin{equation}
\label{6000}
\Biggl.-
\frac{T-t}{2}
\int\limits_{z(v)}^{z(x)}((P_{j+1}(y)-P_{j-1}(y))
{\psi}'(u(y))dy\Biggr),
\end{equation}

\vspace{5mm}
\noindent
where $x, v\in (t, T),$ $j\ge p+1,$ 
$u(y)$ and $z(x)$ are defined by the following relations

\vspace{1mm}
$$
u(y)=\frac{T-t}{2}y+\frac{T+t}{2},\ \ \
z(x)=\left(x-\frac{T+t}{2}\right)\frac{2}{T-t},
$$

\vspace{4mm}
\noindent
and ${\psi}'$ is a derivative of the function $\psi(\theta)$
with respect to the variable $u(y).$

Note that in (\ref{6000}) we used the following well known property
of the Legendre polynomials

\vspace{2mm}
$$
\frac{dP_{j+1}}{dx}(x)-\frac{dP_{j-1}}{dx}(x)=(2j+1)P_j(x),\ \ \ 
j=1, 2,\ldots
$$

\vspace{5mm}

From (\ref{6000}) and the well known estimate for the Legendre
polynomials

\vspace{1mm}
\begin{equation}
\label{200}
|P_j(y)| <\frac{K}{\sqrt{j+1}(1-y^2)^{1/4}},\ \ \ 
y\in (-1, 1),\ \ \ j\in \mathbb{N},
\end{equation}

\vspace{5mm}
\noindent
where constant $K$ does not depend on $y$ and $j$, it follows that

\vspace{2mm}
\begin{equation}
\label{101}
\left|
\int\limits_v^x\phi_{j}(\theta)\psi(\theta)d\theta
\right| <
\frac{C}{j}\Biggl(\frac{1}{(1-(z(x))^2)^{1/4}}+
\frac{1}{(1-(z(v))^2)^{1/4}}+C_1\Biggr),
\end{equation}

\vspace{5mm}
\noindent
where $j\in\mathbb{N},$  
$z(x), z(v)\in (-1, 1),$ $x, v\in (t, T),$
constants $C, C_1$ do not depend on $j$.

From (\ref{101}) we obtain

\vspace{1mm}
\begin{equation}
\label{102}
\left(
\int\limits_v^x\phi_{j}(\theta)\psi(\theta)d\theta
\right)^2 <
\frac{C_2}{j^2}\Biggl(\frac{1}{(1-(z(x))^2)^{1/2}}+
\frac{1}{(1-(z(v))^2)^{1/2}}+C_3\Biggr),
\end{equation}

\vspace{4mm}
\noindent
where $j\in\mathbb{N},$ constants $C_2, C_3$ do not depend on $j$.

Let us apply (\ref{102}) for estimating of the right-hand side
of (\ref{la}). We have

\vspace{1mm}
$$
\int\limits_t^T\int\limits_t^{\tau}
\left(~\int\limits_{u}^{\tau}\phi_{j_s}(\theta)
\psi_s(\theta)d\theta\right)^2 du d\tau\le
$$

$$
\le \frac{K_1}{j_s^2}
\left(
\int\limits_{-1}^1
\frac{dy}{\left(1-y^2\right)^{1/2}}+
\int\limits_{-1}^1\int\limits_{-1}^x
\frac{dy}{\left(1-y^2\right)^{1/2}}dx + K_2\right)\le
$$

\begin{equation}
\label{103}
\le \frac{K_3}{j_s^2},
\end{equation}

\vspace{4mm}
\noindent
where $j_s\in\mathbb{N},$ constants $K_1, K_2, K_3$ are independent of $j_s.$

Now, consider the trigonometric case.
The complete orthonormal system of trigonometric functions
in the space $L_2([t, T])$ has the following form

\begin{equation}
\label{trig11}
\phi_j(\theta)=\frac{1}{\sqrt{T-t}}
\left\{
\begin{matrix}
1,\ & j=0\cr \cr
\sqrt{2}{\rm sin} \left(2\pi r(\theta-t)/(T-t)\right),\ & j=2r-1\cr \cr
\sqrt{2}{\rm cos} \left(2\pi r(\theta-t)/(T-t)\right),\ & j=2r
\end{matrix}
,\right.
\end{equation}

\vspace{3mm}
\noindent
where $r=1, 2,\ldots $

Using the system of functions 
(\ref{trig11}), we have

\vspace{1mm}
$$
\int\limits_v^x\phi_{2r-1}(\theta)\psi(\theta)d\theta=
\sqrt{\frac{2}{T-t}}\int\limits_v^x
{\rm sin} \frac{2\pi r(\theta-t)}{T-t}\psi(\theta)d\theta=
$$

\vspace{1mm}
$$
=-\sqrt{\frac{T-t}{2}}\frac{1}{\pi r}\Biggl(
\psi(x){\rm cos}\frac{2\pi r(x-t)}{T-t}-
\psi(v){\rm cos}\frac{2\pi r(v-t)}{T-t}-\Biggr.
$$

\vspace{1mm}
\begin{equation}
\label{201}
\Biggl.-
\int\limits_v^x
{\rm cos} \frac{2\pi r(\theta-t)}{T-t}\psi'(\theta)d\theta\Biggr),
\end{equation}

\vspace{4mm}

$$
\int\limits_v^x\phi_{2r}(\theta)\psi(\theta)d\theta=
\sqrt{\frac{2}{T-t}}\int\limits_v^x
{\rm cos} \frac{2\pi r(\theta-t)}{T-t}\psi(\theta)d\theta=
$$

\vspace{1mm}
$$
=\sqrt{\frac{T-t}{2}}\frac{1}{\pi r}\Biggl(
\psi(x){\rm sin}\frac{2\pi r(x-t)}{T-t}-
\psi(v){\rm sin}\frac{2\pi r(v-t)}{T-t}-\Biggr.
$$

\vspace{1mm}
\begin{equation}
\label{202}
\Biggl.-
\int\limits_v^x
{\rm sin} \frac{2\pi r(\theta-t)}{T-t}\psi'(\theta)d\theta\Biggr),
\end{equation}

\vspace{5mm}
\noindent
where $\psi'(\theta)$ is a derivative of the function $\psi(\theta)$
with respect to the variable $\theta.$

Combining (\ref{201}) and (\ref{202}), we obtain for the
trigonometric case

\begin{equation}
\label{203}
\left(
\int\limits_v^x\phi_{j}(\theta)\psi(\theta)d\theta
\right)^2 \le 
\frac{C_4}{j^2},
\end{equation}

\vspace{4mm}
\noindent
where $j\in\mathbb{N},$ constant $C_4$ is independent of $j.$

From (\ref{203}) we finally have

\begin{equation}
\label{103x}
\int\limits_t^T\int\limits_t^{\tau}
\left(~\int\limits_{u}^{\tau}\phi_{j_s}(\theta)
\psi_s(\theta)d\theta\right)^2 du d\tau
\le \frac{K_4}{j_s^2},
\end{equation}

\vspace{4mm}
\noindent
where $j_s\in\mathbb{N},$ constant $K_4$ is independent of $j_s.$

Combibing (\ref{la}), (\ref{103}), and (\ref{103x}), we obtain

\vspace{2mm}
$$
\sum_{j_1=0}^{\infty}\ldots
\sum_{j_{s-1}=0}^{\infty}
\sum_{j_s=p+1}^{\infty}\sum_{j_{s+1}=0}^{\infty}\ldots \sum_{j_k=0}^{\infty}
C_{j_k\ldots j_1}^2\le
$$

\vspace{2mm}
\begin{equation}
\label{fff}
\le L_k
\sum_{j_s=p+1}^{\infty}\frac{1}{j_s^2} \le 
\frac{L_k}{p},
\end{equation}

\vspace{5mm}
\noindent
where constant $L_k$ depends on $k$ and $T-t.$

Obviously, the case $s=k$ can be considered absolutely analogously to the
case $s<k$. Then from (\ref{aaap}) and (\ref{fff})
we obtain

\begin{equation}
\label{ddd1}
\int\limits_{[t,T]^k}K^2(t_1,\ldots,t_k)dt_1\ldots dt_k-
\sum_{j_1=0}^{p}\ldots \sum_{j_k=0}^{p}
C_{j_k\ldots j_1}^2\le \frac{G_k}{p},
\end{equation}

\vspace{3mm}
\noindent
where constant $G_k$ depends on $k$ and $T-t.$

For the further consideration we will use the estimate (\ref{99999}).
Using (\ref{ddd1}) and the estimate (\ref{99999})
for the case $p_1=\ldots=p_k=p$ and $n=2$, 
we get

$$
{\sf M}\left\{\biggl(J[\psi^{(k)}]_{T,t}-
J[\psi^{(k)}]_{T,t}^{p,\ldots,p}\biggr)^{4}\right\}\le
$$

\vspace{1mm}
\begin{equation}
\label{fff5}
\le C_{2,k}
\left(
\int\limits_{[t,T]^k}
K^2(t_1,\ldots,t_k)
dt_1\ldots dt_k -\sum_{j_1=0}^{p}\ldots
\sum_{j_k=0}^{p}C^2_{j_k\ldots j_1}
\right)^2\le 
\frac{H_{2,k}}{p^2},
\end{equation}

\vspace{3mm}
\noindent
where 

\vspace{-2mm}
$$
C_{n,k}=(k!)^{n} (2n-1)^{nk}
$$

\vspace{5mm}
\noindent
and $H_{2,k}=G_k^2{C}_{2,k}.$

Let us consider Lemma 1 and put

\vspace{1mm}
$$
\xi_p=\biggl|J[\psi^{(k)}]_{T,t}-
J[\psi^{(k)}]_{T,t}^{p,\ldots,p}\biggr|
$$

\vspace{3mm}
\noindent
and $\alpha=4.$

Then from (\ref{fff5}) we obtain

\vspace{1mm}
$$
\sum\limits_{p=1}^{\infty}
{\sf M}\left\{\biggl(J[\psi^{(k)}]_{T,t}-
J[\psi^{(k)}]_{T,t}^{p,\ldots,p}\biggr)^{4}\right\}\le
$$

\vspace{1mm}
\begin{equation}
\label{qqq1}
\le H_{2,k}\sum\limits_{p=1}^{\infty}\frac{1}{p^2}<\infty.
\end{equation}

\vspace{4mm}

Using Lemma 1, from (\ref{qqq1}) we have

\vspace{2mm}
$$
J[\psi^{(k)}]_{T,t}^{p,\ldots,p}\ \to \ J[\psi^{(k)}]_{T,t}\ \ \ 
\hbox{if}\ \ \ p\to \infty
$$

\vspace{5mm}
\noindent
w.\ p.\ 1, where (see Theorem 1)

\vspace{2mm}
$$
J[\psi^{(k)}]_{T,t}^{p,\ldots,p}=
\sum_{j_1=0}^{p}\ldots\sum_{j_k=0}^{p}
C_{j_k\ldots j_1}\Biggl(
\prod_{l=1}^k\zeta_{j_l}^{(i_l)}\ -
\Biggr.
$$

\vspace{2mm}
\begin{equation}
\label{kk0}
-\ \Biggl.
\hbox{\vtop{\offinterlineskip\halign{
\hfil#\hfil\cr
{\rm l.i.m.}\cr
$\stackrel{}{{}_{N\to \infty}}$\cr
}} }\sum_{(l_1,\ldots,l_k)\in {\rm G}_k}
\phi_{j_{1}}(\tau_{l_1})
\Delta{\bf w}_{\tau_{l_1}}^{(i_1)}\ldots
\phi_{j_{k}}(\tau_{l_k})
\Delta{\bf w}_{\tau_{l_k}}^{(i_k)}\Biggr),
\end{equation}

\vspace{6mm}
\noindent
where $i_1,\ldots,i_k=1,\ldots,m$ in (\ref{kk0}). Theorem 9 is proved.

\vspace{2mm}

Taking into account (\ref{star00011}) and
(\ref{ddd1}), we obtain the following
inequality

\begin{equation}
\label{zsel1}
{\sf M}\left\{\left(
J[\psi^{(k)}]_{T,t}-J[\psi^{(k)}]_{T,t}^{p,\ldots,p}
\right)^2\right\}\le \frac{k! P_k (T-t)^k}{p},
\end{equation}

\vspace{4mm}
\noindent
where constant $P_k$ depends only on $k$.

The estimates
(\ref{99999}) and 
(\ref{ddd1}) imply the following 
inequality

\vspace{1mm}
$$
{\sf M}\left\{\left(
J[\psi^{(k)}]_{T,t}-J[\psi^{(k)}]_{T,t}^{p,\ldots,p}
\right)^{2n}\right\}\le 
$$

\vspace{3mm}
\begin{equation}
\label{xyzyx1}
\le (k!)^{n} (2n-1)^{nk}\
\frac{\left(P_k\right)^n (T-t)^{nk}}{p^n},
\end{equation}

\vspace{5mm}
\noindent
where $n\in\mathbb{N}$ and 
constant $P_k$ depends only on $k$.

Consider the question
on the rate of convergence w.~p.~1 in Theorem 9.
Using the inequality (\ref{xyzyx1}), we obtain

\begin{equation}
\label{xyzyx11}
\left({\sf M}\left\{\left(
J[\psi^{(k)}]_{T,t}-J[\psi^{(k)}]_{T,t}^{p,\ldots,p}
\right)^{2n}\right\}\right)^{1/2n}\le \frac{Q_{n,k}}{\sqrt{p}},
\end{equation}

\vspace{4mm}
\noindent
where $n\in \mathbb{N}$ and 

\vspace{1mm}
$$
Q_{n,k}=(2n-1)^{k/2}\ \sqrt{k!}\
\sqrt{P_k}\ (T-t)^{k/2}.
$$

\vspace{5mm}

According to the Lyapunov inequality and (\ref{xyzyx11}), we have

\vspace{1mm}
\begin{equation}
\label{xyzyx12}
\left({\sf M}\biggl\{\left(
J[\psi^{(k)}]_{T,t}-J[\psi^{(k)}]_{T,t}^{p,\ldots,p}
\right)^{n}\biggr\}\right)^{1/n}\le \frac{Q_{n,k}}{\sqrt{p}}
\end{equation}

\vspace{5mm}
\noindent
for all $n\in \mathbb{N}$. Following \cite{xyz1001} (Lemma 2.1), we get

\vspace{3mm}
$$
\biggl|J[\psi^{(k)}]_{T,t}-
J[\psi^{(k)}]_{T,t}^{p,\ldots,p}\biggr|=
\frac{p^{1/2 - \varepsilon}}{p^{1/2 - \varepsilon}}\biggl|J[\psi^{(k)}]_{T,t}-
J[\psi^{(k)}]_{T,t}^{p,\ldots,p}\biggr|\le
$$

\vspace{4mm}
\begin{equation}
\label{xyzyx13}
\le 
\frac{1}{p^{1/2 - \varepsilon}}
\sup\limits_{p\in \mathbb{N}}\left(p^{1/2 - \varepsilon}
\biggl|J[\psi^{(k)}]_{T,t}-
J[\psi^{(k)}]_{T,t}^{p,\ldots,p}\biggr|\right)=\frac{\eta_{\varepsilon}}{p^{1/2 - \varepsilon}}
\end{equation}

\vspace{6mm}
\noindent
w. p. 1, where
$$
\eta_{\varepsilon}=
\sup\limits_{p\in \mathbb{N}}\left(p^{1/2 - \varepsilon}
\biggl|J[\psi^{(k)}]_{T,t}-
J[\psi^{(k)}]_{T,t}^{p,\ldots,p}\biggr|\right)
$$

\vspace{5mm}
\noindent
and $\varepsilon>0$ is fixed.

For $q>1/\varepsilon,$ $q\in \mathbb{N}$ we obtain \cite{xyz1001} (see (\ref{xyzyx12}))

\vspace{3mm}
$$
{\sf M}\left\{\left|\eta_{\varepsilon}\right|^q\right\}=
{\sf M}\left\{\left(\sup\limits_{p\in \mathbb{N}}\left(p^{1/2 - \varepsilon}
\biggl|J[\psi^{(k)}]_{T,t}-
J[\psi^{(k)}]_{T,t}^{p,\ldots,p}\biggr|\right)\right)^q\right\}=
$$

\vspace{6mm}
$$
=
{\sf M}\left\{\sup\limits_{p\in \mathbb{N}}\left(p^{(1/2 - \varepsilon)q}
\biggl|J[\psi^{(k)}]_{T,t}-
J[\psi^{(k)}]_{T,t}^{p,\ldots,p}\biggr|^q\right)\right\}\le
$$

\vspace{5mm}
$$
\le {\sf M}\left\{\sum\limits_{p=1}^{\infty}p^{(1/2 - \varepsilon)q}
\biggl|J[\psi^{(k)}]_{T,t}-
J[\psi^{(k)}]_{T,t}^{p,\ldots,p}\biggr|^q\right\}=
$$

\vspace{5mm}
$$
= \sum\limits_{p=1}^{\infty}p^{(1/2 - \varepsilon)q}
{\sf M}\left\{\biggl|J[\psi^{(k)}]_{T,t}-
J[\psi^{(k)}]_{T,t}^{p,\ldots,p}\biggr|^q\right\}\le
$$

\vspace{5mm}
\begin{equation}
\label{xyzyx14}
\le
\sum\limits_{p=1}^{\infty}p^{(1/2 - \varepsilon)q}
\frac{\left(Q_{q,k}\right)^q}{p^{q/2}}=
\left(Q_{q,k}\right)^q\sum\limits_{p=1}^{\infty}\frac{1}{p^{\varepsilon q}}<\infty.
\end{equation}

\vspace{6mm}

From (\ref{xyzyx13}) we have that for all $\varepsilon>0$
there exists a random variable $\eta_{\varepsilon}$ such that 
the inequality (\ref{xyzyx13}) is fulfilled w.~p.~1 for all $p\in \mathbb{N}.$
Moreover, from the Lyapunov inequality and (\ref{xyzyx14}) we obtain
${\sf M}\left\{\left|\eta_{\varepsilon}\right|^q\right\}<\infty$
for all $q\ge 1.$

\vspace{5mm}

\section{About the Structure of Functions $K(t_1,\ldots,t_k)$
Used in Applications}

\vspace{5mm}

The systems of iterated  stochastic integrals 
(\ref{ito}), (\ref{str}), (\ref{k1000}), (\ref{k1001})
are part of the 
stochastic Taylor--Ito and Taylor--Stratonovich expansions 
(classical \cite{1995}, \cite{1988}
and unified \cite{3}-\cite{7}, \cite{9a}-\cite{10axx1}).

The function $K(t_1,\ldots,t_k)$ from Theorems 1, 2 for
the family (\ref{k1000}) looks as follows

\vspace{1mm}
\begin{equation}
\label{leto7000}
K(t_1,\ldots,t_k)=
(t-t_k)^{l_k}\ldots (t-t_1)^{l_1}\ {\bf 1}_{\{t_1<\ldots<t_k\}},\ \ \
t_1,\ldots,t_k\in[t, T],
\end{equation}

\vspace{4mm}
\noindent
where ${\bf 1}_A$ is the indicator of the set $A$.

In particular, for the stochastic integrals 

\vspace{1mm}
$$
I_{(1)T,t}^{(i_1)},\ I_{(2)T,t}^{(i_1)},\ I_{(00)T,t}^{(i_1i_2)},\
I_{(000)T,t}^{(i_1i_2i_3)},\ I_{(01)T,t}^{(i_1i_2)},\
I_{(10)T,t}^{(i_1i_2)},\
I_{(0000)T,t}^{(i_1 i_2 i_3 i_4)},
$$

\vspace{2mm}
$$
I_{(20)T,t}^{(i_1i_2)},\ I_{(11)T,t}^{(i_1i_2)},\
I_{(02)T,t}^{(i_1i_2)}\ \ (i_1,\ldots, i_4=1,\ldots,m)
$$

\vspace{5mm}
\noindent
the functions $K(t_1,\ldots,t_k)$ (see (\ref{leto7000}))
correspondently look as follows

\begin{equation}
\label{aaa1}
K_1(t_1)=t-t_1,\ \ \ K_2(t_1)=(t-t_1)^2,\ \ \
K_{00}(t_1,t_2)={\bf 1}_{\{t_1<t_2\}},
\end{equation}

\vspace{1mm}
\begin{equation}
\label{aaa2}
K_{000}(t_1,t_2,t_3)={\bf 1}_{\{t_1<t_2<t_3\}},\ \ \
K_{01}(t_1,t_2)=(t-t_2){\bf 1}_{\{t_1<t_2\}},
\end{equation}

\vspace{1mm}
\begin{equation}
\label{aaa3}
K_{10}(t_1,t_2)=(t-t_1){\bf 1}_{\{t_1<t_2\}},\ \ \
K_{0000}(t_1,t_2)={\bf 1}_{\{t_1<t_2<t_3<t_4\}},
\end{equation}

\vspace{1mm}
\begin{equation}
\label{aaa4}
K_{20}(t_1,t_2)=(t-t_1)^2{\bf 1}_{\{t_1<t_2\}},\ \ \
K_{11}(t_1,t_2)=(t-t_1)(t-t_2){\bf 1}_{\{t_1<t_2\}},
\end{equation}

\vspace{1mm}
\begin{equation}
\label{aaa5}
K_{02}(t_1,t_2)=(t-t_2)^2{\bf 1}_{\{t_1<t_2\}},
\end{equation}

\vspace{3mm}
\noindent
where $t_1,\ldots, t_4\in [t, T].$

It is obviously that the most simple 
expansion 
for the polynomial of a finite degree
into the Fourier series using the 
complete orthonormal system of functions in the space $L_2([t, T])$
will be its Fourier--Legendre 
expansion (finite sum). The polynomial functions 
are included in the functions (\ref{aaa1})--(\ref{aaa5}) 
as their components if $l_1^2+\ldots+l_k^2>0$.
So, it is logical to expect that the most simple expansions 
for the functions
(\ref{aaa1})--(\ref{aaa5}) into multiple Fourier series will 
be their Fourier--Legendre expansions
when $l_1^2+\ldots+l_k^2>0$.

Note that the given assumption is confirmed completely 
(compare the formulas (\ref{4002}), (\ref{4006}) with the formulas
(\ref{444}), (\ref{t1}) correspondently).
So, the usage of Legendre polynomials in the considered scientific
field is 
an obvious step forward.

\vspace{5mm}

\section{Theorems 1--7 from Point
of View of the Wong--Zakai Approximation}

\vspace{5mm}

The iterated Ito stochastic integrals and solutions
of Ito SDEs are complex and important functionals
from the independent components ${\bf f}_{s}^{(i)},$
$i=1,\ldots,m$ of the multidimensional
Wiener process ${\bf f}_{s},$ $s\in[0, T].$
Let ${\bf f}_{s}^{(i)p},$ $p\in\mathbb{N}$ 
be some approximation of
${\bf f}_{s}^{(i)},$
$i=1,\ldots,m$.
Suppose that 
${\bf f}_{s}^{(i)p}$
converges to
${\bf f}_{s}^{(i)},$
$i=1,\ldots,m$ if $p\to\infty$ in some sense and has
differentiable sample trajectories.

A natural question arises: if we replace 
${\bf f}_{s}^{(i)}$
by ${\bf f}_{s}^{(i)p},$
$i=1,\ldots,m$ in the functionals
mentioned above, will the resulting
functionals converge to the original
functionals from the components 
${\bf f}_{s}^{(i)},$
$i=1,\ldots,m$ of the multidimentional
Wiener process ${\bf f}_{s}$?
The answere to this question is negative 
in the general case. However, 
in the pioneering works of Wong E. and Zakai M. \cite{W-Z-1},
\cite{W-Z-2},
it was shown that under the special conditions and 
for some types of approximations 
of the Wiener process the answere is affirmative
with one peculiarity: the convergence takes place 
to the iterated Stratonovich stochastic integrals
and solutions of Stratonovich SDEs and not to iterated 
Ito stochastic integrals and solutions
of Ito SDEs.
The piecewise 
linear approximation 
as well as the regularization by convolution 
\cite{W-Z-1}-\cite{Watanabe} relate the 
mentioned types of approximations
of the Wiener process. The above approximation 
of stochastic integrals and solutions of SDEs 
is often called the Wong--Zakai approximation.

Let ${\bf w}_{\tau},$ $\tau\in[0, T]$ is a random vector with 
an $m+1$ components: ${\bf w}_{\tau}^{(i)}={\bf f}_{\tau}^{(i)}$ 
for $i=1,\ldots,m$ and 
${\bf w}_{\tau}^{(0)}=\tau,$\ 
${\bf f}_{\tau}^{(i)}$ $(i=1,\ldots,m)$
are independent standard Wiener processes.

It is well known that the following representation 
takes place \cite{Lipt}, \cite{7e}

\begin{equation}
\label{um1x}
{\bf w}_{\tau}^{(i)}-{\bf w}_{t}^{(i)}=
\sum_{j=0}^{\infty}\int\limits_t^{\tau}
\phi_j(s)ds\ \zeta_j^{(i)},\ \ \ \zeta_j^{(i)}=
\int\limits_t^T \phi_j(s)d{\bf w}_s^{(i)},
\end{equation}

\vspace{4mm}
\noindent
where $\tau\in[t, T],$ $t\ge 0,$
$\{\phi_j(x)\}_{j=0}^{\infty}$ is an arbitrary complete 
orthonormal system of functions in the space $L_2([t, T]),$ and
$\zeta_j^{(i)}$ are independent standard Gaussian 
random variables for various $i$ or $j.$
Moreover, the series (\ref{um1x}) converges for any $\tau\in [t, T]$
in the mean-square sense.

Let ${\bf w}_{\tau}^{(i)p}-{\bf w}_{t}^{(i)p}$ be 
the mean-square approximation of the process
${\bf w}_{\tau}^{(i)}-{\bf w}_{t}^{(i)},$
which has the following form

\vspace{-3mm}
\begin{equation}
\label{um1xx}
{\bf w}_{\tau}^{(i)p}-{\bf w}_{t}^{(i)p}=
\sum_{j=0}^{p}\int\limits_t^{\tau}
\phi_j(s)ds\ \zeta_j^{(i)}.
\end{equation}

\vspace{3mm}

From (\ref{um1xx}) we obtain

\vspace{-4mm}
\begin{equation}
\label{um1xxx}
d{\bf w}_{\tau}^{(i)p}=
\sum_{j=0}^{p}
\phi_j(\tau)\zeta_j^{(i)} d\tau.
\end{equation}

\vspace{4mm}

Consider the following iterated Riemann--Stieltjes
integral

\begin{equation}
\label{um1xxxx}
\int\limits_t^T
\psi_k(t_k)\ldots \int\limits_t^{t_2}\psi_1(t_1)
d{\bf w}_{t_1}^{(i_1)p_1}\ldots d{\bf w}_{t_k}^{(i_k)p_k},
\end{equation}

\vspace{4mm}
\noindent
where $i_1,\ldots,i_k=0,1,\ldots,m,$\ \ $p_1,\ldots,p_k\in\mathbb{N},$

\begin{equation}
\label{um1xxx1}
d{\bf w}_{\tau}^{(i)p}=
\left\{\begin{matrix}
d{\bf f}_{\tau}^{(i)p}\ &\hbox{\rm for}\ \ \ i=1,\ldots,m\cr\cr\cr
d\tau^p\ &\hbox{\rm for}\ \ \ i=0
\end{matrix}
,\right.
\end{equation}

\vspace{4mm}
\noindent
and $d{\bf f}_{\tau}^{(i)p},$ $d\tau^p$ are defined by the relation (\ref{um1xxx}).

Let us substitute (\ref{um1xxx}) into (\ref{um1xxxx})

\begin{equation}
\label{um1xxxx1}
\int\limits_t^T
\psi_k(t_k)\ldots \int\limits_t^{t_2}\psi_1(t_1)
d{\bf w}_{t_1}^{(i_1)p_1}\ldots d{\bf w}_{t_k}^{(i_k)p_k}=
\sum\limits_{j_1=0}^{p_1}\ldots \sum\limits_{j_k=0}^{p_k}
C_{j_k \ldots j_1}\prod\limits_{l=1}^k \zeta_{j_l}^{(i_l)},
\end{equation}

\vspace{4mm}
\noindent
where 
$$
\zeta_j^{(i)}=\int\limits_t^T \phi_j(s)d{\bf w}_s^{(i)}
$$ 

\vspace{2mm}
\noindent
are independent standard Gaussian random variables for various 
$i$ or $j$ (in the case when $i\ne 0$),
${\bf w}_{s}^{(i)}={\bf f}_{s}^{(i)}$ for
$i=1,\ldots,m$ and 
${\bf w}_{s}^{(0)}=s,$

$$
C_{j_k \ldots j_1}=\int\limits_t^T\psi_k(t_k)\phi_{j_k}(t_k)\ldots
\int\limits_t^{t_2}
\psi_1(t_1)\phi_{j_1}(t_1)
dt_1\ldots dt_k
$$

\vspace{4mm}
\noindent
is the Fourier coefficient.

To best of our knowledge \cite{W-Z-1}-\cite{Watanabe}
the approximations of the Wiener process
in the Wong--Zakai approximation must satisfy fairly strong
restrictions
\cite{Watanabe}
(see Definition 7.1, pp.~480--481).
Moreover, approximations of the Wiener process that are
similar to (\ref{um1xx})
were not considered in \cite{W-Z-1}, \cite{W-Z-2}
(also see \cite{Watanabe}, Theorems 7.1, 7.2).
Therefore, the proof of analogs of Theorems 7.1 and 7.2 \cite{Watanabe}
for approximations of the Wiener 
process based on its series expansion (\ref{um1x})
should be carried out separately.
Thus, the mean-square convergence of the right-hand side
of (\ref{um1xxxx1}) to the iterated Stratonovich stochastic integral 
(\ref{str})
does not follow from the results of the papers
\cite{W-Z-1}, \cite{W-Z-2} (also see \cite{Watanabe},
Theorems 7.1, 7.2).

From the other hand, Theorems 1--7 from this 
paper can be considered as the proof of the
Wong--Zakai approximation for the iterated 
Stratonovich stochastic integrals (\ref{str}) of multiplicities 1 to 6
based on the the Riemann--Stieltjes integrals (\ref{um1xxxx}) and
approximation (\ref{um1xx}) of the Wiener process.
At that, the Riemann--Stieltjes integrals (\ref{um1xxxx})  converge
(according to Theorems 1--7)
to the appropriate Stratonovich 
stochastic integrals (\ref{str}). Recall that
$\{\phi_j(x)\}_{j=0}^{\infty}$ (see (\ref{um1x}), (\ref{um1xx}), and
Theorems 3--7)
is a complete 
orthonormal system of Legendre polynomials or 
trigonometric functions 
in the space $L_2([t, T])$.

To illustrate the above reasoning, 
consider two examples for the case $k=2,$
$\psi_1(s),$ $\psi_2(s)\equiv 1;$ $i_1, i_2=1,\ldots,m.$

The first example relates to the piecewise linear approximation
of the multidimensional Wiener process (these approximations 
were considered in \cite{W-Z-1}-\cite{Watanabe}).

Let ${\bf b}_{\Delta}^{(i)}(t),$ $t\in[0, T]$ be the piecewise
linear approximation of the $i$th component ${\bf f}_t^{(i)}$
of the multidimensional standard Wiener process ${\bf f}_t,$
$t\in [0, T]$ with independent components
${\bf f}_t^{(i)},$ $i=1,\ldots,m,$ i.e.

$$
{\bf b}_{\Delta}^{(i)}(t)={\bf f}_{k\Delta}^{(i)}+
\frac{t-k\Delta}{\Delta}\Delta{\bf f}_{k\Delta}^{(i)},
$$

\vspace{3mm}
\noindent
where 

\vspace{-2mm}
$$
\Delta{\bf f}_{k\Delta}^{(i)}={\bf f}_{(k+1)\Delta}^{(i)}-
{\bf f}_{k\Delta}^{(i)},\ \ \
t\in[k\Delta, (k+1)\Delta),\ \ \ k=0, 1,\ldots, N-1.
$$

\vspace{4mm}

Note that w.~p.~1

\vspace{-1mm}
\begin{equation}
\label{pridum}
\frac{d{\bf b}_{\Delta}^{(i)}}{dt}(t)=
\frac{\Delta{\bf f}_{k\Delta}^{(i)}}{\Delta},\ \ \
t\in[k\Delta, (k+1)\Delta),\ \ \ k=0, 1,\ldots, N-1.
\end{equation}

\vspace{4mm}

Consider the following iterated Riemann--Stieltjes
integral

\vspace{1mm}
$$
\int\limits_0^T
\int\limits_0^{s}
d{\bf b}_{\Delta}^{(i_1)}(\tau)d{\bf b}_{\Delta}^{(i_2)}(s),\ \ \ 
i_1,i_2=1,\ldots,m.
$$

\vspace{4mm}

Using (\ref{pridum}) and additive property of Riemann--Stieltjes integrals, 
we can write w.~p.~1

\vspace{2mm}
$$
\int\limits_0^T
\int\limits_0^{s}
d{\bf b}_{\Delta}^{(i_1)}(\tau)d{\bf b}_{\Delta}^{(i_2)}(s)=
\int\limits_0^T
\int\limits_0^{s}
\frac{d{\bf b}_{\Delta}^{(i_1)}}{d\tau}(\tau)d\tau
\frac{d {\bf b}_{\Delta}^{(i_2)}}{d s}(s)
ds =
$$

\vspace{4mm}
$$
=
\sum\limits_{l=0}^{N-1}\int\limits_{l\Delta}^{(l+1)\Delta}
\left(
\sum\limits_{q=0}^{l-1}\int\limits_{q\Delta}^{(q+1)\Delta}
\frac{\Delta{\bf f}_{q\Delta}^{(i_1)}}{\Delta}d\tau+
\int\limits_{l\Delta}^{s}
\frac{\Delta{\bf f}_{l\Delta}^{(i_1)}}{\Delta}d\tau\right)
\frac{\Delta{\bf f}_{l\Delta}^{(i_2)}}{\Delta}ds=
$$

\vspace{4mm}
$$
=\sum\limits_{l=0}^{N-1}\sum\limits_{q=0}^{l-1}
\Delta{\bf f}_{q\Delta}^{(i_1)}
\Delta{\bf f}_{l\Delta}^{(i_2)}+
\frac{1}{\Delta^2}\sum\limits_{l=0}^{N-1}
\Delta{\bf f}_{l\Delta}^{(i_1)}
\Delta{\bf f}_{l\Delta}^{(i_2)}
\int\limits_{l\Delta}^{(l+1)\Delta}
\int\limits_{l\Delta}^{s}d\tau ds=
$$

\vspace{4mm}
\begin{equation}
\label{oh-ty}
=\sum\limits_{l=0}^{N-1}\sum\limits_{q=0}^{l-1}
\Delta{\bf f}_{q\Delta}^{(i_1)}
\Delta{\bf f}_{l\Delta}^{(i_2)}+
\frac{1}{2}\sum\limits_{l=0}^{N-1}
\Delta{\bf f}_{l\Delta}^{(i_1)}
\Delta{\bf f}_{l\Delta}^{(i_2)}.
\end{equation}

\vspace{6mm}

Using (\ref{oh-ty}) and the standard relation between Stratonovich
and Ito stochastic integrals, it 
is not difficult to show 
that

\vspace{1mm}
$$
\hbox{\vtop{\offinterlineskip\halign{
\hfil#\hfil\cr
{\rm l.i.m.}\cr
$\stackrel{}{{}_{N\to \infty}}$\cr
}} }
\int\limits_0^T
\int\limits_0^{s}
d{\bf b}_{\Delta}^{(i_1)}(\tau)d{\bf b}_{\Delta}^{(i_2)}(s)=
\int\limits_0^T
\int\limits_0^{s}
d{\bf f}_{\tau}^{(i_1)}d{\bf f}_{s}^{(i_2)}+
\frac{1}{2}{\bf 1}_{\{i_1=i_2\}}\int\limits_0^T ds=
$$

\vspace{4mm}
\begin{equation}
\label{uh-111}
=
\int\limits_0^{*T}
\int\limits_0^{*s}
d{\bf f}_{\tau}^{(i_1)}d{\bf f}_{s}^{(i_2)},
\end{equation}

\vspace{6mm}
\noindent
where $\Delta\to 0$ if $N\to\infty$ ($N\Delta=T$).

Obviously, (\ref{uh-111}) agrees with Theorem 7.1 (see \cite{Watanabe},
p.~486).

The next example relates to the approximation
of the Wiener process based on its series expansion
(\ref{um1x}) for $t=0$, where
$\{\phi_j(x)\}_{j=0}^{\infty}$ 
is a complete 
orthonormal system of Legendre polynomials or 
trigonometric functions 
in the space $L_2([0, T])$.

Consider the following iterated Riemann--Stieltjes
integral

\begin{equation}
\label{abcd1}
\int\limits_0^T
\int\limits_0^{s}
d{\bf f}_{\tau}^{(i_1)p}d{\bf f}_{s}^{(i_2)p},\ \ \ 
i_1,i_2=1,\ldots,m,
\end{equation}

\vspace{4mm}
\noindent
where $d{\bf f}_{\tau}^{(i)p}$ is defined by the
relation
(\ref{um1xxx}).

Let us substitute (\ref{um1xxx}) into (\ref{abcd1}) 

\vspace{-1mm}
\begin{equation}
\label{set18}
\int\limits_0^T
\int\limits_0^{s}
d{\bf f}_{\tau}^{(i_1)p}d{\bf f}_{s}^{(i_2)p}=
\sum\limits_{j_1,j_2=0}^p
C_{j_2 j_1} \zeta_{j_1}^{(i_1)}\zeta_{j_2}^{(i_2)},
\end{equation}

\vspace{3mm}
\noindent
where 
$$
C_{j_2 j_1}=
\int\limits_0^T \phi_{j_2}(s)\int\limits_0^s
\phi_{j_1}(\tau)d\tau ds
$$

\vspace{4mm}
\noindent
is the Fourier coefficient; another notations 
are the same as in (\ref{um1xxxx1}).

As we noted above, approximations of the Wiener process that are
similar to (\ref{um1xx})
were not considered in \cite{W-Z-1}, \cite{W-Z-2}
(also see Theorems 7.1, 7.2 in \cite{Watanabe}).
Furthermore, the extension of the results of Theorems 7.1 and 7.2
\cite{Watanabe} to the case under consideration is
not obvious.

On the other hand, we can apply the theory built in Chapters 1 and 2
of the monographs \cite{10a}-\cite{10axx1}. More precisely, 
using 
Theorem 3 for the case $k=2$,  
we obtain from (\ref{set18}) the desired result

$$
\hbox{\vtop{\offinterlineskip\halign{
\hfil#\hfil\cr
{\rm l.i.m.}\cr
$\stackrel{}{{}_{p\to \infty}}$\cr
}} }
\int\limits_0^T
\int\limits_0^{s}
d{\bf f}_{\tau}^{(i_1)p}d{\bf f}_{s}^{(i_2)p}=
\hbox{\vtop{\offinterlineskip\halign{
\hfil#\hfil\cr
{\rm l.i.m.}\cr
$\stackrel{}{{}_{p\to \infty}}$\cr
}} }
\sum\limits_{j_1,j_2=0}^p
C_{j_2 j_1} \zeta_{j_1}^{(i_1)}\zeta_{j_2}^{(i_2)}=
$$

\vspace{3mm}
\begin{equation}
\label{umen-bl}
=
\int\limits_0^{*T}
\int\limits_0^{*s}
d{\bf f}_{\tau}^{(i_1)}d{\bf f}_{s}^{(i_2)}.
\end{equation}

\vspace{6mm}

From the other hand, by Theorem 1
(see (\ref{leto5001})) for the case
$k=2$ we obtain from (\ref{set18}) the following relation

$$
\hbox{\vtop{\offinterlineskip\halign{
\hfil#\hfil\cr
{\rm l.i.m.}\cr
$\stackrel{}{{}_{p\to \infty}}$\cr
}} }
\int\limits_0^T
\int\limits_0^{s}
d{\bf f}_{\tau}^{(i_1)p}d{\bf f}_{s}^{(i_2)p}=
\hbox{\vtop{\offinterlineskip\halign{
\hfil#\hfil\cr
{\rm l.i.m.}\cr
$\stackrel{}{{}_{p\to \infty}}$\cr
}} }
\sum\limits_{j_1,j_2=0}^p
C_{j_2 j_1} \zeta_{j_1}^{(i_1)}\zeta_{j_2}^{(i_2)}=
$$

\vspace{3mm}
$$
=
\hbox{\vtop{\offinterlineskip\halign{
\hfil#\hfil\cr
{\rm l.i.m.}\cr
$\stackrel{}{{}_{p\to \infty}}$\cr
}} }
\sum\limits_{j_1,j_2=0}^p
C_{j_2 j_1} \biggl(\zeta_{j_1}^{(i_1)}\zeta_{j_2}^{(i_2)}-
{\bf 1}_{\{i_1=i_2\}}{\bf 1}_{\{j_1=j_2\}}\biggr)+
{\bf 1}_{\{i_1=i_2\}}\sum\limits_{j_1=0}^{\infty}
C_{j_1 j_1}=
$$

\vspace{3mm}
\begin{equation}
\label{umen-blx}
=
\int\limits_0^T
\int\limits_0^{s}
d{\bf f}_{\tau}^{(i_1)}d{\bf f}_{s}^{(i_2)}+
{\bf 1}_{\{i_1=i_2\}}\sum\limits_{j_1=0}^{\infty}
C_{j_1 j_1}.
\end{equation}

\vspace{6mm}

Since
$$
\sum\limits_{j_1=0}^{\infty}
C_{j_1 j_1}=\frac{1}{2}\sum\limits_{j_1=0}^{\infty}
\left(\int\limits_0^T \phi_j(\tau)d\tau\right)^2
=\frac{1}{2}
\left(\int\limits_0^T \phi_0(\tau)d\tau\right)^2=\frac{1}{2}
\int\limits_0^T ds,
$$

\vspace{4mm}
\noindent
then from (\ref{umen-blx}) 
and the standard relation between Stratonovich
and Ito stochastic integrals we obtain (\ref{umen-bl}).

\vspace{5mm}

\section{Exact Calculation 
of the Mean-Square Approximation Errors
for Iterated Stratonovich Stochastic Integrals
$I_{(0)T,t}^{*(i_1)},$
$I_{(1)T,t}^{*(i_1)},$
$I_{(00)T,t}^{*(i_1i_2)},$
$I_{(000)T,t}^{*(i_1i_2i_3)},$ $I_{(0000)T,t}^{*(i_1i_2i_3i_4)}$}

\vspace{5mm}

First, consider the question on the exact calculation 
of the mean-square approximation errors
for the following iterated Stratonovich stochastic integrals

\begin{equation}
\label{dest1}
I_{(0)T,t}^{*(i_1)},\ \ \ 
I_{(1)T,t}^{*(i_1)},\ \ \ 
I_{(00)T,t}^{*(i_1i_2)},\ \ \ 
I_{(000)T,t}^{*(i_1i_2i_3)}\ \ \ 
(i_1, i_2, i_3=1,\ldots,m)
\end{equation}

\vspace{2mm}
\noindent
defined by (\ref{k1001}).

We assume that the stochastic integrals (\ref{dest1})
are approximated using Theorems 1, 3 and the Legendre
polynomial system. Since
$I_{(0)T,t}^{(i_1)}=I_{(0)T,t}^{*(i_1)},$
$I_{(1)T,t}^{(i_1)}=I_{(1)T,t}^{*(i_1)}$\ w.~p.~1 (see (\ref{k1000})),
then we can use (\ref{4001}), (\ref{4002}) 
to approximate the stochastic integrals 
$I_{(0)T,t}^{*(i_1)},$
$I_{(1)T,t}^{*(i_1)}.$ In this case, we will have zero
mean-square approximation errors.

To approximate the iterated Stratonovich stochastic integral 
$I_{(00)T,t}^{*(i_1i_2)}$ 
we can use the formula (see (\ref{4004}))
\begin{equation}                                                     
\label{dest2}
I_{(00)T,t}^{*(i_1 i_2)q}=
\frac{T-t}{2}\left(\zeta_0^{(i_1)}\zeta_0^{(i_2)}+\sum_{i=1}^{q}
\frac{1}{\sqrt{4i^2-1}}\left(
\zeta_{i-1}^{(i_1)}\zeta_{i}^{(i_2)}-
\zeta_i^{(i_1)}\zeta_{i-1}^{(i_2)}\right)\right).
\end{equation}

\vspace{4mm}

The mean-square approximation error for (\ref{dest2})
will be determined by the formula (\ref{fff09}) $(i_1\ne i_2)$.
For the case $i_1=i_2$ we can use the well known equality

$$
I_{(00)T,t}^{*(i_1 i_1)}
=\frac{T-t}{2}\left(\zeta_0^{(i_1)}\right)^2\ \ \ \hbox{w.~p.~1.}
$$

\vspace{3mm}

Consider now the iterated Stratonovich stochastic integral
$I_{(000)T,t}^{*(i_1i_2i_3)}$ of multiplicity 3
$(i_1, i_2, i_3=$ $1,\ldots,m)$.
For the case of pairwise different 
$i_1, i_2, i_3$ we have the following relation

$$
I_{(000)T,t}^{*(i_1 i_2 i_3)}=I_{(000)T,t}^{(i_1 i_2 i_3)}\ \ \ \hbox{w.~p.~1}.
$$

\vspace{3mm}

Thus, in this case we can use the formulas (\ref{38}) and (\ref{39}).
For the case $i_1=i_2=i_3,$ to approximate the stochastic integral
$I_{(000)T,t}^{*(i_1i_1i_1)},$ we use the formula (\ref{dest4}).

Thus, it remains to consider the following three cases

\vspace{-2mm}
\begin{equation}
\label{dest5}
i_1=i_2\ne i_3,
\end{equation}
\begin{equation}
\label{dest6}
i_1\ne i_2=i_3,
\end{equation}
\begin{equation}
\label{dest7}
i_1=i_3\ne i_2.
\end{equation}

\vspace{3mm}

Taking into account the standard relations between 
Ito and Stratonovich stochastic integrals and
Theorem 1 (the case $k=3$) together with Theorem 3, we obtain

\vspace{-1mm}
$$
{\sf M}\left\{\left(I_{(000)T,t}^{*(i_1i_2i_3)}-
I_{(000)T,t}^{*(i_1i_2i_3)q}\right)^2\right\}=
$$

$$
={\sf M}\left\{\left(
I_{(000)T,t}^{(i_1i_2i_3)}+
\frac{1}{2}{\bf 1}_{\{i_1=i_2\}}
\int\limits_t^T 
\int\limits_t^{\tau}dsd{\bf f}_{\tau}^{(i_3)}
+\frac{1}{2}{\bf 1}_{\{i_2=i_3\}}
\int\limits_t^T \hspace{-1mm}\int\limits_t^{\tau}d{\bf f}_{s}^{(i_1)}d\tau-
I_{(000)T,t}^{*(i_1i_2i_3)q}\right)^{2}\right\}
=
$$

$$
={\sf M}\left\{\left(
I_{(000)T,t}^{(i_1i_2i_3)}-I_{(000)T,t}^{(i_1i_2i_3)q}+
I_{(000)T,t}^{(i_1i_2i_3)q}+
{\bf 1}_{\{i_1=i_2\}}
\frac{1}{2}\int\limits_t^{T}
\int\limits_{t}^{\tau}
dsd{\bf f}_{\tau}^{(i_3)}
\right.\right.+
$$

\begin{equation}
\label{tango3}
\left.\left.
+{\bf 1}_{\{i_2=i_3\}}
\frac{1}{2}\int\limits_t^T\int\limits_t^{\tau}d{\bf f}_{s}^{(i_1)}d\tau-
I_{(000)T,t}^{*(i_1i_2i_3)q}\right)^2\right\},
\end{equation}

\vspace{4mm}
\noindent
where the approximations $I_{(000)T,t}^{*(i_1i_2i_3)q},$
$I_{(000)T,t}^{(i_1i_2i_3)q}$ are defined by 
the relations (see (\ref{good1}), (\ref{zzz1})) 

\vspace{1mm}
$$
I_{(000)T,t}^{(i_1i_2i_3)q}=
\sum_{j_1,j_2,j_3=0}^{q}
C_{j_3j_2j_1}\Biggl(
\zeta_{j_1}^{(i_1)}\zeta_{j_2}^{(i_2)}\zeta_{j_3}^{(i_3)}
-{\bf 1}_{\{i_1=i_2\}}
{\bf 1}_{\{j_1=j_2\}}
\zeta_{j_3}^{(i_3)}-
\Biggr.
$$
\begin{equation}
\label{tango1}
\Biggl.
-{\bf 1}_{\{i_2=i_3\}}
{\bf 1}_{\{j_2=j_3\}}
\zeta_{j_1}^{(i_1)}-
{\bf 1}_{\{i_1=i_3\}}
{\bf 1}_{\{j_1=j_3\}}
\zeta_{j_2}^{(i_2)}\Biggr),
\end{equation}

\vspace{2mm}
\begin{equation}
\label{tango2}
I_{(000)T,t}^{*(i_1i_2i_3)q}=
\sum_{j_1,j_2,j_3=0}^{q}
C_{j_3j_2j_1}
\zeta_{j_1}^{(i_1)}\zeta_{j_2}^{(i_2)}\zeta_{j_3}^{(i_3)}.
\end{equation}

\vspace{6mm} 

Substituting (\ref{tango1}) and (\ref{tango2}) into (\ref{tango3}) yields

\vspace{1mm}
$$
{\sf M}\left\{\left(I_{(000)T,t}^{*(i_1i_2i_3)}-
I_{(000)T,t}^{*(i_1i_2i_3)q}\right)^2\right\}=
$$

$$
={\sf M}\left\{\left(I_{(000)T,t}^{(i_1i_2i_3)}-
I_{(000)T,t}^{(i_1i_2i_3)q}+
{\bf 1}_{\{i_1=i_2\}}
\left(\frac{1}{2}\int\limits_t^T
\int\limits_t^{\tau}dsd{\bf f}_{\tau}^{(i_3)}-
\sum_{j_1,j_3=0}^{q}
C_{j_3j_1j_1}
\zeta_{j_3}^{(i_3)}\right)+\right.\right.
$$

\begin{equation}
\label{dest10}
\left.\left.+{\bf 1}_{\{i_2=i_3\}}\left(
\frac{1}{2}\int\limits_t^T\int\limits_t^{\tau}d{\bf f}_{s}^{(i_1)}d\tau-
\sum_{j_1,j_3=0}^{q}
C_{j_3j_3j_1}
\zeta_{j_1}^{(i_1)}\right)
-{\bf 1}_{\{i_1=i_3\}}
\sum_{j_1,j_2=0}^{q}
C_{j_1j_2j_1}
\zeta_{j_2}^{(i_2)}\right)^{2}\right\}.\ \ \ \ \ \ \ \ \ 
\end{equation}

\vspace{5mm}

Consider the case (\ref{dest5}). From (\ref{dest10}) we obtain

$$
{\sf M}\left\{\left(I_{(000)T,t}^{*(i_1i_2i_3)}-
I_{(000)T,t}^{*(i_1i_2i_3)q}\right)^2\right\}=
$$

\begin{equation}
\label{dest11}
={\sf M}\left\{\left(I_{(000)T,t}^{(i_1i_2i_3)}-
I_{(000)T,t}^{(i_1i_2i_3)q}+
\frac{1}{2}\int\limits_t^T
\int\limits_t^{\tau}dsd{\bf f}_{\tau}^{(i_3)}-
\sum_{j_1,j_3=0}^{q}
C_{j_3j_1j_1}
\zeta_{j_3}^{(i_3)}\right)^2\right\}.
\end{equation}

\vspace{5mm}

According to the results of Sect.~3 in \cite{15b}
(also see Sect.~1.2.2 in \cite{10a}-\cite{10axx1}), the quantity
                                
$$
I_{(000)T,t}^{(i_1i_2i_3)}-
I_{(000)T,t}^{(i_1i_2i_3)q}
$$

\vspace{3mm}
\noindent
includes only iterated Ito stochastic integrals
of multiplicity 3. At the same time, the quantity

$$
\frac{1}{2}\int\limits_t^T
\int\limits_t^{\tau}dsd{\bf f}_{\tau}^{(i_3)}-
\sum_{j_1,j_3=0}^{q}
C_{j_3j_1j_1}
\zeta_{j_3}^{(i_3)}
$$

\vspace{3mm}
\noindent
contains only iterated Ito stochastic integrals
of multiplicity 1. This means that from (\ref{dest11}) we get

\vspace{2mm}
$$
{\sf M}\left\{\left(I_{(000)T,t}^{*(i_1i_2i_3)}-
I_{(000)T,t}^{*(i_1i_2i_3)q}\right)^2\right\}
={\sf M}\left\{\left(I_{(000)T,t}^{(i_1i_2i_3)}-
I_{(000)T,t}^{(i_1i_2i_3)q}\right)^2\right\}+
$$

\begin{equation}
\label{dest12}
+{\sf M}\left\{\left(\frac{1}{2}\int\limits_t^T
(\tau-t)d{\bf f}_{\tau}^{(i_3)}-
\sum_{j_1,j_3=0}^{q}
C_{j_3j_1j_1}
\zeta_{j_3}^{(i_3)}\right)^2\right\}.
\end{equation}

\vspace{3mm}

We have

\vspace{1mm}
$$
{\sf M}\left\{\left(\frac{1}{2}\int\limits_t^T
(\tau-t)d{\bf f}_{\tau}^{(i_3)}-
\sum_{j_1,j_3=0}^{q}
C_{j_3j_1j_1}
\zeta_{j_3}^{(i_3)}\right)^2\right\}=
\frac{1}{4}\int\limits_t^T
(\tau-t)^2 d\tau-
$$

\begin{equation}
\label{dest15}
-\sum_{j_1,j_3=0}^{q}
C_{j_3j_1j_1}\int\limits_t^T
(\tau-t)\phi_{j_3}(\tau)d\tau+
\sum_{j_3=0}^{q}\left(\sum_{j_1=0}^{q}
C_{j_3j_1j_1}\right)^2,
\end{equation}

\vspace{5mm}
\noindent
where $\phi_{j_3}(\tau)$ is the Legendre polynomial defined by (\ref{4009}). 

According to the properties of Legendre polynomials, we obtain

\begin{equation}
\label{dest16}
\int\limits_t^T
(\tau-t)\phi_{j_3}(\tau)d\tau=\frac{(T-t)^{3/2}}{2}
\left\{
\begin{matrix}
1,\ & j_3=0\cr\cr
1/\sqrt{3},\ & j_3=1\cr\cr
0,\ & j_3\ge 2
\end{matrix}
.\right.
\end{equation}

\vspace{4mm}

Combining (\ref{dest12})--(\ref{dest16}) and (\ref{39c}), we get

\vspace{1mm}
$$
{\sf M}\left\{\left(I_{(000)T,t}^{*(i_1i_2i_3)}-
I_{(000)T,t}^{*(i_1i_2i_3)q}\right)^2\right\}=
\frac{(T-t)^3}{4}
-\sum_{j_1,j_2,j_3=0}^q C_{j_3j_2j_1}^2-
\sum_{j_1,j_2,j_3=0}^q C_{j_3j_1j_2}C_{j_3j_2j_1}-
$$

\begin{equation}
\label{dest80}
-\frac{(T-t)^{3/2}}{2}
\sum_{j_1=0}^{q}
\left(C_{0j_1j_1}+\frac{1}{\sqrt{3}}C_{1j_1j_1}\right)+
\sum_{j_3=0}^{q}\left(\sum_{j_1=0}^{q}
C_{j_3j_1j_1}\right)^2,
\end{equation}

\vspace{5mm}
\noindent
where $i_1=i_2\ne i_3.$

Consider the case (\ref{dest6}). From (\ref{dest10}) we obtain

\vspace{1mm}
$$
{\sf M}\left\{\left(I_{(000)T,t}^{*(i_1i_2i_3)}-
I_{(000)T,t}^{*(i_1i_2i_3)q}\right)^2\right\}=
$$

\vspace{1mm}
$$
={\sf M}\left\{\left(I_{(000)T,t}^{(i_1i_2i_3)}-
I_{(000)T,t}^{(i_1i_2i_3)q}+
\frac{1}{2}\int\limits_t^T
\int\limits_t^{\tau}d{\bf f}_{s}^{(i_1)}d\tau-
\sum_{j_1,j_3=0}^{q}
C_{j_3j_3j_1}
\zeta_{j_1}^{(i_1)}\right)^2\right\}=
$$

\vspace{1mm}
$$
={\sf M}\left\{\left(I_{(000)T,t}^{(i_1i_2i_3)}-
I_{(000)T,t}^{(i_1i_2i_3)q}+
\frac{1}{2}\int\limits_t^T
(T-s)d{\bf f}_{s}^{(i_1)}-
\sum_{j_1,j_3=0}^{q}
C_{j_3j_3j_1}
\zeta_{j_1}^{(i_1)}\right)^2\right\}=
$$

\vspace{1mm}
$$
={\sf M}\left\{\left(I_{(000)T,t}^{(i_1i_2i_3)}-
I_{(000)T,t}^{(i_1i_2i_3)q}\right)^2\right\}+
$$

\vspace{1mm}
$$
+
{\sf M}\left\{\left(\frac{1}{2}\int\limits_t^T
(T-s)d{\bf f}_{s}^{(i_1)}-
\sum_{j_1,j_3=0}^{q}
C_{j_3j_3j_1}
\zeta_{j_1}^{(i_1)}\right)^2\right\}=
$$

\vspace{1mm}
$$
={\sf M}\left\{\left(I_{(000)T,t}^{(i_1i_2i_3)}-
I_{(000)T,t}^{(i_1i_2i_3)q}\right)^2\right\}+
$$

\begin{equation}
\label{dest27}
+\frac{1}{4}\int\limits_t^T
(T-s)^2 ds-
\sum_{j_1,j_3=0}^{q}
C_{j_3j_3j_1}\int\limits_t^T
(T-s)\phi_{j_1}(s)ds+
\sum_{j_1=0}^{q}\left(\sum_{j_3=0}^{q}
C_{j_3j_3j_1}\right)^2,
\end{equation}

\vspace{5mm}
\noindent
where $\phi_{j_1}(\tau)$ is the Legendre polynomial defined by (\ref{4009}).

Moreover,

\begin{equation}
\label{dest32}
\int\limits_t^T
(T-s)\phi_{j_1}(s)ds=\frac{(T-t)^{3/2}}{2}
\left\{
\begin{matrix}
1,\ & j_1=0\cr\cr
-1/\sqrt{3},\ & j_1=1\cr\cr
0,\ & j_1\ge 2
\end{matrix}
.\right.
\end{equation}

\vspace{5mm}

Combining (\ref{dest27})--(\ref{dest32}) and (\ref{39a}), we get

\vspace{1mm}
$$
{\sf M}\left\{\left(I_{(000)T,t}^{*(i_1i_2i_3)}-
I_{(000)T,t}^{*(i_1i_2i_3)q}\right)^2\right\}=
\frac{(T-t)^3}{4}
-\sum_{j_1,j_2,j_3=0}^q C_{j_3j_2j_1}^2-
\sum_{j_1,j_2,j_3=0}^q C_{j_2j_3j_1}C_{j_3j_2j_1}-
$$

\begin{equation}
\label{dest70}
-\frac{(T-t)^{3/2}}{2}
\sum_{j_3=0}^{q}
\left(C_{j_3j_3 0}-\frac{1}{\sqrt{3}}C_{j_3j_3 1}\right)+
\sum_{j_1=0}^{q}\left(\sum_{j_3=0}^{q}
C_{j_3j_3j_1}\right)^2,
\end{equation}

\vspace{5mm}
\noindent
where $i_1\ne i_2=i_3.$

Consider the case (\ref{dest7}). From (\ref{dest10}) we obtain

\vspace{1mm}
$$
{\sf M}\left\{\left(I_{(000)T,t}^{*(i_1i_2i_3)}-
I_{(000)T,t}^{*(i_1i_2i_3)q}\right)^2\right\}=
$$

\vspace{1mm}
$$
={\sf M}\left\{\left(I_{(000)T,t}^{(i_1i_2i_3)}-
I_{(000)T,t}^{(i_1i_2i_3)q}-
\sum_{j_1,j_2=0}^{q}
C_{j_1j_2j_1}
\zeta_{j_2}^{(i_2)}\right)^{2}\right\}=
$$

\vspace{1mm}
$$
={\sf M}\left\{\left(I_{(000)T,t}^{(i_1i_2i_3)}-
I_{(000)T,t}^{(i_1i_2i_3)q}\right)^{2}\right\}+
{\sf M}\left\{\left(\sum_{j_1,j_2=0}^{q}
C_{j_1j_2j_1}
\zeta_{j_2}^{(i_2)}\right)^{2}\right\}=
$$

\vspace{1mm}
\begin{equation}
\label{dest49}
={\sf M}\left\{\left(I_{(000)T,t}^{(i_1i_2i_3)}-
I_{(000)T,t}^{(i_1i_2i_3)q}\right)^{2}\right\}+
\sum_{j_2=0}^{q}
\left(\sum_{j_1=0}^{q}C_{j_1j_2j_1}\right)^2.
\end{equation}

\vspace{5mm}

Combining (\ref{dest49}) and (\ref{39b}), we have

\vspace{1mm}
$$
{\sf M}\left\{\left(I_{(000)T,t}^{*(i_1i_2i_3)}-
I_{(000)T,t}^{*(i_1i_2i_3)q}\right)^2\right\}=
\frac{(T-t)^3}{6}-
\sum_{j_1,j_2,j_3=0}^q C_{j_3j_2j_1}^2-
\sum_{j_1,j_2,j_3=0}^q C_{j_3j_2j_1}C_{j_1j_2j_3}+
$$
\begin{equation}
\label{dest60}
+\sum_{j_2=0}^{q}
\left(\sum_{j_1=0}^{q}C_{j_1j_2j_1}\right)^2,
\end{equation}

\vspace{5mm}
\noindent
where $i_1=i_3\ne i_2.$

Thus, the exact calculaton of the mean-square approximation error
for the iterated Stratonovich stochastic integral 
$I_{(000)T,t}^{*(i_1i_2i_3)}$ $(i_1,i_2,i_3=1,\ldots,m)$
is given by the formulas (\ref{39}),
(\ref{dest80}), (\ref{dest70}), and (\ref{dest60}).

Consider now the iterated Stratonovich stochastic integral
$I_{(0000)T,t}^{*(i_1i_2i_3i_4)}$ of multiplicity 4
$(i_1, i_2, i_3, i_4=$ $1,\ldots,m)$.
For $i_1=i_2=i_3=i_4$ we can use the formula (\ref{ud111ee}).
For the case of pairwise different 
$i_1, i_2, i_3, i_4$ 
we have the following relation

$$
I_{(0000)T,t}^{*(i_1 i_2 i_3 i_4)}=I_{(0000)T,t}^{(i_1 i_2 i_3 i_4)}\ \ \ \hbox{w.~p.~1}.
$$

\vspace{3mm}

Then in this case we can use the formulas (\ref{res100}) (for pairwise different
$i_1, i_2, i_3, i_4$) and (\ref{r7})
to approximate the stochastic integral
$I_{(0000)T,t}^{*(i_1i_1i_1i_1)}.$

Thus, it remains to consider the following 13 cases

\vspace{-2mm}
\begin{equation}
\label{casee1}
i_1=i_2\ne i_3, i_4;\ i_3\ne i_4,
\end{equation}
\begin{equation}
\label{casee2}
i_1=i_3\ne i_2, i_4;\ i_2\ne i_4,
\end{equation}
\begin{equation}
\label{casee3}
i_1=i_4\ne i_2, i_3;\ i_2\ne i_3,
\end{equation}
\begin{equation}
\label{casee4}
i_2=i_3\ne i_1, i_4;\ i_1\ne i_4,
\end{equation}
\begin{equation}
\label{casee5}
i_2=i_4\ne i_1, i_3;\ i_1\ne i_3,
\end{equation}
\begin{equation}
\label{casee6}
i_3=i_4\ne i_1, i_2;\ i_1\ne i_2,
\end{equation}
\begin{equation}
\label{casee7}
i_1=i_2=i_3\ne i_4,
\end{equation}
\begin{equation}
\label{casee8}
i_2=i_3=i_4\ne i_1,
\end{equation}
\begin{equation}
\label{casee9}
i_1=i_2=i_4\ne i_3,
\end{equation}
\begin{equation}
\label{casee10}
i_1=i_3=i_4\ne i_2,
\end{equation}
\begin{equation}
\label{casee11}
i_1=i_2\ne i_3=i_4,
\end{equation}
\begin{equation}
\label{casee12}
i_1=i_3\ne i_2=i_4,
\end{equation}
\begin{equation}
\label{casee13}
i_1=i_4\ne i_2=i_3.
\end{equation}

\vspace{3mm}

By analogy with (\ref{dest10}) and using the standard 
relation between Stratonovich and Ito stochastic integrals (\ref{k1001}), (\ref{k1000})
of multiplicity 4 as well as (\ref{res100}), we obtain

\vspace{1mm}
$$
{\sf M}\left\{\left(I_{(0000)T,t}^{*(i_1i_2i_3i_4)}-
I_{(0000)T,t}^{*(i_1i_2i_3i_4)q}\right)^2\right\}=
$$

\vspace{1mm}
$$
={\sf M}\left\{\left(I_{(0000)T,t}^{(i_1i_2i_3i_4)}+
\frac{1}{2}{\bf 1}_{\{i_1=i_2\ne 0\}}
\int\limits_t^T\int\limits_t^{t_4}\int\limits_t^{t_3}dt_1
d{\bf w}_{t_3}^{(i_3)}
d{\bf w}_{t_4}^{(i_4)}+\right.\right.
$$

\vspace{1mm}
$$
+\frac{1}{2}{\bf 1}_{\{i_2=i_3\ne 0\}}
\int\limits_t^T\int\limits_t^{t_4}\int\limits_t^{t_2}
d{\bf w}_{t_1}^{(i_1)}dt_2
d{\bf w}_{t_4}^{(i_4)}
+\frac{1}{2}{\bf 1}_{\{i_3=i_4\ne 0\}}
\int\limits_t^T\int\limits_t^{t_3}\int\limits_t^{t_2}
d{\bf w}_{t_1}^{(i_1)}
d{\bf w}_{t_2}^{(i_2)}dt_3+
$$

\vspace{1mm}
$$
+\frac{1}{4}{\bf 1}_{\{i_1=i_2\ne 0\}}
{\bf 1}_{\{i_3=i_4\ne 0\}}
\int\limits_t^T\int\limits_t^{t_2}dt_1dt_2-
I_{(0000)T,t}^{(i_1i_2i_3i_4)q}-
$$

\vspace{1mm}
$$
-{\bf 1}_{\{i_1=i_2\ne 0\}}\sum\limits_{j_4,j_3=0}^q
\sum\limits_{j_1=0}^q C_{j_4 j_3 j_1 j_1}\zeta_{j_3}^{(i_3)}\zeta_{j_4}^{(i_4)}-
{\bf 1}_{\{i_1=i_3\ne 0\}}\sum\limits_{j_4,j_2=0}^q
\sum\limits_{j_1=0}^q C_{j_4 j_1 j_2 j_1}\zeta_{j_2}^{(i_2)}\zeta_{j_4}^{(i_4)}-
$$

\vspace{1mm}
$$
-{\bf 1}_{\{i_1=i_4\ne 0\}}\sum\limits_{j_3,j_2=0}^q
\sum\limits_{j_1=0}^q C_{j_1 j_3 j_2 j_1}\zeta_{j_2}^{(i_2)}\zeta_{j_3}^{(i_3)}-
{\bf 1}_{\{i_2=i_3\ne 0\}}\sum\limits_{j_4,j_1=0}^q
\sum\limits_{j_2=0}^q C_{j_4 j_2 j_2 j_1}\zeta_{j_1}^{(i_1)}\zeta_{j_4}^{(i_4)}-
$$

\vspace{1mm}
$$
-{\bf 1}_{\{i_2=i_4\ne 0\}}\sum\limits_{j_3,j_1=0}^q
\sum\limits_{j_2=0}^q C_{j_2 j_3 j_2 j_1}\zeta_{j_1}^{(i_1)}\zeta_{j_3}^{(i_3)}-
{\bf 1}_{\{i_3=i_4\ne 0\}}\sum\limits_{j_2,j_1=0}^q
\sum\limits_{j_3=0}^q C_{j_3 j_3 j_2 j_1}\zeta_{j_1}^{(i_1)}\zeta_{j_2}^{(i_2)}+
$$

\vspace{1mm}
$$
+{\bf 1}_{\{i_1=i_2\ne 0\}}{\bf 1}_{\{i_3=i_4\ne 0\}}
\sum\limits_{j_3,j_1=0}^q C_{j_3 j_3 j_1 j_1}+
{\bf 1}_{\{i_1=i_3\ne 0\}}{\bf 1}_{\{i_2=i_4\ne 0\}}
\sum\limits_{j_2,j_1=0}^q C_{j_2 j_1 j_2 j_1}+
$$

\vspace{1mm}
\begin{equation}
\label{casee400}
\left.\left.+{\bf 1}_{\{i_1=i_4\ne 0\}}{\bf 1}_{\{i_2=i_3\ne 0\}}
\sum\limits_{j_2,j_1=0}^q C_{j_1 j_2 j_2 j_1}\right)^2\right\},
\end{equation}

\vspace{5mm}
\noindent 
where $I_{(0000)T,t}^{(i_1i_2i_3i_4)q}$ is defined by (\ref{res100}).

Consider the case (\ref{casee1}). From (\ref{casee400}) we get

\vspace{1mm}
$$
{\sf M}\left\{\left(I_{(0000)T,t}^{*(i_1i_2i_3i_4)}-
I_{(0000)T,t}^{*(i_1i_2i_3i_4)q}\right)^2\right\}=
$$

\vspace{1mm}
\begin{equation}
\label{casee401}
={\sf M}\left\{\left(I_{(0000)T,t}^{(i_1i_2i_3i_4)}
-I_{(0000)T,t}^{(i_1i_2i_3i_4)q}+
\frac{1}{2}\int\limits_t^T\int\limits_t^{t_4}\int\limits_t^{t_3}dt_1
d{\bf w}_{t_3}^{(i_3)}
d{\bf w}_{t_4}^{(i_4)}-
\sum\limits_{j_4,j_3=0}^q
\sum\limits_{j_1=0}^q C_{j_4 j_3 j_1 j_1}\zeta_{j_3}^{(i_3)}\zeta_{j_4}^{(i_4)}
\right)^2\right\}.
\end{equation}

\vspace{5mm}

Note that

\begin{equation}
\label{casee402}
\zeta_{j_3}^{(i_3)}\zeta_{j_4}^{(i_4)}=
\int\limits_t^T \phi_{j_4}(t_4)\int\limits_t^{t_4}\phi_{j_3}(t_3)
d{\bf w}_{t_3}^{(i_3)}
d{\bf w}_{t_4}^{(i_4)}+
\int\limits_t^T \phi_{j_3}(t_3)\int\limits_t^{t_3}\phi_{j_4}(t_4)
d{\bf w}_{t_4}^{(i_4)}
d{\bf w}_{t_3}^{(i_3)}
\end{equation}

\vspace{5mm}
\noindent
w.~p.~1, where $i_3\ne i_4.$

According to the results of Sect.~3 in \cite{15b}
(also see Sect.~1.2.2 in \cite{10a}-\cite{10axx1}), the quantity

$$
I_{(0000)T,t}^{(i_1i_2i_3i_4)}-
I_{(0000)T,t}^{(i_1i_2i_3i_4)q}
$$

\vspace{3mm}
\noindent
includes only iterated Ito stochastic integrals
of multiplicity 4. At the same time (see (\ref{casee402})), the quantity

$$
\frac{1}{2}\int\limits_t^T\int\limits_t^{t_4}\int\limits_t^{t_3}dt_1
d{\bf w}_{t_3}^{(i_3)}
d{\bf w}_{t_4}^{(i_4)}-
\sum\limits_{j_4,j_3=0}^q
\sum\limits_{j_1=0}^p C_{j_4 j_3 j_1 j_1}\zeta_{j_3}^{(i_3)}\zeta_{j_4}^{(i_4)}
$$

\vspace{3mm}
\noindent
contains only iterated Ito stochastic integrals
of multiplicity 2. This means that from (\ref{casee401}) we get

\vspace{1mm}
$$
{\sf M}\left\{\left(I_{(0000)T,t}^{*(i_1i_2i_3i_4)}-
I_{(0000)T,t}^{*(i_1i_2i_3i_4)q}\right)^2\right\}=
{\sf M}\left\{\left(I_{(0000)T,t}^{(i_1i_2i_3i_4)}
-I_{(0000)T,t}^{(i_1i_2i_3i_4)q}\right)^2\right\}
+
$$

\vspace{2mm}
$$
+{\sf M}\left\{\left(\frac{1}{2}\int\limits_t^T\int\limits_t^{t_4}(t_3-t)
d{\bf w}_{t_3}^{(i_3)}
d{\bf w}_{t_4}^{(i_4)}-
\sum\limits_{j_4,j_3=0}^q
\sum\limits_{j_1=0}^q C_{j_4 j_3 j_1 j_1}\zeta_{j_3}^{(i_3)}\zeta_{j_4}^{(i_4)}
\right)^2\right\}=
$$

\vspace{2mm}
$$
=
{\sf M}\left\{\left(I_{(0000)T,t}^{(i_1i_2i_3i_4)}
-I_{(0000)T,t}^{(i_1i_2i_3i_4)q}\right)^2\right\}
+\frac{1}{4}\int\limits_t^T\int\limits_t^{t_4}(t_3-t)^2
dt_3dt_4+
$$

\vspace{2mm}
$$
+\sum\limits_{j_4,j_3=0}^q
\left(\sum\limits_{j_1=0}^q C_{j_4 j_3 j_1 j_1}\right)^2- 
\sum\limits_{j_4,j_3=0}^q
\sum\limits_{j_1=0}^q C_{j_4 j_3 j_1 j_1}
\int\limits_t^T \phi_{j_4}(t_4)\int\limits_t^{t_4}\phi_{j_3}(t_3)(t_3-t)
dt_3 dt_4=
$$

\vspace{2mm}
$$
=
{\sf M}\left\{\left(I_{(0000)T,t}^{(i_1i_2i_3i_4)}
-I_{(0000)T,t}^{(i_1i_2i_3i_4)q}\right)^2\right\}
+\frac{(T-t)^4}{48}+
\sum\limits_{j_4,j_3=0}^q
\left(\sum\limits_{j_1=0}^q C_{j_4 j_3 j_1 j_1}\right)^2+
$$

\vspace{2mm}
\begin{equation}
\label{casee500}
+
\sum\limits_{j_4,j_3=0}^q
\sum\limits_{j_1=0}^q C_{j_4 j_3 j_1 j_1}C_{j_4 j_3}^{10},
\end{equation}

\vspace{4mm}
\noindent
where 
\begin{equation}
\label{casee700}
C_{j_4 j_3}^{10}=
\int\limits_t^T \phi_{j_4}(t_4)\int\limits_t^{t_4}\phi_{j_3}(t_3)(t-t_3)
dt_3 dt_4.
\end{equation}

\vspace{5mm}

Using (\ref{usl1}) and (\ref{casee500}), we finally obtain

\vspace{1mm}
$$
{\sf M}\left\{\left(I_{(0000)T,t}^{*(i_1i_2i_3i_4)}-
I_{(0000)T,t}^{*(i_1i_2i_3i_4)q}\right)^2\right\}=
\frac{(T-t)^4}{16}-
\sum_{j_1,j_2,j_3,j_4=0}^{q}
C_{j_4j_3j_2j_1}\Biggl(\sum\limits_{(j_1,j_2)}
C_{j_4j_3j_2j_1}\Biggr)
+
$$
\begin{equation}
\label{casee501}
+\sum\limits_{j_4,j_3=0}^q
\left(\sum\limits_{j_1=0}^q C_{j_4 j_3 j_1 j_1}\right)^2+
\sum\limits_{j_4,j_3=0}^q
\sum\limits_{j_1=0}^q C_{j_4 j_3 j_1 j_1}C_{j_4 j_3}^{10},
\end{equation}

\vspace{5mm}
\noindent
where $i_1=i_2\ne i_3, i_4;\ i_3\ne i_4.$

Consider the cases (\ref{casee2}), (\ref{casee3}) by analogy
with the case (\ref{casee1}) using (\ref{usl2}), (\ref{usl3}).
We have

\vspace{1mm}
$$
{\sf M}\left\{\left(I_{(0000)T,t}^{*(i_1i_2i_3i_4)}-
I_{(0000)T,t}^{*(i_1i_2i_3i_4)q}\right)^2\right\}=
\frac{(T-t)^4}{24}-
\sum_{j_1,j_2,j_3,j_4=0}^{q}
C_{j_4j_3j_2j_1}\Biggl(\sum\limits_{(j_1,j_3)}
C_{j_4j_3j_2j_1}\Biggr)
+
$$

\vspace{1mm}
$$
+
\sum\limits_{j_4,j_2=0}^q
\left(\sum\limits_{j_1=0}^q C_{j_4 j_1 j_2 j_1}\right)^2,
$$

\vspace{5mm}
\noindent
where $i_1=i_3\ne i_2, i_4$ and $i_2\ne i_4;$

\vspace{1mm}
$$
{\sf M}\left\{\left(I_{(0000)T,t}^{*(i_1i_2i_3i_4)}-
I_{(0000)T,t}^{*(i_1i_2i_3i_4)q}\right)^2\right\}=
\frac{(T-t)^4}{24}-
\sum_{j_1,j_2,j_3,j_4=0}^{q}
C_{j_4j_3j_2j_1}\Biggl(\sum\limits_{(j_1,j_4)}
C_{j_4j_3j_2j_1}\Biggr)
+
$$

\vspace{1mm}
$$
+\sum\limits_{j_3,j_2=0}^q
\left(\sum\limits_{j_1=0}^q C_{j_1 j_3 j_2 j_1}\right)^2,
$$

\vspace{5mm}
\noindent
where $i_1=i_4\ne i_2, i_3$ and $i_2\ne i_3.$

Consider the case (\ref{casee4}) by analogy
with the case (\ref{casee1}). 
We have

\vspace{1mm}
$$
{\sf M}\left\{\left(I_{(0000)T,t}^{*(i_1i_2i_3i_4)}-
I_{(0000)T,t}^{*(i_1i_2i_3i_4)q}\right)^2\right\}=
{\sf M}\left\{\left(I_{(0000)T,t}^{(i_1i_2i_3i_4)}
-I_{(0000)T,t}^{(i_1i_2i_3i_4)q}\right)^2\right\}
+
$$

\vspace{2mm}
$$
+{\sf M}\left\{\left(\frac{1}{2}\int\limits_t^T\int\limits_t^{t_4}\int\limits_t^{t_2}
d{\bf w}_{t_1}^{(i_1)}dt_2
d{\bf w}_{t_4}^{(i_4)}-
\sum\limits_{j_4,j_1=0}^q
\sum\limits_{j_2=0}^q C_{j_4 j_2 j_2 j_1}\zeta_{j_1}^{(i_1)}\zeta_{j_4}^{(i_4)}
\right)^2\right\}=
$$

\vspace{2mm}
$$
=
{\sf M}\left\{\left(I_{(0000)T,t}^{(i_1i_2i_3i_4)}
-I_{(0000)T,t}^{(i_1i_2i_3i_4)q}\right)^2\right\}
+
$$

\vspace{2mm}
$$
+{\sf M}\left\{\left(\frac{1}{2}\int\limits_t^T\int\limits_t^{t_4}(t_4-t_1)
d{\bf w}_{t_1}^{(i_1)}
d{\bf w}_{t_4}^{(i_4)}-
\sum\limits_{j_4,j_1=0}^q
\sum\limits_{j_2=0}^q C_{j_4 j_2 j_2 j_1}\zeta_{j_1}^{(i_1)}\zeta_{j_4}^{(i_4)}
\right)^2\right\}=
$$

\vspace{2mm}
$$
=
{\sf M}\left\{\left(I_{(0000)T,t}^{(i_1i_2i_3i_4)}
-I_{(0000)T,t}^{(i_1i_2i_3i_4)q}\right)^2\right\}
+\frac{(T-t)^4}{48}+
\sum\limits_{j_4,j_1=0}^q
\left(\sum\limits_{j_2=0}^q C_{j_4 j_2 j_2 j_1}\right)^2-
$$

\vspace{2mm}
$$
-
\sum\limits_{j_4,j_1=0}^q
\sum\limits_{j_2=0}^q C_{j_4 j_2 j_2 j_1}
\int\limits_t^T \phi_{j_4}(t_4)\int\limits_t^{t_4}\phi_{j_1}(t_1)(t_4-t_1)
dt_3 dt_4.
$$

\vspace{6mm}

Then applying (\ref{usl4}), we obtain

$$
{\sf M}\left\{\left(I_{(0000)T,t}^{*(i_1i_2i_3i_4)}-
I_{(0000)T,t}^{*(i_1i_2i_3i_4)q}\right)^2\right\}=\frac{(T-t)^4}{16}-
\sum_{j_1,j_2,j_3,j_4=0}^{q}
C_{j_4j_3j_2j_1}\Biggl(\sum\limits_{(j_2,j_3)}
C_{j_4j_3j_2j_1}\Biggr)+
$$

\vspace{2mm}
$$
+
\sum\limits_{j_4,j_1=0}^q
\left(\sum\limits_{j_2=0}^q C_{j_4 j_2 j_2 j_1}\right)^2-
\sum\limits_{j_4,j_1=0}^q
\sum\limits_{j_2=0}^q C_{j_4 j_2 j_2 j_1}
\left(C_{j_4 j_1}^{10}-C_{j_4 j_1}^{01}\right),
$$

\vspace{7mm}
\noindent
where $i_2=i_3\ne i_1, i_4$ and $i_1\ne i_4;$
$C_{j_4 j_1}^{10}$ is defined by (\ref{casee700}) and

\begin{equation}
\label{casee700a}
C_{j_4 j_1}^{01}=
\int\limits_t^T \phi_{j_4}(t_4)(t-t_4)\int\limits_t^{t_4}\phi_{j_1}(t_1)
dt_1 dt_4.
\end{equation}

\vspace{5mm}

For the case (\ref{casee5}) by analogy
with the case (\ref{casee1}) and using (\ref{usl5}), we get

\vspace{1mm}
$$
{\sf M}\left\{\left(I_{(0000)T,t}^{*(i_1i_2i_3i_4)}-
I_{(0000)T,t}^{*(i_1i_2i_3i_4)q}\right)^2\right\}=
\frac{(T-t)^4}{24}-
\sum_{j_1,j_2,j_3,j_4=0}^{q}
C_{j_4j_3j_2j_1}\Biggl(\sum\limits_{(j_2,j_4)}
C_{j_4j_3j_2j_1}\Biggr)
+
$$

\vspace{1mm}
$$
+
\sum\limits_{j_3,j_1=0}^q
\left(\sum\limits_{j_2=0}^q C_{j_2 j_3 j_2 j_1}\right)^2,
$$

\vspace{5mm}
\noindent
where $i_2=i_4\ne i_1, i_3$ and $i_1\ne i_3.$

Consider the case (\ref{casee6}) by analogy
with the case (\ref{casee1}). 
Note that \cite{10a}-\cite{10axx1} (see Example~3.1 in Sect.~3.6)

\begin{equation}
\label{yyee22}
\int\limits_t^T\int\limits_t^{t_3}\int\limits_t^{t_2}
d{\bf w}_{t_1}^{(i_1)}
d{\bf w}_{t_2}^{(i_2)}dt_3
=\int\limits_t^T(T-t_2)\int\limits_t^{t_2}
d{\bf w}_{t_1}^{(i_1)}
d{\bf w}_{t_2}^{(i_2)}\ \ \ \hbox{w.~p.~1.}
\end{equation}

\vspace{5mm}

Using (\ref{yyee22}), we obtain

\vspace{1mm}
$$
{\sf M}\left\{\left(I_{(0000)T,t}^{*(i_1i_2i_3i_4)}-
I_{(0000)T,t}^{*(i_1i_2i_3i_4)q}\right)^2\right\}=
{\sf M}\left\{\left(I_{(0000)T,t}^{(i_1i_2i_3i_4)}
-I_{(0000)T,t}^{(i_1i_2i_3i_4)q}\right)^2\right\}
+
$$

\vspace{2mm}
$$
+{\sf M}\left\{\left(\frac{1}{2}\int\limits_t^T(T-t_2)\int\limits_t^{t_2}
d{\bf w}_{t_1}^{(i_1)}
d{\bf w}_{t_2}^{(i_2)}-
\sum\limits_{j_2,j_1=0}^q
\sum\limits_{j_3=0}^q C_{j_3 j_3 j_2 j_1}\zeta_{j_1}^{(i_1)}\zeta_{j_2}^{(i_2)}
\right)^2\right\}=
$$

\vspace{2mm}
$$
=
{\sf M}\left\{\left(I_{(0000)T,t}^{(i_1i_2i_3i_4)}
-I_{(0000)T,t}^{(i_1i_2i_3i_4)q}\right)^2\right\}
+\frac{(T-t)^4}{48}+
\sum\limits_{j_2,j_1=0}^q
\left(\sum\limits_{j_3=0}^q C_{j_3 j_3 j_2 j_1}\right)^2-
$$

\vspace{2mm}
$$
-
\sum\limits_{j_2,j_1=0}^q
\sum\limits_{j_3=0}^q C_{j_3 j_3 j_2 j_1}
\int\limits_t^T (T-t_2)\phi_{j_2}(t_2)\int\limits_t^{t_2}\phi_{j_1}(t_1)
dt_1 dt_2.
$$

\vspace{5mm}

Then applying (\ref{usl6}), we get

$$
{\sf M}\left\{\left(I_{(0000)T,t}^{*(i_1i_2i_3i_4)}-
I_{(0000)T,t}^{*(i_1i_2i_3i_4)q}\right)^2\right\}=\frac{(T-t)^4}{16}-
\sum_{j_1,j_2,j_3,j_4=0}^{q}
C_{j_4j_3j_2j_1}\Biggl(\sum\limits_{(j_3,j_4)}
C_{j_4j_3j_2j_1}\Biggr)+
$$

\vspace{1mm}
$$
+
\sum\limits_{j_2,j_1=0}^q
\left(\sum\limits_{j_3=0}^q C_{j_3 j_3 j_2 j_1}\right)^2-
\sum\limits_{j_2,j_1=0}^q
\sum\limits_{j_3=0}^q C_{j_3 j_3 j_2 j_1}
\left((T-t)C_{j_2 j_1}+C_{j_2 j_1}^{01}\right),
$$

\vspace{6mm}
\noindent
where $i_3=i_4\ne i_1, i_2$ and $i_1\ne i_2;$
$C_{j_2 j_1}^{01}$ is defined by (\ref{casee700a}) and

$$
C_{j_2 j_1}=
\int\limits_t^T \phi_{j_2}(t_2)\int\limits_t^{t_2}\phi_{j_1}(t_1)
dt_1 dt_2.
$$

\vspace{5mm}

Consider the case (\ref{casee7}). 
From (\ref{casee400}) we have 

\vspace{1mm}
$$
{\sf M}\left\{\left(I_{(0000)T,t}^{*(i_1i_1i_1i_4)}-
I_{(0000)T,t}^{*(i_1i_1i_1i_4)q}\right)^2\right\}=
{\sf M}\left\{\left(I_{(0000)T,t}^{(i_1i_1i_1i_4)}+
\frac{1}{2}
\int\limits_t^T\int\limits_t^{t_4}\int\limits_t^{t_3}dt_1
d{\bf w}_{t_3}^{(i_1)}
d{\bf w}_{t_4}^{(i_4)}+\right.\right.
$$

\vspace{2mm}
$$
+\frac{1}{2}
\int\limits_t^T\int\limits_t^{t_4}\int\limits_t^{t_2}
d{\bf w}_{t_1}^{(i_1)}dt_2
d{\bf w}_{t_4}^{(i_4)}-
I_{(0000)T,t}^{(i_1i_1i_1i_4)q}-
\sum\limits_{j_4,j_3=0}^q
\sum\limits_{j_1=0}^q C_{j_4 j_3 j_1 j_1}\zeta_{j_3}^{(i_1)}\zeta_{j_4}^{(i_4)}-
$$

\vspace{2mm}
\begin{equation}
\label{casee800}
\left.\left.-
\sum\limits_{j_4,j_2=0}^q
\sum\limits_{j_1=0}^q C_{j_4 j_1 j_2 j_1}\zeta_{j_2}^{(i_1)}\zeta_{j_4}^{(i_4)}-
\sum\limits_{j_4,j_1=0}^q
\sum\limits_{j_2=0}^q C_{j_4 j_2 j_2 j_1}\zeta_{j_1}^{(i_1)}\zeta_{j_4}^{(i_4)}
\right)^2\right\}.
\end{equation}

\vspace{5mm}

Furthermore,

\vspace{-1mm}
$$
\int\limits_t^T\int\limits_t^{t_4}\int\limits_t^{t_3}dt_1
d{\bf w}_{t_3}^{(i_1)}
d{\bf w}_{t_4}^{(i_4)}+
\int\limits_t^T\int\limits_t^{t_4}\int\limits_t^{t_2}
d{\bf w}_{t_1}^{(i_1)}dt_2
d{\bf w}_{t_4}^{(i_4)}=
$$

\vspace{2mm}
$$
=\int\limits_t^T\int\limits_t^{t_4}(t_1-t)
d{\bf w}_{t_1}^{(i_1)}
d{\bf w}_{t_4}^{(i_4)}+
\int\limits_t^T\int\limits_t^{t_4}
(t_4-t_1)d{\bf w}_{t_1}^{(i_1)}
d{\bf w}_{t_4}^{(i_4)}=
$$

\vspace{1mm}
\begin{equation}
\label{casee801}
=\int\limits_t^T(t_4-t)\int\limits_t^{t_4}
d{\bf w}_{t_1}^{(i_1)}
d{\bf w}_{t_4}^{(i_4)}\ \ \ \hbox{w.~p.~1.}
\end{equation}

\vspace{5mm}

From (\ref{casee800}) and (\ref{casee801}) we obtain

\vspace{1mm}
$$
{\sf M}\left\{\left(I_{(0000)T,t}^{*(i_1i_1i_1i_4)}-
I_{(0000)T,t}^{*(i_1i_1i_1i_4)q}\right)^2\right\}=
{\sf M}\left\{\left(I_{(0000)T,t}^{(i_1i_1i_1i_4)}-
I_{(0000)T,t}^{(i_1i_1i_1i_4)q}\right)^2\right\}+
$$

\vspace{2mm}
$$
+{\sf M}\left\{\left(
\frac{1}{2}\int\limits_t^T(t_4-t)\int\limits_t^{t_4}
d{\bf w}_{t_1}^{(i_1)}
d{\bf w}_{t_4}^{(i_4)}-
\sum\limits_{j_4,j_1=0}^q
\sum\limits_{j_2=0}^q \left(C_{j_4 j_1 j_2 j_2}+C_{j_4 j_2 j_1 j_2}+
C_{j_4 j_2 j_2 j_1}\right)\zeta_{j_1}^{(i_1)}\zeta_{j_4}^{(i_4)}
\right)^2\right\}=
$$

\vspace{2mm}
$$
=
{\sf M}\left\{\left(I_{(0000)T,t}^{(i_1i_1i_1i_4)}-
I_{(0000)T,t}^{(i_1i_1i_1i_4)q}\right)^2\right\}+\frac{(T-t)^4}{16}+
$$

\vspace{2mm}
$$
+
\sum\limits_{j_4,j_1=0}^q
\left(\sum\limits_{j_2=0}^q \left(C_{j_4 j_1 j_2 j_2}+C_{j_4 j_2 j_1 j_2}+
C_{j_4 j_2 j_2 j_1}\right)\right)^2-
$$

\vspace{2mm}
\begin{equation}
\label{casee805}
-\sum\limits_{j_4,j_1=0}^q
\sum\limits_{j_2=0}^q \left(C_{j_4 j_1 j_2 j_2}+C_{j_4 j_2 j_1 j_2}+
C_{j_4 j_2 j_2 j_1}\right)
\int\limits_t^T (t_4-t)\phi_{j_4}(t_4)\int\limits_t^{t_4}\phi_{j_1}(t_1)
dt_1 dt_4.
\end{equation}

\vspace{5mm}

Using (\ref{usl7}) and (\ref{casee805}), we finally get 

\vspace{1mm}
$$
{\sf M}\left\{\left(I_{(0000)T,t}^{*(i_1i_2i_3i_4)}-
I_{(0000)T,t}^{*(i_1i_2i_3i_4)q}\right)^2\right\}=
\frac{5(T-t)^4}{48}-
\sum_{j_1,j_2,j_3,j_4=0}^{q}
C_{j_4j_3j_2j_1}\Biggl(\sum\limits_{(j_1,j_2,j_3)}
C_{j_4j_3j_2j_1}\Biggr)+
$$

\vspace{1mm}
$$
+
\sum\limits_{j_4,j_1=0}^q
\left(\sum\limits_{j_2=0}^q \left(C_{j_4 j_1 j_2 j_2}+C_{j_4 j_2 j_1 j_2}+
C_{j_4 j_2 j_2 j_1}\right)\right)^2+
$$

\vspace{1mm}
$$
+\sum\limits_{j_4,j_1=0}^q
\sum\limits_{j_2=0}^q \left(C_{j_4 j_1 j_2 j_2}+C_{j_4 j_2 j_1 j_2}+
C_{j_4 j_2 j_2 j_1}\right)C_{j_4j_2}^{01},
$$

\vspace{5mm}
\noindent
where $i_1=i_2=i_3\ne i_4.$

Consider the case (\ref{casee8}). 
From (\ref{casee400}) we have 

\vspace{1mm}
$$
{\sf M}\left\{\left(I_{(0000)T,t}^{*(i_1i_2i_2i_2)}-
I_{(0000)T,t}^{*(i_1i_2i_2i_2)q}\right)^2\right\}=
{\sf M}\left\{\left(I_{(0000)T,t}^{(i_1i_2i_2i_2)}+
\frac{1}{2}
\int\limits_t^T\int\limits_t^{t_4}\int\limits_t^{t_2}
d{\bf w}_{t_1}^{(i_1)}dt_2
d{\bf w}_{t_4}^{(i_2)}+\right.\right.
$$

\vspace{2mm}
$$
+\frac{1}{2}
\int\limits_t^T\int\limits_t^{t_3}\int\limits_t^{t_2}
d{\bf w}_{t_1}^{(i_1)}
d{\bf w}_{t_2}^{(i_2)}dt_3-
I_{(0000)T,t}^{(i_1i_2i_2i_2)q}-
\sum\limits_{j_4,j_1=0}^q
\sum\limits_{j_2=0}^q C_{j_4 j_2 j_2 j_1}\zeta_{j_1}^{(i_1)}\zeta_{j_4}^{(i_2)}-
$$

\vspace{2mm}
\begin{equation}
\label{casee807}
\left.\left.-
\sum\limits_{j_3,j_1=0}^q
\sum\limits_{j_2=0}^q C_{j_2 j_3 j_2 j_1}\zeta_{j_1}^{(i_1)}\zeta_{j_3}^{(i_2)}-
\sum\limits_{j_2,j_1=0}^q
\sum\limits_{j_3=0}^q C_{j_3 j_3 j_2 j_1}\zeta_{j_1}^{(i_1)}\zeta_{j_2}^{(i_2)}
\right)^2\right\}.
\end{equation}

\vspace{5mm}

Moreover,

\vspace{-2mm}
$$
\int\limits_t^T\int\limits_t^{t_4}\int\limits_t^{t_2}
d{\bf w}_{t_1}^{(i_1)}dt_2
d{\bf w}_{t_4}^{(i_2)}+
\int\limits_t^T\int\limits_t^{t_3}\int\limits_t^{t_2}
d{\bf w}_{t_1}^{(i_1)}
d{\bf w}_{t_2}^{(i_2)}dt_3=
$$

\vspace{1mm}
$$
=\int\limits_t^T\int\limits_t^{t_4}(t_4-t_1)
d{\bf w}_{t_1}^{(i_1)}
d{\bf w}_{t_4}^{(i_2)}+
\int\limits_t^T\int\limits_t^{t_4}
(T-t_4)d{\bf w}_{t_1}^{(i_1)}
d{\bf w}_{t_4}^{(i_2)}=
$$

\vspace{1mm}
\begin{equation}
\label{casee808}
=\int\limits_t^T\int\limits_t^{t_4}
(T-t_1)d{\bf w}_{t_1}^{(i_1)}
d{\bf w}_{t_4}^{(i_2)}\ \ \ \hbox{w.~p.~1.}
\end{equation}

\vspace{5mm}

From (\ref{casee807}) and (\ref{casee808}) we get

\vspace{1mm}
$$
{\sf M}\left\{\left(I_{(0000)T,t}^{*(i_1i_2i_2i_2)}-
I_{(0000)T,t}^{*(i_1i_2i_2i_2)q}\right)^2\right\}=
{\sf M}\left\{\left(I_{(0000)T,t}^{(i_1i_2i_2i_2)}-
I_{(0000)T,t}^{(i_1i_2i_2i_2)q}\right)^2\right\}+
$$

\vspace{2mm}
$$
+{\sf M}\left\{\left(
\frac{1}{2}\int\limits_t^T\int\limits_t^{t_4}(T-t_1)
d{\bf w}_{t_1}^{(i_1)}
d{\bf w}_{t_4}^{(i_2)}-
\sum\limits_{j_4,j_1=0}^q
\sum\limits_{j_2=0}^q \left(C_{j_4 j_2 j_2 j_1}+C_{j_2 j_4 j_2 j_1}+
C_{j_2 j_2 j_4 j_1}\right)\zeta_{j_1}^{(i_1)}\zeta_{j_4}^{(i_2)}
\right)^2\right\}=
$$

\vspace{2mm}
$$
=
{\sf M}\left\{\left(I_{(0000)T,t}^{(i_1i_2i_2i_2)}-
I_{(0000)T,t}^{(i_1i_2i_2i_2)q}\right)^2\right\}+\frac{(T-t)^4}{16}+
$$

\vspace{2mm}
$$
+
\sum\limits_{j_4,j_1=0}^q
\left(\sum\limits_{j_2=0}^q \left(C_{j_4 j_2 j_2 j_1}+C_{j_2 j_4 j_2 j_1}+
C_{j_2 j_2 j_4 j_1}\right)\right)^2-
$$

\vspace{2mm}
\begin{equation}
\label{casee809}
-\sum\limits_{j_4,j_1=0}^q
\sum\limits_{j_2=0}^q \left(C_{j_4 j_2 j_2 j_1}+C_{j_2 j_4 j_2 j_1}+
C_{j_2 j_2 j_4 j_1}\right)
\int\limits_t^T \phi_{j_4}(t_4)\int\limits_t^{t_4}(T-t_1)\phi_{j_1}(t_1)
dt_1 dt_4.
\end{equation}

\vspace{5mm}

Applying (\ref{usl8}) and (\ref{casee809}), we finally obtain

\vspace{1mm}
$$
{\sf M}\left\{\left(I_{(0000)T,t}^{*(i_1i_2i_3i_4)}-
I_{(0000)T,t}^{*(i_1i_2i_3i_4)q}\right)^2\right\}=
\frac{5(T-t)^4}{48}-
\sum_{j_1,j_2,j_3,j_4=0}^{q}
C_{j_4j_3j_2j_1}\Biggl(\sum\limits_{(j_2,j_3,j_4)}
C_{j_4j_3j_2j_1}\Biggr)+
$$

\vspace{1mm}
$$
+
\sum\limits_{j_4,j_1=0}^q
\left(\sum\limits_{j_2=0}^q \left(C_{j_4 j_2 j_2 j_1}+C_{j_2 j_4 j_2 j_1}+
C_{j_2 j_2 j_4 j_1}\right)\right)^2-
$$

\vspace{1mm}
$$
-\sum\limits_{j_4,j_1=0}^q
\sum\limits_{j_2=0}^q \left(C_{j_4 j_2 j_2 j_1}+C_{j_2 j_4 j_2 j_1}+
C_{j_2 j_2 j_4 j_1}\right)\left((T-t)C_{j_4 j_1}+C_{j_4 j_1}^{10}\right),
$$

\vspace{5mm}
\noindent
where $i_2=i_3=i_4\ne i_1.$

For the cases (\ref{casee9}), (\ref{casee10}) by analogy with the case (\ref{casee8}) and using
(\ref{usl9}), (\ref{usl10}), we obtain

\vspace{1mm}
$$
{\sf M}\left\{\left(I_{(0000)T,t}^{*(i_1i_2i_3i_4)}-
I_{(0000)T,t}^{*(i_1i_2i_3i_4)q}\right)^2\right\}=
\frac{(T-t)^4}{16}-
\sum_{j_1,j_2,j_3,j_4=0}^{q}
C_{j_4j_3j_2j_1}\Biggl(\sum\limits_{(j_1,j_2,j_4)}
C_{j_4j_3j_2j_1}\Biggr)+
$$

\vspace{1mm}
$$
+
\sum\limits_{j_4,j_3=0}^q
\left(\sum\limits_{j_1=0}^q \left(C_{j_4 j_3 j_1 j_1}+C_{j_1 j_3 j_4 j_1}+
C_{j_1 j_3 j_1 j_4}\right)\right)^2+
$$

\vspace{1mm}
$$
+\sum\limits_{j_4,j_3=0}^q
\sum\limits_{j_1=0}^q \left(C_{j_4 j_3 j_1 j_1}+C_{j_1 j_3 j_4 j_1}+
C_{j_1 j_3 j_1 j_4}\right)C_{j_4 j_3}^{10},
$$

\vspace{5mm}
\noindent
where $i_1=i_2=i_4\ne i_3;$

\vspace{1mm}
$$
{\sf M}\left\{\left(I_{(0000)T,t}^{*(i_1i_2i_3i_4)}-
I_{(0000)T,t}^{*(i_1i_2i_3i_4)q}\right)^2\right\}=
\frac{(T-t)^4}{16}-
\sum_{j_1,j_2,j_3,j_4=0}^{q}
C_{j_4j_3j_2j_1}\Biggl(\sum\limits_{(j_1,j_3,j_4)}
C_{j_4j_3j_2j_1}\Biggr)+
$$

\vspace{1mm}
$$
+
\sum\limits_{j_4,j_2=0}^q
\left(\sum\limits_{j_1=0}^q \left(C_{j_4 j_1 j_2 j_1}+C_{j_1 j_4 j_2 j_1}+
C_{j_1 j_1 j_2 j_4}\right)\right)^2-
$$

\vspace{1mm}
$$
-\sum\limits_{j_4,j_2=0}^q
\sum\limits_{j_1=0}^q \left(C_{j_4 j_1 j_2 j_1}+C_{j_1 j_4 j_2 j_1}+
C_{j_1 j_1 j_2 j_4}\right)\left((T-t)C_{j_2 j_3}+C_{j_2 j_3}^{01}\right),
$$

\vspace{5mm}
\noindent
where $i_1=i_3=i_4\ne i_2.$

Let us consider the case (\ref{casee11}). Using (\ref{casee400}), we have

\vspace{1mm}
$$
{\sf M}\left\{\left(I_{(0000)T,t}^{*(i_1i_1i_3i_3)}-
I_{(0000)T,t}^{*(i_1i_1i_3i_3)q}\right)^2\right\}=
{\sf M}\left\{\left(I_{(0000)T,t}^{(i_1i_1i_3i_3)}+
\frac{1}{2}
\int\limits_t^T\int\limits_t^{t_4}(t_3-t)
d{\bf w}_{t_3}^{(i_3)}
d{\bf w}_{t_4}^{(i_3)}+\right.\right.
$$

\vspace{2mm}
$$
+\frac{1}{2}
\int\limits_t^T\int\limits_t^{t_3}\int\limits_t^{t_2}
d{\bf w}_{t_1}^{(i_1)}
d{\bf w}_{t_2}^{(i_1)}dt_3+ \frac{(T-t)^2}{8}-
I_{(0000)T,t}^{(i_1i_1i_3i_3)q}-
$$

\vspace{2mm}
$$
-\sum\limits_{j_4,j_3=0}^q
\sum\limits_{j_1=0}^q C_{j_4 j_3 j_1 j_1}\zeta_{j_3}^{(i_3)}\zeta_{j_4}^{(i_3)}-
\sum\limits_{j_2,j_1=0}^q
\sum\limits_{j_3=0}^q C_{j_3 j_3 j_2 j_1}\zeta_{j_1}^{(i_1)}\zeta_{j_2}^{(i_1)}+
\left.\left.
\sum\limits_{j_3,j_1=0}^q C_{j_3 j_3 j_1 j_1}\right)^2\right\}=
$$

\vspace{3mm}
$$
=
{\sf M}\left\{\Biggl(I_{(0000)T,t}^{(i_1i_1i_3i_3)}-
I_{(0000)T,t}^{(i_1i_1i_3i_3)q}+\Biggr.\right.
$$

\vspace{2mm}
$$
+
\frac{1}{2}
\int\limits_t^T\int\limits_t^{t_4}(t_3-t)
d{\bf w}_{t_3}^{(i_3)}
d{\bf w}_{t_4}^{(i_3)}-
\sum\limits_{j_4,j_3=0}^q
\sum\limits_{j_1=0}^q C_{j_4 j_3 j_1 j_1}\left(\zeta_{j_3}^{(i_3)}\zeta_{j_4}^{(i_3)}-
{\bf 1}_{\{j_3=j_4\}}\right)+
$$

\vspace{2mm}
$$
+\frac{1}{2}
\int\limits_t^T\int\limits_t^{t_3}\int\limits_t^{t_2}
d{\bf w}_{t_1}^{(i_1)}
d{\bf w}_{t_2}^{(i_1)}dt_3-
\sum\limits_{j_2,j_1=0}^q
\sum\limits_{j_3=0}^q C_{j_3 j_3 j_2 j_1}\left(\zeta_{j_1}^{(i_1)}\zeta_{j_2}^{(i_1)}-
{\bf 1}_{\{j_1=j_2\}}\right)+
$$

\vspace{2mm}
\begin{equation}
\label{caseu1}
\left.\Biggl.+\frac{(T-t)^2}{8}-
\sum\limits_{j_3,j_1=0}^q C_{j_3 j_3 j_1 j_1}\Biggr)^2\right\}.
\end{equation}

\vspace{5mm}

Note that 

\vspace{-1mm}
\begin{equation}
\label{caseu2}
\zeta_{j_3}^{(i_3)}\zeta_{j_4}^{(i_3)}-
{\bf 1}_{\{j_3=j_4\}}=
\int\limits_t^T \phi_{j_4}(t_4)\int\limits_t^{t_4}\phi_{j_3}(t_3)
d{\bf w}_{t_3}^{(i_3)}
d{\bf w}_{t_4}^{(i_3)}+
\int\limits_t^T \phi_{j_3}(t_3)\int\limits_t^{t_3}\phi_{j_4}(t_4)
d{\bf w}_{t_4}^{(i_3)}
d{\bf w}_{t_3}^{(i_3)},
\end{equation}

\vspace{2mm}
\begin{equation}
\label{caseu3}
\zeta_{j_1}^{(i_1)}\zeta_{j_2}^{(i_1)}-
{\bf 1}_{\{j_1=j_2\}}=
\int\limits_t^T \phi_{j_2}(t_2)\int\limits_t^{t_2}\phi_{j_1}(t_1)
d{\bf w}_{t_1}^{(i_1)}
d{\bf w}_{t_2}^{(i_1)}+
\int\limits_t^T \phi_{j_1}(t_1)\int\limits_t^{t_1}\phi_{j_2}(t_2)
d{\bf w}_{t_2}^{(i_1)}
d{\bf w}_{t_1}^{(i_1)}
\end{equation}

\vspace{5mm}
\noindent
w.~p.~1.

The relations (\ref{caseu1})--(\ref{caseu3}) and (\ref{yyee22})
imply the following 

\vspace{1mm}
$$
{\sf M}\left\{\left(I_{(0000)T,t}^{*(i_1i_1i_3i_3)}-
I_{(0000)T,t}^{*(i_1i_1i_3i_3)q}\right)^2\right\}=
{\sf M}\left\{\left(I_{(0000)T,t}^{(i_1i_1i_3i_3)}-
I_{(0000)T,t}^{(i_1i_1i_3i_3)q}\right)^2\right\}+
$$

\vspace{2mm}
$$
+
{\sf M}\left\{\left(\frac{1}{2}
\int\limits_t^T\int\limits_t^{t_4}(t_3-t)
d{\bf w}_{t_3}^{(i_3)}
d{\bf w}_{t_4}^{(i_3)}-
\sum\limits_{j_4,j_3=0}^q
\sum\limits_{j_1=0}^q C_{j_4 j_3 j_1 j_1}\left(\zeta_{j_3}^{(i_3)}\zeta_{j_4}^{(i_3)}-
{\bf 1}_{\{j_3=j_4\}}\right)\right)^2\right\}+
$$

\vspace{2mm}
$$
+{\sf M}\left\{\left(\frac{1}{2}
\int\limits_t^T\int\limits_t^{t_3}\int\limits_t^{t_2}
d{\bf w}_{t_1}^{(i_1)}
d{\bf w}_{t_2}^{(i_1)}dt_3-
\sum\limits_{j_2,j_1=0}^q
\sum\limits_{j_3=0}^q C_{j_3 j_3 j_2 j_1}\left(\zeta_{j_1}^{(i_1)}\zeta_{j_2}^{(i_1)}-
{\bf 1}_{\{j_1=j_2\}}\right)\right)^2\right\}+
$$

\vspace{2mm}
$$
+\left(\frac{(T-t)^2}{8}-
\sum\limits_{j_3,j_1=0}^q C_{j_3 j_3 j_1 j_1}\right)^2=
{\sf M}\left\{\left(I_{(0000)T,t}^{(i_1i_1i_3i_3)}-
I_{(0000)T,t}^{(i_1i_1i_3i_3)q}\right)^2\right\}+
$$

\vspace{2mm}
$$
+{\sf M}\left\{\left(\frac{1}{2}
\int\limits_t^T\int\limits_t^{t_4}(t_3-t)
d{\bf w}_{t_3}^{(i_3)}
d{\bf w}_{t_4}^{(i_3)}-
\sum\limits_{j_4,j_3=0}^q
\sum\limits_{j_1=0}^q C_{j_4 j_3 j_1 j_1}\left(\zeta_{j_3}^{(i_3)}\zeta_{j_4}^{(i_3)}-
{\bf 1}_{\{j_3=j_4\}}\right)\right)^2\right\}+
$$

\vspace{2mm}
$$
+{\sf M}\left\{\left(\frac{1}{2}
\int\limits_t^T(T-t_2)\int\limits_t^{t_2}
d{\bf w}_{t_1}^{(i_1)}
d{\bf w}_{t_2}^{(i_1)}-
\sum\limits_{j_2,j_1=0}^q
\sum\limits_{j_3=0}^q C_{j_3 j_3 j_2 j_1}\left(\zeta_{j_1}^{(i_1)}\zeta_{j_2}^{(i_1)}-
{\bf 1}_{\{j_1=j_2\}}\right)\right)^2\right\}+
$$

\vspace{2mm}
$$
+\left(\frac{(T-t)^2}{8}-
\sum\limits_{j_3,j_1=0}^q C_{j_3 j_3 j_1 j_1}\right)^2=
$$

\vspace{4mm}
$$
=
{\sf M}\left\{\left(I_{(0000)T,t}^{(i_1i_1i_3i_3)}-
I_{(0000)T,t}^{(i_1i_1i_3i_3)q}\right)^2\right\}+
\frac{(T-t)^4}{48}
+\sum\limits_{j_4,j_3=0}^q
\sum\limits_{j_1=0}^q C_{j_4 j_3 j_1 j_1}\left(C_{j_3j_4}^{10}+C_{j_4j_3}^{10}\right)+
$$

\vspace{2mm}
$$
+
{\sf M}\left\{\left(
\sum\limits_{j_4,j_3=0}^q
\sum\limits_{j_1=0}^q C_{j_4 j_3 j_1 j_1}\left(\zeta_{j_3}^{(i_3)}\zeta_{j_4}^{(i_3)}-
{\bf 1}_{\{j_3=j_4\}}\right)\right)^2\right\}+
$$

\vspace{2mm}
$$
+
\frac{(T-t)^4}{48}
-\sum\limits_{j_2,j_1=0}^q
\sum\limits_{j_3=0}^q C_{j_3 j_3 j_2 j_1}\left((T-t)\left(C_{j_1j_2}+C_{j_2j_1}\right)+
C_{j_1j_2}^{01}+C_{j_2j_1}^{01}\right)+
$$

\vspace{2mm}
$$
+
{\sf M}\left\{\left(
\sum\limits_{j_2,j_1=0}^q
\sum\limits_{j_3=0}^q C_{j_3 j_3 j_2 j_1}\left(\zeta_{j_1}^{(i_1)}\zeta_{j_2}^{(i_1)}-
{\bf 1}_{\{j_1=j_2\}}\right)\right)^2\right\}+
$$

\vspace{2mm}
\begin{equation}
\label{caseu100}
+\left(\frac{(T-t)^2}{8}-
\sum\limits_{j_3,j_1=0}^q C_{j_3 j_3 j_1 j_1}\right)^2.
\end{equation}

\vspace{5mm}

Furthermore,

\vspace{-1mm}
$$
{\sf M}\left\{\left(
\sum\limits_{j_4,j_3=0}^q
\sum\limits_{j_1=0}^q C_{j_4 j_3 j_1 j_1}\left(\zeta_{j_3}^{(i_3)}\zeta_{j_4}^{(i_3)}-
{\bf 1}_{\{j_3=j_4\}}\right)\right)^2\right\}=
$$

\vspace{2mm}
$$
={\sf M}\left\{\left(
\sum\limits_{j_4,j_3=0}^q
\sum\limits_{j_1=0}^q C_{j_4 j_3 j_1 j_1}\zeta_{j_3}^{(i_3)}\zeta_{j_4}^{(i_3)}
\right)^2\right\}-
2\left(\sum\limits_{j_3,j_1=0}^q C_{j_3 j_3 j_1 j_1}\right)^2+
\left(\sum\limits_{j_3,j_1=0}^q C_{j_3 j_3 j_1 j_1}\right)^2=
$$

\vspace{2mm}
\begin{equation}
\label{caseu101}
={\sf M}\left\{\left(
\sum\limits_{j_4,j_3=0}^q
\sum\limits_{j_1=0}^q C_{j_4 j_3 j_1 j_1}\zeta_{j_3}^{(i_3)}\zeta_{j_4}^{(i_3)}
\right)^2\right\}-
\left(\sum\limits_{j_3,j_1=0}^q C_{j_3 j_3 j_1 j_1}\right)^2,
\end{equation}

\vspace{3mm}
$$
{\sf M}\left\{\left(
\sum\limits_{j_2,j_1=0}^q
\sum\limits_{j_3=0}^q C_{j_3 j_3 j_2 j_1}\left(\zeta_{j_1}^{(i_1)}\zeta_{j_2}^{(i_1)}-
{\bf 1}_{\{j_1=j_2\}}\right)\right)^2\right\}=
$$
\begin{equation}
\label{caseu102}
={\sf M}\left\{\left(
\sum\limits_{j_2,j_1=0}^q
\sum\limits_{j_3=0}^q C_{j_3 j_3 j_2 j_1}\zeta_{j_1}^{(i_1)}\zeta_{j_2}^{(i_1)}
\right)^2\right\}-
\left(
\sum\limits_{j_1,j_3=0}^q C_{j_3 j_3 j_1 j_1}\right)^2.
\end{equation}

\vspace{5mm}

We have \cite{15c}, p.~71 (also see \cite{10a}, Sect.~2.3)

\vspace{1mm}
\begin{equation}
\label{otit321}
{\sf M}\left\{\left(\sum\limits_{j_3, j_4=0}^{q}
a_{j_4 j_3} \zeta_{j_3}^{(i)}
\zeta_{j_4}^{(i)}\right)^2\right\}=
\left(\sum_{j_3=0}^q a_{j_3j_3}\right)^2+
\sum\limits_{j_4=0}^q\sum\limits_{j_3=0}^{j_4-1}
\left(a_{j_3j_4} + a_{j_4j_3}\right)^2
+2\sum\limits_{j_4=0}^q\left(a_{j_4j_4}\right)^2,
\end{equation}

\vspace{5mm}
\noindent
where $i=1,\ldots,m$ and 
$a_{j_4 j_3}$ $(j_3,j_4=0,1,\ldots,q)$ are scalar nonrandom coefficients.

Applying (\ref{otit321}), we obtain

\vspace{1mm}
$$
{\sf M}\left\{\left(
\sum\limits_{j_4,j_3=0}^q
\sum\limits_{j_1=0}^q C_{j_4 j_3 j_1 j_1}\zeta_{j_3}^{(i_3)}\zeta_{j_4}^{(i_3)}
\right)^2\right\}=
$$

\vspace{2mm}
\begin{equation}
\label{caseu103}
=\left(\sum\limits_{j_3,j_1=0}^q C_{j_3 j_3 j_1 j_1}\right)^2+
\sum\limits_{j_4=0}^q \sum\limits_{j_3=0}^{j_4-1}
\left(\sum\limits_{j_1=0}^q C_{j_3 j_4 j_1 j_1}+ \sum\limits_{j_1=0}^q C_{j_4 j_3 j_1 j_1}
\right)^2
+2\sum\limits_{j_4=0}^q\left(\sum\limits_{j_1=0}^q C_{j_4 j_4 j_1 j_1}\right)^2.
\end{equation}

\vspace{5mm}

From (\ref{caseu101}) and (\ref{caseu103}) we get 

\vspace{1mm}
$$
{\sf M}\left\{\left(
\sum\limits_{j_4,j_3=0}^q
\sum\limits_{j_1=0}^q C_{j_4 j_3 j_1 j_1}\left(\zeta_{j_3}^{(i_3)}\zeta_{j_4}^{(i_3)}-
{\bf 1}_{\{j_3=j_4\}}\right)\right)^2\right\}=
$$

\vspace{2mm}
\begin{equation}
\label{caseu104}
=\sum\limits_{j_4=0}^q \sum\limits_{j_3=0}^{j_4-1}
\left(\sum\limits_{j_1=0}^q C_{j_3 j_4 j_1 j_1}+ 
\sum\limits_{j_1=0}^q C_{j_4 j_3 j_1 j_1}\right)^2
+2\sum\limits_{j_4=0}^q\left(\sum\limits_{j_1=0}^q C_{j_4 j_4 j_1 j_1}\right)^2.
\end{equation}

\vspace{5mm}

By analogy with (\ref{caseu104}) we obtain

\vspace{1mm}
$$
{\sf M}\left\{\left(
\sum\limits_{j_2,j_1=0}^q
\sum\limits_{j_3=0}^q C_{j_3 j_3 j_2 j_1}\left(\zeta_{j_1}^{(i_1)}\zeta_{j_2}^{(i_1)}-
{\bf 1}_{\{j_1=j_2\}}\right)\right)^2\right\}=
$$

\vspace{2mm}
\begin{equation}
\label{caseu105}
=\sum\limits_{j_2=0}^q \sum\limits_{j_1=0}^{j_2-1}
\left(\sum\limits_{j_3=0}^q C_{j_3 j_3 j_1 j_2}+ 
\sum\limits_{j_3=0}^q C_{j_3 j_3 j_2 j_1}\right)^2
+2\sum\limits_{j_2=0}^q\left(\sum\limits_{j_3=0}^q C_{j_3 j_3 j_2 j_2}\right)^2.
\end{equation}

\vspace{5mm}

Combining (\ref{usl11}), (\ref{caseu100}), (\ref{caseu104}),
and (\ref{caseu105}), we finally have 

\vspace{2mm}
$$
{\sf M}\left\{\left(I_{(0000)T,t}^{*(i_1i_2i_3i_4)}-
I_{(0000)T,t}^{*(i_1i_2i_3i_4)q}\right)^2\right\}=
\frac{(T-t)^4}{12}-
$$

\vspace{2mm}
$$
- \sum_{j_1,j_2,j_3,j_4=0}^{q}
C_{j_4j_3j_2j_1}\Biggl(\sum\limits_{(j_1,j_2)}\Biggl(
\sum\limits_{(j_3,j_4)}
C_{j_4j_3j_2j_1}\Biggr)\Biggr)
+\sum\limits_{j_4,j_3=0}^q
\sum\limits_{j_1=0}^q C_{j_4 j_3 j_1 j_1}\left(C_{j_3j_4}^{10}+C_{j_4j_3}^{10}\right)+
$$

\vspace{2mm}
$$
+\sum\limits_{j_4=0}^q \sum\limits_{j_3=0}^{j_4-1}
\left(\sum\limits_{j_1=0}^q C_{j_3 j_4 j_1 j_1}+ 
\sum\limits_{j_1=0}^q C_{j_4 j_3 j_1 j_1}\right)^2
+2\sum\limits_{j_4=0}^q\left(\sum\limits_{j_1=0}^q C_{j_4 j_4 j_1 j_1}\right)^2
-
$$

\vspace{2mm}
$$
-\sum\limits_{j_2,j_1=0}^q
\sum\limits_{j_3=0}^q C_{j_3 j_3 j_2 j_1}\left((T-t)C_{j_1}C_{j_2}+
C_{j_1j_2}^{01}+C_{j_2j_1}^{01}\right)+
$$

\vspace{2mm}
$$
+\sum\limits_{j_2=0}^q \sum\limits_{j_1=0}^{j_2-1}
\left(\sum\limits_{j_3=0}^q C_{j_3 j_3 j_1 j_2}+ 
\sum\limits_{j_3=0}^q C_{j_3 j_3 j_2 j_1}\right)^2
+2\sum\limits_{j_2=0}^q\left(\sum\limits_{j_3=0}^q C_{j_3 j_3 j_2 j_2}\right)^2+
$$

\vspace{2mm}
$$
+\left(\frac{(T-t)^2}{8}-
\sum\limits_{j_3,j_1=0}^q C_{j_3 j_3 j_1 j_1}\right)^2,
$$

\vspace{5mm}
\noindent
where $i_1=i_2\ne i_3=i_4$ and

\vspace{-2mm}
$$
C_{j}=\int\limits_t^T
\phi_{j}(\tau)d\tau
=\left\{
\begin{matrix}
\sqrt{T-t},\ & j=0\cr\cr
0,\ & j\ne 0
\end{matrix}
.\right.
$$

\vspace{5mm}

Consider the case (\ref{casee12}) by analogy with the case (\ref{casee11}).
Using (\ref{casee400}), we obtain 

\vspace{1mm}
$$
{\sf M}\left\{\left(I_{(0000)T,t}^{*(i_1i_2i_1i_2)}-
I_{(0000)T,t}^{*(i_1i_2i_1i_2)q}\right)^2\right\}=
{\sf M}\left\{\Biggl(I_{(0000)T,t}^{(i_1i_2i_1i_2)}-
I_{(0000)T,t}^{(i_1i_2i_1i_2)q}\Biggr.\right.-
$$

\vspace{2mm}
$$
-\sum\limits_{j_4,j_2=0}^q
\sum\limits_{j_1=0}^q C_{j_4 j_1 j_2 j_1}\zeta_{j_2}^{(i_2)}\zeta_{j_4}^{(i_2)}
-\sum\limits_{j_3,j_1=0}^q
\sum\limits_{j_2=0}^q C_{j_2 j_3 j_2 j_1}\zeta_{j_1}^{(i_1)}\zeta_{j_3}^{(i_1)}
\left.\left.+\sum\limits_{j_2,j_1=0}^q C_{j_2 j_1 j_2 j_1}\right)^2\right\}=
$$

\vspace{2mm}
$$
={\sf M}\left\{\Biggl(I_{(0000)T,t}^{(i_1i_2i_1i_2)}-
I_{(0000)T,t}^{(i_1i_2i_1i_2)q}\Biggr.\right.
-\sum\limits_{j_4,j_2=0}^q
\sum\limits_{j_1=0}^q C_{j_4 j_1 j_2 j_1}\left(\zeta_{j_2}^{(i_2)}\zeta_{j_4}^{(i_2)}-
{\bf 1}_{\{j_2=j_4\}}\right)-
$$

\vspace{2mm}
$$
-\sum\limits_{j_3,j_1=0}^q
\sum\limits_{j_2=0}^q C_{j_2 j_3 j_2 j_1}\left(\zeta_{j_1}^{(i_1)}\zeta_{j_3}^{(i_1)}-
{\bf 1}_{\{j_1=j_3\}}\right)
\left.\left.-\sum\limits_{j_2,j_1=0}^q C_{j_2 j_1 j_2 j_1}\right)^2\right\}=
$$

\vspace{2mm}
$$
={\sf M}\left\{\left(I_{(0000)T,t}^{(i_1i_2i_1i_2)}-
I_{(0000)T,t}^{(i_1i_2i_1i_2)q}\right)^2\right\}+
$$

\vspace{2mm}
$$
+{\sf M}\left\{\left(\sum\limits_{j_4,j_2=0}^q
\sum\limits_{j_1=0}^q C_{j_4 j_1 j_2 j_1}\left(\zeta_{j_2}^{(i_2)}\zeta_{j_4}^{(i_2)}-
{\bf 1}_{\{j_2=j_4\}}\right)\right)^2\right\}+
$$

\vspace{2mm}
$$
+{\sf M}\left\{\left(\sum\limits_{j_3,j_1=0}^q
\sum\limits_{j_2=0}^q C_{j_2 j_3 j_2 j_1}\left(\zeta_{j_1}^{(i_1)}\zeta_{j_3}^{(i_1)}-
{\bf 1}_{\{j_1=j_3\}}\right)\right)^2\right\}+
$$

\vspace{2mm}
\begin{equation}
\label{caseu700}
+\left(\sum\limits_{j_2,j_1=0}^q C_{j_2 j_1 j_2 j_1}\right)^2.
\end{equation}

\vspace{5mm}

Applying (\ref{usl12}) and (\ref{caseu700}), we finally get 

\vspace{1mm}
$$
{\sf M}\left\{\left(I_{(0000)T,t}^{*(i_1i_2i_3i_4)}-
I_{(0000)T,t}^{*(i_1i_2i_3i_4)q}\right)^2\right\}=
\frac{(T-t)^4}{24}-\sum_{j_1,j_2,j_3,j_4=0}^{q}
C_{j_4j_3j_2j_1}\Biggl(\sum\limits_{(j_1,j_3)}\Biggl(
\sum\limits_{(j_2,j_4)}
C_{j_4j_3j_2j_1}\Biggr)\Biggr)+
$$

\vspace{2mm}
$$
+\sum\limits_{j_4=0}^q \sum\limits_{j_2=0}^{j_4-1}
\left(\sum\limits_{j_1=0}^q C_{j_2 j_1 j_4 j_1}+ \sum\limits_{j_1=0}^q C_{j_4 j_1 j_2 j_1}
\right)^2
+2\sum\limits_{j_4=0}^q\left(\sum\limits_{j_1=0}^q C_{j_4 j_1 j_4 j_1}\right)^2+
$$

\vspace{2mm}
$$
+\sum\limits_{j_3=0}^q \sum\limits_{j_1=0}^{j_3-1}
\left(\sum\limits_{j_2=0}^q C_{j_2 j_1 j_2 j_3}+\sum\limits_{j_2=0}^q C_{j_2 j_3 j_2 j_1}
\right)^2
+2\sum\limits_{j_3=0}^q\left(\sum\limits_{j_2=0}^q C_{j_2 j_3 j_2 j_3}\right)^2+
$$

\vspace{2mm}
$$
+
\left(\sum\limits_{j_2,j_1=0}^q C_{j_2 j_1 j_2 j_1}\right)^2,
$$

\vspace{5mm}
\noindent
where $i_1=i_3\ne i_2=i_4$.

Consider the case (\ref{casee13}) by analogy with the cases (\ref{casee11}) and (\ref{casee12}).
Applying (\ref{casee400}), we obtain

\vspace{1mm}
$$
{\sf M}\left\{\left(I_{(0000)T,t}^{*(i_1i_2i_2i_1)}-
I_{(0000)T,t}^{*(i_1i_2i_2i_1)q}\right)^2\right\}=
$$

\vspace{2mm}
$$
={\sf M}\left\{\left(I_{(0000)T,t}^{(i_1i_2i_2i_1)}
+\frac{1}{2}
\int\limits_t^T\int\limits_t^{t_4}\int\limits_t^{t_2}
d{\bf w}_{t_1}^{(i_1)}dt_2
d{\bf w}_{t_4}^{(i_1)}
-
I_{(0000)T,t}^{(i_1i_2i_2i_1)q}-\right.\right.
$$

\vspace{2mm}
$$
-\sum\limits_{j_3,j_2=0}^q
\sum\limits_{j_1=0}^q C_{j_1 j_3 j_2 j_1}\zeta_{j_2}^{(i_2)}\zeta_{j_3}^{(i_2)}-
\sum\limits_{j_4,j_1=0}^q
\sum\limits_{j_2=0}^q C_{j_4 j_2 j_2 j_1}\zeta_{j_1}^{(i_1)}\zeta_{j_4}^{(i_1)}+
\left.\left.
\sum\limits_{j_2,j_1=0}^q C_{j_1 j_2 j_2 j_1}\right)^2\right\}=
$$

\vspace{3mm}
$$
={\sf M}\left\{\Biggl(I_{(0000)T,t}^{(i_1i_2i_2i_1)}-
I_{(0000)T,t}^{(i_1i_2i_2i_1)q}+\Biggr.\right.
$$

\vspace{2mm}
$$
+\frac{1}{2}
\int\limits_t^T\int\limits_t^{t_4}\int\limits_t^{t_2}
d{\bf w}_{t_1}^{(i_1)}dt_2
d{\bf w}_{t_4}^{(i_1)}
-\sum\limits_{j_4,j_1=0}^q
\sum\limits_{j_2=0}^q C_{j_4 j_2 j_2 j_1}\left(\zeta_{j_1}^{(i_1)}\zeta_{j_4}^{(i_1)}-
{\bf 1}_{\{j_1=j_4\}}\right)-
$$

\vspace{2mm}
$$
\left.\left.-\sum\limits_{j_3,j_2=0}^q
\sum\limits_{j_1=0}^q C_{j_1 j_3 j_2 j_1}\left(\zeta_{j_2}^{(i_2)}\zeta_{j_3}^{(i_2)}-
{\bf 1}_{\{j_2=j_3\}}\right)
-
\sum\limits_{j_2,j_1=0}^q C_{j_1 j_2 j_2 j_1}\right)^2\right\}=
$$

\vspace{3mm}
$$
={\sf M}\left\{\left(I_{(0000)T,t}^{(i_1i_2i_2i_1)}-
I_{(0000)T,t}^{(i_1i_2i_2i_1)q}\right)^2\right\}+
$$

\vspace{2mm}
$$
+
{\sf M}\left\{\left(\frac{1}{2}
\int\limits_t^T\int\limits_t^{t_4}(t_4-t_1)
d{\bf w}_{t_1}^{(i_1)}
d{\bf w}_{t_4}^{(i_1)}
-\sum\limits_{j_4,j_1=0}^q
\sum\limits_{j_2=0}^q C_{j_4 j_2 j_2 j_1}\left(\zeta_{j_1}^{(i_1)}\zeta_{j_4}^{(i_1)}-
{\bf 1}_{\{j_1=j_4\}}\right)\right)^2\right\}+
$$

\vspace{2mm}
$$
+{\sf M}\left\{\left(\sum\limits_{j_3,j_2=0}^q
\sum\limits_{j_1=0}^q C_{j_1 j_3 j_2 j_1}\left(\zeta_{j_2}^{(i_2)}\zeta_{j_3}^{(i_2)}-
{\bf 1}_{\{j_2=j_3\}}\right)\right)^2\right\}+
$$

\vspace{2mm}
$$
+
\left(\sum\limits_{j_2,j_1=0}^q C_{j_1 j_2 j_2 j_1}\right)^2=
$$

\vspace{3mm}
$$
={\sf M}\left\{\left(I_{(0000)T,t}^{(i_1i_2i_2i_1)}-
I_{(0000)T,t}^{(i_1i_2i_2i_1)q}\right)^2\right\}+\frac{(T-t)^4}{48}-
$$

\vspace{2mm}
$$
-
\sum\limits_{j_4,j_1=0}^q
\sum\limits_{j_2=0}^q C_{j_4 j_2 j_2 j_1}\left(
\int\limits_t^T \phi_{j_4}(t_4)\int\limits_t^{t_4}(t_4-t_1)\phi_{j_1}(t_1)
dt_1 dt_4+\int\limits_t^T \phi_{j_1}(t_4)\int\limits_t^{t_4}(t_4-t_1)\phi_{j_4}(t_1)
dt_1 dt_4\right)+
$$

\vspace{2mm}
$$
+{\sf M}\left\{\left(\sum\limits_{j_4,j_1=0}^q
\sum\limits_{j_2=0}^q C_{j_4 j_2 j_2 j_1}\left(\zeta_{j_1}^{(i_1)}\zeta_{j_4}^{(i_1)}-
{\bf 1}_{\{j_1=j_4\}}\right)\right)^2\right\}+
$$

\vspace{2mm}
$$
+{\sf M}\left\{\left(\sum\limits_{j_3,j_2=0}^q
\sum\limits_{j_1=0}^q C_{j_1 j_3 j_2 j_1}\left(\zeta_{j_2}^{(i_2)}\zeta_{j_3}^{(i_2)}-
{\bf 1}_{\{j_2=j_3\}}\right)\right)^2\right\}+
$$

\vspace{2mm}
\begin{equation}
\label{caseu800}
+
\left(\sum\limits_{j_2,j_1=0}^q C_{j_1 j_2 j_2 j_1}\right)^2.
\end{equation}

\vspace{5mm}

Applying (\ref{usl13}) and (\ref{caseu800}), we finally obtain 

\vspace{1mm}
$$
{\sf M}\left\{\left(I_{(0000)T,t}^{*(i_1i_2i_3i_4)}-
I_{(0000)T,t}^{*(i_1i_2i_3i_4)q}\right)^2\right\}=\frac{(T-t)^4}{16}-
\sum_{j_1,j_2,j_3,j_4=0}^{q}
C_{j_4j_3j_2j_1}\Biggl(\sum\limits_{(j_1,j_4)}\Biggl(
\sum\limits_{(j_2,j_3)}
C_{j_4j_3j_2j_1}\Biggr)\Biggr)-
$$

\vspace{2mm}
$$
-\sum\limits_{j_4,j_1=0}^q
\sum\limits_{j_2=0}^q C_{j_4 j_2 j_2 j_1}\left(C_{j_4j_1}^{10}+
C_{j_1j_4}^{10}-C_{j_4j_1}^{01}-C_{j_1j_4}^{01}\right)+
$$

\vspace{2mm}
$$
+\sum\limits_{j_4=0}^q \sum\limits_{j_1=0}^{j_4-1}
\left(\sum\limits_{j_2=0}^q C_{j_1 j_2 j_2 j_4} + \sum\limits_{j_2=0}^q C_{j_4 j_2 j_2 j_1}
\right)^2+
2\sum\limits_{j_4=0}^q\left(\sum\limits_{j_2=0}^q C_{j_4 j_2 j_2 j_4}\right)^2+
$$

\vspace{2mm}
$$
+\sum\limits_{j_3=0}^q \sum\limits_{j_2=0}^{j_3-1}
\left(\sum\limits_{j_1=0}^q C_{j_1 j_2 j_3 j_1} + \sum\limits_{j_1=0}^q C_{j_1 j_3 j_2 j_1}
\right)^2
+2\sum\limits_{j_3=0}^q\left(\sum\limits_{j_1=0}^q C_{j_1 j_3 j_3 j_1}\right)^2+
$$

\vspace{2mm}
$$
+
\left(\sum\limits_{j_2,j_1=0}^q C_{j_1 j_2 j_2 j_1}\right)^2,
$$

\vspace{5mm}
\noindent
where $i_1=i_4\ne i_2=i_3$.


\vspace{8mm}


\begin{thebibliography}{199}

\vspace{5mm}


\bibitem{1982}
Gihman I.I., Skorohod A.V. Stochastic Differential Equations and 
its Applications.
Naukova Dumka, Kiev, 1982, 612 pp.



\bibitem{1995}
Kloeden P.E., Platen E. Numerical Solution of Stochastic
Differential Equations. Springer, Berlin, 1992, 632 pp.


\bibitem{1988}
Milstein G.N. Numerical Integration of Stochastic Differential 
Equations. Ural University Press, Sverdlovsk, 1988, 225 pp. 


\bibitem{2004}
Milstein G.N., Tretyakov M.V. 
Stochastic Numerics for Mathematical Physics. 
Springer, Berlin, 2004, 616 pp.


\bibitem{1}
Kuznetsov D.F. A method of expansion and approximation of repeated 
stochastic Stratonovich integrals based on multiple Fourier series 
on full orthonormal systems. [In Russian].
Electronic Journal "Differential Equations and Control Processes"
ISSN 1817-2172 (online),
1 (1997), 18-77.
Available at:\\
{\color{blue}http://diffjournal.spbu.ru/EN/numbers/1997.1/article.1.2.html}


\bibitem{1a}
Kuznetsov D.F. Problems of the numerical analysis of Ito stochastic 
differential equations. 
[In Russian].
Electronic Journal "Differential Equations and Control Processes"
ISSN 1817-2172 (online),
1 (1998), 66-367.
Available at:\
{\color{blue}http://diffjournal.spbu.ru/EN/numbers/1998.1/article.1.3.html}\
Hard Cover Edition: SPbGTU, Saint-Petersburg, 1998, 204 pp. (ISBN 5-7422-0045-5)


\bibitem{2}
Kuznetsov D.F. Mean Square Approximation of Solutions 
of Stochastic Differential 
Equations Using Legendres Polynomials. [In English]. Journal of 
Automation and 
Information Sciences (Begell House), 32, Issue 12, (2000), 69-86.
DOI: {\color{blue}http://doi.org/10.1615/JAutomatInfScien.v32.i12.80} 




\bibitem{3}
Kuznetsov D.F. Numerical Integration of Stochastic Differential Equations. 2.
[In Russian]. Polytechnical University Publishing House, 
Saint-Petersburg, 2006, 764 pp.
DOI: {\color{blue}http://doi.org/10.18720/SPBPU/2/s17-227}\
Available at:\ {\color{blue}http://www.sde-kuznetsov.spb.ru/06.pdf}\
(ISBN 5-7422-1191-0)



\bibitem{3a}
Kuznetsov D.F. Stochastic Differential Equations: Theory and Practice 
of Numerical Solution. With MatLab programs, 1st Edition. [In Russian]. 
Polytechnical University Publishing House, Saint-Petersburg, 2007, 778 pp.
DOI: {\color{blue}http://doi.org/10.18720/SPBPU/2/s17-228}\
Available at:\ {\color{blue}http://www.sde-kuznetsov.spb.ru/07b.pdf}\
(ISBN 5-7422-1394-8)




\bibitem{4}
Kuznetsov D.F. Stochastic Differential Equations: Theory and Practice 
of Numerical Solution. With MatLab programs, 2nd Edition. [In Russian]. 
Polytechnical University 
Publishing House, Saint-Petersburg, 2007, XXXII+770 pp.
DOI: {\color{blue}http://doi.org/10.18720/SPBPU/2/s17-229}\
Available at:\
{\color{blue}http://www.sde-kuznetsov.spb.ru/07a.pdf}\
(ISBN 5-7422-1439-1)




\bibitem{4a}
Kuznetsov D.F. Stochastic Differential Equations: Theory and Practice 
of Numerical Solution. With MatLab programs, 3rd Edition. [In Russian]. 
Polytechnical 
University Publishing House, Saint-Petersburg, 2009, XXXIV+768 pp.
DOI: {\color{blue}http://doi.org/10.18720/SPBPU/2/s17-230}\
Available at:\
{\color{blue}http://www.sde-kuznetsov.spb.ru/09.pdf}\
(ISBN 978-5-7422-2132-6)


                                        

\bibitem{5}
Kuznetsov D.F. Stochastic Differential Equations: Theory and Practice 
of Numerical Solution. With MatLab programs. 4th Edition. [In Russian].
Polytechnical University Publishing House, Saint-Petersburg, 2010, 
XXX+786 pp. DOI: {\color{blue}http://doi.org/10.18720/SPBPU/2/s17-231}\
Available at:\
{\color{blue}http://www.sde-kuznetsov.spb.ru/10.pdf}\
(ISBN 978-5-7422-2448-8)



\bibitem{5a}
Kuznetsov D.F. Multiple stochastic Ito and Stratonovich integrals 
and multiple Fourier series.
[In Russian].
Electronic Journal "Differential Equations and Control Processes"
ISSN 1817-2172 (online),
3 (2010), A.1-A.257. DOI: {\color{blue}http://doi.org/10.18720/SPBPU/2/z17-7}\
Available at:\\
{\color{blue}http://diffjournal.spbu.ru/EN/numbers/2010.3/article.2.1.html}





\bibitem{5b}
Kuznetsov D.F. Strong Approximation of Multiple Ito and 
Stratonovich Stochastic Integrals: Multiple Fourier Series Approach.
1st Edition. [In English]. 
Polytechnical University Publishing House, Saint-Petersburg, 
2011, 250 pp. DOI: {\color{blue}http://doi.org/10.18720/SPBPU/2/s17-232}\
Available at:\\
{\color{blue}http://www.sde-kuznetsov.spb.ru/11b.pdf}\
(ISBN 978-5-7422-2988-9)



\bibitem{6}
Kuznetsov D.F. Strong Approximation of Multiple Ito and 
Stratonovich Stochastic Integrals: Multiple Fourier Series Approach.
2nd Edition. [In English]. 
Polytechnical University Publishing House, Saint-Petersburg, 
2011, 284 pp. DOI: {\color{blue}http://doi.org/10.18720/SPBPU/2/s17-233}\
Available at:\\
{\color{blue}http://www.sde-kuznetsov.spb.ru/11a.pdf}\
(ISBN 978-5-7422-3162-2)



\bibitem{6a}
Kuznetsov D.F. Multiple Ito and Stratonovich stochastic 
integrals: approximations, properties, formulas. [In English].
Polytechnical University Publishing House, Saint-Petersburg,
2013, 382 pp.\\
DOI: {\color{blue}http://doi.org/10.18720/SPBPU/2/s17-234}\\
Available at:\
{\color{blue}http://www.sde-kuznetsov.spb.ru/13.pdf}\
(ISBN 978-5-7422-3973-4)




\bibitem{7}
Kuznetsov D.F. Multiple Ito and Stratonovich stochastic integrals: 
Fourier-Legendre and trigonometric expansions, approximations, formulas.
[In English].
Electronic Journal "Differential Equations and Control Processes"
ISSN 1817-2172 (online),
1 (2017), A.1--A.385. DOI: {\color{blue}http://doi.org/10.18720/SPBPU/2/z17-3}\
Available at:\
{\color{blue}http://diffjournal.spbu.ru/EN/numbers/2017.1/article.2.1.html} 


\bibitem{8}
Kuznetsov D.F. Development and Application of the Fourier 
Method for the Numerical Solution of Ito Stochastic Differential 
Equations. [In English]. Computational Mathematics and 
Mathematical Physics, 58, 7 (2018), 1058-1070.
DOI: {\color{blue}http://doi.org/10.1134/S0965542518070096}



\bibitem{9}
Kuznetsov D.F. On Numerical Modeling of the Multidimensional 
Dynamic Systems Under Random Perturbations With the 1.5 and 2.0 
Orders of Strong Convergence [In English]. Automation and Remote Control, 
79, 7 (2018), 1240-1254.
DOI: {\color{blue}http://doi.org/10.1134/S0005117918070056}


\bibitem{9a}
Kuznetsov D.F. Stochastic Differential Equations: Theory and Practice of 
Numerical Solution. With Programs on MATLAB, 5th Edition. [In Russian].
Electronic Journal "Differential Equations and Control Processes"
ISSN 1817-2172 (online), 2 (2017), A.1-A.1000.
DOI: {\color{blue}http://doi.org/10.18720/SPBPU/2/z17-4}\
Available at:\\
{\color{blue}http://diffjournal.spbu.ru/EN/numbers/2017.2/article.2.1.html}





\bibitem{10}
Kuznetsov D.F. Stochastic Differential Equations: Theory and Practice of 
Numerical Solution. With MATLAB Programs, 6th Edition. [In Russian].
Electronic Journal "Differential Equations and Control Processes"
ISSN 1817-2172 (online), 4 (2018), A.1-A.1073.\
Available at:\\
{\color{blue}http://diffjournal.spbu.ru/EN/numbers/2018.4/article.2.1.html}


\bibitem{10a}
Kuznetsov D.F.
Strong approximation of iterated Ito and Stratonovich stochastic 
integrals based on generalized multiple Fourier series. 
Application to numerical solution of Ito SDEs and semilinear SPDEs.
[In English].
\href{http://doi.org/10.48550/arXiv.2003.14184}{arXiv:2003.14184}{arXiv:2003.14184} [math.PR], 2026, 1246 pp.


\bibitem{10axx}
Kuznetsov D.F.
Strong approximation of iterated Ito and Stratonovich stochastic 
integrals based on generalized multiple Fourier series. 
Application to numerical solution of Ito SDEs and semilinear SPDEs.
Electronic Journal "Differential Equations and Control Processes"
ISSN 1817-2172 (online), 4 (2020), A.1-A.606.\
Available at:\
{\color{blue}http://diffjournal.spbu.ru/EN/numbers/2020.4/article.1.8.html}


\bibitem{10axx1}
Kuznetsov D.F.
Mean-Square Approximation of Iterated It\^{o} and Stratonovich Stochastic 
Integrals Based on Generalized Multiple Fourier Series. 
Application to Numerical Integration of It\^{o} SDEs and Semilinear SPDEs.
[In English].
Electronic Journal "Differential Equations and Control Processes"
ISSN 1817-2172 (online),
4 (2021), A.1-A.788.\ Available at:\
{\color{blue}http://diffjournal.spbu.ru/EN/numbers/2021.4/article.1.9.html}



\bibitem{11}
Kuznetsov D.F. 
Expansion of iterated Ito stochastic integrals of arbitrary multiplicity
based on generalized multiple Fourier series converging in the mean. 
[In English].
\href{http://doi.org/10.48550/arXiv.1712.09746}{arXiv:1712.09746} [math.PR]. 2026, 151 pp. 


\bibitem{12}
Kuznetsov D.F. 
Expansions of Iterated Stratonovich stochastic integrals 
based on generalized multiple Fourier series:
multiplicities 1 to 8 and beyond. [in English].
\href{http://doi.org/10.48550/arXiv.1712.09516}{arXiv:1712.09516} [math.PR]. 2026, 392 pp. 



\bibitem{13}
Kuznetsov D.F. Mean-square approximation of iterated Ito and 
Stratonovich stochastic 
integrals of multiplicities 1 to 6 
from the Taylor-Ito and Taylor-Stratonovich expansions
using Legendre polynomials. [In English]. 
\href{http://doi.org/10.48550/arXiv.1801.00231}{arXiv:1801.00231} [math.PR]. 2022, 106 pp. 




\bibitem{14}
Kuznetsov D.F.
Expansion of iterated Stratonovich stochastic integrals of arbitrary 
multiplicity
based on generalized iterated Fourier series converging pointwise. 
[In English].
\href{http://doi.org/10.48550/arXiv.1801.00784}{arXiv:1801.00784} [math.PR]. 2023, 80 pp. 



\bibitem{15}
Kuznetsov D.F. Strong numerical methods of orders 2.0, 2.5, and 3.0 for 
Ito stochastic differential equations based on the unified 
stochastic Taylor expansions and multiple Fourier-Legendre series.
[In English]. 
\href{http://doi.org/10.48550/arXiv.1807.02190}{arXiv:1807.02190} [math.PR]. 2022, 44 pp.



\bibitem{15a}
Kuznetsov D.F.
Expansion of iterated Stratonovich stochastic integrals of fifth, sixth,
seventh and eighth
multiplicities based on generalized multiple Fourier series.
[In English]. 
\href{http://doi.org/10.48550/arXiv.1802.00643}{arXiv:1802.00643} [math.PR]. 2026, 304 pp. 



\bibitem{15b}
Kuznetsov D.F.
Exact calculation of mean-square error in the method of approximation of 
iterated Ito stochastic 
integrals based on the generalized multiple Fourier series. [In English].
\href{http://doi.org/10.48550/arXiv.1801.01079}{arXiv:1801.01079} [math.PR]. 2023, 71 pp. 




\bibitem{15c}
Kuznetsov D.F. Expansion of iterated Stratonovich stochastic integrals
based on generalized multiple Fourier series. 
[In English]. Ufa Mathematical Journal, 
11, 4 (2019), 49-77. 
DOI: {\color{blue}http://doi.org/10.13108/2019-11-4-49}\\
Available at:\
{\color{blue}http://matem.anrb.ru/en/article?art\_id=604}



\bibitem{15d}
Kuznetsov D.F. On Numerical Modeling of the Multidimentional Dynamic 
Systems Under Random Perturbations With the 2.5 Order of Strong 
Convergence. [In English]. Automation and Remote Control, 
80, 5 (2019), 867-881. DOI: {\color{blue}http://doi.org/10.1134/S0005117919050060}


\bibitem{15f}
Kuznetsov D.F. Application of the method of approximation of iterated 
Ito stochastic integrals based on generalized multiple Fourier 
series to the high-order strong numerical methods for non-commutative 
semilinear stochastic partial differential equations. [In English].
\href{http://doi.org/10.48550/arXiv.1905.03724}{arXiv:1905.03724} [math.GM], 2022, 41 pp.



\bibitem{17}
Kuznetsov D.F. Comparative analysis of the efficiency of application 
of Legendre polynomials and trigonometric functions to the numerical 
integration of It\^{o} stochastic differential equations.
[In English]. Computational Mathematics and Mathematical Physics, 
59, 8 (2019),  1236-1250. DOI: {\color{blue}http://doi.org/10.1134/S0965542519080116}



\bibitem{18}
Kuznetsov D.F. New representations of the Taylor-Stratonovich expansion.
Journal of Mathematical Sciences (N.Y.), 118, 6 (2003), 5586-5596.\	
DOI: {\color{blue}http://doi.org/10.1023/A:1026138522239}



\bibitem{xxx}
Kuznetsov D.F. Expansion of iterated stochastic integrals with respect
to martingale Poisson measures and with respect to martingales
based on generalized multiple Fourier series. 
[In English]. 
\href{http://doi.org/10.48550/arXiv.1801.06501}{arXiv:1801.06501} [math.PR].
2023, 40 pp.


\bibitem{new-1}
Kuznetsov D.F. Application of the method of approximation 
of iterated stochastic It\^{o} integrals based on generalized multiple 
Fourier series to the high-order strong numeri\-cal methods 
for non-commutative semilinear stochastic partial 
differential equations. Electronic Journal "Differential Equations and Control Processes"
ISSN 1817-2172 (online),
3 (2019), 18-62.
Available at:\\
{\color{blue}http://diffjournal.spbu.ru/EN/numbers/2019.3/article.1.2.html}



\bibitem{new-3}
Kuznetsov D.F. Application of multiple Fourier-Legendre series 
to implementation of strong exponential Milstein and 
Wagner-Platen methods for non-commutative semilinear stochastic 
partial differential equations. [In English].
\href{http://doi.org/10.48550/arXiv.1912.02612}{arXiv:1912.02612} [math.PR], 2022, 32 pp.


\bibitem{new-4}
Kuznetsov D.F. Application of multiple Fourier-Legendre series 
to strong exponential Milstein and 
Wagner-Platen methods for non-commutative semilinear stochastic 
partial differential equations. [In English]. 
Electronic Journal "Differential Equations and Control Processes"
ISSN 1817-2172 (online),
3 (2020), 129-162.
Available at:\
{\color{blue}http://diffjournal.spbu.ru/EN/numbers/2020.3/article.1.6.html}




\bibitem{arxiv-4}
Kuznetsov D.F.
The hypotheses on expansions of iterated Stratonovich stochastic 
integrals of arbitrary multiplicity and their partial proof. [in English].
\href{http://doi.org/10.48550/arXiv.1801.03195}{arXiv:1801.03195} [math.PR]. 
2026, 318 pp. 





\bibitem{9999}
Kuznetsov D.F. The proof of convergence with probability 1 
in the method of expansion 
of iterated Ito stochastic integrals based on generalized multiple 
Fourier series.
[In English]. Electronic Journal "Differential Equations and Control Processes"
ISSN 1817-2172 (online),
2 (2020), 89-117. Available at:\\ 
{\color{blue}http://diffjournal.spbu.ru/RU/numbers/2020.2/article.1.6.html}



\bibitem{Kuz-Kuz}
Kuznetsov M.D., Kuznetsov D.F.
SDE-MATH: A software package for the implementation of strong high-order 
numerical methods for Ito SDEs with multidimensional non-commutative noise 
based on multiple Fourier-Legendre series. [In English].
Electronic Journal "Differential Equations and Control Processes"
ISSN 1817-2172 (online),
1 (2021), 93-422. Available at:\
{\color{blue}http://diffjournal.spbu.ru/EN/numbers/2021.1/article.1.5.html}



\bibitem{Mikh-1}
Kuznetsov M.D., Kuznetsov D.F.
Implementation of strong numerical methods 
of orders 0.5, 1.0, 1.5, 2.0, 2.5, and 3.0 for Ito SDEs with non-commutative 
noise based on the unified Taylor-Ito and Taylor-Stratonovich 
expansions and multiple Fourier-Legendre series. [In English].
\href{http://doi.org/10.48550/arXiv.2009.14011}{arXiv:2009.14011} [math.PR], 2025, 347 pp. 


\bibitem{Mikh-2}
Kuznetsov M.D., Kuznetsov D.F.
Optimization of the mean-square approximation procedures for 
iterated Ito stochastic integrals of multiplicities 1 to 5 from 
the unified Taylor-Ito expansion based on multiple Fourier-Legendre series.
[In English].
\href{http://doi.org/10.48550/arXiv.2010.13564}{arXiv:2010.13564} [math.PR], 2022, 63 pp. 



\bibitem{Kuzh-1}
Kuznetsov D.F. 
Application of multiple Fourier-Legendre series to the implementation 
of strong exponential Milstein and Wagner-Platen methods for 
non-commutative semilinear SPDEs. 
Proceedings of the XIII International Conference on Applied 
Mathematics and Mechanics in the Aerospace Industry (AMMAI-2020). 
MAI, Moscow, 2020, pp. 451-453.\
Available at:\ {\color{blue}http://www.sde-kuznetsov.spb.ru/20e.pdf}


\bibitem{new-new-1}
Kuznetsov D.F., Kuznetsov M.D. Mean-square approximation of iterated 
stochastic integrals from strong exponential Milstein and Wagner--Platen 
methods for non-commutative semilinear SPDEs based on multiple 
Fourier--Legendre series. Recent Developments in Stochastic Methods and 
Applications. ICSM-5 2020. 
Springer Proceedings in Mathematics \& Statistics, vol 371, Eds. 
Shiryaev, A.N., Samouylov, K.E., Kozyrev, D.V.
Springer, Cham, 2021, pp. 17-32.\ DOI: {\color{blue}http://doi.org/10.1007/978-3-030-83266-7\_2}


\bibitem{new-art-1-xxy}
Kuznetsov D.F. A new approach to the series expansion of iterated 
Stratonovich stochastic integrals of arbitrary multiplicity with respect 
to components of the multidimensional Wiener process. [In English].
Electronic Journal "Differential Equations and Control Processes"
ISSN 1817-2172 (online), 2 (2022), 83-186. 
Available at:\\
{\color{blue}http://diffjournal.spbu.ru/EN/numbers/2022.2/article.1.6.html}



\bibitem{new-new-2}
Kuznetsov, D.F. The proof of convergence with probability 1 in the method of expansion 
of iterated Ito stochastic integrals based on generalized multiple Fourier series.
\href{http://doi.org/10.48550/arXiv.2006.16040}{arXiv:2006.16040} [math.PR], 2026, 33 pp. [In English].



\bibitem{KPS}
Kloeden P.E., Platen E., Schurz H. Numerical solution
of SDE through computer experiments. Berlin: Springer, 1994, 292 pp.



\bibitem{Zapad-2}
Kloeden P.E., Platen E., Wright I.W. The approximation of multiple 
stochastic integrals. Stoch. Anal. Appl.,
10, 4 (1992), 431-441. 




\bibitem{Zapad-9}
Platen E., Bruti-Liberati N. Numerical Solution of Stochastic 
Differential Equations
with Jumps in Finance. Springer, Berlin-Heidelberg, 2010, 868 pp.


\bibitem{Rybakov1000}
Rybakov K.A. Orthogonal expansion of multiple It\^{o} stochastic integrals.
Electronic Journal "Differential Equations and Control Processes"
ISSN 1817-2172 (online),
3 (2021), 109-140. Available at:\\
{\color{blue}http://diffjournal.spbu.ru/EN/numbers/2021.3/article.1.8.html}




\bibitem{xyz1001}
Kloeden P.E., Neuenkirch A.
The pathwise convergence of approximation schemes for stochastic differential equations.
LMS Journal of Computation and Mathematics. 10 (2007), 235-253.


\bibitem{W-Z-1}
Wong E., Zakai M. On the convergence of ordinary integrals to 
stochastic integrals. Ann. Math. Stat.,
5, 36 (1965), 1560-1564.


\bibitem{W-Z-2}
Wong E., Zakai M. On the relation between ordinary and stochastic 
differential equations. Int. J. Eng. Sci., 3 (1965), 213-229.


\bibitem{Watanabe}
Ikeda N., Watanabe S. Stochastic
Differential Equations and Diffusion Processes.
2nd Edition. North-Holland Publishing Company,
Amsterdam, Oxford, New-York, 1989, 555 pp.


\bibitem{Lipt}
Liptser R.Sh., Shirjaev A.N. 
Statistics of Stochastic Processes: Nonlinear 
Filtering and Related Problems. [In Russian]. Moscow, Nauka, 1974, 696 pp.

\bibitem{7e}
Luo W. Wiener chaos expansion and numerical solutions of stochastic
partial differential equations. PhD thesis, 
California Inst. of Technology,
2006, 225 pp.




\end{thebibliography}
\end{document}